\documentclass[12pt]{article}
\usepackage{amssymb,amsmath,amsfonts}
\usepackage{relsize}
\usepackage{cite}
\usepackage{esvect}
\usepackage{txfonts}
\newtheorem{theorem}{Theorem}[section]
\newtheorem{definition}[theorem]{Definition}
\newtheorem{corollary}[theorem]{Corollary}
\newtheorem{proposition}[theorem]{Proposition}
\newtheorem{remark}[theorem]{Remark}
\newtheorem{lemma}[theorem]{Lemma}
\newtheorem{problem}[theorem]{Problem}
\newcommand {\Ac}      {{\mathcal A}}
\newcommand {\Cc}      {{\mathcal C}}
\newcommand {\Dc}      {{\mathcal D}}
\newcommand {\Ec}      {{\mathcal E}}
\newcommand {\Fc}      {{\mathcal F}}
\newcommand {\Jc}      {{\mathcal J}}
\newcommand {\Kc}      {{\mathcal K}}
\newcommand {\Lc}      {{\mathcal L}}
\newcommand {\Mc}      {{\mathcal M}}
\newcommand {\Nc}      {{\mathcal N}}
\newcommand {\Pc}      {{\mathcal P}}
\newcommand {\Qc}      {{\mathcal Q}}
\newcommand {\Tc}      {{\mathcal T}}
\newcommand {\Uc}      {{\mathcal U}}
\newcommand {\mZ}      {{\mathbb Z}}
\newcommand {\R}       {{\mathbb R}}
\newcommand {\N}       {{\mathbb N}}
\newcommand {\tA}      {\widetilde{A}}
\newcommand {\tK}      {\widetilde{K}}
\newcommand {\tQ}      {\widetilde{Q}}
\newcommand {\tF}      {\widetilde{F}}
\newcommand {\tB}      {\widetilde{B}}
\newcommand {\tP}      {\widetilde{P}}
\newcommand {\tL}      {\widetilde{{\mathcal L}}}
\newcommand {\tQc}     {\widetilde{\Qc}}
\newcommand {\Fcw}     {\widetilde{\Fc}}
\newcommand {\tx}      {\tilde{x}}
\newcommand {\teps}    {\tilde{\ve}}
\newcommand {\tgm}     {\widetilde{\gamma}}
\newcommand {\RN}      {\R^n}
\newcommand {\ve}      {\varepsilon}
\newcommand {\veh}     {\hat{\ve}}
\newcommand {\gmh}     {\hat{\gamma}}
\newcommand {\SP}      {Sobolev-Poincar\'{e}~}
\newcommand {\LOP}     {L^1_p(\RN)}
\newcommand {\LMOP}    {L_{p}^{m+1}(\RN)}
\newcommand {\LMP}     {L_{p}^{m}(\RN)}
\newcommand {\WMP}     {W_{p}^{m}(\RN)}
\newcommand {\LPRN}    {L_p(\RN)}
\newcommand {\CMD}     {C^{m,(d)}(\RN)}
\newcommand {\CMLD}    {C^{\,m-1,(d)}(\RN)}
\newcommand {\CMO}     {C^{m,\omega}(\RN)}
\newcommand {\CM}      {C^{m}(\RN)}
\newcommand {\CMON}    {C^{m-1}(\RN)}
\newcommand {\LIPD}    {\Lip(\RN;d)}
\newcommand {\intl}    {\int\limits}
\newcommand {\emp}     {\emptyset}
\newcommand {\MC}      {\Mc\Cc}
\newcommand {\phq}     {\varphi_q}
\newcommand {\dq}      {\delta_q}
\newcommand {\DPQ}     {\Dc_{p,q}(\RN)}
\newcommand {\PMRN}    {\Pc_{m-1}(\RN)}
\newcommand {\PMP}     {\Pc_{m}(\RN)}
\newcommand {\fs}      {f^{\sharp}_{E}}
\newcommand {\prs}     {\gamma}
\newcommand {\VP}      {{\bf P}}
\newcommand {\JV}      {{\bf \Jc}}
\newcommand {\LX}      {L^{(x)}}
\newcommand {\VST}     {\vspace*{1mm}}
\newcommand {\smed}    {\mathlarger{\sum}}
\newcommand {\sbig}    {\mathlarger{\mathlarger{\sum}}}
\newcommand {\q}       {90}
\newcommand {\LE}      {\Lc_E}
\newcommand {\hL}      {\hat{{\mathcal L}}}
\newcommand {\QL}      {Q^{(L)}}
\newcommand {\OV}      {\hat{Q}}
\newcommand {\PRL}     {\mathbb{PR}_{\mathbb{E}}}
\newcommand {\PRC}     {\mathbb{PR}_{\mathbb{E}}^{-1}}
\newcommand {\GE}      {\Gamma_E}
\newcommand {\vs}      {\theta}
\newcommand {\NCP}     {\|\VSH\|_{\LPRN}}
\newcommand {\PME}     {\|\VP\|_{m,p,E}}
\newcommand {\PMEW}    {\|\VP\|^*_{m,p,E}}
\newcommand {\VSH}     {\VP^\sharp_{m,E}}
\newcommand {\ld}      {\lambda\hspace{0.3mm}_d}
\newcommand {\lr}      {\leftrightarrow}
\newcommand {\lcr}     {\Leftrightarrow}
\newcommand {\TC}      {T^{-1}}
\newcommand {\hx}      {\bar{x}}
\newcommand {\hy}      {\bar{y}}
\newcommand {\DGE}     {d_{\GE}}
\newcommand {\VE}      {V}
\newcommand {\lw}      {\widetilde{\lambda}}
\newcommand {\lcap}    {\lambda}
\newcommand {\KM}      {\hat{K}}
\newcommand {\dw}      {\eta}
\newcommand {\td}      {10^{-5}}
\newcommand {\NW}      {\widetilde{\Nc}}
\DeclareMathOperator{\esssup}{ess\,sup}
\newcommand {\Lip}     {\operatorname{Lip}}
\newcommand {\diam}    {\operatorname{diam}}
\newcommand {\dist}    {\operatorname{dist}}
\newcommand {\supp}    {\operatorname{supp}}
\newcommand {\MP}      {\operatorname{M}}
\newcommand {\ED}      {\operatorname{Ed}}
\newcommand {\AD}      {\operatorname{Ad}}
\newcommand {\intr}    {\operatorname{int}}
\newcommand {\bx}      {\hspace{10mm}$\Box$}
\newcommand {\rbx}     {\hspace{10mm}$\vartriangleleft$}
\newcommand {\nn}      {\nonumber}
\newcommand {\rf}[1]    {(\ref{#1})}      %no references
\newcommand {\reff}[1] {\ref{#1}}         %no references
\newcommand{\lbl}[1]      {\label{#1}}       %no ref
\newcommand{\be}          {\begin{eqnarray}}
\newcommand{\bel}[1]      {\begin{eqnarray} \label{#1}}
\newcommand{\ee}           {\end{eqnarray}}
\newcommand {\SECT}[2] {\section*{\centerline{\normalsize
{\bf #1}}} \setcounter{section}{#2}
\setcounter{theorem}{0}\setcounter{equation}{0}}
\begin{document}
%----------------------------------------------------------
\parindent 1em
%----------------------------------------------------------
\parskip 0mm
%----------------------------------------------------------
\medskip
%@@@@@@@@@@@@@@@@@@@@@@@@@@@@@@@@@@@@@@@@@@@@@@@@@@@@@@@@@@
%@@@@@@@@@@@@@@@@@@@@@@@@@@@@@@@@@@@@@@@@@@@@@@@@@@@@@@@@@@
%----------------------------------------------------------
\centerline{\large{\bf Whitney-type extension theorems for jets}}\vspace*{5mm}
%----------------------------------------------------------
\centerline{\large{\bf generated by Sobolev functions}}
%----------------------------------------------------------
\vspace*{10mm} \centerline{By~  {\sc Pavel Shvartsman}}\vspace*{5 mm}
%----------------------------------------------------------
\centerline {\it Department of Mathematics, Technion - Israel Institute of Technology}\vspace*{2 mm}
%----------------------------------------------------------
\centerline{\it 32000 Haifa, Israel}\vspace*{2 mm}
%----------------------------------------------------------
\centerline{\it e-mail: pshv@tx.technion.ac.il}
%----------------------------------------------------------
\vspace*{10 mm}
%----------------------------------------------------------
\renewcommand{\thefootnote}{ }
%----------------------------------------------------------
\footnotetext[1]{{\it\hspace{-6mm}Math Subject
Classification} 46E35\smallskip
%----------------------------------------------------------
\par This research was supported by Grant No 2014055 from the United States-Israel Binational Science Foundation (BSF).} %----------------------------------------------------------
%@@@@@@@@@@@@@@@@@@@@@@@@@@@@@@@@@@@@@@@@@@@@@@@@@@@@@@@@@@
%@@@@@@@@@@@@@@@@@@@@@@@@@@@@@@@@@@@@@@@@@@@@@@@@@@@@@@@@@@
%----------------------------------------------------------
\begin{abstract} Let $\LMP$, $1\le p\le\infty$, be the homogeneous Sobolev space, and let $E\subset\RN$ be a closed set. For each $p>n$ and each non-negative integer $m$ we give an intrinsic characterization of the restrictions $\{D^\alpha F|_E:|\alpha|\le m\}$ to $E$ of $m$-jets generated by functions $F\in\LMOP$. Our trace criterion is expressed in terms of variations of corresponding Taylor remainders of $m$-jets evaluated on a certain family of ``well separated'' two point subsets of $E$. For $p=\infty$ this result coincides with the classical Whitney-Glaeser extension theorem for $m$-jets.
%----------------------------------------------------------
\par Our approach is based on a representation of the Sobolev space $\LMOP$, $p>n,$ as a union of $\CMD$-spaces where $d$ belongs to a family of metrics on $\RN$ with certain ``nice'' properties. Here $\CMD$ is the space of $C^m$-functions on $\RN$ whose partial derivatives of order $m$ are Lipschitz functions with respect to $d$.
This enables us to show that, for every  non-negative integer $m$ and every $p\in (n,\infty)$, the very same classical linear Whitney extension operator as in \cite{W1} provides an almost optimal extension of $m$-jets generated by $L^{m+1}_p$-functions.
%----------------------------------------------------------
\end{abstract}
%----------------------------------------------------------
\renewcommand{\contentsname}{ }
\tableofcontents
%----------------------------------------------------------
\addtocontents{toc}{{\centerline{\sc{Contents}}}
\vspace*{10mm}\par}
%----------------------------------------------------------
%@@@@@@@@@@@@@@@@@@@@@@@@@@@@@@@@@@@@@@@@@@@@@@@@@@@@@@@@@@
%@@@@@@@@@@@@@@@@@@@@@@@@@@@@@@@@@@@@@@@@@@@@@@@@@@@@@@@@@@
%@@@@@@@@@@@@@@@@@@@@@@@@@@@@@@@@@@@@@@@@@@@@@@@@@@@@@@@@@@
%@@@@@@@@@@@@@@@@@@@@@@@@@      @@@@@@@@@@@@@@@@@@@@@@@@@@@
%@@@@@@@@@@@@@@@@@@@@@@@          @@@@@@@@@@@@@@@@@@@@@@@@@
%@@@@@@@@@@@@@@@@@@@@@              @@@@@@@@@@@@@@@@@@@@@@@
%@@@@@@@@@@@@@@@@@@@     SECTION 1    @@@@@@@@@@@@@@@@@@@@@
%@@@@@@@@@@@@@@@@@@@@@              @@@@@@@@@@@@@@@@@@@@@@@
%@@@@@@@@@@@@@@@@@@@@@@@          @@@@@@@@@@@@@@@@@@@@@@@@@
%@@@@@@@@@@@@@@@@@@@@@@@@@      @@@@@@@@@@@@@@@@@@@@@@@@@@@
%@@@@@@@@@@@@@@@@@@@@@@@@@@@@@@@@@@@@@@@@@@@@@@@@@@@@@@@@@@
%@@@@@@@@@@@@@@@@@@@@@@@@@@@@@@@@@@@@@@@@@@@@@@@@@@@@@@@@@@
%@@@@@@@@@@@@@@@@@@@@@@@@@@@@@@@@@@@@@@@@@@@@@@@@@@@@@@@@@@
%----------------------------------------------------------
\SECT{1. Introduction.}{1}
%----------------------------------------------------------
%@@@@@@@@@@@@@@@@@@@@@@@@@@@@@@@@@@@@@@@@@@@@@@@@@@@@@@@@@@
%----------------------------------------------------------
\addtocontents{toc}{~~~1. Introduction.\hfill \thepage\par\VST}
%----------------------------------------------------------
%@@@@@@@@@@@@@@@@@@@@@@@@@@@@@@@@@@@@@@@@@@@@@@@@@@@@@@@@@@
%----------------------------------------------------------
\par {\bf 1.1. Main definitions and main results.}
%----------------------------------------------------------
\addtocontents{toc}{~~~~1.1. Main definitions and main results. \hfill \thepage\par}\medskip
%----------------------------------------------------------
\indent\par Given $m\in\N$ and $p\in[1,\infty]$, let $\LMP$ be the homogeneous Sobolev space of all (equi\-valence classes of) real valued functions $f\in L_{1,loc}(\RN)$ whose distributional partial derivatives on $\RN$ {\it of order $m$} belong to $L_p(\RN)$. $\LMP$ is seminormed by
%@@@@@@@@@@@@@@@@@@@@@@@@@@@@@@@@@@@@@@@@@@@@@@@@@@@@@@@@@@
%----------------------------------------------------------
\bel{DEF-LMP}
\|f\|_{\LMP}:= \smed_{|\alpha|=m}\,\, \|D^\alpha
f\|_{L_p(\RN)}\,.
\ee
%----------------------------------------------------------
\par As usual, we let $\WMP$ denote the corresponding Sobolev space of all functions $f\in \LPRN$ whose distributional partial derivatives on $\RN$ of {\it all orders up to $m$} belong to $\LPRN$. This space is normed by
%----------------------------------------------------------
$$
\|f\|_{\WMP}:=\smed_{|\alpha|\le m}\, \|D^\alpha f\|_{\LPRN}.
$$
%----------------------------------------------------------
%@@@@@@@@@@@@@@@@@@@@@@@@@@@@@@@@@@@@@@@@@@@@@@@@@@@@@@@@@@
%@@@@@@@@@@@@@@@@@@@@@@@@@@@@@@@@@@@@@@@@@@@@@@@@@@@@@@@@@@
%@@@@@@@@@@@@@@@@@@@@@@@@@@@@@@@@@@@@@@@@@@@@@@@@@@@@@@@@@@
%----------------------------------------------------------
%@@@@@@@@@@@@@@@@@@@@@@@@@@@@@@@@@@@@@@@@@@@@@@@@@@@@@@@@@@
%----------------------------------------------------------
\par We recall that, by the Sobolev imbedding theorem, every function $f\in \LMP$, $p>n$, can be redefined, if necessary, on a set of Lebesgue measure zero so that it belongs to the space $C^{m-1}(\RN)$. (See e.g., \cite{M}, p. 73.) Thus, for $p>n$, we can identify each element $f\in \LMP$ with its unique $C^{m-1}$-representative on $\RN$. This will allow us to restrict our attention to the case of Sobolev $C^{m-1}$-functions. \medskip
%----------------------------------------------------------
%@@@@@@@@@@@@@@@@@@@@@@@@@@@@@@@@@@@@@@@@@@@@@@@@@@@@@@@@@@ %----------------------------------------------------------
\par Given a function $F\in C^m(\RN)$ and $x\in\RN$, we let %----------------------------------------------------------
$$
T^m_x[F](y)=\smed_{|\alpha|\le m}\,\,\frac{1}{\alpha !}\, D^{\alpha}
F(x)(y-x)^{\alpha},~~~~y\in\RN,
$$
%----------------------------------------------------------
denote the Taylor polynomial of $F$ of degree $m$ at $x$.
By $\PMP$ we denote the space of all polynomials of degree at most $m$ defined on $\RN$.
%---------------------------------------------------------
\par  In this paper we study the following
%----------------------------------------------------------
\begin{problem}\lbl{PR-MAIN} {\em Let $p\in[1,\infty]$, $m\in\N$, and let $E$ be a closed subset of $\RN$. Suppose that, for each point $x\in E$, we are given a polynomial $P_x\in\PMRN$. We ask two questions:\smallskip
%---------------------------------------------------------
\par {\it 1. How can we decide whether there exists a function $F\in \LMP$ such that $T^{m-1}_x[F]=P_x$ for all $x\in E$\,?}\smallskip
%---------------------------------------------------------
%@@@@@@@@@@@@@@@@@@@@@@@@@@@@@@@@@@@@@@@@@@@@@@@@@@@@@@@@@
%---------------------------------------------------------
\par 2. Consider the $\LMP$-norms of all functions $F\in \LMP$ such that $T^{m-1}_x[F]=P_x$ on $E$.  {\it How small can these norms be?}}
%----------------------------------------------------------
\end{problem}
%----------------------------------------------------------
\smallskip
%----------------------------------------------------------
\par This problem is a variant of a classical extension problem posed by H. Whitney in 1934 in his pioneering papers \cite{W1,W2}, namely: {\it How can one tell whether a given function $f$ defined on an arbitrary subset $E\subset\RN$ extends to a $C^m$-function on all of $\RN$?} Over the years since 1934 this problem, often called the Whitney Extension Problem, has attracted a lot of attention, and there is an extensive literature devoted to different aspects of this problem and its analogues for various spaces of smooth functions. Among the multitude of results known so far we mention those in the papers \cite{BMP1,BS1,BS2,BS3,F2,F-J,
F4,F3,F-IO,F-Bl,F9,FIL,G,Is,Sh1,Sh2,Sh3,Sh4,Z1,Z2}. We refer the reader to all of these  papers and references therein, for numerous results and techniques concerning this topic.
%---------------------------------------------------------
%@@@@@@@@@@@@@@@@@@@@@@@@@@@@@@@@@@@@@@@@@@@@@@@@@@@@@@@@@@
%@@@@@@@@@@@@@@@@@@@@@@@@@@@@@@@@@@@@@@@@@@@@@@@@@@@@@@@@@@
%@@@@@@@@@@@@@@@@@@@@@@@@@@@@@@@@@@@@@@@@@@@@@@@@@@@@@@@@@@
%---------------------------------------------------------
\par In \cite{Sh2} we solved Problem \reff{PR-MAIN} for $m=1$, $n\in\N$ and $p>n$. In the present paper we give a complete solution to Problem \reff{PR-MAIN} for arbitrary $m,n\in\N$ and $p>n$.
%---------------------------------------------------------
\par Note that, for the case $p=\infty$, Problem \reff{PR-MAIN} was solved by Whitney \cite{W1} and Glaeser \cite{G}. In that case the space $L^m_\infty(\RN)$ can be identified with the space $C^{m-1,1}(\RN)$ of all $C^{m-1}$-functions on $\RN$ whose partial derivatives of
order $m-1$ all satisfy Lipschitz conditions.
%---------------------------------------------------------
\par We recall the statement of the classical Whitney \cite{W1}-Glaeser \cite{G} extension theorem: {\it Let $E$ be an arbitrary closed subset of $\RN$. There exists a $C^{m-1}$-function $F\in L^m_\infty(\RN)$ such that $T^{m-1}_x[F]=P_x$ for every $x\in E$ if and only if}
%----------------------------------------------------------
\bel{W-GL}
\sup_{x,y\in E,\, x\ne y}\,\,\smed_{|\alpha|\le m-1}\frac{|D^\alpha P_{x}(x)-D^\alpha P_{y}(x)|}
{\|x-y\|^{m-|\alpha|}} <\infty\,.
\ee
%----------------------------------------------------------
\par Theorem \reff{EX-TK}, our main contribution in this  paper, generalizes this result to the case $n<p<\infty$. Let us prepare the ingredients that are needed to formulate this theorem.
We shall always use the word ``cube'' to mean a closed cube in $\RN$ whose sides are parallel to the coordinate axes. $Q(c,r)$ will denote the cube in $\RN$ centered at $c$ with side length $2r$. Given $\lambda >0$ and a cube $Q$, $\lambda Q$ will denote the dilation of $Q$ with respect to its center by a factor of $\lambda$. (Thus $\lambda\,Q(c,r)=Q(c,\lambda r)$.) The Lebesgue measure of a measurable set $A\subset \RN$ will be denoted by $\left|A\right|$.
%----------------------------------------------------------
\par Let $\VP=\{P_x: x\in E\}$ be a family of polynomials of degree at most $m$ indexed by points of a given closed subset $E$ of $\RN$. (Thus $P_x\in \PMP$ for every $x\in E$.) Following \cite{F9} we refer to $\VP$ as {\it a Whitney $m$-field defined on $E$}.
%----------------------------------------------------------
\par We say that a function $F\in \CM$ {\it agrees with the  Whitney $m$-field $\VP=\{P_x: x\in E\}$ on $E$}, if $T^{m}_x[F]=P_x$ for each $x\in E$. In that case we also refer to $\VP$ as the Whitney $m$-field on $E$ {\it generated by $F$} or as {\it the $m$-jet generated by $F$.}
We define the $L^m_p$-``norm'' of the $m$-jet $\VP=\{P_x: x\in E\}$ by
%---------------------------------------------------------
\bel{N-VP}
\PME:=\inf\left\{\|F\|_{L^m_p(\RN)}:F\in \LMP,\, T^{m-1}_x[F]=P_x~~\text{for every}~~x\in E\right\}.
\ee
%---------------------------------------------------------
%@@@@@@@@@@@@@@@@@@@@@@@@@@@@@@@@@@@@@@@@@@@@@@@@@@@@@@@@@
%----------------------------------------------------------
\par We also need the following notion:
%----------------------------------------------------------
%@@@@@@@@@@@@@@@@@@@@@@@@@@@@@@@@@@@@@@@@@@@@@@@@@@@@@@@@@@
%@@@@@@@@@@@@@@@@@@@@@@@@@@@@@@@@@@@@@@@@@@@@@@@@@@@@@@@@@@
%@@@@@@@@@@@@@@@@@@@@@@@@@@@@@@@@@@@@@@@@@@@@@@@@@@@@@@@@@@
%----------------------------------------------------------
%@@@@@@@@@@@@@@@@@@@@@@@@@@@@@@@@@@@@@@@@@@@@@@@@@@@@@@@@@@
\begin{definition}\lbl{DF-PR} {\em Let $\prs\ge 1$ and let $\Ac=\{\{x_i,y_i\}: i\in I\}$ be a family of two point subsets of $\RN$. We say that the family $\Ac$ is {\it $\prs$-sparse} if there exists a collection $\{Q_i: i\in I\}$ of pairwise disjoint cubes in $\RN$ such that
%---------------------------------------------------------
\bel{XY-Q}
x_i, y_i\in \prs\, Q_i~~~\text{and}~~~
\diam Q_i\le \prs \|x_i-y_i\|~~~
\text{for all}~~i\in I.
\ee
%----------------------------------------------------------
}
%----------------------------------------------------------
\end{definition}
%----------------------------------------------------------
\par We can now explicitly formulate the above mentioned main result of this paper.
%----------------------------------------------------------
%@@@@@@@@@@@@@@@@@@@@@@@@@@@@@@@@@@@@@@@@@@@@@@@@@@@@@@@@@@
%@@@@@@@@@@@@@@@@@@@@@@@@@@@@@@@@@@@@@@@@@@@@@@@@@@@@@@@@@@
%@@@@@@@@@@@@@@@@@@@@@@@@@@@@@@@@@@@@@@@@@@@@@@@@@@@@@@@@@@
%@@@@@@@@@@@@@@@@@@@@@@@@@@@@@@@@@@@@@@@@@@@@@@@@@@@@@@@@@@
%@@@@@@@@@@@@@@@@@@@@@@@@@@@@@@@@@@@@@@@@@@@@@@@@@@@@@@@@@@
%@@@@@@@@@@@@@@@@@@@@@@@@@@@@@@@@@@@@@@@@@@@@@@@@@@@@@@@@@@
%----------------------------------------------------------
\begin{theorem}\lbl{EX-TK} Let $m\in\N$, $p\in(n,\infty)$, and let $E$ be a closed subset of $\RN$.
There exists an absolute constant $\gamma \ge 1$ for which the following result holds:
%----------------------------------------------------------
\par Suppose we are given a family $\VP=\{P_x: x\in E\}$ of polynomials of degree at most $m-1$ indexed by points of $E$. There exists a $C^{m-1}$-function $F\in\LMP$ such that %----------------------------------------------------------
\bel{J-PE}
T_{x}^{m-1}[F]=P_{x}~~~~\text{for every}~~~~x\in E
\ee
%----------------------------------------------------------
if and only if the following quantity
%----------------------------------------------------------
\bel{N-P-NEW}
\Nc_{m,p,E}(\VP):=\sup\left\{\,\smed_{i=1}^k\,\,
\smed_{|\alpha|\le m-1}\frac{|D^\alpha P_{x_i}(x_i)-D^\alpha P_{y_i}(x_i)|^p}
{\|x_i-y_i\|^{(m-|\alpha|)p-n}}\right\}^{1/p}
\ee
%----------------------------------------------------------
is finite.  Here the supremum is taken over all finite $\gamma$-sparse collections $\{\{x_i, y_i\} : i = 1, ..., k\}$ of two point subsets of $E$.
%----------------------------------------------------------
\par Furthermore,
%---------------------------------------------------------
\bel{EQV-JET}
\PME\sim \Nc_{m,p,E}(\VP)\,.
\ee
%---------------------------------------------------------
The constants of equivalence in \rf{EQV-JET} depend only on $m,n$ and $p$.
%----------------------------------------------------------
\end{theorem}
%----------------------------------------------------------
%@@@@@@@@@@@@@@@@@@@@@@@@@@@@@@@@@@@@@@@@@@@@@@@@@@@@@@@@@@
%@@@@@@@@@@@@@@@@@@@@@@@@@@@@@@@@@@@@@@@@@@@@@@@@@@@@@@@@@@
%@@@@@@@@@@@@@@@@@@@@@@@@@@@@@@@@@@@@@@@@@@@@@@@@@@@@@@@@@@
%@@@@@@@@@@@@@@@@@@@@@@@@@@@@@@@@@@@@@@@@@@@@@@@@@@@@@@@@@@
%@@@@@@@@@@@@@@@@@@@@@@@@@@@@@@@@@@@@@@@@@@@@@@@@@@@@@@@@@@
%@@@@@@@@@@@@@@@@@@@@@@@@@@@@@@@@@@@@@@@@@@@@@@@@@@@@@@@@@@
%----------------------------------------------------------
\par We refer to this result as {\it a variational criterion for the traces of $L^m_p$-jets}.\medskip
%----------------------------------------------------------
%@@@@@@@@@@@@@@@@@@@@@@@@@@@@@@@@@@@@@@@@@@@@@@@@@@@@@@@@@@
%@@@@@@@@@@@@@@@@@@@@@@@@@@@@@@@@@@@@@@@@@@@@@@@@@@@@@@@@@@
%@@@@@@@@@@@@@@@@@@                           @@@@@@@@@@@@
%@@@@@@@@@@@@@@@@@@ A REFINEMENT FOR FINITE E @@@@@@@@@@@@@ @@@@@@@@@@@@@@@@@@@                           @@@@@@@@@@@@
%@@@@@@@@@@@@@@@@@@@@@@@@@@@@@@@@@@@@@@@@@@@@@@@@@@@@@@@@@@
%@@@@@@@@@@@@@@@@@@@@@@@@@@@@@@@@@@@@@@@@@@@@@@@@@@@@@@@@@@
%@@@@@@@@@@@@@@@@@@@@@@@@@@@@@@@@@@@@@@@@@@@@@@@@@@@@@@@@@@
%@@@@@@@@@@@@@@@@@@@@@@@@@@@@@@@@@@@@@@@@@@@@@@@@@@@@@@@@@@
%----------------------------------------------------------
\par The variational criterion describes the structure of the linear space
%---------------------------------------------------------
$$
\JV(\LMP)|_E:=\{\VP=\{T^{m-1}_x[F]: x\in E\}: F\in \LMP\}
$$
%----------------------------------------------------------
of all Whitney $(m-1)$-fields on $E$ generated by $C^{m-1}$-functions belonging to $\LMP$. We refer to the space $\JV(\LMP)|_E$ as the trace jet-space of $\LMP$ to $E$.
The functional $\|\cdot\|_{m,p,E}$ defined above is
a seminorm on $\JV(\LMP)|_E$.
\par When $E=\RN$ we simply write $\JV(\LMP)$ instead of $\JV(\LMP)|_E$. \medskip
%----------------------------------------------------------
\par The criterion \rf{N-P-NEW} and equivalence \rf{EQV-JET} show which properties of $\VP$ on $E$ control its almost optimal extension to a jet generated by a function from $\LMP$. At first sight, this criterion seems to be extremely difficult to check ``in practice'' even
when $E$ is a finite set. The statement of Theorem \reff{EX-TK} indicates that this check requires one to verify a very large number of conditions, a number which is very much larger than the cardinality of $E$. (It is the number of {\it all} $\gamma$-sparse families of two point subsets of $E$.)
%----------------------------------------------------------
\par However an examination of our proof for a finite set $E$ shows that, after all, it is only necessary to deal with just one of these families of two point subsets of $E$. This is a certain special $\gamma$-sparse family which is constructed by a particular step in the proof. {\it It is enough to examine the behavior of given field $\VP$ only on this particular family}. Furthermore, this special family of two point subsets of $E$ generates a certain graph structure on $E$ with rather nice properties.
%----------------------------------------------------------
%@@@@@@@@@@@@@@@@@@@@@@@@@@@@@@@@@@@@@@@@@@@@@@@@@@@@@@@@@@
%@@@@@@@@@@@@@@@@@@@@@@@@@@@@@@@@@@@@@@@@@@@@@@@@@@@@@@@@@@
%@@@@@@@@@@@@@@@@@@@@@@@@@@@@@@@@@@@@@@@@@@@@@@@@@@@@@@@@@@
%@@@@@@@@@@@@@@@@@@@@@@@@@@@@@@@@@@@@@@@@@@@@@@@@@@@@@@@@@@
%----------------------------------------------------------
\begin{definition}\lbl{DF-GRPR} {\em Let $\prs\ge 1$ and let $\Gamma$ be a graph whose set of vertices belongs to $\RN$. The set of edges of $\Gamma$ defines a family $\Ac$ of two point subsets of $E$ consisting of all pairs of points which are joined by an edge. We say that the graph $\Gamma$ is {\it $\prs$-sparse} if the family $\Ac$ defined in this way is $\gamma$-sparse, in the sense of Definition \reff{DF-PR}.
\smallskip
%----------------------------------------------------------
\par In other words, $\Gamma$ is $\gamma$-sparse if there exists a family
%----------------------------------------------------------
$$
\{Q_{xy}: x,y~ \text{are vertices of}~\Gamma~ \text{joined by an edge}\}
$$
%----------------------------------------------------------
of {\it pairwise disjoint} cubes $Q_{xy}$ such that $x,y\in\gamma\, Q_{xy}$ and $\diam Q_{xy}\le \gamma\,\|x-y\|$ whenever  $x$ and $y$ are two arbitrary vertices joined by an edge in the graph $\Gamma$.}
%----------------------------------------------------------
\end{definition}
%----------------------------------------------------------
\par The following theorem provides a refinement of the variational criterion given in Theorem \reff{EX-TK} for {\it finite subsets} of $\RN$.
%----------------------------------------------------------
%@@@@@@@@@@@@@@@@@@@@@@@@@@@@@@@@@@@@@@@@@@@@@@@@@@@@@@@@@@
%@@@@@@@@@@@@@@@@@@@@@@@@@@@@@@@@@@@@@@@@@@@@@@@@@@@@@@@@@@
%@@@@@@@@@@@@@@@@@@@@@@@@@@@@@@@@@@@@@@@@@@@@@@@@@@@@@@@@@@
%@@@@@@@@@@@@@@@@@@@@@@@@@@@@@@@@@@@@@@@@@@@@@@@@@@@@@@@@@@
%----------------------------------------------------------
\begin{theorem}\lbl{EX-REFTK} Let $E$ be a finite subset of $\RN$. There exists a constant $\gamma=\gamma(n)\ge 1$ and a connected $\gamma$-sparse graph $\GE$ whose set of vertices coincides with $E$ which has the following properties:\smallskip
%----------------------------------------------------------
\par (i). The degree of each vertex of $\GE$ is bounded by a constant $C=C(n)$;
\smallskip
%----------------------------------------------------------
\par (ii). Let $m\in\N$, $n<p<\infty$, and let  $\VP=\{P_x: x\in E\}$ be a Whitney $(m-1)$-field on $E$. Then
%---------------------------------------------------------
\bel{GE-NORM}
\PME\,\sim \left\{\sbig_{x,y\in E,\,
x\,\underset{\Gamma_E}{\lr}\,y}\,\,\,\,\sbig_{|\alpha|\le m-1}\frac{|D^\alpha P_{x}(x)-D^\alpha P_{y}(x)|^p}
{\|x-y\|^{(m-|\alpha|)p-n}}\right\}^{\frac1p}
\ee
%---------------------------------------------------------
where the first sum is taken over all points $x,y\in E$ joined by an edge in the graph $\GE$ ($x\underset{\Gamma_E}{\lr}y$).
%----------------------------------------------------------
\par The constants in this equivalence depend only on $m,n$ and $p$.
%----------------------------------------------------------
\end{theorem}
%----------------------------------------------------------
%@@@@@@@@@@@@@@@@@@@@@@@@@@@@@@@@@@@@@@@@@@@@@@@@@@@@@@@@@@
%@@@@@@@@@@@@@@@@@@@@@@@@@@@@@@@@@@@@@@@@@@@@@@@@@@@@@@@@@@
%@@@@@@@@@@@@@@@@@@@@@@@@@@@@@@@@@@@@@@@@@@@@@@@@@@@@@@@@@@
%@@@@@@@@@@@@@@@@@@@@@@@@@@@@@@@@@@@@@@@@@@@@@@@@@@@@@@@@@@
%----------------------------------------------------------
\begin{remark}{\em Using some obvious modifications of the proof of Theorem \reff{EX-REFTK} we can also obtain an analogue of \rf{GE-NORM} for $p=\infty$, namely that the equivalence
%---------------------------------------------------------
\bel{P-INF1}
\|\VP\|_{m,\infty,E}\,\sim \sup_{x,y\in E,\, x\lr y}\,\sbig_{|\alpha|\le m-1}\frac{|D^\alpha P_{x}(x)-D^\alpha P_{y}(x)|}
{\|x-y\|^{m-|\alpha|}}
\ee
%---------------------------------------------------------
holds, for every finite set $E$, and every
Whitney $(m-1)$-field $\VP$ on $E$. As the notation indicates, here the supremum is taken over all points $x,y\in E$ joined by an edge in $\GE$.
The
constants in this equivalence depend only on $m,n$ and $p$. C.f., \rf{W-GL}.\rbx
%----------------------------------------------------------
}
%----------------------------------------------------------
\end{remark}
%----------------------------------------------------------
%@@@@@@@@@@@@@@@@@@@@@@@@@@@@@@@@@@@@@@@@@@@@@@@@@@@@@@@@@@
%@@@@@@@@@@@@@@@@@@@@@@@@@@@@@@@@@@@@@@@@@@@@@@@@@@@@@@@@@@
%@@@@@@@@@@@@@@@@@@@@@@@@@@@@@@@@@@@@@@@@@@@@@@@@@@@@@@@@@@
%@@@@@@@@@@@@@@@@@@@@@@@@@@@@@@@@@@@@@@@@@@@@@@@@@@@@@@@@@@
%----------------------------------------------------------
\par Combining the result of Theorem \reff{EX-REFTK} with Theorem \reff{EX-TK} we obtain the following
%----------------------------------------------------------
%@@@@@@@@@@@@@@@@@@@@@@@@@@@@@@@@@@@@@@@@@@@@@@@@@@@@@@@@@@
%@@@@@@@@@@@@@@@@@@@@@@@@@@@@@@@@@@@@@@@@@@@@@@@@@@@@@@@@@@
%@@@@@@@@@@@@@@@@@@@@@@@@@@@@@@@@@@@@@@@@@@@@@@@@@@@@@@@@@@
%@@@@@@@@@@@@@@@@@@@@@@@@@@@@@@@@@@@@@@@@@@@@@@@@@@@@@@@@@@
%@@@@@@@@@@@@@@@@@@@@@@@@@@@@@@@@@@@@@@@@@@@@@@@@@@@@@@@@@@
%@@@@@@@@@@@@@@@@@@@@@@@@@@@@@@@@@@@@@@@@@@@@@@@@@@@@@@@@@@
%----------------------------------------------------------
\begin{theorem}\lbl{CR-TREE} Let $m\in\N$, $p\in(n,\infty)$, and let $E$ be a closed subset of $\RN$. Let $\gamma=\gamma(n)\ge 1$ and $C=C(n)\ge 1$ be the same as in Theorem \reff{EX-REFTK}. Then  the following result holds:
%----------------------------------------------------------
\par Suppose we are given a family $\VP=\{P_x: x\in E\}$ of polynomials of degree at most $m-1$ indexed by points of $E$. There exists a $C^{m-1}$-function $F\in\LMP$ such that $T_{x}^{m-1}[F]=P_{x}$ for every $x\in E$ if and only if the following quantity
%---------------------------------------------------------
\bel{NP-WIGL}
\NW_{m,p,E}(\VP):=\sup_{\Gamma} \left\{\sbig_{x\underset{\Gamma}{\lr}y}
\,\,\,\,\sbig_{|\alpha|\le m-1}\frac{|D^\alpha P_{x}(x)-D^\alpha P_{y}(x)|^p}
{\|x-y\|^{(m-|\alpha|)p-n}}\right\}^{\frac1p}
\ee
%---------------------------------------------------------
is finite. Here the supremum is taken over all finite connected $\gamma$-sparse graphs $\Gamma$ with vertices in $E$ and with the degree of each vertex bounded by $C$. The first sum in \rf{NP-WIGL} is taken over all vertices $x,y$ of a graph $\Gamma$ joined by an edge in $\Gamma$ ($x\underset{\Gamma}{\lr}y$).
%----------------------------------------------------------
%@@@@@@@@@@@@@@@@@@@@@@@@@@@@@@@@@@@@@@@@@@@@@@@@@@@@@@@@@@
%----------------------------------------------------------
\par Furthermore, $\PME\sim \NW_{m,p,E}(\VP)$ with the constants of equivalence depending only on $m,n$ and $p$.
%----------------------------------------------------------
\end{theorem}
%----------------------------------------------------------
%@@@@@@@@@@@@@@@@@@@@@@@@@@@@@@@@@@@@@@@@@@@@@@@@@@@@@@@@@@
%@@@@@@@@@@@@@@@@@@@@@@@@@@@@@@@@@@@@@@@@@@@@@@@@@@@@@@@@@@
%@@@@@@@@@@@@@@@@@@@@@@@@@@@@@@@@@@@@@@@@@@@@@@@@@@@@@@@@@@
%@@@@@@@@@@@@@@@@@@@@@@@@@@@@@@@@@@@@@@@@@@@@@@@@@@@@@@@@@@
%@@@@@@@@@@@@@@@@@@@@@@@@@@@@@@@@@@@@@@@@@@@@@@@@@@@@@@@@@@
%@@@@@@@@@@@@@@@@@@@@@@@@@@@@@@@@@@@@@@@@@@@@@@@@@@@@@@@@@@
%----------------------------------------------------------
\par See Remark \reff{CR-TREE-PR}.\medskip
%----------------------------------------------------------
%@@@@@@@@@@@@@@@@@@@@@@@@@@@@@@@@@@@@@@@@@@@@@@@@@@@@@@@@@@
%@@@@@@@@@@@@@@@@@@@@@@@@@@@@@@@@@@@@@@@@@@@@@@@@@@@@@@@@@@
%----------------------------------------------------------
\par The next theorem gives another characterization  of $L^m_p$-jets on $E$ expressed in terms of $L_p$-norms of certain maximal functions. For each family $\VP=\{P_x\in\PMRN: x\in E\}$ of polynomials we let $\VSH$ denote a certain kind of ``sharp maximal function'' associated with $\VP$ which is defined by
%----------------------------------------------------------
\bel{IEF}
\VSH(x):=\sup_{y,z\in E,\,\, y\ne z}\,\, \frac{|P_y(x)-P_z(x)|}
{\|x-y\|^{m}+\|x-z\|^{m}},~~~~~x\in\RN.
%----------------------------------------------------------
\ee
%----------------------------------------------------------
%@@@@@@@@@@@@@@@@@@@@@@@@@@@@@@@@@@@@@@@@@@@@@@@@@@@@@@@@@@
%@@@@@@@@@@@@@@@@@@@@@@@@@@@@@@@@@@@@@@@@@@@@@@@@@@@@@@@@@@
%@@@@@@@@@@@@@@@@@@@@@@@@@@@@@@@@@@@@@@@@@@@@@@@@@@@@@@@@@@
%@@@@@@@@@@@@@@@@@@@@@@@@@@@@@@@@@@@@@@@@@@@@@@@@@@@@@@@@@@
%----------------------------------------------------------
\begin{theorem} \lbl{JET-S} Let $m$, $p$, $E$ and
$\VP=\{P_x\in\PMRN: x\in E\}$ be as in the statement of Theorem \reff{EX-TK}. Then there exists a $C^{m-1}$-function $F\in\LMP$ such that $T_{x}^{m-1}[F]=P_{x}$ \ for every $x\in E$
if and only if\, $\VSH\in\LPRN$.
%----------------------------------------------------------
\par Furthermore,
%---------------------------------------------------------
\bel{PME-VS}
\PME\sim \|\VSH\|_{\LPRN}
\ee
%---------------------------------------------------------
with the constants in this equivalence depending only on $m,n,$ and $p$.
%----------------------------------------------------------
\end{theorem}
%----------------------------------------------------------
%\medskip
%----------------------------------------------------------
%@@@@@@@@@@@@@@@@@@@@@@@@@@@@@@@@@@@@@@@@@@@@@@@@@@@@@@@@@@
%@@@@@@@@@@@@@@@@@@@@@@@@@@@@@@@@@@@@@@@@@@@@@@@@@@@@@@@@@@
%@@@@@@@@@@@@@@@@@@@@@@@@@@@@@@@@@@@@@@@@@@@@@@@@@@@@@@@@@@
%@@@@@@@@@@@@@@@@@@@@@@@@@@@@@@@@@@@@@@@@@@@@@@@@@@@@@@@@@@
%----------------------------------------------------------
\begin{remark}{\em Let us note two interesting results related to the equivalence \rf{PME-VS}.
%----------------------------------------------------------
\par First, one can show that \rf{PME-VS} becomes an equality whenever $m=1$ and $p=\infty$, i.e,
$$\|\VP\|_{1,\infty,E}=\|\VP^\sharp_{1,E}\|_{L_\infty(\RN)}.$$
This easily follows from the McShane extension theorem \cite{McS} which states that every function $f\in \Lip(E)$ extends to a function $F\in\Lip(\RN)$ such that $\|F\|_{\Lip(\RN)}=\|f\|_{\Lip(E)}$. (Recall that $L^1_\infty(\RN)=\Lip(\RN)$.)
%----------------------------------------------------------
\par Secondly, we can express a deep and interesting result proved by Le Gruyer \cite{LG} in our notation here by the formula
%----------------------------------------------------------
$$
\|\VP\|_{2,\infty,E}=\tfrac12\,
\|\VP^\sharp_{2,E}\|_{L_\infty(\RN)}.
$$
%----------------------------------------------------------
\par For further results related to the equivence \rf{PME-VS} for the case $m=1,2$ and $p=\infty$ we refer  the reader to papers \cite{HHC,HL,WEL} and references therein.}\rbx
%----------------------------------------------------------
\end{remark}
%----------------------------------------------------------
%@@@@@@@@@@@@@@@@@@@@@@@@@@@@@@@@@@@@@@@@@@@@@@@@@@@@@@@@@@
%@@@@@@@@@@@@@@@@@@@@@@@@@@@@@@@@@@@@@@@@@@@@@@@@@@@@@@@@@@
%@@@@@@@@@@@@@@@@@@@@@@@@@@@@@@@@@@@@@@@@@@@@@@@@@@@@@@@@@@
%@@@@@@@@@@@@@@@@@@@@@@@@@@@@@@@@@@@@@@@@@@@@@@@@@@@@@@@@@@
%----------------------------------------------------------
\par The case $m=1$ merits particular attention. Theorems \reff{EX-TK} and \reff{JET-S} immediately imply the following characterization of the trace space $\LOP|_E$.
%----------------------------------------------------------
%@@@@@@@@@@@@@@@@@@@@@@@@@@@@@@@@@@@@@@@@@@@@@@@@@@@@@@@@@@
%@@@@@@@@@@@@@@@@@@@@@@@@@@@@@@@@@@@@@@@@@@@@@@@@@@@@@@@@@@
%@@@@@@@@@@@@@@@@@@@@@@@@@@@@@@@@@@@@@@@@@@@@@@@@@@@@@@@@@@
%@@@@@@@@@@@@@@@@@@@@@@@@@@@@@@@@@@@@@@@@@@@@@@@@@@@@@@@@@@
%----------------------------------------------------------
\begin{theorem}\lbl{LONE-TR} Let $p\in(n,\infty)$, and let $E$ be a closed subset of $\RN$. Let $f$ be a function defined on $E$. Then the following three statements are equivalent:\medskip
%----------------------------------------------------------
%@@@@@@@@@@@@@@@@@@@@@@@@@@@@@@@@@@@@@@@@@@@@@@@@@@@@@@@@@@
%@@@@@@@@@@@@@@@@@@@@@@@@@@@@@@@@@@@@@@@@@@@@@@@@@@@@@@@@@@
%----------------------------------------------------------
\par (i). $f$ extends to a continuous function $F\in\LOP$;\smallskip
%----------------------------------------------------------
%@@@@@@@@@@@@@@@@@@@@@@@@@@@@@@@@@@@@@@@@@@@@@@@@@@@@@@@@@@
%@@@@@@@@@@@@@@@@@@@@@@@@@@@@@@@@@@@@@@@@@@@@@@@@@@@@@@@@@@
%----------------------------------------------------------
\par (ii). The following quantity
%----------------------------------------------------------
$$
\Phi_{p,E}(f):=\sup
\left\{\,\smed_{i=1}^k\,\,\frac{|f(x_i)-f(y_i)|^p}
{\|x_i-y_i\|^{p-n}}\right\}^{1/p}
$$
%----------------------------------------------------------
is finite. Here the supremum is taken over all finite $\gamma$-sparse collections $\{\{x_i, y_i\} : i = 1, ..., k\}$ of two point subsets of $E$, and $\gamma\ge 1$ is a certain absolute constant;\smallskip
%----------------------------------------------------------
%@@@@@@@@@@@@@@@@@@@@@@@@@@@@@@@@@@@@@@@@@@@@@@@@@@@@@@@@@@
%@@@@@@@@@@@@@@@@@@@@@@@@@@@@@@@@@@@@@@@@@@@@@@@@@@@@@@@@@@
%----------------------------------------------------------
\par (iii). The following quantity
%---------------------------------------------------------
$$
\Psi_{p,E}(f):=\left(\,\,\intl\limits_{\RN}\,
\sup_{y,z\in E,\,y\ne z}\,\, \frac{|f(y)-f(z)|^p}{\|x-y\|^p+\|x-z\|^p}\,\,dx
\right)^{1/p}<\infty\,.
$$
%----------------------------------------------------------
%@@@@@@@@@@@@@@@@@@@@@@@@@@@@@@@@@@@@@@@@@@@@@@@@@@@@@@@@@@
%@@@@@@@@@@@@@@@@@@@@@@@@@@@@@@@@@@@@@@@@@@@@@@@@@@@@@@@@@@
%----------------------------------------------------------
\par Furthermore,
%---------------------------------------------------------
\bel{EQ-L1}
\|f\|_{\LOP|_E}\sim \Phi_{p,E}(f)\sim \Psi_{p,E}(f)
\ee
%---------------------------------------------------------
where
$$
\|f\|_{\LOP|_E}:=\inf\left\{\|F\|_{\LOP}:F\in\LOP\cap C(\RN), F|_E=f\right\}.
$$
%---------------------------------------------------------
\smallskip
The constants of equivalences in \rf{EQ-L1} depend only on $n$ and $p$.
%----------------------------------------------------------
\end{theorem}
%----------------------------------------------------------
%@@@@@@@@@@@@@@@@@@@@@@@@@@@@@@@@@@@@@@@@@@@@@@@@@@@@@@@@@@
%@@@@@@@@@@@@@@@@@@@@@@@@@@@@@@@@@@@@@@@@@@@@@@@@@@@@@@@@@@
%----------------------------------------------------------
\par As we have mentioned above, the restrictions of $L^1_p(\RN)$-functions to subsets of $\RN$, $p>n$, have been studied in \cite{Sh2}. The trace criteria given in part (ii) and (iii) of Theorem \reff{LONE-TR} are improvements of Theorem 1.1 and Theorem 1.4, part (i), of \cite{Sh2}. \smallskip
%----------------------------------------------------------
%@@@@@@@@@@@@@@@@@@@@@@@@@@@@@@@@@@@@@@@@@@@@@@@@@@@@@@@@@@
%@@@@@@@@@@@@@@@@@@@@@@@@@@@@@@@@@@@@@@@@@@@@@@@@@@@@@@@@@@
%----------------------------------------------------------
\par Our next result, Theorem \reff{LIN-OP}, states that there is a solution to Problem \reff{PR-MAIN} which depends {\it linearly} on the initial data, i.e., the families of polynomials $\{P_x\in\PMRN: x\in E\}$.
%----------------------------------------------------------
%@@@@@@@@@@@@@@@@@@@@@@@@@@@@@@@@@@@@@@@@@@@@@@@@@@@@@@@@@@
%@@@@@@@@@@@@@@@@@@@@@@@@@@@@@@@@@@@@@@@@@@@@@@@@@@@@@@@@@@
%@@@@@@@@@@@@@@@@@@@@@@@@@@@@@@@@@@@@@@@@@@@@@@@@@@@@@@@@@@
%@@@@@@@@@@@@@@@@@@@@@@@@@@@@@@@@@@@@@@@@@@@@@@@@@@@@@@@@@@
%@@@@@@@@@@@@@@@@@@@@@@@@@@@@@@@@@@@@@@@@@@@@@@@@@@@@@@@@@@
%@@@@@@@@@@@@@@@@@@@@@@@@@@@@@@@@@@@@@@@@@@@@@@@@@@@@@@@@@@
%----------------------------------------------------------
\begin{theorem} \lbl{LIN-OP} For every closed subset $E\subset\RN$ and every $p>n$ there exists a continuous linear operator
%----------------------------------------------------------
$$
\Fc:\,\JV(\LMP)|_E\,\to\, \LMP
$$
%----------------------------------------------------------
such that for every Whitney $(m-1)$-field $\VP=\{P_x:x\in E\}\in \JV(\LMP)|_E$ the function $\Fc(\VP)$ agrees with $\VP$ on $E$, i.e.,
%----------------------------------------------------------
$$
T^{m-1}_x[\Fc(\VP)]=P_x~~~\text{for every}~~~x\in E.
$$
%----------------------------------------------------------
%@@@@@@@@@@@@@@@@@@@@@@@@@@@@@@@@@@@@@@@@@@@@@@@@@@@@@@@@@@
%----------------------------------------------------------
\par Furthermore, the operator norm of $\Fc$ is bounded by a constant depending only on $m,n,$ and $p$.
%----------------------------------------------------------
\end{theorem}
%----------------------------------------------------------
%@@@@@@@@@@@@@@@@@@@@@@@@@@@@@@@@@@@@@@@@@@@@@@@@@@@@@@@@@@
%@@@@@@@@@@@@@@@@@@@@@@@@@@@@@@@@@@@@@@@@@@@@@@@@@@@@@@@@@@
%@@@@@@@@@@@@@@@@@@@@@@@@@@@@@@@@@@@@@@@@@@@@@@@@@@@@@@@@@@
%@@@@@@@@@@@@@@@@@@@@@@@@@@@@@@@@@@@@@@@@@@@@@@@@@@@@@@@@@@
%----------------------------------------------------------
\begin{remark} {\em Actually we show, perhaps surprisingly, that the very same classical linear Whitney extension operator
%----------------------------------------------------------
$$
\Fc_W:\,\JV(\CMON)|_E\,\to\, \CMON \,,
$$
%----------------------------------------------------------
which was introduced in  \cite{W1} for the space $\JV(\CMON)|_E$ of Whitney $(m-1)$-fields on $E$ generated by $C^{m-1}$-functions, has the properties described in Theorem \reff{LIN-OP}. See Remark \reff{VM-E}. %----------------------------------------------------------
\par In fact, this ``universality'' of the Whitney extension operator for the scale of $(m-1)$-jets ge\-ne\-rated by the spaces $\LMP$ for all $p>n$
is the consequence of another result which will be formulated below, namely that it is possible to represent the space $\LMP$ as a union (see \rf{DM-PR} below) of  $C^{m-1}$-spaces with certain specific Lipschitz properties of higher partial derivatives of order $m-1$.}\rbx
%----------------------------------------------------------
\end{remark}
%----------------------------------------------------------
\medskip
%----------------------------------------------------------
\par Our results mentioned so far deal only with homogeneous Sobolev spaces. But we also wish to treat the non-homogeneous (normed) case.
We defer doing this to Section 8. There we present analogues of Theorems \reff{EX-TK} and \reff{LIN-OP} for spaces of $(m-1)$-jets generated by functions from the normed Sobolev space $\WMP$. See Theorems \reff{EX-WTK}, \reff{WP-PR}, and \reff{LIN-WOP}.
%----------------------------------------------------------
\bigskip\smallskip
%----------------------------------------------------------
%@@@@@@@@@@@@@@@@@@@@@@@@@@@@@@@@@@@@@@@@@@@@@@@@@@@@@@@@@@
%@@@@@@@@@@@@@@@@@@@@@@@@@@@@@@@@@@@@@@@@@@@@@@@@@@@@@@@@@@
%@@@@@@@@@@@@@@@@@@@@@@@@@@@@@@@@@@@@@@@@@@@@@@@@@@@@@@@@@@
%@@@@@@@@@@@@@@@@@@@@@@@@@@@@@@@@@@@@@@@@@@@@@@@@@@@@@@@@@@
%@@@@@@@@@@@@@@@@@@@@@@@@@@@@@@@@@@@@@@@@@@@@@@@@@@@@@@@@@@
%@@@@@@@@@@@@@@@@@@@@@@@@@@@@@@@@@@@@@@@@@@@@@@@@@@@@@@@@@@
%---------------------------------------------------------- %@@@@@@@@@@@@@@@@@@@@@@@@@@@@@@@@@@@@@@@@@@@@@@@@@@@@@@@@@@
%----------------------------------------------------------
\par {\bf 1.2. Our approach: Sobolev-Poincar\'e inequality and $\CMD$-spaces.}\medskip
%----------------------------------------------------------
\addtocontents{toc}{~~~~1.2. Our approach: Sobolev-Poincar\'e inequality and $\CMD$-spaces. \hfill \thepage\par\VST}
%----------------------------------------------------------
\indent\par Let us briefly describe the main ideas of our approach. Let $d$ be a metric on $\RN$ and let $\LIPD$ be the space of functions on $\RN$ satisfying the Lipschitz condition with respect to the metric $d$. $\LIPD$ is equipped with the standard seminorm
%----------------------------------------------------------
$$
\|F\|_{\LIPD}:=\sup_{x,y\in\RN,\,x\ne y} \,\,\frac{|F(x)-F(y)|}{d(x,y)}\,.
$$
%----------------------------------------------------------
%@@@@@@@@@@@@@@@@@@@@@@@@@@@@@@@@@@@@@@@@@@@@@@@@@@@@@@@@@@
%@@@@@@@@@@@@@@@@@@@@@@@@@@@@@@@@@@@@@@@@@@@@@@@@@@@@@@@@@@
%@@@@@@@@@@@@@@@@@@@@@@@@@@@@@@@@@@@@@@@@@@@@@@@@@@@@@@@@@@
%----------------------------------------------------------
%@@@@@@@@@@@@@@@@@@@@@@@@@@@@@@@@@@@@@@@@@@@@@@@@@@@@@@@@@@
\begin{definition}\lbl{CMD} {\em Let $\CMD$ be a space of $C^m$-functions on $\RN$ whose partial derivatives of order $m$ are Lipschitz continuous on $\RN$ with respect to $d$. This space is seminormed by
%----------------------------------------------------------
$$
\|F\|_{\CMD}:=\smed_{|\alpha|=m}\,\|D^\alpha F\|_{\Lip(\RN;d)}.
$$
%----------------------------------------------------------
}
%----------------------------------------------------------
\end{definition}
%----------------------------------------------------------
%@@@@@@@@@@@@@@@@@@@@@@@@@@@@@@@@@@@@@@@@@@@@@@@@@@@@@@@@@@
%----------------------------------------------------------
\par The main ingredient of our approach is a representation of the Sobolev space $\LMP$,  $p\in(n,\infty)$, as a union of {\it $C^{m-1,(d)}$-spaces}  where $d$ belongs to a certain family $\Dc$ of metrics on $\RN$:
%----------------------------------------------------------
\bel{LMP-UN}
\LMP=\bigcup_{d\in\Dc}\,\,\CMLD.
\ee
%----------------------------------------------------------
%@@@@@@@@@@@@@@@@@@@@@@@@@@@@@@@@@@@@@@@@@@@@@@@@@@@@@@@@@@
%----------------------------------------------------------
See \rf{DM-PR}.
%----------------------------------------------------------
\par We obtain this representation using a slight modification of the classical Sobolev-Poincar\'{e} inequality for $L^m_p$-functions. More specifically, our aim is to reformulate the \SP inequality for $p>n$ in the form of a certain {\it Lipschitz condition} for partial derivatives of order $m-1$ {\it with respect to a certain metric on $\RN$.}
%----------------------------------------------------------
%@@@@@@@@@@@@@@@@@@@@@@@@@@@@@@@@@@@@@@@@@@@@@@@@@@@@@@@@@@
%----------------------------------------------------------
\par For functions $F\in \LMP$ it is convenient to use the notation
%----------------------------------------------------------
$$
\nabla^m F(x):=\left(\smed_{|\alpha|= m}(D^\alpha F(x))^2\right)^{\frac{1}{2}},~~~~x\in\RN,
$$
%----------------------------------------------------------
so that
%----------------------------------------------------------
$$
\|F\|_{\LMP}\sim\|\nabla^mF\|_{\LPRN}
$$
%----------------------------------------------------------
with constants depending only on $n,m$ and $p$. See \rf{DEF-LMP}.
%----------------------------------------------------------
\par We recall a variant of the \SP inequality for $\LMP$-functions which holds whenever $p>n$:
%----------------------------------------------------------
\par Let $q\in(n,p)$ and let $F\in \LMP$. Then for every cube $Q\subset\RN$, every  $x,y\in Q$ and every multiindex $\beta$, $|\beta|\le m-1,$ the following inequality
%----------------------------------------------------------
\bel{SP-BTQ}
|D^{\beta}F(x)-D^{\beta}(T^{m-1}_y[F])(x)|\le C\,(\diam Q)^{m-|\beta|}\left(\frac{1}{|Q|}
\intl_{Q}(\nabla^m F(u))^qdu\right)^{\frac{1}{q}}
\ee
%----------------------------------------------------------
holds. Here $C>0$ is a constant depending only on $n,m$ and $q$. See, e.g. \cite{M}, p. 61, or \cite{MP}, p. 55.
%----------------------------------------------------------
\par In particular, for every $\alpha$ which satisfies $|\alpha|=m-1,$
%----------------------------------------------------------
$$
|D^{\alpha}F(x)-D^{\alpha}F(y)|\le C\,\|x-y\|\left(\frac{1}{|Q_{xy}|}\intl_{Q_{xy}}
(\nabla^m F(u))^q\,du\right)^{\frac{1}{q}},~~~x,y\in\RN,
$$
%----------------------------------------------------------
where
%----------------------------------------------------------
$$
Q_{xy}:=Q(x,\|x-y\|).
$$
%----------------------------------------------------------
Hence
%----------------------------------------------------------
\bel{F-XY-N}
|D^{\alpha}F(x)-D^{\alpha}F(y)|\le \|x-y\|
\sup_{Q\ni x,y} \left(\frac{1}{|Q|}\intl_{Q} h^q(u)\,du\right)^{\frac{1}{q}},~~~x,y\in\RN,
\ee
%----------------------------------------------------------
where $h=C(n,m,q)\|\nabla^m F\|$.
%----------------------------------------------------------
%@@@@@@@@@@@@@@@@@@@@@@@@@@@@@@@@@@@@@@@@@@@@@@@@@@@@@@@@@@
%----------------------------------------------------------
\par Let $h\in L_p(\RN)$ be an arbitrary non-negative function. Inequality \rf{F-XY-N} motivates us to introduce the function
%---------------------------------------------------------
\bel{DQ}
\dq(x,y:h)=\|x-y\| \sup_{Q\ni x,y} \left(\frac{1}{|Q|}\intl_{Q} h^q(u)\,du\right)^{\frac{1}{q}}~,~~~~~x,y\in \RN.
\ee
%----------------------------------------------------------
\par By \rf{F-XY-N}, for each $p\in(n,\infty)$, $q\in(n,p)$, and every $F\in \LMP$ there exists a non-negative function $h\in \LPRN$ such that $\|h\|_{\LPRN}\le C(n,p,q)\|F\|_{\LMP}$ and
%---------------------------------------------------------
\bel{F-SP}
|D^{\alpha}F(x)-D^{\alpha}F(y)|\le \dq(x,y:h)~~~\text{for all}~\alpha,|\alpha|=m-1,~\text{and}~x,y\in\RN.
\ee
%---------------------------------------------------------
\par As we prove below, see Theorem \reff{SM-MAIN}, the converse statement is also true: Let $n<q<p<\infty$ and let $F$ be a $C^{m-1}$-function on $\RN$. Suppose that there exists a non-negative function $h\in\LPRN$ such that \rf{F-SP} holds. Then $F\in \LMP$ and $\|F\|_{\LMP}\le C\|h\|_{\LPRN}$ with $C$ depending only on $n,m,p$ and $q$. %----------------------------------------------------------
\smallskip
%----------------------------------------------------------
%@@@@@@@@@@@@@@@@@@@@@@@@@@@@@@@@@@@@@@@@@@@@@@@@@@@@@@@@@@
%@@@@@@@@@@@@@@@@@@@@@@@@@@@@@@@@@@@@@@@@@@@@@@@@@@@@@@@@@@
%@@@@@@@@@@@@@@@@@@@@@@@@@@@@@@@@@@@@@@@@@@@@@@@@@@@@@@@@@@
%@@@@@@@@@@@@@@@@@@@@@@@@@@@@@@@@@@@@@@@@@@@@@@@@@@@@@@@@@@
%@@@@@@@@@@@@@@@@@@@@@@@@@@@@@@@@@@@@@@@@@@@@@@@@@@@@@@@@@@
%---------------------------------------------------------
\par Thus we have an alternative equivalent definition of the homogeneous Sobolev space $\LMP$ in terms of the ``Lipschitz-like'' conditions \rf{F-SP} with respect to the functions $\dq(h)$ whenever $h\in \LPRN$.
%---------------------------------------------------------
\par Of course, these observations lead us to the desired representation \rf{LMP-UN}, provided $\dq(h)$ is a {\it metric} on $\RN$.
%----------------------------------------------------------
\par However, in general, $\dq(h)$ {\it is not a metric.}
Nevertheless, we prove that {\it the geodesic distance} $d_q(h)$ associated with the function $\dq(h)$ is equivalent to $\dq(h)$. Recall that, given $x,y\in\RN$ this distance is defined by the formula
%----------------------------------------------------------
\bel{GD-M}
d_q(x,y:h):=\inf\,\smed_{i=0}^{m-1}\,\dq(x_i,x_{i+1}:h)
\ee
%---------------------------------------------------------
where the infimum is taken over all finite sequences of points $\{x_0,x_1,...,x_m\}$ in $\RN$ such that $x_0=x$  and $x_m=y$. In Section 2 we prove the following
%----------------------------------------------------------
%@@@@@@@@@@@@@@@@@@@@@@@@@@@@@@@@@@@@@@@@@@@@@@@@@@@@@@@@@@
%@@@@@@@@@@@@@@@@@@@@@@@@@@@@@@@@@@@@@@@@@@@@@@@@@@@@@@@@@@
%@@@@@@@@@@@@@@@@@@@@@@@@@@@@@@@@@@@@@@@@@@@@@@@@@@@@@@@@@@
%@@@@@@@@@@@@@@@@@@@@@@@@@@@@@@@@@@@@@@@@@@@@@@@@@@@@@@@@@@
%@@@@@@@@@@@@@@@@@@@@@@@@@@@@@@@@@@@@@@@@@@@@@@@@@@@@@@@@@@
%@@@@@@@@@@@@@@@@@@@@@@@@@@@@@@@@@@@@@@@@@@@@@@@@@@@@@@@@@@
%@@@@@@@@@@@@@@@@@@@@@@@@@@@@@@@@@@@@@@@@@@@@@@@@@@@@@@@@@@
%----------------------------------------------------------
\begin{theorem}\lbl{MAIN-MT} Let $n\le q<\infty$ and let $h\in L_{q,loc}(\RN)$ be a non-negative function. Then for every $x,y\in\RN$ the following inequality
%----------------------------------------------------------
\bel{W-G}
d_q(x,y:h)\le \dq(x,y:h)\le 16\, d_q(x,y:h)
\ee
%----------------------------------------------------------
holds.
%----------------------------------------------------------
%@@@@@@@@@@@@@@@@@@@@@@@@@@@@@@@@@@@@@@@@@@@@@@@@@@@@@@@@@@
%@@@@@@@@@@@@@@@@@@@@@@@@@@@@@@@@@@@@@@@@@@@@@@@@@@@@@@@@@@
%----------------------------------------------------------
\end{theorem}
%----------------------------------------------------------
%@@@@@@@@@@@@@@@@@@@@@@@@@@@@@@@@@@@@@@@@@@@@@@@@@@@@@@@@@@
%@@@@@@@@@@@@@@@@@@@@@@@@@@@@@@@@@@@@@@@@@@@@@@@@@@@@@@@@@@
%@@@@@@@@@@@@@@@@@@@@@@@@@@@@@@@@@@@@@@@@@@@@@@@@@@@@@@@@@@
%@@@@@@@@@@@@@@@@@@@@@@@@@@@@@@@@@@@@@@@@@@@@@@@@@@@@@@@@@@
%@@@@@@@@@@@@@@@@@@@@@@@@@@@@@@@@@@@@@@@@@@@@@@@@@@@@@@@@@@
%----------------------------------------------------------
\par This leads us to the following result which enables us to explicitly describe the family of metrics $\Dc$ which provides the representation in \rf{LMP-UN}.
%----------------------------------------------------------
%@@@@@@@@@@@@@@@@@@@@@@@@@@@@@@@@@@@@@@@@@@@@@@@@@@@@@@@@@@
%@@@@@@@@@@@@@@@@@@@@@@@@@@@@@@@@@@@@@@@@@@@@@@@@@@@@@@@@@@
%@@@@@@@@@@@@@@@@@@@@@@@@@@@@@@@@@@@@@@@@@@@@@@@@@@@@@@@@@@
%@@@@@@@@@@@@@@@@@@@@@@@@@@@@@@@@@@@@@@@@@@@@@@@@@@@@@@@@@@
%@@@@@@@@@@@@@@@@@@@@@@@@@@@@@@@@@@@@@@@@@@@@@@@@@@@@@@@@@@
%----------------------------------------------------------
\begin{theorem}\lbl{SM-MAIN} Let $m\in\N$ and let $n<q< p<\infty$. A $C^{m-1}$-function $F$ belongs to $\LMP$ if and only if there exists a non-negative function $h\in\LPRN$ such that for every multiindex $\alpha$, $|\alpha|=m-1$, and every $x,y\in\RN$ the following inequality
%----------------------------------------------------------
\bel{LIP2}
|D^\alpha F(x)-D^\alpha F(y)|\le d_q(x,y:h)
\ee
%----------------------------------------------------------
holds. Furthermore,
%----------------------------------------------------------
$$
\|F\|_{\LMP}\sim \inf \|h\|_{\LPRN}.
$$
%----------------------------------------------------------
The constants of this equivalence depend only on $m,n,q,$ and $p$.
%----------------------------------------------------------
%@@@@@@@@@@@@@@@@@@@@@@@@@@@@@@@@@@@@@@@@@@@@@@@@@@@@@@@@@@
%@@@@@@@@@@@@@@@@@@@@@@@@@@@@@@@@@@@@@@@@@@@@@@@@@@@@@@@@@@
%----------------------------------------------------------
\end{theorem}
%----------------------------------------------------------
%@@@@@@@@@@@@@@@@@@@@@@@@@@@@@@@@@@@@@@@@@@@@@@@@@@@@@@@@@@
%@@@@@@@@@@@@@@@@@@@@@@@@@@@@@@@@@@@@@@@@@@@@@@@@@@@@@@@@@@
%@@@@@@@@@@@@@@@@@@@@@@@@@@@@@@@@@@@@@@@@@@@@@@@@@@@@@@@@@@
%@@@@@@@@@@@@@@@@@@@@@@@@@@@@@@@@@@@@@@@@@@@@@@@@@@@@@@@@@@
%@@@@@@@@@@@@@@@@@@@@@@@@@@@@@@@@@@@@@@@@@@@@@@@@@@@@@@@@@@
%----------------------------------------------------------
\par We prove this theorem in Section 2. By this result, for each $p\in(n,\infty)$ and each $q\in (n,p)$, the space $\LMP$ can be represented in the form
%----------------------------------------------------------
\bel{DM-PR}
\LMP=\bigcup_{d\in\DPQ}\CMLD
\ee
%----------------------------------------------------------
where
%----------------------------------------------------------
\bel{D-LIP}
\DPQ:=\{d_q(h): h\in\LPRN, h\ge 0\}.
\ee
%----------------------------------------------------------
%@@@@@@@@@@@@@@@@@@@@@@@@@@@@@@@@@@@@@@@@@@@@@@@@@@@@@@@@@@
%@@@@@@@@@@@@@@@@@@@@@@@@@@@@@@@@@@@@@@@@@@@@@@@@@@@@@@@@@@
%@@@@@@@@@@@@@@@@@@@@@@@@@@@@@@@@@@@@@@@@@@@@@@@@@@@@@@@@@@
%@@@@@@@@@@@@@@@@@@@@@@@@@@@@@@@@@@@@@@@@@@@@@@@@@@@@@@@@@@
%@@@@@@@@@@@@@@@@@@@@@@@@@@@@@@@@@@@@@@@@@@@@@@@@@@@@@@@@@@
%@@@@@@@@@@@@@@@@@@@@@@@@@@@@@@@@@@@@@@@@@@@@@@@@@@@@@@@@@@
%@@@@@@@@@@@@@@@@@@@@@@@@@@@@@@@@@@@@@@@@@@@@@@@@@@@@@@@@@@
%----------------------------------------------------------
\par In particular, this implies the following representation of the space $\LOP$:
%----------------------------------------------------------
$$
\LOP=\bigcup_{d\in\,\DPQ}\LIPD\,.
$$
%----------------------------------------------------------
\par Representation \rf{DM-PR} motivates us to study an analog of Problem \reff{PR-MAIN} for the spaces $\CMD$ whenever $d\in\DPQ$. The following theorem provides a solution to this problem.
%----------------------------------------------------------
%@@@@@@@@@@@@@@@@@@@@@@@@@@@@@@@@@@@@@@@@@@@@@@@@@@@@@@@@@@
%@@@@@@@@@@@@@@@@@@@@@@@@@@@@@@@@@@@@@@@@@@@@@@@@@@@@@@@@@@
%@@@@@@@@@@@@@@@@@@@@@@@@@@@@@@@@@@@@@@@@@@@@@@@@@@@@@@@@@@
%@@@@@@@@@@@@@@@@@@@@@@@@@@@@@@@@@@@@@@@@@@@@@@@@@@@@@@@@@@
%@@@@@@@@@@@@@@@@@@@@@@@@@@@@@@@@@@@@@@@@@@@@@@@@@@@@@@@@@@
%@@@@@@@@@@@@@@@@@@@@@@@@@@@@@@@@@@@@@@@@@@@@@@@@@@@@@@@@@@
%@@@@@@@@@@@@@@@@@@@@@@@@@@@@@@@@@@@@@@@@@@@@@@@@@@@@@@@@@@
%----------------------------------------------------------
\begin{theorem} \lbl{CMD-E-NEW} Let $n<q<p<\infty$, $m\in\N$, and let $d\in\DPQ$. Let
%----------------------------------------------------------
$$
\VP=\{P_x\in\PMP: x\in E\}
$$
%----------------------------------------------------------
be a Whitney $m$-field defined on a closed set $E\subset\RN$. There exists a function $F\in\CMD$ which agrees with $\VP$ on $E$ if and only if the following quantity
%---------------------------------------------------------
$$
\Lc_{m,d}(\VP):=\sup_{x,y\in E,\,x\ne y}\,\, \sbig_{|\alpha|\le m}\,\,\, \frac{|D^{\alpha}P_x(x)-D^{\alpha}P_y(x)|}
{\|x-y\|^{m-|\alpha|}\,d(x,y)}
$$
%---------------------------------------------------------
is finite. Furthermore,
%---------------------------------------------------------
\bel{LMD-INF}
\Lc_{m,d}(\VP)\sim\inf\left\{\|F\|_{\CMD}:F\in \CMD, F~\text{agrees with}~\VP~\text{on}~E\right\}
\ee
%---------------------------------------------------------
with constants of equivalence depending only on $n,p,q,$ and $m$.
%----------------------------------------------------------
\end{theorem}
%----------------------------------------------------------
%@@@@@@@@@@@@@@@@@@@@@@@@@@@@@@@@@@@@@@@@@@@@@@@@@@@@@@@@@@
%@@@@@@@@@@@@@@@@@@@@@@@@@@@@@@@@@@@@@@@@@@@@@@@@@@@@@@@@@@
%@@@@@@@@@@@@@@@@@@@@@@@@@@@@@@@@@@@@@@@@@@@@@@@@@@@@@@@@@@
%@@@@@@@@@@@@@@@@@@@@@@@@@@@@@@@@@@@@@@@@@@@@@@@@@@@@@@@@@@
%@@@@@@@@@@@@@@@@@@@@@@@@@@@@@@@@@@@@@@@@@@@@@@@@@@@@@@@@@@
%----------------------------------------------------------
\begin{remark} \lbl{VM-E}{\em When $d$ is the Euclidean metric this result of course coincides with the Whitney-Glaeser extension theorem. See \rf{W-GL}.
%----------------------------------------------------------
\par Theorem \reff{CMD-E-NEW} is also well-known for the space $\CMO:= C^{m,(d_\omega)}(\RN)$ where
%----------------------------------------------------------
$$
d_{\omega}(x,y)=\omega(\|x-y\|),~~~~ x,y\in\RN,
$$
%----------------------------------------------------------
and $\omega$ is a {\it concave non-decreasing continuous} function on $[0,\infty)$ such that $\omega(0)=0$. See
\cite{St}, Ch. 6, \S 2.2.3 and \S 4.6.\medskip
%----------------------------------------------------------
\par We refer to any function $\omega$ with these properties as a ``modulus of continuity'', and we let $\MC$ denote the family of all ``moduli of continuity''. It should be noted that, for each $\omega \in \MC$, the metric space $(\RN,d_\omega)$ is known in the literature as {\it the metric transform of $\RN$ by $\omega$} or {\it the $\omega$-metric transform} of $\RN$.\medskip
%----------------------------------------------------------
\par In \cite{St} it is shown that the classical Whitney extension operator $\Fc_W$ constructed in \cite{W1} maps the space $\Tc(\CMO)|_E$ into the space $\CMO$ in such a way that, for every Whitney $m$-field
$$\VP=\{P_x\in\PMP:x\in E\}$$ belonging to $\Tc(\CMO)|_E$,
the function $\Fc_W(\VP)$ agrees with $\VP$. Furthermore, the operator norm of $\Fc_W$ satisfies the inequality  $$\|\Fc_W\|\le C(n,m).\vspace*{3mm}$$
%----------------------------------------------------------
%@@@@@@@@@@@@@@@@@@@@@@@@@@@@@@@@@@@@@@@@@@@@@@@@@@@@@@@@@@
%----------------------------------------------------------
\par We prove that the same statement is true for {\it an arbitrary metric $d\in\DPQ$}. Our proof of this property of $\Fc_W$ mainly follows the classical scheme used in \cite{St} for the metric $d_\omega$ with $\omega\in\MC$. As we show below, see Proposition \reff{VMT}, the fact that  the proofs of Theorem \reff{CMD-E-NEW} for $d_\omega$ and for arbitrary $d\in\DPQ$ are similar to each other
can be explained by the following important property of metrics in $\DPQ$: for each $d\in\DPQ$ there exists a mapping $\RN\ni x\mapsto\omega_x\in\MC$ such that
%----------------------------------------------------------
\bel{DXY-1}
d(x,y)\sim \omega _x(\|x-y\|)\,, ~~~~x,y\in\RN,
\ee
%---------------------------------------------------------
with constants of equivalence depending only on $n,p,$ and $q$. This equivalence motivates us to refer to the metric space $(\RN,d)$ as a \textit{variable metric transform} of $\RN$.
\smallskip
%----------------------------------------------------------
\par Property \rf{DXY-1} and the representation \rf{LMP-UN} also explain the above-mentioned phenomenon of the universality of the Whitney extension operator for the scale of $\LMP$-spaces when $p>n$. See Remark \reff{MT-RM}.\rbx}
%---------------------------------------------------------
%@@@@@@@@@@@@@@@@@@@@@@@@@@@@@@@@@@@@@@@@@@@@@@@@@@@@@@@@@@
%----------------------------------------------------------
\end{remark}
%----------------------------------------------------------
%@@@@@@@@@@@@@@@@@@@@@@@@@@@@@@@@@@@@@@@@@@@@@@@@@@@@@@@@@@
%----------------------------------------------------------
\par
The proof of Theorem \reff{JET-S} is given in Section 5. It is a consequence of Theorem \reff{CMD-E-NEW} and the representation \rf{LMP-UN}.
\smallskip
%----------------------------------------------------------
\par Theorem \reff{EX-TK} is proven in Sections 6 and 7. The difficult part of its proof is the sufficiency, which relies on a modification of the Whitney extension method \cite{W1} used in the author's paper \cite{Sh3}.
The main feature of that modification is that, instead of treating each Whitney cube separately, as is done in the original extension method in \cite{W1}, we deal simultaneously with all members of certain {\it families} of Whitney cubes.
We refer to these families of Whitney cubes as {\it lacunae}. See Section 6.\smallskip
%----------------------------------------------------------
\par Finally, direct calculations of the $L^m_p$-norm of the extension obtained by the lacunary extension method lead us to criterion \rf{N-P-NEW} whiich proves the sufficiency part of Theorem \reff{EX-TK}.\bigskip
%----------------------------------------------------------
%@@@@@@@@@@@@@@@@@@@@@@@@@@@@@@@@@@@@@@@@@@@@@@@@@@@@@@@@@@
%@@@@@@@@@@@@@@@@@@@@@@@@@@@@@@@@@@@@@@@@@@@@@@@@@@@@@@@@@@
%@@@@@@@@@@@@@@@@@@@@@@@@@@@@@@@@@@@@@@@@@@@@@@@@@@@@@@@@@@
%@@@@@@@@@@@@@@@@@@@@@@@@@@@@@@@@@@@@@@@@@@@@@@@@@@@@@@@@@@
%@@@@@@@@@@@@@@@@@@@@@@@@@@@@@@@@@@@@@@@@@@@@@@@@@@@@@@@@@@
%@@@@@@@@@@@@@@@@@@@@@@@@@@@@@@@@@@@@@@@@@@@@@@@@@@@@@@@@@@
%@@@@@@@@@@@@@@@@@@@@@@@@@@@@@@@@@@@@@@@@@@@@@@@@@@@@@@@@@@
%----------------------------------------------------------
\par {\bf Acknowledgements.} I am very thankful to M. Cwikel for useful suggestions and remarks.
\par The results of this paper were presented at the ``Whitney Problems Workshop" in Luminy, France, October 2015. I am very grateful to all participants of this conference for stimulating discussions and valuable advice.
%---------------------------------------------------------
\par I am also pleased to thank P. Haj{\l}asz who kindly drew my attention to several important papers and results related to Theorems \reff{MAIN-MT} and \reff{SM-MAIN}.
%----------------------------------------------------------
%@@@@@@@@@@@@@@@@@@@@@@@@@@@@@@@@@@@@@@@@@@@@@@@@@@@@@@@@@@
%@@@@@@@@@@@@@@@@@@@@@@@@@@@@@@@@@@@@@@@@@@@@@@@@@@@@@@@@@@
%@@@@@@@@@@@@@@@@@@@@@@@@@@@@@@@@@@@@@@@@@@@@@@@@@@@@@@@@@@
%@@@@@@@@@@@@@@@@@@@@@@@@@      @@@@@@@@@@@@@@@@@@@@@@@@@@@
%@@@@@@@@@@@@@@@@@@@@@@@          @@@@@@@@@@@@@@@@@@@@@@@@@
%@@@@@@@@@@@@@@@@@@@@@              @@@@@@@@@@@@@@@@@@@@@@@
%@@@@@@@@@@@@@@@@@@@     SECTION 2    @@@@@@@@@@@@@@@@@@@@@
%@@@@@@@@@@@@@@@@@@@@@              @@@@@@@@@@@@@@@@@@@@@@@
%@@@@@@@@@@@@@@@@@@@@@@@          @@@@@@@@@@@@@@@@@@@@@@@@@
%@@@@@@@@@@@@@@@@@@@@@@@@@      @@@@@@@@@@@@@@@@@@@@@@@@@@@
%@@@@@@@@@@@@@@@@@@@@@@@@@@@@@@@@@@@@@@@@@@@@@@@@@@@@@@@@@@
%@@@@@@@@@@@@@@@@@@@@@@@@@@@@@@@@@@@@@@@@@@@@@@@@@@@@@@@@@@
%@@@@@@@@@@@@@@@@@@@@@@@@@@@@@@@@@@@@@@@@@@@@@@@@@@@@@@@@@@
%----------------------------------------------------------
\SECT{2. Metrics on $\RN$ generated by the Sobolev-Poincar\'e inequality.}{2}
%----------------------------------------------------------
\addtocontents{toc}{2. Metrics on $\RN$ generated by the Sobolev-Poincar\'e inequality. \hfill\thepage\par\VST}
%----------------------------------------------------------
%@@@@@@@@@@@@@@@@@@@@@@@@@@@@@@@@@@@@@@@@@@@@@@@@@@@@@@@@@@@
%----------------------------------------------------------
\indent
%@@@@@@@@@@@@@@@@@@@@@@@@@@@@@@@@@@@@@@@@@@@@@@@@@@@@@@@@@@
%----------------------------------------------------------
\par Let us fix additional notation. Throughout the paper $\gamma,\gamma_1,\gamma_2...,$ and $C,C_1,C_2,...$ will be generic positive constants which depend only on parameters determining function spaces ($m,n,p,q,$ etc). These constants can change even in a single string of estimates. The dependence of a constant on certain parameters is expressed, for example, by the notation $C=C(n,p,q)$. We write $A\sim B$ if there is a constant $C\ge 1$ such that $A/C\le B\le CA$.
%----------------------------------------------------------
%@@@@@@@@@@@@@@@@@@@@@@@@@@@@@@@@@@@@@@@@@@@@@@@@@@@@@@@@@@
%----------------------------------------------------------
\par Given $x=(x_1,x_2,...,x_n)\in\RN$ by $\|x\|:=\max \{|x_1|,|x_2|,...,|x_n|\}$ we denote the uniform norm in $\RN$. Let $A,B\subset \RN$. We put
$\diam A:=\sup\{\|a-a'\|:~a,a'\in A\}$ and
%----------------------------------------------------------
$$
\dist(A,B):=\inf\{\|a-b\|:~a\in A, b\in B\}.
$$
%----------------------------------------------------------
For $x\in \RN$ we also set $\dist(x,A):=\dist(\{x\},A)$.
We put  $\dist(A,B)=+\infty$ and $\dist({x},B)=+\infty$ whenever $B=\emp.$ For each pair of points $z_1$ and $z_2$ in $\RN$ we let $(z_1,z_2)$ denote the open line segment joining them. Given a cube $Q$ in $\RN$ by $c_Q$ we denote its center, and by $r_Q$ a half of its side length. (Thus $Q=Q(c_Q,r_Q).$)
%---------------------------------------------------------
\par Finally, given a function $g\in L_{1,loc}(\RN)$ we let $\Mc[g]$ denote its Hardy-Littlewood maximal function:
%----------------------------------------------------------
$$
\Mc[g](x):=\sup_{Q\ni x}\frac{1}{|Q|}\intl_Q|g|dx,~~~~x\in\RN.
$$
%----------------------------------------------------------
Here the supremum is taken over all cubes $Q$ in $\RN$ containing $x$.
%----------------------------------------------------------
\bigskip
%----------------------------------------------------------
%@@@@@@@@@@@@@@@@@@@@@@@@@@@@@@@@@@@@@@@@@@@@@@@@@@@@@@@@@@
%@@@@@@@@@@@@@@@@@@@@@@@@@@@@@@@@@@@@@@@@@@@@@@@@@@@@@@@@@@
%@@@@@@@@@@@@@@@@@@@@@@@@@@@@@@@@@@@@@@@@@@@@@@@@@@@@@@@@@@
%@@@@@@@@@@@@@@@@@@@@@@@@@@@@@@@@@@@@@@@@@@@@@@@@@@@@@@@@@@
%@@@@@@@@@@@@@@@@@@@@@@@@@@@@@@@@@@@@@@@@@@@@@@@@@@@@@@@@@@
%----------------------------------------------------------
%@@@@@@@@@@@@@@@@@@@@@@@@@@@@@@@@@@@@@@@@@@@@@@@@@@@@@@@@@@
%@@@@@@@@@@@@@@@@@@@@@@@@@@@@@@@@@@@@@@@@@@@@@@@@@@@@@@@@@@
%@@@@@@@@@@@@@@@@@@@@@@@@@@@@@@@@@@@@@@@@@@@@@@@@@@@@@@@@@@
%@@@@@@@@@@@@@@@@@@@@@@@@@@@@@@@@@@@@@@@@@@@@@@@@@@@@@@@@@@
%@@@@@@@@@@@@@@@@@@@@@@@@@@@@@@@@@@@@@@@@@@@@@@@@@@@@@@@@@@
%----------------------------------------------------------
\par {\bf 2.1. The geodesic distance $d_q(h)$: a proof of Theorem \reff{MAIN-MT}.}\medskip
%----------------------------------------------------------
\addtocontents{toc}{~~~~2.1. The geodesic distance $d_q(h)$: a proof of Theorem \reff{MAIN-MT}. \hfill \thepage\par}
%----------------------------------------------------------
\indent\par Clearly, by definition \rf{GD-M}, $d_q(x,y:h)\le \dq(x,y:h)$. Prove that
%----------------------------------------------------------
\bel{G-HF}
\dq(x,y:h)\le 16\, d_q(x,y:h).
\ee
%----------------------------------------------------------
\par Let $w$ be a weight on $\RN$, i.e., a non-negative locally integrable function. Let $q>0$ and let $\phq(w):\RN\times\RN\to\R_+$ be a function defined by the formula
%----------------------------------------------------------
\bel{DF-PHI}
\phq(x,y:w):=\|x-y\| \sup_{Q\ni x,y} \left(\frac{1}{|Q|}\intl_{Q} w(u)\,du\right)^{\frac{1}{q}},~~~x,y\in\RN.
\ee
%---------------------------------------------------------
\begin{proposition} Let $n\le q<\infty$. Then for every  $x,y\in \RN$ and every finite family of points
$x_0=x,x_1,...x_{m-1},x_m=y$ in $\RN$ the following inequality
%---------------------------------------------------------
\bel{3.4'}
\phq(x,y:w)\leq 16\, \smed^{m-1}_{i=0}\,\phq(x_i,x_{i+1}:w)
\ee
%---------------------------------------------------------
holds.
\end{proposition}
%----------------------------------------------------------
\par \textit{Proof.}  Let $K$ be an arbitrary cube in $\RN$ and let $x,y\in K$. Prove that
%----------------------------------------------------------
\bel{K-I}
\|x-y\| \left(\frac{1}{|K|}\intl_{K} w(u)\,du\right)^{\frac{1}{q}}\leq 16\,\smed^{m-1}_{i=0}\,\phq(x_i,x_{i+1}:w).
\ee
%---------------------------------------------------------
Let
%----------------------------------------------------------
$$
\tQ=Q(x,2\|x-y\|)=2\,Q_{xy}.
$$
%----------------------------------------------------------
\par Consider two cases.
%----------------------------------------------------------
\par \textit{The first case.} Suppose that {\it there exists $j\in\{0,...,m-1\}$ such that $x_j\in 2\,Q_{xy}$, but}
%----------------------------------------------------------
$$
x_{j+1}\notin 4\,Q_{xy}.
$$
%----------------------------------------------------------
Hence $\|x-x_j\|\le 2\|x-y\|$ and $\|x-x_{j+1}\|\ge 4\|x-y\|$ so that
%---------------------------------------------------------
\bel{3.6}
\|x_i-x_{j+1}\|\ge 2\|x-y\|.
\ee
%---------------------------------------------------------
Since $x,y\in K$, we have  $\diam K\geq \|x-y\|$ which easily implies the inclusion
%---------------------------------------------------------
\bel{QXY-K}
2\,Q_{xy}\subset 5K.
\ee
%----------------------------------------------------------
\par Hence $x_j\in 2\,Q_{xy}\subset 5K.$\medskip
%----------------------------------------------------------
\par Consider the following two subcases. First let us assume that
%----------------------------------------------------------
$$
x_{j+1}\in 8K.
$$
%----------------------------------------------------------
Since $x_j\in 2\,Q_{xy}\subset 5K$, we have $x_j,x_{j+1}\in 8K$. Hence, by \rf{3.6} and definition \rf{DF-PHI} of $\phq$, we obtain
%----------------------------------------------------------
\be
\|x-y\|\left(\frac{1}{|K|}\intl_Kw(u)du\right) ^{\frac{1}{q}}&\leq& \|x_j-x_{j+1}\|\left(\frac{1}{|K|}\intl_Kw(u)du\right) ^{\frac{1}{q}}\nn\\
&\leq& 8^{\frac{n}{q}}\|x_j-x_{j+1}\|\left(\frac{1}{|8K|} \intl_{8K}w(u)du\right)^{\frac{1}{q}}\nn\\&\leq& 8^{\frac{n}{q}}\phq(x_j,x_{j+1}:w).\nn
\ee
%----------------------------------------------------------
Since $n\le q$, we have
%----------------------------------------------------------
$$
\|x-y\|\left(\frac{1}{|K|}\intl_Kw(u)du\right) ^{\frac{1}{q}}\le 8\phq(x_j,x_{j+1}:w)
$$
%----------------------------------------------------------
proving \rf{K-I} in the case under consideration.
%----------------------------------------------------------
\par Now consider the second subcase where
%----------------------------------------------------------
$$
x_{j+1}\notin 8K.
$$
%----------------------------------------------------------
Since $x_j\in 5K$, we conclude that
$\|x_j-x_{j+1}\|\ge 3r_K.$
%----------------------------------------------------------
\par Let $\overline{Q}:=Q(x_j,2\|x_j-x_{j+1}\|).$ Since $x_j\in 5K$, for every $z\in K$ we have
%----------------------------------------------------------
$$
\|x_j-z\|\le \|x_j-c_K\|+\|c_K-z\|\le 5r_K+r_K=6r_K \le 2\|x_j-x_{j+1}\|=r_{\overline{Q}}
$$
%----------------------------------------------------------
proving that $\overline{Q}\supset K$. Clearly,
$r_K\le r_{\overline{Q}}$\, and $\overline{Q}\ni x_j,x_{j+1}$.
%----------------------------------------------------------
\par Hence
%----------------------------------------------------------
\be
\|x-y\|\left(\frac{1}{|K|}\intl_Kw(u)\,du\right)
^{\frac{1}{q}}&=&
2^{\frac{n}{q}}\|x-y\|r_K^{-\frac{n}{q}}
\left(\intl_Kw(u)\,du\right)
^{\frac{1}{q}}\nn\\
&\le& 2^{\frac{n}{q}}\|x-y\|r_K^{-\frac{n}{q}}
\left(\intl_{\overline{Q}} w(u)\,du\right)
^{\frac{1}{q}}.\nn
\ee
%----------------------------------------------------------
Since $x,y\in K$, we have $\|x-y\|\le \diam K=2r_K$. Combining this with inequality $r_K\le r_{\overline{Q}}$\,, we obtain
%----------------------------------------------------------
\be
\|x-y\|\left(\frac{1}{|K|}\intl_K w(u)\,du\right)
^{\frac{1}{q}}&\le&
2^{\frac{n}{q}+1}r_K^{1-\frac{n}{q}}
\left(\intl_{\overline{Q}}
w(u)\,du\right)^{\frac{1}{q}}
\le 2^{\frac{n}{q}+1}r_{\overline{Q}}^{1-\frac{n}{q}}
\left(\intl_{\overline{Q}}
w(u)\,du\right)^{\frac{1}{q}}\nn\\
&=&2^{\frac{2n}{q}+1}r_{\overline{Q}}
\left(\frac{1}{|\overline{Q}|}
\intl_{\overline{Q}}w(u)\,du\right)
^{\frac{1}{q}}.\nn
\ee
%----------------------------------------------------------
Since $r_{\overline{Q}}=2\,\|x_j-x_{j+1}\|$ and $\frac{n}{q}\le 1$, we have
%----------------------------------------------------------
$$
\|x-y\|\left(\frac{1}{|K|}
\intl_{K}w(u)\,du\right)^{\frac{1}{q}}\le 2^4\|x_j-x_{j+1}\|\left(\frac{1}{|\overline{Q}|}
\intl_{\overline{Q}}w(u)\,du\right)
^{\frac{1}{q}}.
$$
%----------------------------------------------------------
But $x_j,x_{j+1}\in \overline{Q}$ so that, by definition \rf{3.4'},
%----------------------------------------------------------
$$
\|x-y\|\left(\frac{1}{|K|}\intl_{K}w(u)\,du\right)
^{\frac{1}{q}}\le 2^4\phq(x_j,x_{j+1}:w)
$$
%----------------------------------------------------------
proving \rf{K-I}.\medskip
%----------------------------------------------------------
%@@@@@@@@@@@@@@@@@@@@@@@@@@@@@@@@@@@@@@@@@@@@@@@@@@@@@@@@@@
%@@@@@@@@@@@@@@@@@@@@@@@@@@@@@@@@@@@@@@@@@@@@@@@@@@@@@@@@@@
%@@@@@@@@@@@@@@@@@@@@@@@@@@@@@@@@@@@@@@@@@@@@@@@@@@@@@@@@@@
%@@@@@@@@@@@@@@@@@@@@@@@@@@@@@@@@@@@@@@@@@@@@@@@@@@@@@@@@@@
%@@@@@@@@@@@@@@@@@@@@@@@@@@@@@@@@@@@@@@@@@@@@@@@@@@@@@@@@@@
%@@@@@@@@@@@@@@@@@@@@@@@@@@@@@@@@@@@@@@@@@@@@@@@@@@@@@@@@@@
%----------------------------------------------------------
\par \textit{The second case}. Suppose that the assumption of the first case is not satisfied, i.e., {\it for each $i\in\{0,...,m-1\}$ such that $x_i\in 2\, Q_{xy}$ we have  $x_{i+1}\in 4\, Q_{xy}$}.
%----------------------------------------------------------
\par Let us define a number $j\in \{0,1,...,m-1\}$ as follows. If
%----------------------------------------------------------
$$\{x_0,x_1,...,x_m\}\subset 2\,Q_{xy},$$
%----------------------------------------------------------
we put $j=m-1$. If
%----------------------------------------------------------
$$\{x_0,x_1,...,x_m\}\nsubseteq 2\,Q_{xy},$$
%----------------------------------------------------------
then there exists $j\in {0,1,...,m-1}$, such that
%----------------------------------------------------------
$$
\{x_0,x_1,...,x_j\}\subset 2\,Q_{xy}~~~\text{but}~~~
x_{j+1}\notin 2\,Q_{xy}.
$$
%----------------------------------------------------------
%@@@@@@@@@@@@@@@@@@@@@@@@@@@@@@@@@@@@@@@@@@@@@@@@@@@@@@@@@@
%----------------------------------------------------------
Note that, by the assumption, $x_{j+1}\in 4\,Q_{xy}$.
%----------------------------------------------------------
\par Prove that
%----------------------------------------------------------
$$
\|x-y\|\left(\frac{1}{|K|}\intl_{K}w(u)\,du\right)
^{\frac{1}{q}}\le 16\,\smed_{i=0}^j\,\phq(x_i,x_{i+1}:w).
$$
%----------------------------------------------------------
In fact, since $x_{j+1}\notin 2\,Q_{xy}=Q(x,2\|x-y\|)$ we have
%----------------------------------------------------------
$$
\|x_0-x_{j+1}\|=\|x-x_{j+1}\|\ge 2\|x-y\|
$$
%----------------------------------------------------------
so that
%----------------------------------------------------------
\bel{3.11}
\smed_{i=0}^j\,\|x_i-x_{i+1}\|\ge \|x_0-x_{j+1}\|\ge 2\,\|x-y\|.
\ee
%---------------------------------------------------------
\par Recall that, by \rf{QXY-K}, $2\,Q_{xy}\subset 5K$ so that $10K\supset 4\,Q_{xy}$. Hence
%----------------------------------------------------------
\bel{3.12}
x_0,x_1,...,x_j,x_{j+1}\in 10K.
\ee
%----------------------------------------------------------
By \rf{3.11},
%----------------------------------------------------------
\be
\|x-y\|\left(\frac{1}{|K|}
\intl_{K}w(u)\,du\right)^{\frac{1}{q}}&\le& \frac{1}{2}\left(\smed_{i=0}^j\,\|x_i-x_{i+1}\|\right)
\left(\frac{1}{|K|}
\intl_{K}w(u)\,du\right)^{\frac{1}{q}}\nn\\
&\le&\frac{(10)^{\frac{n}{q}}}{2}
\left(\smed_{i=0}^j\,\|x_i-x_{i+1}\|\right)
\left(\frac{1}{|10K|}
\intl_{10K}w(u)\,du\right)^{\frac{1}{q}}\nn.
\ee
%----------------------------------------------------------
Since $10^{\frac{n}{q}}\le 10$, we obtain
%----------------------------------------------------------
$$
\|x-y\|\left(\frac{1}{|K|}\intl_{K}w(u)\,du\right)
^{\frac{1}{q}}\le
5\left(\smed_{i=0}^j\,\|x_i-x_{i+1}\|\right)
\left(\frac{1}{|10K|}\intl_{10K}w(u)\,du\right)^{\frac{1}{q}}.
$$
%----------------------------------------------------------
But, by \rf{3.12}, $x_i,x_{i+1}\in 10K$ for every $0\le i\le j$, so that
%----------------------------------------------------------
$$
\|x_i-x_{i+1}\|\left(\frac{1}{|10K|}
\intl_{10K}w(u)\,du\right)^{\frac{1}{q}}\le
\phq(x_{i},x_{i+1}:w).
$$
%----------------------------------------------------------
Hence
%----------------------------------------------------------
$$
\|x-y\|\left(\frac{1}{|K|}
\intl_{K}w(u)\,du\right)^{\frac{1}{q}}\le
5\,\smed_{i=0}^j\,\phq(x_{i},x_{i+1}:w)\le
5\,\smed_{i=0}^{m-1}\,\phq(x_{i},x_{i+1}:w)
$$
%----------------------------------------------------------
proving inequality \rf{K-I}.\smallskip
%----------------------------------------------------------
\par Thus \rf{K-I} is proven for an arbitrary family of points $x_0=x,x_1,...,x_m=y$ in $\RN$. Taking the supremum in this inequality over all cubes $K\ni x,y$, we finally obtain the statement of the proposition.\bx
%----------------------------------------------------------
%@@@@@@@@@@@@@@@@@@@@@@@@@@@@@@@@@@@@@@@@@@@@@@@@@@@@@@@@@@
%@@@@@@@@@@@@@@@@@@@@@@@@@@@@@@@@@@@@@@@@@@@@@@@@@@@@@@@@@@
%@@@@@@@@@@@@@@@@@@@@@@@@@@@@@@@@@@@@@@@@@@@@@@@@@@@@@@@@@@
%@@@@@@@@@@@@@@@@@@@@@@@@@@@@@@@@@@@@@@@@@@@@@@@@@@@@@@@@@@
%@@@@@@@@@@@@@@@@@@@@@@@@@@@@@@@@@@@@@@@@@@@@@@@@@@@@@@@@@@
%@@@@@@@@@@@@@@@@@@@@@@@@@@@@@@@@@@@@@@@@@@@@@@@@@@@@@@@@@@
%@@@@@@@@@@@@@@@@@@@@@@@@@@@@@@@@@@@@@@@@@@@@@@@@@@@@@@@@@@
%@@@@@@@@@@@@@@@@@@@@@@@@@@@@@@@@@@@@@@@@@@@@@@@@@@@@@@@@@@
%----------------------------------------------------------
\par Clearly, applying this proposition to $w:=h^q$ we obtain the required inequality \rf{G-HF}. The proof of Theorem \reff{MAIN-MT} is complete.\bx
%----------------------------------------------------------
\bigskip
%----------------------------------------------------------
%@@@@@@@@@@@@@@@@@@@@@@@@@@@@@@@@@@@@@@@@@@@@@@@@@@@@@@@@@@
%@@@@@@@@@@@@@@@@@@@@@@@@@@@@@@@@@@@@@@@@@@@@@@@@@@@@@@@@@@
%@@@@@@@@@@@@@@@@@@@@@@@@@@@@@@@@@@@@@@@@@@@@@@@@@@@@@@@@@@
%@@@@@@@@@@@@@@@@@@@@@@@@@@@@@@@@@@@@@@@@@@@@@@@@@@@@@@@@@@
%@@@@@@@@@@@@@@@@@@@@@@@@@@@@@@@@@@@@@@@@@@@@@@@@@@@@@@@@@@
%@@@@@@@@@@@@@@@@@@@@@@@@@@@@@@@@@@@@@@@@@@@@@@@@@@@@@@@@@@
%@@@@@@@@@@@@@@@@@@@@@@@@@@@@@@@@@@@@@@@@@@@@@@@@@@@@@@@@@@
%@@@@@@@@@@@@@@@@@@@@@@@@@@@@@@@@@@@@@@@@@@@@@@@@@@@@@@@@@@
%----------------------------------------------------------
\par {\bf 2.2. Variable metric transforms.}\medskip
%----------------------------------------------------------
\addtocontents{toc}{~~~~2.2. Variable metric transforms. \hfill \thepage\par\VST}
%----------------------------------------------------------
\indent\par In this subsection we present several important properties of metrics from the family $\DPQ$. See \rf{D-LIP}.
%----------------------------------------------------------
%@@@@@@@@@@@@@@@@@@@@@@@@@@@@@@@@@@@@@@@@@@@@@@@@@@@@@@@@@@
%@@@@@@@@@@@@@@@@@@@@@@@@@@@@@@@@@@@@@@@@@@@@@@@@@@@@@@@@@@
%@@@@@@@@@@@@@@@@@@@@@@@@@@@@@@@@@@@@@@@@@@@@@@@@@@@@@@@@@@
%@@@@@@@@@@@@@@@@@@@@@@@@@@@@@@@@@@@@@@@@@@@@@@@@@@@@@@@@@@
%@@@@@@@@@@@@@@@@@@@@@@@@@@@@@@@@@@@@@@@@@@@@@@@@@@@@@@@@@@
%----------------------------------------------------------
\begin{proposition}\lbl{VMT} Let $n\le q\le p<\infty$ and let a metric $d\in\DPQ$. Then:\medskip
%----------------------------------------------------------
\par (a). There exists a mapping $\RN\in x\to \omega_x\in \MC$ such that
%---------------------------------------------------------
\bel{D-VM}
\frac12\, d(x,y)\le \omega_x(\|x-y\|)\le 32\,d(x,y)~~\text{for every}~x,y\in\RN;
\ee
%----------------------------------------------------------
%@@@@@@@@@@@@@@@@@@@@@@@@@@@@@@@@@@@@@@@@@@@@@@@@@@@@@@@@@@
%@@@@@@@@@@@@@@@@@@@@@@@@@@@@@@@@@@@@@@@@@@@@@@@@@@@@@@@@@@
%@@@@@@@@@@@@@@@@@@@@@@@@@@@@@@@@@@@@@@@@@@@@@@@@@@@@@@@@@@
%@@@@@@@@@@@@@@@@@@@@@@@@@@@@@@@@@@@@@@@@@@@@@@@@@@@@@@@@@@
%@@@@@@@@@@@@@@@@@@@@@@@@@@@@@@@@@@@@@@@@@@@@@@@@@@@@@@@@@@
%----------------------------------------------------------
\par (b). Let $x,y,z\in\RN$ and let $\lambda\ge 1$ be a constant such that $\|y-z\|\le\lambda\|x-z\|$. Then
%----------------------------------------------------------
\bel{O-1}
d(y,z)\le 32\,\lambda\,d(x,z)
\ee
%----------------------------------------------------------
and
%----------------------------------------------------------
\bel{O-2}
\frac{d(x,z)}{\|x-z\|}\le 32\lambda\,\frac{d(y,z)}{\|y-z\|};
\ee
%----------------------------------------------------------
\medskip
%----------------------------------------------------------
%@@@@@@@@@@@@@@@@@@@@@@@@@@@@@@@@@@@@@@@@@@@@@@@@@@@@@@@@@@
%@@@@@@@@@@@@@@@@@@@@@@@@@@@@@@@@@@@@@@@@@@@@@@@@@@@@@@@@@@
%@@@@@@@@@@@@@@@@@@@@@@@@@@@@@@@@@@@@@@@@@@@@@@@@@@@@@@@@@@
%@@@@@@@@@@@@@@@@@@@@@@@@@@@@@@@@@@@@@@@@@@@@@@@@@@@@@@@@@@
%@@@@@@@@@@@@@@@@@@@@@@@@@@@@@@@@@@@@@@@@@@@@@@@@@@@@@@@@@@
%----------------------------------------------------------
\par (c). For every $x,y\in\RN$ and $z\in(x,y)$ the following inequality
%---------------------------------------------------------
$$
d(x,z)+d(y,z)\le 64\, d(x,y)
$$
%---------------------------------------------------------
holds.
%---------------------------------------------------------
\end{proposition}
%---------------------------------------------------------
\par {\it Proof.} $(a).$ Since $d\in\DPQ$, there exists a non-negative function $h\in\LPRN$  such that $d=d_q(h)$.
Given $x\in\RN$ we let $v_x$ denote a function on $\R_+$ defined by the formula
%---------------------------------------------------------
\bel{VX}
v_x(t):=t\sup_{s\ge t}\left(\frac{1}{|Q(x,s)|}\intl_{Q(x,s)}
h^q(u)\,du\right)^{\frac{1}{q}},~~~~t\ge 0.
\ee
%---------------------------------------------------------
\par Prove that $d(x,y)\sim v_x(\|x-y\|)$. In fact, since $y\in Q(x,s)$ whenever $s\ge \|x-y\|$, by \rf{DQ}, $v_x(\|x-y\|)\le \dq(x,y:h)$.
%---------------------------------------------------------
\par On the other hand, for every cube $Q=Q(a,r)\ni x,y$ we have $\diam Q=2r\ge \|x-y\|$. Furthermore, $Q(a,r)\subset \tQ=Q(x,2r)$ and $|\tQ|=2^n|Q|$. Hence
%---------------------------------------------------------
$$
\|x-y\|\left(\,\frac{1}{|Q|}\intl_{Q}
h^q(u)\,du\right)^{\frac{1}{q}}\le 2^{\frac{n}{q}}
\|x-y\|\left(\,\frac{1}{|\tQ|}\intl_{\tQ}
h^q(u)\,du\right)^{\frac{1}{q}}\le 2 v_x(\|x-y\|)
$$
%---------------------------------------------------------
proving that $\dq(x,y:h)\le 2v_x(\|x-y\|)$.
%---------------------------------------------------------
\par Combining these inequalities with inequality \rf{W-G} of Theorem \reff{MAIN-MT}, we obtain:
%---------------------------------------------------------
\bel{D-VX}
\frac12\, d(x,y)\le v_x(\|x-y\|)\le 16\,d(x,y).
\ee
%----------------------------------------------------------
\par Prove that $v_x$ is equivalent to a ``modulus of  continuity'' $\omega_x\in\MC$. Clearly, by \rf{VX}, the function $v_x(t)/t$ is non-increasing. Let us show that $v_x$ is a non-decreasing function.
%----------------------------------------------------------
\par Let $0\le t_1\le t_2$. Then
$v_x(t_1)=\max\{I_1,I_2\}$
%---------------------------------------------------------
where
%---------------------------------------------------------
$$
I_1:=
t_1\sup_{s\ge t_2}\left(\frac{1}{|Q(x,s)|}\intl_{Q(x,s)}
h^q(u)\,du\right)^{\frac{1}{q}}~~\text{and}~~
I_2:=t_1\sup_{t_1\le s<t_2} \left(\frac{1}{|Q(x,s)|}\intl_{Q(x,s)}
h^q(u)\,du\right)^{\frac{1}{q}}.
$$
%---------------------------------------------------------
Since $t_1\le t_2$, by \rf{VX}, $I_1\le v_x(t_2)$. On the other hand,
%---------------------------------------------------------
$$
I_2\le t_1\left(\frac{1}{|Q(x,t_1)|}\intl_{Q(x,t_2)}
h^q(u)\,du\right)^{\frac{1}{q}}=
2^{-\frac{n}{q}}t_1^{1-\frac{n}{q}}\left(\intl_{Q(x,t_2)}
h^q(u)\,du\right)^{\frac{1}{q}}.
$$
%---------------------------------------------------------
Hence
%---------------------------------------------------------
$$
I_2\le
2^{-\frac{n}{q}}t_2^{1-\frac{n}{q}}\left(\,\intl_{Q(x,t_2)}
h^q(u)\,du\right)^{\frac{1}{q}}=
t_2\left(\frac{1}{|Q(x,t_2)|}\intl_{Q(x,t_2)}
h^q(u)\,du\right)^{\frac{1}{q}}\le v_x(t_2)
$$
%---------------------------------------------------------
so that $$v_x(t_1)=\max\{I_1,I_2\}\le v_x(t_2).$$
%---------------------------------------------------------
\par Let $\omega_x$ be the least concave majorant of $v_x$. Since $v_x$ and $t/v_x(t)$ are non-decreasing, for every $t>0$ the following inequality
%---------------------------------------------------------
$$
v_x(t)\le \omega_x(t)\le 2 v_x(t)
$$
%---------------------------------------------------------
holds. Combining this inequality with \rf{D-VX} we obtain the required inequality \rf{D-VM}.\medskip
%----------------------------------------------------------
%@@@@@@@@@@@@@@@@@@@@@@@@@@@@@@@@@@@@@@@@@@@@@@@@@@@@@@@@@@
%@@@@@@@@@@@@@@@@@@@@@@@@@@@@@@@@@@@@@@@@@@@@@@@@@@@@@@@@@@
%@@@@@@@@@@@@@@@@@@@@@@@@@@@@@@@@@@@@@@@@@@@@@@@@@@@@@@@@@@
%@@@@@@@@@@@@@@@@@@@@@@@@@@@@@@@@@@@@@@@@@@@@@@@@@@@@@@@@@@
%@@@@@@@@@@@@@@@@@@@@@@@@@@@@@@@@@@@@@@@@@@@@@@@@@@@@@@@@@@
%----------------------------------------------------------
\par $(b).$ In part $(a)$ we have proved that for each $z\in\RN$ the functions $v_z$ is non-decreasing and the function $v_z(t)/t$ is non-increasing. Hence $v_z(\lambda t)\le\lambda v_z(t)$ provided $\lambda\ge 1$. Therefore, by \rf{D-VX},
%---------------------------------------------------------
$$
d(y,z)\le 2v_z(\|y-z\|)\le 2v_z(\lambda\|x-z\|)\le 2\lambda v_z(\|x-z\|)\le 32\,\lambda d(x,z)
$$
%---------------------------------------------------------
proving \rf{O-1}. On the other hand, by the first inequality in \rf{D-VX},
%----------------------------------------------------------
$$
\frac{d(x,z)}{\|x-z\|}\le \frac{2v_z(\|x-z\|)}{\|x-z\|}
\le 2\lambda \frac{v_z(\|y-z\|/\lambda)}{\|y-z\|}\le
2\lambda \frac{v_z(\|y-z\|)}{\|y-z\|}
$$
%----------------------------------------------------------
so that, by the second inequality in \rf{D-VX},
%----------------------------------------------------------
$$
\frac{d(x,z)}{\|x-z\|}\le
32\lambda\,\frac{d(y,z)}{\|y-z\|}
$$
%----------------------------------------------------------
proving \rf{O-2}.
%----------------------------------------------------------
\medskip
%----------------------------------------------------------
%@@@@@@@@@@@@@@@@@@@@@@@@@@@@@@@@@@@@@@@@@@@@@@@@@@@@@@@@@@
%@@@@@@@@@@@@@@@@@@@@@@@@@@@@@@@@@@@@@@@@@@@@@@@@@@@@@@@@@@
%@@@@@@@@@@@@@@@@@@@@@@@@@@@@@@@@@@@@@@@@@@@@@@@@@@@@@@@@@@
%@@@@@@@@@@@@@@@@@@@@@@@@@@@@@@@@@@@@@@@@@@@@@@@@@@@@@@@@@@
%@@@@@@@@@@@@@@@@@@@@@@@@@@@@@@@@@@@@@@@@@@@@@@@@@@@@@@@@@@
%----------------------------------------------------------
\par $(c).$ Since $z\in(x,y)$, we have $\|y-z\|,\|x-z\|\le\|x-y\|$ so that by part $(b)$
%----------------------------------------------------------
$$
d(x,z)\le 32\,d(x,y)~~~\text{and}~~~d(y,z)\le 32\,d(x,y)
$$
%----------------------------------------------------------
proving the statement $(c)$ and the proposition.\bx
%----------------------------------------------------------
%@@@@@@@@@@@@@@@@@@@@@@@@@@@@@@@@@@@@@@@@@@@@@@@@@@@@@@@@@@
%@@@@@@@@@@@@@@@@@@@@@@@@@@@@@@@@@@@@@@@@@@@@@@@@@@@@@@@@@@
%@@@@@@@@@@@@@@@@@@@@@@@@@@@@@@@@@@@@@@@@@@@@@@@@@@@@@@@@@@
%@@@@@@@@@@@@@@@@@@@@@@@@@@@@@@@@@@@@@@@@@@@@@@@@@@@@@@@@@@
%@@@@@@@@@@@@@@@@@@@@@@@@@@@@@@@@@@@@@@@@@@@@@@@@@@@@@@@@@@
%----------------------------------------------------------
\begin{remark}\lbl{MT-RM} {\em Equivalence \rf{D-VM} motivates us to refer to the metric space $(\RN,d)$ where $d\in\DPQ$ as a {\it variable metric transform} of $\RN$; see Remark \reff{VM-E}. This equivalence shows that given $x\in\RN$ the local behavior of the metric $d$ is similar to the behaviour of a certain regular metric transform $d_{\omega_x}:=\omega_x(\|\cdot\|)$ where $\omega_x\in\MC$ is a ``modulus of continuity''. The function $\omega_x$ varies from point to point, and this is the main difference between a regular metric transform (where $\omega_x$ is ``constant'', i.e., the same ``modulus of continuity'' $\omega$ for all $x\in\RN$) and a variable metric transform.
%---------------------------------------------------------
\par Nevertheless, in spite of $\omega_x$ changes, the metric $d\in\DPQ$ preserves several important pro\-perties of regular metric transforms.
%---------------------------------------------------------
\par For instance, let $E\subset \RN$ be a closed set and let $x\in \RN\setminus E$. Let $\tx\in E$ be an almost nearest point to $x$ on $E$ with respect to the
Euclidean distance, i.e., $\|x-\tx\|\sim \dist(x,E)$. Then $\tx$ is an almost nearest to $x$ point with respect to the variable majorant $d$ as well.
%----------------------------------------------------------
\par Another example is the standard Whitney covering of $\RN\setminus E$ by a family of Whitney's cubes. See, e.g. \cite{St}. It is well known that this covering is universal with respect to the family
$$
\Mc\Tc=\{(\RN,d_\omega),~~\omega\in\MC\}
$$
of all metric transforms, i.e., it provides an almost optimal Whitney type extension construction for the family of Lipschitz spaces with respect to metric transforms. As we shall see below, the same property also holds for variable metric transforms.
%---------------------------------------------------------
\par Thus there exists a more or less complete analogy between extension methods for regular and variable metric transforms. In the next sections we present several applications of this approach to extensions of jets generated by Sobolev functions.\rbx}
%---------------------------------------------------------
\end{remark}
%----------------------------------------------------------
%@@@@@@@@@@@@@@@@@@@@@@@@@@@@@@@@@@@@@@@@@@@@@@@@@@@@@@@@@@
%@@@@@@@@@@@@@@@@@@@@@@@@@@@@@@@@@@@@@@@@@@@@@@@@@@@@@@@@@@
%@@@@@@@@@@@@@@@@@@@@@@@@@@@@@@@@@@@@@@@@@@@@@@@@@@@@@@@@@@
%@@@@@@@@@@@@@@@@@@@@@@@@@      @@@@@@@@@@@@@@@@@@@@@@@@@@@
%@@@@@@@@@@@@@@@@@@@@@@@          @@@@@@@@@@@@@@@@@@@@@@@@@
%@@@@@@@@@@@@@@@@@@@@@              @@@@@@@@@@@@@@@@@@@@@@@
%@@@@@@@@@@@@@@@@@@@     SECTION 3    @@@@@@@@@@@@@@@@@@@@@
%@@@@@@@@@@@@@@@@@@@@@              @@@@@@@@@@@@@@@@@@@@@@@
%@@@@@@@@@@@@@@@@@@@@@@@          @@@@@@@@@@@@@@@@@@@@@@@@@
%@@@@@@@@@@@@@@@@@@@@@@@@@      @@@@@@@@@@@@@@@@@@@@@@@@@@@
%@@@@@@@@@@@@@@@@@@@@@@@@@@@@@@@@@@@@@@@@@@@@@@@@@@@@@@@@@@
%@@@@@@@@@@@@@@@@@@@@@@@@@@@@@@@@@@@@@@@@@@@@@@@@@@@@@@@@@@
%@@@@@@@@@@@@@@@@@@@@@@@@@@@@@@@@@@@@@@@@@@@@@@@@@@@@@@@@@@
%----------------------------------------------------------
\SECT{3. Sobolev $L^{m}_p$-space as a union of $C^{m-1,(d)}$-spaces.}{3}
%----------------------------------------------------------
\addtocontents{toc}{3. Sobolev $L^{m}_p$-space as a union of $C^{m-1,(d)}$-spaces. \hfill\thepage\par\VST}
%----------------------------------------------------------
\indent\par In this subsection we prove Theorem \reff{SM-MAIN}.\medskip
%---------------------------------------------------------- %@@@@@@@@@@@@@@@@@@@@@@@@@@@@@@@@@@@@@@@@@@@@@@@@@@@@@@@@@@
%----------------------------------------------------------
\par {\it (Necessity).} The necessity part directly follows from inequality \rf{F-SP} and inequality \rf{W-G} with  $h=C(n,p,q)\,\|\nabla^m F\|$.\medskip
%---------------------------------------------------------- %@@@@@@@@@@@@@@@@@@@@@@@@@@@@@@@@@@@@@@@@@@@@@@@@@@@@@@@@@@
%----------------------------------------------------------
\par {\it (Sufficiency).} Let $F$ be a $C^{m-1}$-function. Suppose there exists a non-negative function $h\in\LPRN$ such that inequality \rf{LIP2} holds for every multiindex $\alpha, |\alpha|=m-1$. Prove that $F\in\LMP$ and $\|F\|_{\LMP}\le C(n,p,q)\|h\|_{\LPRN}$.
%----------------------------------------------------------
\par Our proof relies on a result of Calder\'{o}n \cite{C1} which characterizes Sobolev functions in terms of certain sharp maximal functions. See also \cite{CS}. Let us recall this characterization of Sobolev spaces.
%----------------------------------------------------------
\par Given a cube $Q\subset\RN$ and a function $f\in
L_{q}(Q)$, $0<q\le\infty,$ we let $\Ec_m(f;Q)_{L_q}$ denote the {\it normalized local best approximation} of $f$ on $Q$ in $L_q$-norm by polynomials of degree at most $m-1$. More explicitly, we define
%----------------------------------------------------------
$$
\Ec_m(f;Q)_{L_q}:=\inf_{P\in\PMRN}\left(\frac{1}{|Q|}
\intl_Q|f-P|^q\,dx\right)^{\frac{1}{q}}.
$$
%----------------------------------------------------------
%@@@@@@@@@@@@@@@@@@@@@@@@@@@@@@@@@@@@@@@@@@@@@@@@@@@@@@@@@@
%@@@@@@@@@@@@@@@@@@@@@@@@@@@@@@@@@@@@@@@@@@@@@@@@@@@@@@@@@@
%@@@@@@@@@@@@@@@@@@@@@@@@@@@@@@@@@@@@@@@@@@@@@@@@@@@@@@@@@@
%----------------------------------------------------------
\par Given a locally integrable function $f$ on $\RN$, we define its {\it sharp maximal function} $f^{\sharp}_m$ by letting
%----------------------------------------------------------
$$
f^{\sharp}_m(x):=
\sup_{r>0} r^{-m}\,\Ec_{m}(f;Q(x,r))_{L_1}.
$$
%----------------------------------------------------------
\par Calder\'{o}n \cite{C1} proved that a locally integrable function $f\in\LMP$, $1<p<\infty$, if and only if  $f^{\sharp}_m$ is in $\LPRN$. Furthermore, the following equivalence
%----------------------------------------------------------
\bel{AN}
\|f\|_{\LMP}\sim \|f_{m}^{\sharp}\|_{\LPRN}
\ee
%----------------------------------------------------------
%@@@@@@@@@@@@@@@@@@@@@@@@@@@@@@@@@@@@@@@@@@@@@@@@@@@@@@@@@@
holds with constants depending only on $n,m$ and $p$.\smallskip
%----------------------------------------------------------
\par Let us show that $F^{\sharp}_m\in\LPRN$. Our proof of this fact relies on series of auxiliary lemmas. To formulate the first of them we introduce the following notion.
%---------------------------------------------------------
%@@@@@@@@@@@@@@@@@@@@@@@@@@@@@@@@@@@@@@@@@@@@@@@@@@@@@@@@@
%@@@@@@@@@@@@@@@@@@@@@@@@@@@@@@@@@@@@@@@@@@@@@@@@@@@@@@@@@
%@@@@@@@@@@@@@@@@@@@@@@@@@@@@@@@@@@@@@@@@@@@@@@@@@@@@@@@@@
%@@@@@@@@@@@@@@@@@@@@@@@@@@@@@@@@@@@@@@@@@@@@@@@@@@@@@@@@@
%---------------------------------------------------------
\par We say that a metric $d$ on $\RN$ is \textit{pseudoconvex} if there exists a constant $\ld\ge 1$ such that for every $x,y,z\in \RN$, $z\in (x,y)$, the following inequality
%----------------------------------------------------------
\bel{B-3}
d(x,z)+d(z,y)\le \ld\,d(x,y)
\ee
%---------------------------------------------------------
holds.
%---------------------------------------------------------
%@@@@@@@@@@@@@@@@@@@@@@@@@@@@@@@@@@@@@@@@@@@@@@@@@@@@@@@@@@
%@@@@@@@@@@@@@@@@@@@@@@@@@@@@@@@@@@@@@@@@@@@@@@@@@@@@@@@@@@
%@@@@@@@@@@@@@@@@@@@@@@@@@@@@@@@@@@@@@@@@@@@@@@@@@@@@@@@@@@
%@@@@@@@@@@@@@@@@@@@@@@@@@@@@@@@@@@@@@@@@@@@@@@@@@@@@@@@@@@
%@@@@@@@@@@@@@@@@@@@@@@@@@@@@@@@@@@@@@@@@@@@@@@@@@@@@@@@@@@
%@@@@@@@@@@@@@@@@@@@@@@@@@@@@@@@@@@@@@@@@@@@@@@@@@@@@@@@@@@
%---------------------------------------------------------
\begin{lemma}\lbl{Cl7.1} Let $d$ be a pseudoconvex metric on $\RN$. Then for every $F\in \CMD$ and every multiindex $\beta$ with $|\beta|\le m$ the following inequality
%----------------------------------------------------------
$$
|D^{\beta}F(x)-D^{\beta}\left(T^m_y[F]\right)(x)|\le C\|F\|_{\CMD}\|x-y\|^{m-|\beta|}d(x,y),~~~~x,y\in \RN,
$$
%---------------------------------------------------------
holds. Here $C=C(m,\ld)$.
%----------------------------------------------------------
\end{lemma}
%----------------------------------------------------------
\par {\it Proof.} For $|\beta|=m$ the lemma follows from Definition \reff{CMD}, so we can assume that $|\beta|<m$.
%----------------------------------------------------------
\par  Our proof for this case relies on the following well known identity: Let $m>0$ and let $F$ be a $C^m$-function on $\RN$. Then for every $x,y\in\RN$ the following equality
%----------------------------------------------------------
$$
F(x)=T^m_y[F](x)+m\,\sbig_{|\alpha|=m}\,\frac{1}{\alpha !}(x-y)^{\alpha}\int_0^1(1-t)^{m-1}(D^{\alpha}F(x+t(x-y))-
D^{\alpha}F(y))dt
$$
%---------------------------------------------------------
holds.
%---------------------------------------------------------
\par Let us apply this identity to $D^{\beta}F$. We obtain
%----------------------------------------------------------
\be
D^{\beta}F(x)&=&
T^{m-|\beta|}_y[D^{\beta}F](x)
+(m-|\beta|)\smed_{|\alpha|=m-|\beta|}\,
\frac{1}{\alpha!}(x-y)^{\alpha}\nn\\
&\cdot&
\int^1_0(1-t)^{m-|\beta|-1}
\left\{D^{\alpha}(D^{\beta}F)(x+t(x-y))
-D^{\alpha}(D^{\beta}F)(y)\right\}dt\,.\nn
\ee
%----------------------------------------------------------
Since
%----------------------------------------------------------
$$
T^{m-|\beta|}_y[D^{\beta}F](x)=D^{\beta}
\left(T^m_y[F]\right)(x)~~~\text{for every}~~~x,y\in \RN,
$$
%----------------------------------------------------------
we have
%----------------------------------------------------------
\be
&&D^{\beta}F(x)-D^{\beta}(T^m_y[F])(x)
=(m-|\beta|)\smed_{|\alpha|=m-|\beta|}
\frac{1}{\alpha!}\,(x-y)^{\alpha}\nn\\
&\cdot&\int^1_0(1-t)^{m-|\beta|-1}
\left\{D^{\alpha+\beta}F(x+t(x-y))
-D^{\alpha+\beta}F(y)\right\}dt.\nn
\ee
%---------------------------------------------------------
Hence
%----------------------------------------------------------
$$
|D^{\beta}F(x)-D^{\beta}(T^m_y[F])(x)|\le
m\smed_{|\alpha|+|\beta|=m}
\,\|x-y\|^{|\alpha|}
\sup_{z\in [x,y]}|D^{\alpha+\beta}F)(z)
-D^{\alpha+\beta}F)(y)|.
$$
%---------------------------------------------------------
Combining this inequality with Definition \reff{CMD} we obtain
%----------------------------------------------------------
$$
|D^{\beta}F(x)-D^{\beta}(T^m_y[F])(x)|\le
m\|x-y\|^{m-|\beta|}\|F\|_{\CMD}
\sup_{z\in [x,y]}d(z,y).
$$
%---------------------------------------------------------
Since $d$ is pseudoconvex, by \rf{B-3},
%----------------------------------------------------------
$$
|D^{\beta}F(x)-D^{\beta}(T^m_y[F])(x)|\le
\ld\,m\, \|x-y\|^{m-|\beta|}\,\|F\|_{\CMD}\,d(x,y)
$$
%---------------------------------------------------------
proving the lemma. \bx\medskip
%---------------------------------------------------------
%@@@@@@@@@@@@@@@@@@@@@@@@@@@@@@@@@@@@@@@@@@@@@@@@@@@@@@@@@
%@@@@@@@@@@@@@@@@@@@@@@@@@@@@@@@@@@@@@@@@@@@@@@@@@@@@@@@@@
%@@@@@@@@@@@@@@@@@@@@@@@@@@@@@@@@@@@@@@@@@@@@@@@@@@@@@@@@@
%@@@@@@@@@@@@@@@@@@@@@@@@@@@@@@@@@@@@@@@@@@@@@@@@@@@@@@@@@
%@@@@@@@@@@@@@@@@@@@@@@@@@@@@@@@@@@@@@@@@@@@@@@@@@@@@@@@@@
%@@@@@@@@@@@@@@@@@@@@@@@@@@@@@@@@@@@@@@@@@@@@@@@@@@@@@@@@@
%---------------------------------------------------------
\par Let us apply this result to metrics from the family  $\DPQ$ whenever $q\in(n,p)$, see \rf{D-LIP}. Note that, by part (c) of Proposition \reff{VMT}, every metric $d\in\DPQ$ is pseudoconvex. This property of $d$ and Lemma \reff{Cl7.1} imply the following
%----------------------------------------------------------
%@@@@@@@@@@@@@@@@@@@@@@@@@@@@@@@@@@@@@@@@@@@@@@@@@@@@@@@@@
%@@@@@@@@@@@@@@@@@@@@@@@@@@@@@@@@@@@@@@@@@@@@@@@@@@@@@@@@@
%@@@@@@@@@@@@@@@@@@@@@@@@@@@@@@@@@@@@@@@@@@@@@@@@@@@@@@@@@
%---------------------------------------------------------
\begin{proposition}\lbl{Prop-FT} Let $d\in\DPQ$ where $q\in(n,p)$. Then for every $F\in \CMD$, every $x,y\in \RN$ and every $\beta$, $|\beta|\le m$, the following inequality
%----------------------------------------------------------
$$
|D^{\beta}F(x)-D^{\beta}(T^m_y[F])(x)|\le
C\|F\|_{\CMD}\|x-y\|^{m-|\beta|}d(x,y)
$$
%---------------------------------------------------------
holds. Here $C=C(m,n,q,p)$.
\end{proposition}
%---------------------------------------------------------
%@@@@@@@@@@@@@@@@@@@@@@@@@@@@@@@@@@@@@@@@@@@@@@@@@@@@@@@@@
%@@@@@@@@@@@@@@@@@@@@@@@@@@@@@@@@@@@@@@@@@@@@@@@@@@@@@@@@@
%@@@@@@@@@@@@@@@@@@@@@@@@@@@@@@@@@@@@@@@@@@@@@@@@@@@@@@@@@
%@@@@@@@@@@@@@@@@@@@@@@@@@@@@@@@@@@@@@@@@@@@@@@@@@@@@@@@@@
%@@@@@@@@@@@@@@@@@@@@@@@@@@@@@@@@@@@@@@@@@@@@@@@@@@@@@@@@@
%@@@@@@@@@@@@@@@@@@@@@@@@@@@@@@@@@@@@@@@@@@@@@@@@@@@@@@@@@
%---------------------------------------------------------
\medskip
%---------------------------------------------------------
\par We are in a position to finish the proof of the sufficiency part of Theorem \reff{SM-MAIN}. First recall that, by \rf{LIP2}, for every $\alpha,|\alpha|=m-1,$
%----------------------------------------------------------
$$
|D^\alpha F(x)-D^\alpha F(y)|\le d(x,y),~~~x,y\in\RN,
$$
%----------------------------------------------------------
where $d=d_q(h)$. Thus, by Definition \reff{CMD},
%----------------------------------------------------------
$$
F\in \CMLD~~~\text{and}~~~\|F\|_{\CMLD}\le C(m,n).
$$
%---------------------------------------------------------
Furthermore, since $h\ge 0$ and $h\in\LPRN$, the metric $d\in\DPQ$, see \rf{D-LIP}, so that, by Proposition \reff{Prop-FT}, for every $x,y\in\RN$ the following inequality
%----------------------------------------------------------
$$
|F(x)-(T^m_y[F])(x)|\le C\|x-y\|^{m-1}\,d(x,y)
$$
%---------------------------------------------------------
holds. Here $C=C(m,n,q,p)$.
%---------------------------------------------------------
\par Hence, by Theorem \reff{MAIN-MT} and definition \rf{DQ},
%----------------------------------------------------------
$$
|F(x)-(T^m_y[F])(x)|\le C\|x-y\|^m\,
\sup_{Q\ni x,y} \left(\frac{1}{|Q|}\intl_{Q} h^q(u)\,du\right)^{\frac{1}{q}},~~~x,y\in\RN,
$$
%---------------------------------------------------------
which leads to the following inequality
%----------------------------------------------------------
$$
|F(x)-(T^m_y[F])(x)|\le C\|x-y\|^m\,
\left(\Mc[h^q]\right)^{\frac{1}{q}}(y),~~~y\in\RN.
$$
%----------------------------------------------------------
\par Let $Q=Q(y,r)$ be a cube in $\RN$ centered in $y$. Then, by the latter inequality, for every $y\in\RN$,
%----------------------------------------------------------
$$
\sup_Q|F-T^m_y[F]|\le C\,r^m\,\left(\Mc[h^q]\right)^{\frac{1}{q}}(y),
$$
%---------------------------------------------------------
so that
%----------------------------------------------------------
\be
r^{-m}\Ec_m(F;Q)_{L_1}&=&r^{-m}\,
\inf_{P\in\PMRN}\,\frac{1}{|Q|}
\intl_Q|F-P|\,du
\le
r^{-m}\,\frac{1}{|Q|}
\intl_Q|F-T^m_y[F]|\,du\nn\\
&\le&
C\,\left(\Mc[h^q]\right)^{\frac{1}{q}}(y).\nn
\ee
%---------------------------------------------------------
Hence
%----------------------------------------------------------
$$
F^{\sharp}_m(y)=\sup_{r>0}
r^{-m}\Ec_m(F;Q(y,r))_{L_1}\le
C\left(\Mc[h^q]\right)^{\frac{1}{q}}(y).
$$
%---------------------------------------------------------
Since $p>q$, by the Hardy-Littlewood maximal theorem,
%----------------------------------------------------------
$$
\|F^{\sharp}_m\|_{\LPRN}\le C\left(\,\,\intl_{\RN}\Mc[h^q]^{\frac{p}{q}}\,du
\right)^{\frac1p}
\le C\,\left(\,\,\intl_{\RN}(h^q)^{\frac{p}{q}}\,du
\right)^{\frac1p}= C\,\|h\|_{\LPRN},
$$
%---------------------------------------------------------
so that, by \rf{AN}, $F\in\LMP$ and $\|F\|_{\LMP}\le C\|h\|_{\LPRN}$.
%----------------------------------------------------------
\par The proof of Theorem \reff{SM-MAIN} is complete.\bx
%----------------------------------------------------------
%@@@@@@@@@@@@@@@@@@@@@@@@@@@@@@@@@@@@@@@@@@@@@@@@@@@@@@@@@@
%@@@@@@@@@@@@@@@@@@@@@@@@@@@@@@@@@@@@@@@@@@@@@@@@@@@@@@@@@@
%----------------------------------------------------------
%@@@@@@@@@@@@@@@@@@@@@@@@@@@@@@@@@@@@@@@@@@@@@@@@@@@@@@@@@@
%@@@@@@@@@@@@@@@@@@@@@@@@@@@@@@@@@@@@@@@@@@@@@@@@@@@@@@@@@@
%@@@@@@@@@@@@@@@@@@@@@@@@@@@@@@@@@@@@@@@@@@@@@@@@@@@@@@@@@@
%@@@@@@@@@@@@@@@@@@@@@@@@@      @@@@@@@@@@@@@@@@@@@@@@@@@@@
%@@@@@@@@@@@@@@@@@@@@@@@          @@@@@@@@@@@@@@@@@@@@@@@@@
%@@@@@@@@@@@@@@@@@@@@@              @@@@@@@@@@@@@@@@@@@@@@@
%@@@@@@@@@@@@@@@@@@@     SECTION 4    @@@@@@@@@@@@@@@@@@@@@
%@@@@@@@@@@@@@@@@@@@@@              @@@@@@@@@@@@@@@@@@@@@@@
%@@@@@@@@@@@@@@@@@@@@@@@          @@@@@@@@@@@@@@@@@@@@@@@@@
%@@@@@@@@@@@@@@@@@@@@@@@@@      @@@@@@@@@@@@@@@@@@@@@@@@@@@
%@@@@@@@@@@@@@@@@@@@@@@@@@@@@@@@@@@@@@@@@@@@@@@@@@@@@@@@@@@
%@@@@@@@@@@@@@@@@@@@@@@@@@@@@@@@@@@@@@@@@@@@@@@@@@@@@@@@@@@
%@@@@@@@@@@@@@@@@@@@@@@@@@@@@@@@@@@@@@@@@@@@@@@@@@@@@@@@@@@
%----------------------------------------------------------
\SECT{4. Extensions of jets generated by  $\CMD$-functions.}{4}
%----------------------------------------------------------
\addtocontents{toc}{4. Extensions of jets generated by  $\CMD$-functions. \hfill\thepage\par\VST}
%----------------------------------------------------------
%@@@@@@@@@@@@@@@@@@@@@@@@@@@@@@@@@@@@@@@@@@@@@@@@@@@@@@@@@@
%----------------------------------------------------------
\indent\par In this section we prove Theorem \reff{CMD-E-NEW}.
%---------------------------------------------------------
%@@@@@@@@@@@@@@@@@@@@@@@@@@@@@@@@@@@@@@@@@@@@@@@@@@@@@@@@@@
%@@@@@@@@@@@@@@@@@@@@@@@@@@@@@@@@@@@@@@@@@@@@@@@@@@@@@@@@@@
%@@@@@@@@@@@@@@@@@@@@@@@@@@@@@@@@@@@@@@@@@@@@@@@@@@@@@@@@@@
%@@@@@@@@@@@@@@@@@@@@@@@@@@@@@@@@@@@@@@@@@@@@@@@@@@@@@@@@@@
%----------------------------------------------------------
\par {\it Necessity.} Suppose that there exists a function $F\in\CMD$ which agrees with the Whitney $m$-field $\VP=\{P_x\in\PMP: x\in E\}$, i.e., $T_x^m[F]=P_x$ for every $x\in E$. Then, by Proposition \reff{Prop-FT},
%---------------------------------------------------------
$$
\Lc_{m,d}(\VP):=\sbig_{|\alpha|\le m}\,\,\,
\sup_{x,y\in E,\,x\ne y}\,\, \frac{|D^{\alpha}P_x(x)-D^{\alpha}P_y(x)|}
{\|x-y\|^{m-|\alpha|}\,d(x,y)}\le C\,\|F\|_{\CMD}
$$
%---------------------------------------------------------
where $C=C(m,n,p,q)$. This proves the necessity and the inequality
%---------------------------------------------------------
$$
\Lc_{m,d}(\VP)\le C\,\inf\left\{\|F\|_{\CMD}:F\in \CMD, F~\text{agrees with}~\VP~\text{on}~E\right\}.
$$
%---------------------------------------------------------
%@@@@@@@@@@@@@@@@@@@@@@@@@@@@@@@@@@@@@@@@@@@@@@@@@@@@@@@@@@
%@@@@@@@@@@@@@@@@@@@@@@@@@@@@@@@@@@@@@@@@@@@@@@@@@@@@@@@@@@
%@@@@@@@@@@@@@@@@@@@@@@@@@@@@@@@@@@@@@@@@@@@@@@@@@@@@@@@@@@
%@@@@@@@@@@@@@@@@@@@@@@@@@@@@@@@@@@@@@@@@@@@@@@@@@@@@@@@@@@
%----------------------------------------------------------
\medskip
\par {\it (Sufficiency.)} Let $\VP=\{P_x\in\PMP: x\in E\}$ be a Whitney $m$-field on $E$ and let $\lambda:=\Lc_{m,d}(\VP)$. Suppose that $\lambda<\infty$. Thus for every multiindex $\alpha$, $|\alpha|\le m$, and every $x,y\in E$ the following inequality
%---------------------------------------------------------
\bel{W-P}
|D^{\alpha}P_x(x)-D^{\alpha}P_y(x)|\le \,\|x-y\|^{m-|\alpha|}\,d(x,y)\,\lambda
\ee
%---------------------------------------------------------
holds.
%---------------------------------------------------------
\par Prove the existence of a function $F\in\CMD$ such that $T^m_x[F]=P_x$ for every $x\in E$, and
$\|F\|_{\CMD}\le C\,\lambda$ where $C$ is a positive constant depending only on $m,n,p,$ and $q$.\medskip
%---------------------------------------------------------
%@@@@@@@@@@@@@@@@@@@@@@@@@@@@@@@@@@@@@@@@@@@@@@@@@@@@@@@@@
%---------------------------------------------------------
\par We construct $F$ with the help of a slight modification of the classical Whitney extension method which we present below.
%---------------------------------------------------------
\par Since $E$ is a closed set, the set $\RN\setminus E$ is open so that it admits a Whitney covering by a family $W_E$ of non-overlapping cubes. See, e.g. \cite{St}, or \cite{Guz}. These cubes have the following properties:
%----------------------------------------------------------
\medskip
%----------------------------------------------------------
%@@@@@@@@@@@@@@@@@@@@@@@@@@@@@@@@@@@@@@@@@@@@@@@@@@@@@@@@@@
%----------------------------------------------------------
\par (i). $\RN\setminus E=\cup\{Q:Q\in W_E\}$;\smallskip
%----------------------------------------------------------
\par (ii). For every cube $Q\in W_E$ we have
%----------------------------------------------------------
\bel{DQ-E}
\diam Q\le \dist(Q,E)\le 4\diam Q.
\ee
%----------------------------------------------------------
%@@@@@@@@@@@@@@@@@@@@@@@@@@@@@@@@@@@@@@@@@@@@@@@@@@@@@@@@@@
%----------------------------------------------------------
\par We are also needed certain additional properties of
Whitney's cubes which we present in the next lemma. These
properties easily follow from constructions of the Whitney covering given in \cite{St} and \cite{Guz}.
%----------------------------------------------------------
%@@@@@@@@@@@@@@@@@@@@@@@@@@@@@@@@@@@@@@@@@@@@@@@@@@@@@@@@@@
\par Given a cube $Q\subset\RN$ let $Q^*:=\frac{9}{8}Q$.
%----------------------------------------------------------
%@@@@@@@@@@@@@@@@@@@@@@@@@@@@@@@@@@@@@@@@@@@@@@@@@@@@@@@@@@
%@@@@@@@@@@@@@@@@@@@@@@@@@@@@@@@@@@@@@@@@@@@@@@@@@@@@@@@@@@
%@@@@@@@@@@@@@@@@@@@@@@@@@@@@@@@@@@@@@@@@@@@@@@@@@@@@@@@@@@
%@@@@@@@@@@@@@@@@@@@@@@@@@@@@@@@@@@@@@@@@@@@@@@@@@@@@@@@@@@
%----------------------------------------------------------
\begin{lemma}\lbl{Wadd}
%----------------------------------------------------------
%@@@@@@@@@@@@@@@@@@@@@@@@@@@@@@@@@@@@@@@@@@@@@@@@@@@@@@@@@@
%----------------------------------------------------------
(1). If $Q,K\in W_E$ and $Q^*\cap K^*\ne\emptyset$, then
%----------------------------------------------------------
$$
\frac{1}{4}\diam Q\le \diam K\le 4\diam Q\,;
$$
%----------------------------------------------------------
%@@@@@@@@@@@@@@@@@@@@@@@@@@@@@@@@@@@@@@@@@@@@@@@@@@@@@@@@@@
%----------------------------------------------------------
%@@@@@@@@@@@@@@@@@@@@@@@@@@@@@@@@@@@@@@@@@@@@@@@@@@@@@@@@@@
\par (2). For every $K\in W_E$ there are at most
$N=N(n)$ cubes from the family
%----------------------------------------------------------
$W_E^*:=\{Q^*:Q\in W_E\}$
%----------------------------------------------------------
which intersect $K^*$;
%----------------------------------------------------------
\medskip
%----------------------------------------------------------
\par (3). If $Q,K\in W_E$, then $Q^*\cap K^*\ne\emptyset$
if and only if  $Q\cap K\ne\emptyset$.
%----------------------------------------------------------
%@@@@@@@@@@@@@@@@@@@@@@@@@@@@@@@@@@@@@@@@@@@@@@@@@@@@@@@@@@
\end{lemma}
%----------------------------------------------------------
%@@@@@@@@@@@@@@@@@@@@@@@@@@@@@@@@@@@@@@@@@@@@@@@@@@@@@@@@@@
%@@@@@@@@@@@@@@@@@@@@@@@@@@@@@@@@@@@@@@@@@@@@@@@@@@@@@@@@@@
%@@@@@@@@@@@@@@@@@@@@@@@@@@@@@@@@@@@@@@@@@@@@@@@@@@@@@@@@@@
%----------------------------------------------------------
\par Let $\Phi_E:=\{\varphi_Q:Q\in W_E\}$ be a smooth partition of unity subordinated to the Whitney decomposition $W_E$. Recall the main properties of this partition.
%@@@@@@@@@@@@@@@@@@@@@@@@@@@@@@@@@@@@@@@@@@@@@@@@@@@@@@@@@@
\begin{lemma}\lbl{P-U} The family of functions $\Phi_E$ has the following properties:
%----------------------------------------------------------
%@@@@@@@@@@@@@@@@@@@@@@@@@@@@@@@@@@@@@@@@@@@@@@@@@@@@@@@@@@
%----------------------------------------------------------
\medskip
%----------------------------------------------------------
\par (a). $\varphi_Q\in C^\infty(\RN)$ and
$0\le\varphi_Q\le 1$ for every $Q\in W_E$;\smallskip
%----------------------------------------------------------
\medskip
%----------------------------------------------------------
\par (b). $\supp \varphi_Q\subset Q^*(:=\frac{9}{8}Q),$
$Q\in W_E$;\smallskip
%----------------------------------------------------------
\medskip
%----------------------------------------------------------
\par (c). $\smed\,\{\varphi_Q(x):Q\in W_E\}=1$ for every
$x\in\RN\setminus S$;\smallskip
%----------------------------------------------------------
\medskip
%----------------------------------------------------------
\par (d). For every cube $Q\in W_E$, every $x\in\RN$ and every multiindex $\beta, |\beta|\le m,$ the following inequality
%----------------------------------------------------------
$$
|D^\beta\varphi_Q(x)| \le C(n,m)\,(\diam Q)^{-|\beta|}
$$
%----------------------------------------------------------
holds.
%----------------------------------------------------------
%@@@@@@@@@@@@@@@@@@@@@@@@@@@@@@@@@@@@@@@@@@@@@@@@@@@@@@@@@@
\end{lemma}
%----------------------------------------------------------
%@@@@@@@@@@@@@@@@@@@@@@@@@@@@@@@@@@@@@@@@@@@@@@@@@@@@@@@@@@
%@@@@@@@@@@@@@@@@@@@@@@@@@@@@@@@@@@@@@@@@@@@@@@@@@@@@@@@@@@
%----------------------------------------------------------
\medskip
%----------------------------------------------------------
%@@@@@@@@@@@@@@@@@@@@@@@@@@@@@@@@@@@@@@@@@@@@@@@@@@@@@@@@@@
%@@@@@@@@@@@@@@@@@@@@@@@@@@@@@@@@@@@@@@@@@@@@@@@@@@@@@@@@@@
%@@@@@@@@@@@@@@@@@@@@@@@@@@@@@@@@@@@@@@@@@@@@@@@@@@@@@@@@@@
%@@@@@@@@@@@@@@@@@@@@@@@@@@@@@@@@@@@@@@@@@@@@@@@@@@@@@@@@@@
%@@@@@@@@@@@@@@@@@@@@@@@@@@@@@@@@@@@@@@@@@@@@@@@@@@@@@@@@@@
%----------------------------------------------------------
\par Let $\vs\ge 1$ be a constant. Let $Q\in W_E$ be a Whitney cube and let $A\in E$ be a point such that
%---------------------------------------------------------
\bel{A-NP}
\dist(A,Q)\le\vs \dist(Q,E).
\ee
%---------------------------------------------------------
We refer to $A$ as a {\it $\vs$-nearest point to the cube $Q$}.
%---------------------------------------------------------
\par Clearly, a point $A\in E$ is a nearest point to $Q$ on $E$ if and only if $A$ is a $1$-nearest point to $Q$. Also it can be readily seen that
%---------------------------------------------------------
\bel{IN-ANP}
A\in (8\vs+1)\,Q
\ee
%---------------------------------------------------------
provided $A\in E$ is a $\vs$-nearest point to a Whitney cube $Q\in W_E$. Conversely, if $A\in (\gamma Q)\cap E$ where $Q\in W_E$ and $\gamma>0$ is a constant, then
%---------------------------------------------------------
\bel{IN-GQ}
A~~~\text{is a}~~\frac{\gamma+1}{2}-\text{nearest point to}~~ Q.
\ee
%--------------------------------------------------------- %@@@@@@@@@@@@@@@@@@@@@@@@@@@@@@@@@@@@@@@@@@@@@@@@@@@@@@@@@
%---------------------------------------------------------
\smallskip
%---------------------------------------------------------
\par Suppose that to every cube $Q\in W_E$ we have assigned a $\vs$-nearest point $A_{Q,\vs}\in E$. For the sake of brevity we denote this point by $a_Q$; thus $a_Q=A_{Q,\vs}$. In particular, by \rf{A-NP} and \rf{IN-ANP},
%---------------------------------------------------------
\bel{AQ-TNP}
\dist(a_Q,Q)\le\vs\dist(Q,E)
~~\text{and}~~a_Q\in(8\vs+1)\,Q~~\text{for every}~~Q\in W_E.
\ee
%--------------------------------------------------------- %@@@@@@@@@@@@@@@@@@@@@@@@@@@@@@@@@@@@@@@@@@@@@@@@@@@@@@@@@
%---------------------------------------------------------
\par By $P^{(Q)}$ we denote the polynomial $P_{a_Q}$. Finally, we define the extension $F$ by the Whitney extension formula:
%----------------------------------------------------------
\bel{DEF-F}
F(x):=\left \{
%----------------------------------------------------------
\begin{array}{ll}
P_x(x),& x\in E,\smallskip\\
\sbig\limits_{Q\in W_E}\,\,
\varphi_Q(x)P^{(Q)}(x),& x\in\RN\setminus E.
\end{array}
%----------------------------------------------------------
\right.
\ee
%----------------------------------------------------------
\medskip
%----------------------------------------------------------
%@@@@@@@@@@@@@@@@@@@@@@@@@@@@@@@@@@@@@@@@@@@@@@@@@@@@@@@@@@
%@@@@@@@@@@@@@@@@@@@@@@@@@@@@@@@@@@@@@@@@@@@@@@@@@@@@@@@@@@
%@@@@@@@@@@@@@@@@@@@@@@@@@@@@@@@@@@@@@@@@@@@@@@@@@@@@@@@@@@
%@@@@@@@@@@@@@@@@@@@@@@@@@@@@@@@@@@@@@@@@@@@@@@@@@@@@@@@@@@
%@@@@@@@@@@@@@@@@@@@@@@@@@@@@@@@@@@@@@@@@@@@@@@@@@@@@@@@@@@
%----------------------------------------------------------
\par Let us note that the metric $d$ is continuous with respect to the Euclidean distance, i.e., for every $x\in\RN$
%---------------------------------------------------------
\bel{C-M}
d(x,y)\to 0 ~~~\text{as}~~~\|x-y\|\to 0.
\ee
%---------------------------------------------------------
\par In fact, since $d\in\DPQ$, by definition \rf{D-LIP}, $d=d_q(h)$ for some non-negative function $h\in\LPRN$.
Then, by Theorem \reff{MAIN-MT},
%---------------------------------------------------------
$$
d(x,y)\sim\dq(x,y:h)=\|x-y\| \sup_{Q\ni x,y} \left(\frac{1}{|Q|}\intl_{Q} h^q(u)\,du\right)^{\frac{1}{q}}.
$$
%----------------------------------------------------------
Since $n\le q$, $h\in L_{q,loc}(\RN)$ and $\diam Q_{xy}=2\,\|x-y\|\to 0$ as $\|x-y\|\to 0$, we have $d(x,y)\to 0$ proving \rf{C-M}.\medskip
%----------------------------------------------------------
%@@@@@@@@@@@@@@@@@@@@@@@@@@@@@@@@@@@@@@@@@@@@@@@@@@@@@@@@@@
%@@@@@@@@@@@@@@@@@@@@@@@@@@@@@@@@@@@@@@@@@@@@@@@@@@@@@@@@@@
%@@@@@@@@@@@@@@@@@@@@@@@@@@@@@@@@@@@@@@@@@@@@@@@@@@@@@@@@@@
%@@@@@@@@@@@@@@@@@@@@@@@@@@@@@@@@@@@@@@@@@@@@@@@@@@@@@@@@@@
%@@@@@@@@@@@@@@@@@@@@@@@@@@@@@@@@@@@@@@@@@@@@@@@@@@@@@@@@@@
%----------------------------------------------------------
\par Clearly, if a cube $Q\ni x,y$, then $|Q|\ge \|x-y\|^n$. (Recall that we measure the distance in the uniform metric.) Since $n<q<p$, we obtain
%---------------------------------------------------------
\be
d(x,y)&\le& C\,\|x-y\| \sup_{Q\ni x,y} \left(\frac{1}{|Q|}\intl_{Q} h^q(u)\,du\right)^{\frac{1}{q}}
\le C\,\|x-y\| \sup_{Q\ni x,y} \left(\frac{1}{|Q|}\intl_{Q} h^p(u)\,du\right)^{\frac{1}{p}}\nn\\
&\le&
C\,\|x-y\|\,\|x-y\|^{-\frac{n}{p}}\,\|h\|_{\LPRN}=
C\,\|x-y\|^{1-\frac{n}{p}}\,\|h\|_{\LPRN}\to 0
\nn
\ee
%----------------------------------------------------------
as $\|x-y\|\to 0$.
%----------------------------------------------------------
\par Hence, by \rf{W-P}, for every multiindex $\beta$, $|\beta|\le m$, we have
%---------------------------------------------------------
$$
D^{\beta}P_x(x)-D^{\beta}P_y(x)
=o( \|x-y\|^{m-|\beta|}),~~~~~x,y\in E,
$$
%---------------------------------------------------------
so that the $m$-jet $\VP=\{P_x\in\PMP: x\in E\}$ satisfies the hypothesis of the Whitney extension theorem \cite{W1}. In this paper Whitney proved that, whenever $\vs=1$ (i.e., $a_Q$ is a point nearest to $Q$ on $E$) the extension $F:\RN\to\R$ defined by formula \rf{DEF-F} is a $C^m$-function which agrees with the Whitney $m$-field $\VP$, i.e.,
%---------------------------------------------------------
\bel{W-EF}
D^{\beta}F(x)=D^{\beta}P_x(x)~~\text{for every}~~x\in E~~\text{and every}~~\beta, |\beta|\le m.
\ee
%----------------------------------------------------------
%@@@@@@@@@@@@@@@@@@@@@@@@@@@@@@@@@@@@@@@@@@@@@@@@@@@@@@@@@@
%@@@@@@@@@@@@@@@@@@@@@@@@@@@@@@@@@@@@@@@@@@@@@@@@@@@@@@@@@@
%----------------------------------------------------------
\par Stein \cite{St}, p. 172, noticed that the approach suggested by Whitney in \cite{W1} works for any mapping
%---------------------------------------------------------
$$
W_E\ni Q\mapsto a_Q\in E
$$
%---------------------------------------------------------
provided $a_Q$ is a $\vs$-nearest point to $Q\in W_E$. In particular, the property \rf{W-EF} holds for such a choice of the point $a_Q$.
%----------------------------------------------------------
%@@@@@@@@@@@@@@@@@@@@@@@@@@@@@@@@@@@@@@@@@@@@@@@@@@@@@@@@@@
%@@@@@@@@@@@@@@@@@@@@@@@@@@@@@@@@@@@@@@@@@@@@@@@@@@@@@@@@@@
%@@@@@@@@@@@@@@@@@@@@@@@@@@@@@@@@@@@@@@@@@@@@@@@@@@@@@@@@@@
%@@@@@@@@@@@@@@@@@@@@@@@@@@@@@@@@@@@@@@@@@@@@@@@@@@@@@@@@@@
%----------------------------------------------------------
\smallskip
%----------------------------------------------------------
\par Prove that $\|F\|_{\CMD}\le C\lambda$, i.e., for every multiindex $\alpha$, $|\alpha|=m$, and every $x,y\in \RN$ the following inequality
%---------------------------------------------------------
\bel{AIM}
|D^{\alpha}F(x)-D^{\alpha}F(y)|\le C\lambda\,d(x,y)
\ee
%---------------------------------------------------------
holds. Here $C=C(m,n,p,\vs)$.\medskip
%----------------------------------------------------------
%@@@@@@@@@@@@@@@@@@@@@@@@@@@@@@@@@@@@@@@@@@@@@@@@@@@@@@@@@@
%@@@@@@@@@@@@@@@@@@@@@@@@@@@@@@@@@@@@@@@@@@@@@@@@@@@@@@@@@@
%----------------------------------------------------------
\par Consider four cases.
%----------------------------------------------------------
%@@@@@@@@@@@@@@@@@@@@@@@@@@@@@@@@@@@@@@@@@@@@@@@@@@@@@@@@@@
%@@@@@@@@@@@@@@@@@@@@@@@@@@@@@@@@@@@@@@@@@@@@@@@@@@@@@@@@@@
%@@@@@@@@@@@@@@@@@@@@@@@@@@@@@@@@@@@@@@@@@@@@@@@@@@@@@@@@@@
%@@@@@@@@@@@@@@@@@@@@@@@@@@@@@@@@@@@@@@@@@@@@@@@@@@@@@@@@@@
%----------------------------------------------------------
\par {\it The first case: $x,y\in E$.} Since $P_y\in\PMP$, for every multiindex $\alpha$ of order $m$ the function $D^\alpha P_y$ is a constant function. In particular, $D^{\alpha}P_y(x)=D^{\alpha}P_y(y)$. Hence, by \rf{W-EF},
%---------------------------------------------------------
$$
|D^{\alpha}F(x)-D^{\alpha}F(y)|=
|D^{\alpha}P_x(x)-D^{\alpha}P_y(y)|
=|D^{\alpha}P_x(x)-D^{\alpha}P_y(x)|
$$
%---------------------------------------------------------
so that, by \rf{W-P},
%---------------------------------------------------------
$$
|D^{\alpha}F(x)-D^{\alpha}F(y)|\le \lambda\,d(x,y)
$$
%---------------------------------------------------------
proving \rf{AIM} in the case under consideration.\bigskip
%----------------------------------------------------------
%@@@@@@@@@@@@@@@@@@@@@@@@@@@@@@@@@@@@@@@@@@@@@@@@@@@@@@@@@@
%@@@@@@@@@@@@@@@@@@@@@@@@@@@@@@@@@@@@@@@@@@@@@@@@@@@@@@@@@@
%@@@@@@@@@@@@@@@@@@@@@@@@@@@@@@@@@@@@@@@@@@@@@@@@@@@@@@@@@@
%@@@@@@@@@@@@@@@@@@@@@@@@@@@@@@@@@@@@@@@@@@@@@@@@@@@@@@@@@@
%@@@@@@@@@@@@@@@@@@@@@@@@@@@@@@@@@@@@@@@@@@@@@@@@@@@@@@@@@@
%@@@@@@@@@@@@@@@@@@@@@@@@@@@@@@@@@@@@@@@@@@@@@@@@@@@@@@@@@@
%----------------------------------------------------------
\par {\it The second case: $x\in E$, $y\in\RN\setminus E$.} Given a Whitney cube $K\in W_E$ let
%---------------------------------------------------------
\bel{TK}
T(K):=\{Q\in W_E: Q\cap K\ne\emp\}
\ee
%---------------------------------------------------------
be a family of all Whitney cubes touching $K$.
%----------------------------------------------------------
%@@@@@@@@@@@@@@@@@@@@@@@@@@@@@@@@@@@@@@@@@@@@@@@@@@@@@@@@@@
%@@@@@@@@@@@@@@@@@@@@@@@@@@@@@@@@@@@@@@@@@@@@@@@@@@@@@@@@@@
%@@@@@@@@@@@@@@@@@@@@@@@@@@@@@@@@@@@@@@@@@@@@@@@@@@@@@@@@@@
%@@@@@@@@@@@@@@@@@@@@@@@@@@@@@@@@@@@@@@@@@@@@@@@@@@@@@@@@@@
%@@@@@@@@@@@@@@@@@@@@@@@@@@@@@@@@@@@@@@@@@@@@@@@@@@@@@@@@@@
%@@@@@@@@@@@@@@@@@@@@@@@@@@@@@@@@@@@@@@@@@@@@@@@@@@@@@@@@@@
%----------------------------------------------------------
\begin{lemma}\lbl{C-A} Let $K\in W_E$ be a Whitney cube and let $y\in K^*=\tfrac98 K.$ Then for every multiindex $\alpha$ the following inequality
%----------------------------------------------------------
%@@@@@@@@@@@@@@@@@@@@@@@@@@@@@@@@@@@@@@@@@@@@@@@@@@@@@@@@@@
%----------------------------------------------------------
$$
|D^{\alpha}F(y)-D^{\alpha}P_{a_K}(y)|\le C\,
\smed\limits_{Q\in T(K)}\,\,\smed\limits_{|\xi|\le m}\,\,
(\diam K)^{|\xi|-|\alpha|} \,|D^{\xi}P_{a_Q}(a_K)-D^{\xi}P_{a_K}(a_K)|
$$
%----------------------------------------------------------
holds. Here $C$ is a constant depending only on $n,m,\alpha$ and $\vs$.
%@@@@@@@@@@@@@@@@@@@@@@@@@@@@@@@@@@@@@@@@@@@@@@@@@@@@@@@@@@
\end{lemma}
%----------------------------------------------------------
%@@@@@@@@@@@@@@@@@@@@@@@@@@@@@@@@@@@@@@@@@@@@@@@@@@@@@@@@@@
%@@@@@@@@@@@@@@@@@@@@@@@@@@@@@@@@@@@@@@@@@@@@@@@@@@@@@@@@@@
%@@@@@@@@@@@@@@@@@@@@@@@@@@@@@@@@@@@@@@@@@@@@@@@@@@@@@@@@@@
%----------------------------------------------------------
\par {\it Proof.} Note that, by part (2) of Lemma \reff{Wadd}, $\#T(K)\le N(n)$, and, by part (3) of this lemma,
%---------------------------------------------------------
$$
T(K)=\{Q\in W_E: Q^*\cap K^*\ne\emp\}.
$$
%---------------------------------------------------------
\par Recall that $P^{(Q)}=P_{a_Q}$ and $a_Q\in 9Q$ for every $Q\in W_E$. Let us estimate the quantity
%---------------------------------------------------------
$$
I:=|D^{\alpha}F(y)-D^{\alpha}P_{a_K}(y)|.
$$
%---------------------------------------------------------
By formula \rf{DEF-F} and by part (c) of Lemma \reff{Wadd}, %---------------------------------------------------------
$$
F(y)-P_{a_K}(y)=\smed\limits_{Q\in W_E}
\varphi_Q(y)(P^{(Q)}(y)-P_{a_K}(y))
$$
%---------------------------------------------------------
so that, by part (b) of Lemma \reff{Wadd} and by definition \rf{TK},
%---------------------------------------------------------
$$
F(y)-P_{a_K}(y)=\smed\limits_{Q\in T(K)}
\varphi_Q(y)(P^{(Q)}(y)-P_{a_K}(y)).
$$
%---------------------------------------------------------
Hence
%---------------------------------------------------------
$$
I:=|D^\alpha F(y)-D^\alpha P_{a_K}(y)|\le \smed\limits_{Q\in T(K)}
|D^\alpha(\varphi_Q(y)(P^{(Q)}(y)-P_{a_K}(y)))|
$$
%---------------------------------------------------------
so that
%---------------------------------------------------------
\bel{I2}
I\le \smed\limits_{Q\in T(K)} A_Q(y;\alpha)
\ee
%---------------------------------------------------------
where
%---------------------------------------------------------
$$
A_Q(y;\alpha):=
|D^\alpha(\varphi_Q(y)(P^{(Q)}(y)-P_{a_K}(y)))|.
$$
%---------------------------------------------------------
\par Let $Q\in T(K)$. Then
%---------------------------------------------------------
$$
A_Q(y;\alpha)\le C\smed\limits_{|\beta|+|\gamma|=|\alpha|}
|D^\beta\varphi_Q(y)|\,|D^\gamma(P^{(Q)}(y)-P_{a_K}(y))|
$$
%---------------------------------------------------------
so that, by part (d) of Lemma \reff{P-U},
%---------------------------------------------------------
$$
A_Q(y;\alpha)\le C\smed\limits_{|\beta|+|\gamma|=|\alpha|}
(\diam Q)^{-|\beta|} \,|D^\gamma(P^{(Q)}(y)-P_{a_K}(y))|.
$$
%---------------------------------------------------------
Since $Q\cap K\ne\emp$, by part (1) of Lemma \reff{Wadd}, $\diam Q\sim\diam K$ so that
%---------------------------------------------------------
\bel{AQA}
A_Q(y;\alpha)\le C\smed\limits_{|\beta|+|\gamma|=|\alpha|}
(\diam K)^{-|\beta|} \,|D^\gamma(P^{(Q)}(y)-P_{a_K}(y))|.
\ee
%---------------------------------------------------------
%@@@@@@@@@@@@@@@@@@@@@@@@@@@@@@@@@@@@@@@@@@@@@@@@@@@@@@@@@
%---------------------------------------------------------
\par Let us estimate the distance between $a_Q$ and $y$. By \rf{DQ-E} and \rf{AQ-TNP},
%---------------------------------------------------------
\be
\|a_Q-y\|&\le& \diam K^*+\diam Q+\dist(a_Q,Q)\nn\\
&\le& 2\diam K+\diam Q+\vs\dist(Q,E)\nn\\
&\le&2\diam K+(1+4\vs)\diam Q.\nn
\ee
%---------------------------------------------------------
Since $Q\cap K\ne\emp$, by part (1) of Lemma \reff{Wadd},
%---------------------------------------------------------
\bel{D-AQ}
\|a_Q-y\|\le 2\diam K+4(1+4\vs)\diam K\le 22\,\vs\diam K.
\ee
%---------------------------------------------------------
%@@@@@@@@@@@@@@@@@@@@@@@@@@@@@@@@@@@@@@@@@@@@@@@@@@@@@@@@@
%---------------------------------------------------------
\par Let
%---------------------------------------------------------
$$
\tP_Q:=P^{(Q)}-P_{a_K}=P_{a_Q}-P_{a_K}.
$$
%---------------------------------------------------------
Let us estimate the quantity $|D^\gamma\tP_Q(y)|$. Since $\tP_Q\in\PMP$, we can represent this polynomial in the form
%---------------------------------------------------------
$$
\tP_Q(z)=\smed\limits_{|\xi|\le m}\,\,\frac{1}{\xi!}\, D^{\xi}\tP_Q(a_K)\,(z-a_K)^{\xi}.
$$
%---------------------------------------------------------
Hence
%---------------------------------------------------------
$$
D^\gamma\tP_Q(z)=\smed\limits_{|\gamma|\le |\xi|\le m}\,\,\frac{1}{(\xi-\gamma)!}\, D^{\xi}\tP_Q(a_K)\,(z-a_K)^{\xi-\gamma}
$$
%---------------------------------------------------------
so that
%---------------------------------------------------------
$$
|D^\gamma\tP_Q(y)|\le C\smed\limits_{|\gamma|\le |\xi|\le m} \,
|D^{\xi}\tP_Q(a_K)|\,\|y-a_K\|^{|\xi|-|\gamma|}\le
C\smed\limits_{|\gamma|\le |\xi|\le m} \,\,
(\diam K)^{|\xi|-|\gamma|}\,|D^{\xi}\tP_Q(a_K)|.
$$
%---------------------------------------------------------
Combining this inequality with \rf{AQA} we obtain
%---------------------------------------------------------
\be
A_Q(y;\alpha)&\le& C\smed\limits_{|\beta|+|\gamma|=|\alpha|}
(\diam K)^{-|\beta|} \,|D^\gamma\tP_Q(y)|\nn\\
&\le& C\smed\limits_{|\beta|+|\gamma|=|\alpha|}
(\diam K)^{-|\beta|} \,\smed\limits_{|\gamma|\le |\xi|\le m}(\diam K)^{|\xi|-|\gamma|}\,|D^{\xi}\tP_Q(a_K)|\nn\\
&\le&
C \,\smed\limits_{|\xi|\le m}(\diam K)^{|\xi|-|\alpha|}\,|D^{\xi}\tP_Q(a_K)|.\nn
\ee
%---------------------------------------------------------
Hence, by \rf{I2},
%---------------------------------------------------------
$$
I\le C\,\smed\limits_{Q\in T(K)}\,\,
\smed\limits_{|\xi|\le m}
(\diam K)^{|\xi|-|\alpha|}\,|D^{\xi}\tP_Q(a_K)|
$$
%---------------------------------------------------------
proving the lemma.\bx
%----------------------------------------------------------
\smallskip
\par Let us apply Lemma \reff{C-A} to an arbitrary multiindex $\alpha$ of order $m+1$. Since $D^\alpha P=0$ for every polynomial $P\in\PMP$, we obtain the following statement.
%----------------------------------------------------------
%@@@@@@@@@@@@@@@@@@@@@@@@@@@@@@@@@@@@@@@@@@@@@@@@@@@@@@@@@@
%@@@@@@@@@@@@@@@@@@@@@@@@@@@@@@@@@@@@@@@@@@@@@@@@@@@@@@@@@@
%@@@@@@@@@@@@@@@@@@@@@@@@@@@@@@@@@@@@@@@@@@@@@@@@@@@@@@@@@@
%@@@@@@@@@@@@@@@@@@@@@@@@@@@@@@@@@@@@@@@@@@@@@@@@@@@@@@@@@@
%@@@@@@@@@@@@@@@@@@@@@@@@@@@@@@@@@@@@@@@@@@@@@@@@@@@@@@@@@@
%@@@@@@@@@@@@@@@@@@@@@@@@@@@@@@@@@@@@@@@@@@@@@@@@@@@@@@@@@@
%----------------------------------------------------------
\begin{lemma}\lbl{MPO-C} Let $K\in W_E$ be a Whitney cube and let $y\in K^*=\tfrac98 K.$ Then for every multiindex $\alpha$, $|\alpha|=m+1$, the following inequality
%----------------------------------------------------------
%@@@@@@@@@@@@@@@@@@@@@@@@@@@@@@@@@@@@@@@@@@@@@@@@@@@@@@@@@@
%----------------------------------------------------------
$$
|D^{\alpha}F(y)|\le C\,
\smed\limits_{Q\in T(K)}\,\,\smed\limits_{|\xi|\le m}\,\,
(\diam K)^{|\xi|-m-1} \,|D^{\xi}P_{a_Q}(a_K)-D^{\xi}P_{a_K}(a_K)|
$$
%----------------------------------------------------------
holds. Here $C$ is a constant depending only on $m,n,$ and $\vs$.
%@@@@@@@@@@@@@@@@@@@@@@@@@@@@@@@@@@@@@@@@@@@@@@@@@@@@@@@@@@
\end{lemma}
%----------------------------------------------------------
%@@@@@@@@@@@@@@@@@@@@@@@@@@@@@@@@@@@@@@@@@@@@@@@@@@@@@@@@@@
%@@@@@@@@@@@@@@@@@@@@@@@@@@@@@@@@@@@@@@@@@@@@@@@@@@@@@@@@@@
%@@@@@@@@@@@@@@@@@@@@@@@@@@@@@@@@@@@@@@@@@@@@@@@@@@@@@@@@@@
%@@@@@@@@@@@@@@@@@@@@@@@@@@@@@@@@@@@@@@@@@@@@@@@@@@@@@@@@@@
%----------------------------------------------------------
\begin{lemma}\lbl{C-2} Let $x\in E$ and let $K\in W_E$ be a Whitney cube. Then for every $y\in K$ and every $\alpha,|\alpha|=m,$ the following inequality
%----------------------------------------------------------
%@@@@@@@@@@@@@@@@@@@@@@@@@@@@@@@@@@@@@@@@@@@@@@@@@@@@@@@@@@
%----------------------------------------------------------
\be
|D^{\alpha}F(x)-D^{\alpha}F(y)|&\le& C\,\{
|D^{\alpha}P_x(x)-D^{\alpha}P_{a_K}(x)|\nn\\
&+&
\smed\limits_{Q\in T(K)}\,\,\smed\limits_{|\xi|\le m}
(\diam K)^{|\xi|-m} \,|D^{\xi}P_{a_Q}(a_K)-D^{\xi}P_{a_K}(a_K)|\}\nn
\ee
%----------------------------------------------------------
holds. Here $C$ is a constant depending only on $n,m$ and $\vs$.
%@@@@@@@@@@@@@@@@@@@@@@@@@@@@@@@@@@@@@@@@@@@@@@@@@@@@@@@@@@
\end{lemma}
%----------------------------------------------------------
%@@@@@@@@@@@@@@@@@@@@@@@@@@@@@@@@@@@@@@@@@@@@@@@@@@@@@@@@@@
%@@@@@@@@@@@@@@@@@@@@@@@@@@@@@@@@@@@@@@@@@@@@@@@@@@@@@@@@@@
%@@@@@@@@@@@@@@@@@@@@@@@@@@@@@@@@@@@@@@@@@@@@@@@@@@@@@@@@@@
%----------------------------------------------------------
\par {\it Proof.} We have
%---------------------------------------------------------
$$
|D^{\alpha}F(x)-D^{\alpha}F(y)|
\le |D^{\alpha}F(x)-D^{\alpha}P_{a_K}(y)|
+|D^{\alpha}P_{a_K}(y)-D^{\alpha}F(y)|=I_1+I_2.
$$
%---------------------------------------------------------
Since $P_{a_K}\in\PMP$ and $|\alpha|=m$, the function $D^\alpha P_y$ is a constant function so that
%---------------------------------------------------------
$$
D^{\alpha}P_{a_K}(y)=D^{\alpha}P_{a_K}(x).
$$
%---------------------------------------------------------
Since $x\in E$, by \rf{W-EF}, $D^{\alpha}F(x)=D^{\alpha}P_{x}(x)$ so that, by \rf{W-P},  %---------------------------------------------------------
$$
I_1:=|D^{\alpha}F(x)-D^{\alpha}P_{a_K}(y)|=
|D^{\alpha}P_{x}(x)-D^{\alpha}P_{a_K}(x)|.
$$
%---------------------------------------------------------
It remains to apply Lemma \reff{C-A} to
$I_2:=|D^{\alpha}F(y)-D^{\alpha}P_{a_K}(y)|$,
and the lemma follows.\bx
%----------------------------------------------------------
%@@@@@@@@@@@@@@@@@@@@@@@@@@@@@@@@@@@@@@@@@@@@@@@@@@@@@@@@@@
%@@@@@@@@@@@@@@@@@@@@@@@@@@@@@@@@@@@@@@@@@@@@@@@@@@@@@@@@@@
%@@@@@@@@@@@@@@@@@@@@@@@@@@@@@@@@@@@@@@@@@@@@@@@@@@@@@@@@@@
%----------------------------------------------------------
\bigskip
%@@@@@@@@@@@@@@@@@@@@@@@@@@@@@@@@@@@@@@@@@@@@@@@@@@@@@@@@@@
%@@@@@@@@@@@@@@@@@@@@@@@@@@@@@@@@@@@@@@@@@@@@@@@@@@@@@@@@@@
%@@@@@@@@@@@@@@@@@@@@@@@@@@@@@@@@@@@@@@@@@@@@@@@@@@@@@@@@@@
%@@@@@@@@@@@@@@@@@@@@@@@@@@@@@@@@@@@@@@@@@@@@@@@@@@@@@@@@@@
%@@@@@@@@@@@@@@@@@@@@@@@@@@@@@@@@@@@@@@@@@@@@@@@@@@@@@@@@@@
%@@@@@@@@@@@@@@@@@@@@@@@@@@@@@@@@@@@@@@@@@@@@@@@@@@@@@@@@@@
%----------------------------------------------------------
\par Let us prove inequality \rf{AIM} for arbitrary $x\in E$, $y\in\RN\setminus E$ and $\alpha,\,|\alpha|=m$. Let $y\in K$ for some $K\in W_E$. By Lemma \reff{C-2},
%----------------------------------------------------------
$$
|D^{\alpha}F(x)-D^{\alpha}F(y)|\le C\{J_1+J_2\}
$$
%----------------------------------------------------------
where
%----------------------------------------------------------
$$
J_1:=|D^{\alpha}P_x(x)-D^{\alpha}P_{a_K}(x)|
$$
%----------------------------------------------------------
and
%----------------------------------------------------------
$$
J_2:=
\smed\limits_{Q\in T(K)}\,\,\smed\limits_{|\xi|\le m}
(\diam K)^{|\xi|-m} \,|D^{\xi}P_{a_Q}(a_K)-D^{\xi}P_{a_K}(a_K)|
$$
%----------------------------------------------------------
\par First let us estimate $J_2$. By \rf{W-P}, for every $\xi, |\xi|\le m,$ and every $Q\in T(K)$, we have
%---------------------------------------------------------
$$
|D^{\xi}P_{a_Q}(a_K)-D^{\xi}P_{a_K}(a_K)|
\le
\lambda\,\|a_Q-a_K\|^{m-|\xi|}d(a_Q,a_K).
$$
%---------------------------------------------------------
Hence
%---------------------------------------------------------
$$
J_2\le\lambda
\smed\limits_{Q\in T(K)}\,\,\smed\limits_{|\xi|\le m}
(\diam K)^{|\xi|-m}\|a_Q-a_K\|^{m-|\xi|}d(a_Q,a_K).
$$
%---------------------------------------------------------
But, by \rf{D-AQ},
%---------------------------------------------------------
\bel{A-4}
\|a_Q-a_K\|\le \|a_Q-y\|+\|y-a_K\|\le (22\vs+1)\diam K
\ee
%---------------------------------------------------------
proving that
%---------------------------------------------------------
\bel{E-J2}
J_2\le C(n,m,\vs)\,\lambda
\smed\limits_{Q\in T(K)}d(a_Q,a_K).
\ee
%---------------------------------------------------------
\medskip
\par Now prove that for some constant $C=C(n,p,q,\vs)$
%---------------------------------------------------------
\bel{D-AY}
d(a_Q,y)\le C\,d(x,y)~~~\text{for every}~~~Q\in T(K).
\ee
%---------------------------------------------------------
In fact, by \rf{D-AQ}, $\|a_Q-y\|\le 22\vs\diam K$. Since $x\in E$ and $y\in K$, by \rf{DQ-E},
%---------------------------------------------------------
$$
\diam K\le 4\dist(K,E)\le 4\|x-y\|
$$
%---------------------------------------------------------
so that $\|a_Q-y\|\le 88\vs\|x-y\|$. Hence, by part (b) of Proposition \reff{VMT}, see inequality \rf{O-1},
$d(a_Q,y)\le C\,d(x,y)$ proving \rf{D-AY}.
%---------------------------------------------------------
\medskip
\par Now we have
%---------------------------------------------------------
$$
d(a_Q,a_K)\le d(a_Q,y)+d(y,a_K)\le C\,d(x,y)
$$
%---------------------------------------------------------
so that, by \rf{E-J2},
%---------------------------------------------------------
$$
J_2\le C\,\lambda\, \#T(K)\,d(x,y)\le C\,\lambda\, d(x,y).
$$
%---------------------------------------------------------
See part (2) of Lemma \reff{Wadd}.
%---------------------------------------------------------
%@@@@@@@@@@@@@@@@@@@@@@@@@@@@@@@@@@@@@@@@@@@@@@@@@@@@@@@@@@
%@@@@@@@@@@@@@@@@@@@@@@@@@@@@@@@@@@@@@@@@@@@@@@@@@@@@@@@@@@
%@@@@@@@@@@@@@@@@@@@@@@@@@@@@@@@@@@@@@@@@@@@@@@@@@@@@@@@@@@
%@@@@@@@@@@@@@@@@@@@@@@@@@@@@@@@@@@@@@@@@@@@@@@@@@@@@@@@@@@
%---------------------------------------------------------
\par On the other hand, by \rf{W-P} and \rf{D-AY},
%---------------------------------------------------------
$$
J_1:=|D^{\alpha}P_x(x)-D^{\alpha}P_{a_K}(x)|\le \lambda \,d(x,a_K)\le \lambda (d(x,y)+d(y,a_K))\le C\,\lambda\,d(x,y).
$$
%---------------------------------------------------------
Finally,
%---------------------------------------------------------
$$
|D^{\alpha}F(x)-D^{\alpha}F(y)|\le C\,\{J_1+J_2\}\le C\lambda\,d(x,y).
$$
%----------------------------------------------------------
\bigskip
%@@@@@@@@@@@@@@@@@@@@@@@@@@@@@@@@@@@@@@@@@@@@@@@@@@@@@@@@@@
%@@@@@@@@@@@@@@@@@@@@@@@@@@@@@@@@@@@@@@@@@@@@@@@@@@@@@@@@@@
%@@@@@@@@@@@@@@@@@@@@@@@@@@@@@@@@@@@@@@@@@@@@@@@@@@@@@@@@@@
%@@@@@@@@@@@@@@@@@@@@@@@@@@@@@@@@@@@@@@@@@@@@@@@@@@@@@@@@@@
%@@@@@@@@@@@@@@@@@@@@@@@@@@@@@@@@@@@@@@@@@@@@@@@@@@@@@@@@@@
%@@@@@@@@@@@@@@@@@@@@@@@@@@@@@@@@@@@@@@@@@@@@@@@@@@@@@@@@@@
%----------------------------------------------------------
\par {\it The third case: $y\in K$, $K\in W_E$ and $x\in\RN\setminus K^*$.}
%----------------------------------------------------------
\par Since $K^*=\tfrac98 K$ and $x\notin K^*$, we have
%---------------------------------------------------------
$$
\|x-y\|\ge \tfrac{1}{16}\diam K.
$$
%---------------------------------------------------------
Let $a\in E$ be a point nearest to $x$ on $E$. Then
%---------------------------------------------------------
$$
\|a-x\|=\dist(x,E)\le\dist(y,E)+\|x-y\|\le \dist(K,E)+
\diam K+\|x-y\|
$$
%---------------------------------------------------------
so that, by \rf{DQ-E},
%---------------------------------------------------------
$$
\|a-x\|\le 4\diam K+\diam K+\|x-y\|\le 81\|x-y\|.
$$
%---------------------------------------------------------
Hence, by part (b) of Claim \reff{VMT}, see \rf{O-1}, $d(a,x)\le C\,d(x,y)$.
%---------------------------------------------------------
\par We have
%---------------------------------------------------------
$$
\|y-a\|\le \|x-y\|+\|x-a\|\le 82\|x-y\|
$$
%---------------------------------------------------------
so that again, by \rf{O-1}, $d(y,a)\le C\,d(x,y)$. We obtain
%---------------------------------------------------------
$$
|D^{\alpha}F(x)-D^{\alpha}F(y)|\le |D^{\alpha}F(x)-D^{\alpha}F(a)|+
|D^{\alpha}F(a)-D^{\alpha}F(y)|
$$
%----------------------------------------------------------
so that, by the result proven in the second case,
%---------------------------------------------------------
$$
|D^{\alpha}F(x)-D^{\alpha}F(y)|\le C\,\lambda\,(d(x,a)+d(y,a))\le C\,\lambda\,d(x,y).
$$
%----------------------------------------------------------
\bigskip
%@@@@@@@@@@@@@@@@@@@@@@@@@@@@@@@@@@@@@@@@@@@@@@@@@@@@@@@@@@
%@@@@@@@@@@@@@@@@@@@@@@@@@@@@@@@@@@@@@@@@@@@@@@@@@@@@@@@@@@
%@@@@@@@@@@@@@@@@@@@@@@@@@@@@@@@@@@@@@@@@@@@@@@@@@@@@@@@@@@
%@@@@@@@@@@@@@@@@@@@@@@@@@@@@@@@@@@@@@@@@@@@@@@@@@@@@@@@@@@
%@@@@@@@@@@@@@@@@@@@@@@@@@@@@@@@@@@@@@@@@@@@@@@@@@@@@@@@@@@
%@@@@@@@@@@@@@@@@@@@@@@@@@@@@@@@@@@@@@@@@@@@@@@@@@@@@@@@@@@
%----------------------------------------------------------
\par {\it The fourth case: $y\in K$, $x\in K^*$ where $K\in W_E$.} The proof of inequality \rf{AIM} in this case is based on the next
%----------------------------------------------------------
%@@@@@@@@@@@@@@@@@@@@@@@@@@@@@@@@@@@@@@@@@@@@@@@@@@@@@@@@@@
%@@@@@@@@@@@@@@@@@@@@@@@@@@@@@@@@@@@@@@@@@@@@@@@@@@@@@@@@@@
%@@@@@@@@@@@@@@@@@@@@@@@@@@@@@@@@@@@@@@@@@@@@@@@@@@@@@@@@@@
%@@@@@@@@@@@@@@@@@@@@@@@@@@@@@@@@@@@@@@@@@@@@@@@@@@@@@@@@@@
%----------------------------------------------------------
\begin{lemma}\lbl{C-4} Let $K\in W_E$ be a Whitney cube and let $x,y\in K^*$. Then for every multiindex  $\alpha,|\alpha|=m,$ the following inequality
%----------------------------------------------------------
%@@@@@@@@@@@@@@@@@@@@@@@@@@@@@@@@@@@@@@@@@@@@@@@@@@@@@@@@@@
%----------------------------------------------------------
$$
|D^{\alpha}F(x)-D^{\alpha}F(y)|\le C\,
\|x-y\|\,
\smed\limits_{Q\in\, T(K)}\,\,\smed\limits_{|\xi|\le m}
(\diam K)^{|\xi|-m-1} \,|D^{\xi}P_{a_Q}(a_K)-D^{\xi}P_{a_K}(a_K)|
$$
%----------------------------------------------------------
holds. Here $C$ is a constant depending only on $n,m$ and $\vs$.
%@@@@@@@@@@@@@@@@@@@@@@@@@@@@@@@@@@@@@@@@@@@@@@@@@@@@@@@@@@
\end{lemma}
%----------------------------------------------------------
%@@@@@@@@@@@@@@@@@@@@@@@@@@@@@@@@@@@@@@@@@@@@@@@@@@@@@@@@@@
%@@@@@@@@@@@@@@@@@@@@@@@@@@@@@@@@@@@@@@@@@@@@@@@@@@@@@@@@@@
%@@@@@@@@@@@@@@@@@@@@@@@@@@@@@@@@@@@@@@@@@@@@@@@@@@@@@@@@@@
%----------------------------------------------------------
\par {\it Proof.} Note that the function $F|_{\RN\setminus E}\in C^{\infty}(\RN\setminus E)$, see formula \rf{DEF-F}, so that, by the Lagrange theorem, for every $\alpha$, $|\alpha|=m$, there exists $z\in[x,y]$ such that
%---------------------------------------------------------
\bel{L-V}
|D^{\alpha}F(x)-D^{\alpha}F(y)|\le C\|x-y\|
\smed\limits_{|\beta|=m+1}\,\, |D^{\beta}F(z)|.
\ee
%----------------------------------------------------------
Since $x,y\in K^*$, the point $z\in K^*$ as well.
%----------------------------------------------------------
%@@@@@@@@@@@@@@@@@@@@@@@@@@@@@@@@@@@@@@@@@@@@@@@@@@@@@@@@@@
%----------------------------------------------------------
\par Combining this inequality with Lemma \reff{MPO-C} we obtain the statement of the lemma.\bx
%----------------------------------------------------------
%@@@@@@@@@@@@@@@@@@@@@@@@@@@@@@@@@@@@@@@@@@@@@@@@@@@@@@@@@@
%@@@@@@@@@@@@@@@@@@@@@@@@@@@@@@@@@@@@@@@@@@@@@@@@@@@@@@@@@@
%@@@@@@@@@@@@@@@@@@@@@@@@@@@@@@@@@@@@@@@@@@@@@@@@@@@@@@@@@@
%@@@@@@@@@@@@@@@@@@@@@@@@@@@@@@@@@@@@@@@@@@@@@@@@@@@@@@@@@@
%----------------------------------------------------------
\bigskip
\par We are in a position to prove inequality \rf{AIM} for arbitrary $y\in K$ and $x\in K^*$. By inequality \rf{W-P}, for every cube $Q\in T(K)$ and every $\xi, |\xi|\le m$,
%---------------------------------------------------------
$$
|D^{\xi}P_{a_Q}(a_K)-D^{\xi}P_{a_K}(a_K)|\le \lambda\,d(a_Q,a_K)(\diam K)^{m-|\xi|}
$$
%----------------------------------------------------------
so that, by Lemma \reff{L-V},
%----------------------------------------------------------
\be
I&:=&|D^{\alpha}F(x)-D^{\alpha}F(y)|\nn\\
&\le&
C\frac{\|x-y\|}{\diam K}
\smed\limits_{Q\in\, T(K)}\,\,\smed\limits_{|\xi|\le m}
(\diam K)^{|\xi|-m} \,|D^{\xi}P_{a_Q}(a_K)-D^{\xi}P_{a_K}(a_K)|\nn\\
&\le&
C\frac{\|x-y\|}{\diam K}
\smed\limits_{Q\in\, T(K)}\,\,\smed\limits_{|\xi|\le m}
(\diam K)^{|\xi|-m} \,(\lambda\,d(a_Q,a_K)(\diam K)^{m-|\xi|})\nn\\
&\le&
C\lambda\frac{\|x-y\|}{\diam K}
\smed\limits_{Q\in\, T(K)}\,d(a_Q,a_K).\nn
\ee
%----------------------------------------------------------
\par Note that, by \rf{A-4} and \rf{DQ-E},
%---------------------------------------------------------
\be
\|a_Q-y\|&\le& \|a_Q-a_K\|+\|a_K-y\|\le 23\diam K+\diam K+\dist(K,E)\nn\\&\le& 23\cdot 4\dist(K,E)+4\dist(K,E)+\dist(K,E)\nn\\&=&97\dist(K,E)\le 97\,\|y-a_K\|.\nn
\ee
%---------------------------------------------------------
Hence, by part (b) of Claim \reff{VMT}, see \rf{O-1}, $d(a_Q,y)\le C\,d(y,a_K)$ so that
%---------------------------------------------------------
$$
d(a_Q,a_K)\le d(a_Q,y)+d(y,a_K)\le C\,d(y,a_K).
$$
%---------------------------------------------------------
\par This implies the following inequality
%---------------------------------------------------------
$$
I\le C\lambda\,(\#T(K))\,\frac{\|x-y\|}{\diam K}\,d(y,a_K).
$$
%---------------------------------------------------------
But, by part (2) of Lemma \reff{Wadd}, $\#T(K)\le N(n)$ so that
%---------------------------------------------------------
$$
I\le C\lambda\,\frac{\|x-y\|}{\diam K}\,d(y,a_K).
$$
%---------------------------------------------------------
\par Since $x,y\in K^*$, the distance
$\|x-y\|\le\diam K^*=\tfrac98\diam K.$ But
%---------------------------------------------------------
$$
\diam K\le 4\dist(K,E)\le 4\|y-a_K\|
$$
%----------------------------------------------------------
so that $\|x-y\|\le 5\|y-a_K\|$. Therefore, by part (b) of Proposition \reff{VMT}, see \rf{O-2},
%---------------------------------------------------------
$$
\frac{\|x-y\|}{\|y-a_K\|}\,d(y,a_K)\le C\,d(x,y).
$$
%----------------------------------------------------------
\par On the other hand,
%---------------------------------------------------------
\be
\|y-a_K\|&\le& \diam K+\dist(a_K,K)=\diam K+\dist(K,E)\nn\\
&\le&
\diam K+4\diam K=5\diam K.\nn
\ee
%----------------------------------------------------------
Finally,
%----------------------------------------------------------
$$
I\le C\lambda\,\frac{\|x-y\|}{\diam K}\,d(y,a_K)\le
5\,C\lambda\,\frac{\|x-y\|}{\|y-a_K\|}\,d(y,a_K)
\le C\lambda\,d(x,y).
$$
%----------------------------------------------------------
%@@@@@@@@@@@@@@@@@@@@@@@@@@@@@@@@@@@@@@@@@@@@@@@@@@@@@@@@@@
%----------------------------------------------------------
\par The proof of Theorem \reff{CMD-E-NEW} is complete.\bx
%----------------------------------------------------------
%@@@@@@@@@@@@@@@@@@@@@@@@@@@@@@@@@@@@@@@@@@@@@@@@@@@@@@@@@@
%@@@@@@@@@@@@@@@@@@@@@@@@@@@@@@@@@@@@@@@@@@@@@@@@@@@@@@@@@@
%@@@@@@@@@@@@@@@@@@@@@@@@@@@@@@@@@@@@@@@@@@@@@@@@@@@@@@@@@@
%@@@@@@@@@@@@@@@@@@@@@@@@@@@@@@@@@@@@@@@@@@@@@@@@@@@@@@@@@@
%----------------------------------------------------------
\begin{remark} \lbl{MC-EXT-REM} {\em For $m=0$ the statement of Theorem \reff{CMD-E-NEW} is true for an {\it arbitrary} metric $d$ on $\RN$. In fact, in this case $C^{0,(d)}(\RN)=\Lip(\RN;d)$. By the McShane extension theorem \cite{McS}, every function $f\in \Lip(E;d)$ extends to a function $F\in\Lip(\RN;d)$ such that $\|F\|_{\Lip(\RN;d)}=\|f\|_{\Lip(E;d)}$. This extension property of Lipschitz functions coincides with the statement of Theorem \reff{CMD-E-NEW} for $m=0$. (Furthermore, in this case the equivalence \rf{LMD-INF} is actually an equality).
%----------------------------------------------------------
\par We also note that, by the McShane extension formula, the function $F$ can be chosen in the form
%---------------------------------------------------------
\bel{Mc-FORM}
F(x)=\inf_{y\in E}\{f(y)+d(x,y)\},~~~x\in\RN.
\ee
%---------------------------------------------------------
\par This observation enables us to simplify considerably an almost optimal algorithm for extension of functions from the Sobolev space $L^1_p(\RN)$. See Remark \reff{ALG-MC}.}\rbx
%----------------------------------------------------------
\end{remark}
%----------------------------------------------------------
%@@@@@@@@@@@@@@@@@@@@@@@@@@@@@@@@@@@@@@@@@@@@@@@@@@@@@@@@@@
%@@@@@@@@@@@@@@@@@@@@@@@@@@@@@@@@@@@@@@@@@@@@@@@@@@@@@@@@@@
%@@@@@@@@@@@@@@@@@@@@@@@@@@@@@@@@@@@@@@@@@@@@@@@@@@@@@@@@@@
%@@@@@@@@@@@@@@@@@@@@@@@@@      @@@@@@@@@@@@@@@@@@@@@@@@@@@
%@@@@@@@@@@@@@@@@@@@@@@@          @@@@@@@@@@@@@@@@@@@@@@@@@
%@@@@@@@@@@@@@@@@@@@@@              @@@@@@@@@@@@@@@@@@@@@@@
%@@@@@@@@@@@@@@@@@@@     SECTION 5    @@@@@@@@@@@@@@@@@@@@@
%@@@@@@@@@@@@@@@@@@@@@              @@@@@@@@@@@@@@@@@@@@@@@
%@@@@@@@@@@@@@@@@@@@@@@@          @@@@@@@@@@@@@@@@@@@@@@@@@
%@@@@@@@@@@@@@@@@@@@@@@@@@      @@@@@@@@@@@@@@@@@@@@@@@@@@@
%@@@@@@@@@@@@@@@@@@@@@@@@@@@@@@@@@@@@@@@@@@@@@@@@@@@@@@@@@@
%@@@@@@@@@@@@@@@@@@@@@@@@@@@@@@@@@@@@@@@@@@@@@@@@@@@@@@@@@@
%@@@@@@@@@@@@@@@@@@@@@@@@@@@@@@@@@@@@@@@@@@@@@@@@@@@@@@@@@@
%----------------------------------------------------------
\SECT{5. Extensions of $\LMP$\,-jets: a proof of Theorem \reff{JET-S}.}{5}
%----------------------------------------------------------
\addtocontents{toc}{5. Extensions of $\LMP$\,-jets: a proof of Theorem \reff{JET-S}. \hfill\thepage\par\VST}
%----------------------------------------------------------
%@@@@@@@@@@@@@@@@@@@@@@@@@@@@@@@@@@@@@@@@@@@@@@@@@@@@@@@@@@
%----------------------------------------------------------
\indent\par {\it (Necessity.)} Let $\VP=\{P_x:x\in E\}$ be a Whitney $(m-1)$-field so that $P_x\in\PMRN$ for every $x\in E$. Let $F\in\LMP$ be a $C^{m-1}$-function such that $T^{m-1}_x[F]=P_x$ for all $x\in E$.
Prove that $\NCP\le C\|F\|_{\LMP}$ with $C=C(m,n,p)$.
%----------------------------------------------------------
\par Let $q=(p+n)/2$, and let $Q$ be a cube in $\RN$. Let $y,z\in Q\cap E$ and let $\tP:=P_y-P_z$. Let
%----------------------------------------------------------
$$
I_Q:=\,\left(\frac{1}{|Q|}
\intl_{Q}(\nabla^m F(u))^qdu\right)^{\frac{1}{q}}.
$$
%----------------------------------------------------------
Then, by inequality \rf{SP-BTQ}, for every multiindex $\alpha$, $|\alpha|\le m-1$, the following inequality
%----------------------------------------------------------
\bel{AL-P}
|D^{\alpha}\tP(y)|\,(\diam Q)^{|\alpha|-m}\le\, C\,I_Q
%----------------------------------------------------------
\ee
%----------------------------------------------------------
holds. Here $C=C(m,n,p)$.
%----------------------------------------------------------
\par Prove that for every $x\in Q$
%----------------------------------------------------------
\bel{AD-1}
|\tP(x)|\,(\diam Q)^{-m}\le\, C\,I_Q.
\ee
%----------------------------------------------------------
\par In fact, since $\tP\in \PMRN$,
%----------------------------------------------------------
$$
\tP(x)=\,\smed_{|\alpha|\le m-1}\,\,\frac{1}{\alpha!}\,D^\alpha\tP(y)\,(x-y)^\alpha
$$
%----------------------------------------------------------
so that
%----------------------------------------------------------
$$
|\tP(x)|\le C\,\smed_{|\alpha|\le m-1} \,\,|D^\alpha\tP(y)|\,\|x-y\|^{|\alpha|}\le
C\,\smed_{|\alpha|\le m-1} \,\,|D^\alpha\tP(y)|\,(\diam Q)^{|\alpha|}.
$$
%----------------------------------------------------------
\par Combining this inequality with \rf{AL-P}, we obtain the required inequality \rf{AD-1}.
%----------------------------------------------------------
\par  Let us apply this inequality to the cube $Q=Q(x,R)$ where $R:=\|x-y\|+\|x-z\|$. (Clearly, $Q\ni y,z$.) We obtain that for every $y,z\in E$ and every $x\in\RN$
%----------------------------------------------------------
\be
\frac{|P_y(x)-P_z(x)|}
{\|x-y\|^m+\|x-z\|^m}
&\le&\, C\,\left(\frac{1}{|Q(x,R)|}
\intl_{Q(x,R)}(\nabla^m F(u))^qdu\right)^{\frac{1}{q}}\nn\\
&\le&
\,C\,(\Mc[(\nabla^m F)^q])^\frac1q(x).\nn
%----------------------------------------------------------
\ee
%----------------------------------------------------------
Hence
%----------------------------------------------------------
$$
\VSH(x)\le\,C\,(\Mc[(\nabla^m F)^q])^\frac1q(x),
$$
%----------------------------------------------------------
see \rf{IEF}, so that
%----------------------------------------------------------
$$
\NCP\le\,C\,\|(\Mc[(\nabla^m F)^q])^\frac1q\|_{\LPRN}.
$$
%----------------------------------------------------------
Since $1<q<p$, by the Hardy-Littlewood maximal theorem,
%----------------------------------------------------------
$$
\NCP\le\,C\,\left(\,\,\intl_{\RN}\,
[(\nabla^m F)^q]) ^\frac{p}{q}\,du\right)^{\frac1p}
=\,C\,\|\nabla^m F\|_{\LPRN}\sim \|F\|_{\LMP}
$$
%----------------------------------------------------------
proving the necessity.\medskip
%----------------------------------------------------------
%@@@@@@@@@@@@@@@@@@@@@@@@@@@@@@@@@@@@@@@@@@@@@@@@@@@@@@@@@@
%@@@@@@@@@@@@@@@@@@@@@@@@@@@@@@@@@@@@@@@@@@@@@@@@@@@@@@@@@@
%@@@@@@@@@@@@@@@@@@@@@@@@@@@@@@@@@@@@@@@@@@@@@@@@@@@@@@@@@@
%@@@@@@@@@@@@@@@@@@@@@@@@@@@@@@@@@@@@@@@@@@@@@@@@@@@@@@@@@@
%----------------------------------------------------------
\par {\it(Sufficiency.)} Let $\VP=\{P_x:x\in E\}$ be a Whitney $(m-1)$-field such that $\VSH\in\LPRN$. Prove the existence of a function $F\in\LMP$ such that $T^{m-1}_x[F]=P_x$ for all $x\in E$, and
%----------------------------------------------------------
$$
%---------------------------------------------------------
\|F\|_{\LMP}\le C(m,n,p)\,\NCP.
%----------------------------------------------------------
$$
%----------------------------------------------------------
%@@@@@@@@@@@@@@@@@@@@@@@@@@@@@@@@@@@@@@@@@@@@@@@@@@@@@@@@@@
%----------------------------------------------------------
\par Let $x,y\in E$, and let $q:=(n+p)/2$. Let $h_1:=\VSH$ and let $\tP:=P_x-P_y$. Prove that for every multiindex $\alpha$, $|\alpha|\le m-1$,
%----------------------------------------------------------
\bel{DP-2}
|D^{\alpha}\tP(x)|
\le C\,\|x-y\|^{m-1-|\alpha|}\,\dq(x,y:h_1)
\ee
%---------------------------------------------------------
where $C=C(m,n,p)$. See definition \rf{DQ}.
%----------------------------------------------------------
\par Let $u\in Q_{xy}:=Q(x,\|x-y\|)$. Then $\|u-x\|,\,\|u-y\|\le 2\|x-y\|$, so that, by definition \rf{IEF}, the following inequality
%----------------------------------------------------------
$$
|\tP(u)|\le
(\|u-x\|^m+\|u-y\|^m)\,\VSH(u)
\le 4^m\,\|x-y\|^m\,\VSH(u)=4^m\,\|x-y\|^m\,h_1(u)
%----------------------------------------------------------
$$
%----------------------------------------------------------
holds. Hence
%----------------------------------------------------------
$$
|\tP(u)|^q\le
C\,\|x-y\|^{mq}\,h^q_1(u)~~\text{for every}~~u\in Q_{xy}.
$$
%----------------------------------------------------------
Integrating this inequality over $Q$ with respect to $u$ we obtain the following:
%----------------------------------------------------------
$$
\intl_{Q_{xy}}|\tP(u)|^q\,du\le \, C\,\|x-y\|^{mq}\,
\intl_{Q_{xy}}\,h_1^q(u)\,du.
$$
%----------------------------------------------------------
Hence
%----------------------------------------------------------
\be
\left(\frac{1}{|Q_{xy}|}\intl_{Q_{xy}}|\tP(u)|^q\,du\right)^
{\frac1q}&\le& \, C\,\|x-y\|^{m}\,
\left(\frac{1}{|Q_{xy}|}\intl_{Q_{xy}}\,h_1^q(u)\,du\right)^
{\frac1q}\nn\\
&\le& C\,\|x-y\|^{m}\sup_{Q\ni x,y}\left(\,\,
\frac{1}{|Q|}\intl_{Q}h_1^q(u)\,du\right)^{\frac1q}\nn\\
&=&\, C\|x-y\|^{m-1}\dq(x,y:h_1).\nn
\ee
%----------------------------------------------------------
\medskip
\par Since $\tP\in\PMRN$,
%----------------------------------------------------------
$$
\sup_{u\in Q_{xy}} |\tP(u)|\le C(m,n,q)\, \left(\frac{1}{|Q_{xy}|}\intl_{Q_{xy}}|\tP(u)|^q\,du\right)^
{\frac1q}.
$$
%----------------------------------------------------------
On the other hand, by Markov's inequality,
%----------------------------------------------------------
$$
\sup_{u\in Q_{xy}} |D^\alpha\tP(u)|\le\,
C(m,n)\,(\diam Q_{xy})^{-|\alpha|}\sup_{u\in Q_{xy}} |\tP(u)|\le C(m,n)\,\|x-y\|^{-|\alpha|}\sup_{u\in Q_{xy}} |\tP(u)|.
$$
%----------------------------------------------------------
Hence,
%----------------------------------------------------------
\be
|D^\alpha\tP(x)|&\le& \sup_{u\in Q_{xy}} |D^\alpha\tP(u)|\le\, C(m,n)\,\|x-y\|^{-|\alpha|}\left(\frac{1}{|Q_{xy}|}
\intl_{Q_{xy}}|\tP(u)|^q\,du\right)^
{\frac1q}\nn\\
&\le& C\,\|x-y\|^{m-|\alpha|-1}\dq(x,y:h_1)\nn
\ee
%---------------------------------------------------------
proving \rf{DP-2}.\medskip
%---------------------------------------------------------
\par This inequality and Theorem \reff{MAIN-MT} imply the following: %----------------------------------------------------------
\bel{R-FS}
|D^{\alpha}P_x(x)-D^{\alpha}P_y(x)|
\le C\,\|x-y\|^{m-1-|\alpha|}\, d_q(x,y:h_1).
\ee
%----------------------------------------------------------
\par By definition \rf{D-LIP}, the metric $d=d_q(h_1)$ belongs to the family $\DPQ$. Furthermore, by the latter inequality,
%---------------------------------------------------------
$$
\Lc_{m-1,d}(\VP):=\sbig_{|\alpha|\le m-1}\,\,\,
\sup_{x,y\in E,\,x\ne y}\,\, \frac{|D^{\alpha}P_x(x)-D^{\alpha}P_y(x)|}
{\|x-y\|^{m-1-|\alpha|}\,d(x,y)}\le C=C(m,n,p),
$$
%---------------------------------------------------------
so that, by Theorem \reff{CMD-E-NEW}, there exists a function $F\in C^{m-1,(d)}(\RN)$ which agrees with $\VP$ on $E$ (i.e., $T^{m-1}_x[F]=P_x$ for all $x\in E$) and
%---------------------------------------------------------
$$
\|F\|_{C^{m-1,(d)}(\RN)}\le C(m,n,p)\,\Lc_{m-1,d}(\VP)\le C(m,n,p).
$$
%---------------------------------------------------------
\par Thus,
%----------------------------------------------------------
$$
|D^\alpha F(x)-D^\alpha F(y)|\le d_q(x,y:h_2)~~\text{for all}~~\alpha,|\alpha|=m-1,
$$
%----------------------------------------------------------
where $h_2=C(m,n,p)\,h_1$, so that, by Theorem \reff{SM-MAIN},
%----------------------------------------------------------
\bel{R-FPS}
\|F\|_{\LMP}\le C\,\|h_2\|_{\LPRN}\le C \|h_1\|_{\LPRN}= C\,\|\VSH\|_{\LPRN}.
\ee
%---------------------------------------------------------
\par Theorem \reff{JET-S} is completely proved.\bx
%----------------------------------------------------------
%@@@@@@@@@@@@@@@@@@@@@@@@@@@@@@@@@@@@@@@@@@@@@@@@@@@@@@@@@@
%@@@@@@@@@@@@@@@@@@@@@@@@@@@@@@@@@@@@@@@@@@@@@@@@@@@@@@@@@@
%@@@@@@@@@@@@@@@@@@@@@@@@@@@@@@@@@@@@@@@@@@@@@@@@@@@@@@@@@@
%----------------------------------------------------------
\medskip
%----------------------------------------------------------
%@@@@@@@@@@@@@@@@@@@@@@@@@@@@@@@@@@@@@@@@@@@@@@@@@@@@@@@@@@
%@@@@@@@@@@@@@@@@@@@@@@@@@@@@@@@@@@@@@@@@@@@@@@@@@@@@@@@@@@
%@@@@@@@@@@@@@@@@@@@@@@@@@@@@@@@@@@@@@@@@@@@@@@@@@@@@@@@@@@
%@@@@@@@@@@@@@@@@@@@@@@@@@@@@@@@@@@@@@@@@@@@@@@@@@@@@@@@@@@
%@@@@@@@@@@@@@@@@@@@@@@@@@@@@@@@@@@@@@@@@@@@@@@@@@@@@@@@@@@
%@@@@@@@@@@@@@@@@@@@@@@@@@@@@@@@@@@@@@@@@@@@@@@@@@@@@@@@@@@
%----------------------------------------------------------
\begin{remark}\lbl{ALG-MC} {\em Note that Theorem \reff{CMD-E-NEW} is the main ingredient of the proof of the sufficiency presented above. As we have noted in Remark \reff{MC-EXT-REM}, for the space $\LOP$ we can replace Theorem  \reff{CMD-E-NEW} with the McShane's extension formula \rf{Mc-FORM} which simplifies essentially the proof of the sufficiency part of Theorem \reff{JET-S}.
%---------------------------------------------------------
\par Actually we obtain a new (non-linear) algorithm for extension of $\LOP$-functions whenever $p>n$. Let us describe the main steps of this algorithm.
%----------------------------------------------------------
\par Let $q=(n+p)/2$ and let $f$ be a continuous function defined on $E$.\smallskip
%----------------------------------------------------------
\par \textit{Step 1.} We introduce the ``sharp maximal function'' associated with $f$ by
%----------------------------------------------------------
$$
\fs(x):=\sup_{y,z\in E}\frac{|f(x)-f(y)|}{\|x-y\|+\|y-z\|},~~~x\in\RN.
$$
%----------------------------------------------------------
\smallskip
\par \textit{Step 2.} We introduce a ``pre-metric'' $\dq(\fs)$ associated with the function $\fs$, i.e., a function
%---------------------------------------------------------
$$
\dq(x,y:\fs)=\|x-y\| \sup_{Q\ni x,y} \left(\frac{1}{|Q|}\intl_{Q} (\fs)^q(u)\,du\right)^{\frac{1}{q}}~,~~~~~x,y\in \RN.
$$
%----------------------------------------------------------
See \rf{DQ}.
%----------------------------------------------------------
%@@@@@@@@@@@@@@@@@@@@@@@@@@@@@@@@@@@@@@@@@@@@@@@@@@@@@@@@@@
%----------------------------------------------------------
\smallskip
\par \textit{Step 3.} Using formula \rf{GD-M} we construct the geodesic distance $d=d_q(x,y:\fs)$ associated with the ``pre-metric'' $\dq(\fs)$:
%----------------------------------------------------------
$$
d_q(x,y:\fs):=\inf\,\smed_{i=0}^{m-1}\,\dq(x_i,x_{i+1}:\fs)
$$
%---------------------------------------------------------
where the infimum is taken over all finite sequences of points $\{x_0,x_1,...,x_m\}$ in $\RN$ such that $x_0=x$  and $x_m=y$.
%----------------------------------------------------------
\smallskip
\par \textit{Step 4.} Using the McShane's formula \rf{Mc-FORM} we construct a function
%----------------------------------------------------------
\bel{MC-EXT}
F(x):=\inf_{y\in E}\left\{f(y)+48\,d_q(x,y:\fs)\right\},~~~~x\in\RN.
\ee
%----------------------------------------------------------
%@@@@@@@@@@@@@@@@@@@@@@@@@@@@@@@@@@@@@@@@@@@@@@@@@@@@@@@@@@
%----------------------------------------------------------
\smallskip
\par Repeating the proof of the sufficiency part of Theorem \reff{JET-S} (for $m=1$), we are able to show that the function $F$ provides an almost optimal extension of the function $f$ to a function from the Sobolev space $\LOP$ provided $f\in\LOP|_E$. In fact, by \rf{R-FS},
%----------------------------------------------------------
$$
|f(x)-f(y)|
\le 48\, d_q(x,y:\fs)
$$
%----------------------------------------------------------
for every $x,y\in\RN$, so that, by the McShane's formula, the function $F$ is {\it an extension of $f$} from $E$ on all of $\RN$. Following the proof of the sufficiency part of Theorem \reff{JET-S} we conclude that inequality \rf{R-FPS} holds whenever $m=1$ and $F$ is defined by \rf{MC-EXT}. Thus $\|F\|_{\LOP}\le C\,\|\fs\|_{\LPRN}$.
%----------------------------------------------------------
\par On the other hand, by \rf{EQ-L1},
$\|\fs\|_{\LPRN}\sim \|f\|_{\LOP|_E}$. Hence  $\|F\|_{\LOP}\le C\,\|f\|_{\LOP|_E}$ proving that $F$ provides an almost optimal extension of $f$ to a function from $\LOP$.\rbx}
%----------------------------------------------------------
\end{remark}
%----------------------------------------------------------
%@@@@@@@@@@@@@@@@@@@@@@@@@@@@@@@@@@@@@@@@@@@@@@@@@@@@@@@@@@
%@@@@@@@@@@@@@@@@@@@@@@@@@@@@@@@@@@@@@@@@@@@@@@@@@@@@@@@@@@
%@@@@@@@@@@@@@@@@@@@@@@@@@@@@@@@@@@@@@@@@@@@@@@@@@@@@@@@@@@
%@@@@@@@@@@@@@@@@@@@@@@@@@      @@@@@@@@@@@@@@@@@@@@@@@@@@@
%@@@@@@@@@@@@@@@@@@@@@@@          @@@@@@@@@@@@@@@@@@@@@@@@@
%@@@@@@@@@@@@@@@@@@@@@              @@@@@@@@@@@@@@@@@@@@@@@
%@@@@@@@@@@@@@@@@@@@     SECTION 6    @@@@@@@@@@@@@@@@@@@@@
%@@@@@@@@@@@@@@@@@@@@@              @@@@@@@@@@@@@@@@@@@@@@@
%@@@@@@@@@@@@@@@@@@@@@@@          @@@@@@@@@@@@@@@@@@@@@@@@@
%@@@@@@@@@@@@@@@@@@@@@@@@@      @@@@@@@@@@@@@@@@@@@@@@@@@@@
%@@@@@@@@@@@@@@@@@@@@@@@@@@@@@@@@@@@@@@@@@@@@@@@@@@@@@@@@@@
%@@@@@@@@@@@@@@@@@@@@@@@@@@@@@@@@@@@@@@@@@@@@@@@@@@@@@@@@@@
%@@@@@@@@@@@@@@@@@@@@@@@@@@@@@@@@@@@@@@@@@@@@@@@@@@@@@@@@@@
%----------------------------------------------------------
\SECT{6. Lacunae of Whitney's cubes and a lacunary extension operator.}{6}
%----------------------------------------------------------
\addtocontents{toc}{6. Lacunae of Whitney's cubes and a lacunary extension operator.\hfill\thepage\par\VST}
%----------------------------------------------------------
%@@@@@@@@@@@@@@@@@@@@@@@@@@@@@@@@@@@@@@@@@@@@@@@@@@@@@@@@@@
\indent\par We prove the sufficiency part of Theorem \reff{EX-TK} with the help of a modification of the classical Whitney extension method \cite{W1} used in the author's paper \cite{Sh3}. As we have noted in Introduction, the main idea of this approach is that, instead of treating each Whitney cube separately, as is done in \cite{W1}, we deal simultaneously with all members of certain {\it families} of Whitney cubes.
We refer to these families of Whitney cubes as {\it lacunae}.
%----------------------------------------------------------
\par In Subsection 6.1 we present main definitions related to this notion, and main properties of the lacunae. For the proof of these results we refer the reader to \cite{Sh3}, Sections 4-5.\bigskip
%----------------------------------------------------------
\par {\bf 6.1. Lacunae of Whitney's cubes.}\medskip
%----------------------------------------------------------
\addtocontents{toc}{~~~~6.1. Lacunae of Whitney's cubes. \hfill \thepage\par}
%----------------------------------------------------------
\par Let $E$ be a closed subset of $\RN$ and let $W_E$ be a Whitney covering of its complement $\RN\setminus E$ satisfying inequality \rf{DQ-E}. This inequality implies the following property of Whitney's cubes:
%----------------------------------------------------------
\bel{9Q-INT}
(9Q)\cap E\ne\emp~~~ \text{for every}~~~ Q\in W_E.
\ee
%----------------------------------------------------------
%@@@@@@@@@@@@@@@@@@@@@@@@@@@@@@@@@@@@@@@@@@@@@@@@@@@@@@@@@@
%@@@@@@@@@@@@@@@@@@@@@@@@@@@@@@@@@@@@@@@@@@@@@@@@@@@@@@@@@@
%----------------------------------------------------------
\par By $LW_E$ we denote a subfamily of Whitney cubes satisfying the following condition:
%----------------------------------------------------------
\bel{L-PR}
(10Q)\cap E=(\q Q)\cap E.
\ee
%----------------------------------------------------------
\par Then we introduce a binary relation $\sim$ on $LW_E$: for every $Q_1,Q_2\in LW_E$
%----------------------------------------------------------
\bel{BIN-R}
Q_1\sim Q_2 ~~\Longleftrightarrow~ (10Q_1)\cap E= (10Q_2)\cap E.
\ee
%----------------------------------------------------------
\par It can be easily seen that $\sim$ satisfies the axioms of equivalence relations, i.e., it is reflexive, symmetric and transitive. Given a cube $Q\in LW_E$ by
%----------------------------------------------------------
$$
[Q]:=\{K\in LW_E: K\sim Q\}
$$
%----------------------------------------------------------
we denote the equivalence class of $Q$. We refer to this equivalence class as {\it a true lacuna} with respect to the set $E$.
%----------------------------------------------------------
\par Let
%----------------------------------------------------------
$$
\tL_E=LW_E/\sim\,=\{[Q]: Q\in LW_E\}
$$
%----------------------------------------------------------
be the corresponding quotient set of $LW_E$ by $\sim$\,, i.e., the set of all possible equivalence classes (lacunae) of $LW_E$ by $\sim$\,.
%----------------------------------------------------------
\par Thus for every pair of Whitney cubes $Q_1,Q_2\in W_E$ which belong to a true lacuna $L\in\tL_E$ we have
%----------------------------------------------------------
\bel{I-L}
(10Q_1)\cap E=(\q Q_1)\cap E=(10Q_2)\cap E=(\q Q_2)\cap E.
\ee
%----------------------------------------------------------
By $V_L$ we denote the associated set of the lacuna $L$
%----------------------------------------------------------
\bel{D-VL}
V_L:=(\q Q)\cap E.
\ee
%----------------------------------------------------------
Here $Q$ is an arbitrary cube from $L$. By \rf{I-L}, any choice of a cube $Q\in L$ provides the same set $V_L$ so that $V_L$ is well-defined. We also note that for each cube $Q$ which belong to a true lacuna $L$ we have $V_L=(10Q)\cap E.$\medskip
%----------------------------------------------------------
\par We extend the family $\tL_E$ of true lacunae to a family $\LE$ of {\it all lacunae} in the following way. Suppose that $Q\in W_E\setminus LW_E$, see \rf{L-PR}, i.e.,
%----------------------------------------------------------
$$
(10Q)\cap E\ne(\q Q)\cap E.
$$
%----------------------------------------------------------
In this case to the cube $Q$ we assign a lacuna $L:=\{Q\}$ consisting of a unique cube - the cube $Q$ itself.  We also put $V_L:=(\q Q)\cap E$ as in \rf{D-VL}.
%----------------------------------------------------------
\par We refer to such a lacuna $L:=\{Q\}$ as an {\it elementary lacuna} with respect to the set $E$. By $\hL_E$ we denote the family of all elementary lacunae with respect to $E$:
%----------------------------------------------------------
\bel{EL-LC}
\hL_E:=\{L=\{Q\}:Q\in W_E\setminus LW_E\}.
\ee
%----------------------------------------------------------
\par In \cite{Sh3} we prove that for every elementary lacuna $L=\{Q\}\in \hL_E$ we have
%----------------------------------------------------------
\bel{A-DL}
\diam Q\le 2\diam V_L=2 \diam ((\q Q)\cap E).
\ee
%----------------------------------------------------------
\par Finally, by $\LE$ we denote the family of all lacunae with respect to $E$:
%----------------------------------------------------------
$$
\LE=\tL_E\cup \hL_E.
$$
%----------------------------------------------------------
%\smallskip
%----------------------------------------------------------
%@@@@@@@@@@@@@@@@@@@@@@@@@@@@@@@@@@@@@@@@@@@@@@@@@@@@@@@@@@
%@@@@@@@@@@@@@@@@@@@@@@@@@@@@@@@@@@@@@@@@@@@@@@@@@@@@@@@@@@
%@@@@@@@@@@@@@@@@@@@@@@@@@@@@@@@@@@@@@@@@@@@@@@@@@@@@@@@@@@
%@@@@@@@@@@@@@@@@@@@@@@@@@@@@@@@@@@@@@@@@@@@@@@@@@@@@@@@@@@
\par Let us present several important properties of lacunae.
%----------------------------------------------------------
%@@@@@@@@@@@@@@@@@@@@@@@@@@@@@@@@@@@@@@@@@@@@@@@@@@@@@@@@@@
%@@@@@@@@@@@@@@@@@@@@@@@@@@@@@@@@@@@@@@@@@@@@@@@@@@@@@@@@@@
%@@@@@@@@@@@@@@@@@@@@@@@@@@@@@@@@@@@@@@@@@@@@@@@@@@@@@@@@@@
%@@@@@@@@@@@@@@@@@@@@@@@@@@@@@@@@@@@@@@@@@@@@@@@@@@@@@@@@@@
%@@@@@@@@@@@@@@@@@@@@@@@@@@@@@@@@@@@@@@@@@@@@@@@@@@@@@@@@@@
%@@@@@@@@@@@@@@@@@@@@@@@@@@@@@@@@@@@@@@@@@@@@@@@@@@@@@@@@@@
%@@@@@@@@@@@@@@@@@@@@@@@@@@@@@@@@@@@@@@@@@@@@@@@@@@@@@@@@@@
%----------------------------------------------------------
\par Let $L\in\LE$ and let $U_L:=\cup\{Q:Q\in L\}$. We say that $L$ is a {\it bounded} lacuna if $U_L$ is a bounded set.
%----------------------------------------------------------
\par In \cite{Sh3}, Section 4, we prove that
for every lacuna $L\in\LE$ there exists a cube $Q_L\in L$ such that
%----------------------------------------------------------
$$
\diam Q_L=\inf\,\{\diam Q: Q\in L\}
$$
%----------------------------------------------------------
provided $\diam V_L>0$. We also prove that for every {\it bounded} lacuna $L\in\LE$ there exists a cube  $\QL\in L$ such that
%----------------------------------------------------------
\bel{DF-QLU}
\diam \QL=\sup\,\{\diam K: K\in L\}.
\ee
%----------------------------------------------------------
%@@@@@@@@@@@@@@@@@@@@@@@@@@@@@@@@@@@@@@@@@@@@@@@@@@@@@@@@@@
%@@@@@@@@@@@@@@@@@@@@@@@@@@@@@@@@@@@@@@@@@@@@@@@@@@@@@@@@@@
%@@@@@@@@@@@@@@@@@@@@@@@@@@@@@@@@@@@@@@@@@@@@@@@@@@@@@@@@@@
%@@@@@@@@@@@@@@@@@@@@@@@@@@@@@@@@@@@@@@@@@@@@@@@@@@@@@@@@@@
%----------------------------------------------------------
\par For the proof of the following four propositions we refer the reader to \cite{Sh3}, Sections 4 and 5.
%----------------------------------------------------------
%@@@@@@@@@@@@@@@@@@@@@@@@@@@@@@@@@@@@@@@@@@@@@@@@@@@@@@@@@@
%@@@@@@@@@@@@@@@@@@@@@@@@@@@@@@@@@@@@@@@@@@@@@@@@@@@@@@@@@@
%@@@@@@@@@@@@@@@@@@@@@@@@@@@@@@@@@@@@@@@@@@@@@@@@@@@@@@@@@@
%@@@@@@@@@@@@@@@@@@@@@@@@@@@@@@@@@@@@@@@@@@@@@@@@@@@@@@@@@@
%----------------------------------------------------------
\begin{proposition}\lbl{P-2L} (i). If $E$ is an unbounded set, then every lacuna $L\in\LE$ is bounded;\smallskip
%----------------------------------------------------------
\par (ii). If $E$ is bounded, then there exists the unique unbounded lacuna $L^{\max}\in\LE$. The lacuna $L^{\max}$ is a true lacuna for which $V_{L^{\max}}=E$.
%----------------------------------------------------------
\end{proposition}
%----------------------------------------------------------
%@@@@@@@@@@@@@@@@@@@@@@@@@@@@@@@@@@@@@@@@@@@@@@@@@@@@@@@@@@
%@@@@@@@@@@@@@@@@@@@@@@@@@@@@@@@@@@@@@@@@@@@@@@@@@@@@@@@@@@
%@@@@@@@@@@@@@@@@@@@@@@@@@@@@@@@@@@@@@@@@@@@@@@@@@@@@@@@@@@
%@@@@@@@@@@@@@@@@@@@@@@@@@@@@@@@@@@@@@@@@@@@@@@@@@@@@@@@@@@
%----------------------------------------------------------
\begin{proposition}\lbl{P-SM} \par (i). For every bounded true lacuna $L$ the following inequalities
%----------------------------------------------------------
\bel{YY}
40\,\diam \QL\le \dist(V_L,E\setminus V_L)\le \gamma_1 \diam\QL
\ee
%----------------------------------------------------------
hold;\smallskip
%----------------------------------------------------------
%@@@@@@@@@@@@@@@@@@@@@@@@@@@@@@@@@@@@@@@@@@@@@@@@@@@@@@@@@@
%----------------------------------------------------------
\par (ii). For every lacuna $L$ with $\diam V_L>0$
%----------------------------------------------------------
\bel{YG}
\diam Q_L\le \gamma_1 \diam V_L\,.
\ee
%----------------------------------------------------------
\par Here $\gamma_1>0$ is an absolute constant.
%----------------------------------------------------------
\end{proposition}
%----------------------------------------------------------
%@@@@@@@@@@@@@@@@@@@@@@@@@@@@@@@@@@@@@@@@@@@@@@@@@@@@@@@@@@
%@@@@@@@@@@@@@@@@@@@@@@@@@@@@@@@@@@@@@@@@@@@@@@@@@@@@@@@@@@
%@@@@@@@@@@@@@@@@@@@@@@@@@@@@@@@@@@@@@@@@@@@@@@@@@@@@@@@@@@
%@@@@@@@@@@@@@@@@@@@@@@@@@@@@@@@@@@@@@@@@@@@@@@@@@@@@@@@@@@
%@@@@@@@@@@@@@@@@@@@@@@@@@@@@@@@@@@@@@@@@@@@@@@@@@@@@@@@@@@
%@@@@@@@@@@@@@@@@@@@@@@@@@@@@@@@@@@@@@@@@@@@@@@@@@@@@@@@@@@
%----------------------------------------------------------
\begin{proposition}\lbl{M-LAC} Let $L\in\LE$ be a lacuna and let $Q\in L$. Suppose that there exist a lacuna $L'\in \LE$, $L\ne L'$, and a cube $Q'\in L'$ such that $Q\cap Q'\ne\emp$. Then:\smallskip
%----------------------------------------------------------
\par (i). If $L$ is a true lacuna, then $L'$ is an elementary lacuna\,;\smallskip
%----------------------------------------------------------
\par (ii). Either $\diam \QL\le \tau\,\diam Q$
or $\diam Q\le \tau\,\diam Q_L$ where $\tau$ is a positive absolute constant.
%@@@@@@@@@@@@@@@@@@@@@@@@@@@@@@@@@@@@@@@@@@@@@@@@@@@@@@@@@@
%----------------------------------------------------------
\end{proposition}
%----------------------------------------------------------
\medskip
%----------------------------------------------------------
%@@@@@@@@@@@@@@@@@@@@@@@@@@@@@@@@@@@@@@@@@@@@@@@@@@@@@@@@@@
%@@@@@@@@@@@@@@@@@@@@@@@@@@@@@@@@@@@@@@@@@@@@@@@@@@@@@@@@@@
%@@@@@@@@@@@@@@@@@@@@@@@@@@@@@@@@@@@@@@@@@@@@@@@@@@@@@@@@@@
%----------------------------------------------------------
\par This proposition motivates us to introduce the following
%@@@@@@@@@@@@@@@@@@@@@@@@@@@@@@@@@@@@@@@@@@@@@@@@@@@@@@@@@@
%@@@@@@@@@@@@@@@@@@@@@@@@@@@@@@@@@@@@@@@@@@@@@@@@@@@@@@@@@@
%@@@@@@@@@@@@@@@@@@@@@@@@@@@@@@@@@@@@@@@@@@@@@@@@@@@@@@@@@@
%@@@@@@@@@@@@@@@@@@@@@@@@@@@@@@@@@@@@@@@@@@@@@@@@@@@@@@@@@@
%----------------------------------------------------------
\begin{definition} \lbl{CONT-L}{\em Let $L,L'\in\LE$ be lacunae, $L\ne L'$. We say that $L$ and $L'$ are {\it contacting lacunae} if there exist cubes $Q\in L$ and $Q'\in L'$ such that $Q\cap Q'\ne\emp$. We refer to the pair of such cubes as {\it contacting cubes}.
%----------------------------------------------------------
\par We write
%----------------------------------------------------------
$$
L\lcr L' ~~\text{for contacting lacunae}~~L,L'\in\LE.
$$
%----------------------------------------------------------
}
%----------------------------------------------------------
\end{definition}
%----------------------------------------------------------
%@@@@@@@@@@@@@@@@@@@@@@@@@@@@@@@@@@@@@@@@@@@@@@@@@@@@@@@@@@
%@@@@@@@@@@@@@@@@@@@@@@@@@@@@@@@@@@@@@@@@@@@@@@@@@@@@@@@@@@
%@@@@@@@@@@@@@@@@@@@@@@@@@@@@@@@@@@@@@@@@@@@@@@@@@@@@@@@@@@
%----------------------------------------------------------
\begin{proposition}\lbl{D-LQN} Every lacuna $L\in\LE$ contacts with at most $M$ lacunae. Furthermore, $L$ contains at most $M$ contacting cubes. Here $M=M(n)$ is a positive integer depending only on $n$.
%----------------------------------------------------------
\end{proposition}\medskip
%----------------------------------------------------------
%@@@@@@@@@@@@@@@@@@@@@@@@@@@@@@@@@@@@@@@@@@@@@@@@@@@@@@@@@@
%@@@@@@@@@@@@@@@@@@@@@@@@@@@@@@@@@@@@@@@@@@@@@@@@@@@@@@@@@@
%@@@@@@@@@@@@@@@@@@@@@@@@@@@@@@@@@@@@@@@@@@@@@@@@@@@@@@@@@@
%@@@@@@@@@@@@@@@@@@@@@@@@@@@@@@@@@@@@@@@@@@@@@@@@@@@@@@@@@@
%@@@@@@@@@@@@@@@@@@@@@@@@@@@@@@@@@@@@@@@@@@@@@@@@@@@@@@@@@@
%@@@@@@@@@@@@@@@@@@@@@@@@@@@@@@@@@@@@@@@@@@@@@@@@@@@@@@@@@@
%@@@@@@@@@@@@@@@@@@@@@@@@@@@@@@@@@@@@@@@@@@@@@@@@@@@@@@@@@@
%----------------------------------------------------------
\par {\bf 6.2. A lacunary projector and centers of lacunae. }\medskip
%----------------------------------------------------------
\addtocontents{toc}{~~~~6.2. A lacunary projector and centers of lacunae. \hfill \thepage\par}
%----------------------------------------------------------
%@@@@@@@@@@@@@@@@@@@@@@@@@@@@@@@@@@@@@@@@@@@@@@@@@@@@@@@@@@
%----------------------------------------------------------
\par One of the main ingredient of the lacunary approach is a mapping $\PRL:\LE\to E$ whose properties are described in Theorem \reff{L-PE} below. We refer to this mapping as a ``projector" from $\Lc_E$ into the set $E$. Also given $L\in\LE$ we refer to the point $\PRL(L)\in E$ as {\it a center} of the lacuna $L$.
%----------------------------------------------------------
\par Theorem \reff{L-PE} is a refinement of a result proven in \cite{Sh3}, Section 5.
%----------------------------------------------------------
%@@@@@@@@@@@@@@@@@@@@@@@@@@@@@@@@@@@@@@@@@@@@@@@@@@@@@@@@@@
%@@@@@@@@@@@@@@@@@@@@@@@@@@@@@@@@@@@@@@@@@@@@@@@@@@@@@@@@@@
%@@@@@@@@@@@@@@@@@@@@@@@@@@@@@@@@@@@@@@@@@@@@@@@@@@@@@@@@@@
%@@@@@@@@@@@@@@@@@@@@@@@@@@@@@@@@@@@@@@@@@@@@@@@@@@@@@@@@@@
%@@@@@@@@@@@@@@@@@@@@@@@@@@@@@@@@@@@@@@@@@@@@@@@@@@@@@@@@@@
%@@@@@@@@@@@@@@@@@@@@@@@@@@@@@@@@@@@@@@@@@@@@@@@@@@@@@@@@@@
%----------------------------------------------------------
\begin{theorem}\lbl{L-PE} There exist an absolute constant $\tgm\ge 1$ and a mapping $\PRL:\LE\to E$ such that:
%----------------------------------------------------------
\par (i). For every lacuna $L\in\LE$ and every cube $Q\in L$ we have
%----------------------------------------------------------
\bel{PR-GM}
\PRL(L)\in (\tgm\,Q)\cap E~;
\ee
%----------------------------------------------------------
%\smallskip
%----------------------------------------------------------
%@@@@@@@@@@@@@@@@@@@@@@@@@@@@@@@@@@@@@@@@@@@@@@@@@@@@@@@@@@
%@@@@@@@@@@@@@@@@@@@@@@@@@@@@@@@@@@@@@@@@@@@@@@@@@@@@@@@@@@
%@@@@@@@@@@@@@@@@@@@@@@@@@@@@@@@@@@@@@@@@@@@@@@@@@@@@@@@@@@
%----------------------------------------------------------
\par (ii). Let $L,L'\in\LE$ be two distinct lacunae such that $\PRL(L)\ne \PRL(L')$. Then for every two cubes $Q\in L$ and $Q'\in L'$ such that $Q\cap Q'\ne\emp$ the following inequality
%----------------------------------------------------------
\bel{DM-QQP}
\diam Q+\diam Q'\le \tgm\,\|\PRL(L)-\PRL(L')\|
\ee
%----------------------------------------------------------
holds;
%----------------------------------------------------------
\smallskip
%----------------------------------------------------------
%@@@@@@@@@@@@@@@@@@@@@@@@@@@@@@@@@@@@@@@@@@@@@@@@@@@@@@@@@@
%@@@@@@@@@@@@@@@@@@@@@@@@@@@@@@@@@@@@@@@@@@@@@@@@@@@@@@@@@@
%@@@@@@@@@@@@@@@@@@@@@@@@@@@@@@@@@@@@@@@@@@@@@@@@@@@@@@@@@@
%----------------------------------------------------------
\par (iii). For every point $A\in E$
%----------------------------------------------------------
$$
\#\{L\in\Lc_E:\PRL(L)=A\}\le C
$$
%----------------------------------------------------------
where $C=C(n)$ is a constant depending only on $n$.
%@@@@@@@@@@@@@@@@@@@@@@@@@@@@@@@@@@@@@@@@@@@@@@@@@@@@@@@@@@
%----------------------------------------------------------
\end{theorem}
%----------------------------------------------------------
%@@@@@@@@@@@@@@@@@@@@@@@@@@@@@@@@@@@@@@@@@@@@@@@@@@@@@@@@@@
%@@@@@@@@@@@@@@@@@@@@@@@@@@@@@@@@@@@@@@@@@@@@@@@@@@@@@@@@@@
%----------------------------------------------------------
\par {\it Proof.} We define the projector $\PRL:\LE\to E$ in several steps. First consider a {\it true lacuna} $L\in\LE$ such that $\diam V_L=0$, i.e., $V_L=\{a\}$ for some $a\in E$. In this case we put $\PRL(L):=a$.\medskip
%----------------------------------------------------------
\par Let now $\diam V_L>0$. Fix points $A_L,B_L\in V_L$ such that
%----------------------------------------------------------
$$
\|A_L-B_L\|=\diam V_L.
$$
%----------------------------------------------------------
By $i_L\in\mZ$ we denote an integer such that
%----------------------------------------------------------
\bel{IL-VL}
2^{i_L}<\diam V_L=\|A_L-B_L\|\le 2^{i_L+1}.
\ee
%----------------------------------------------------------
%@@@@@@@@@@@@@@@@@@@@@@@@@@@@@@@@@@@@@@@@@@@@@@@@@@@@@@@@@@
%@@@@@@@@@@@@@@@@@@@@@@@@@@@@@@@@@@@@@@@@@@@@@@@@@@@@@@@@@@
%@@@@@@@@@@@@@@@@@@@@@@@@@@@@@@@@@@@@@@@@@@@@@@@@@@@@@@@@@@
%@@@@@@@@@@@@@@@@@@@@@@@@@@@@@@@@@@@@@@@@@@@@@@@@@@@@@@@@@@
%----------------------------------------------------------
\par In what follows we will be needed a result proven in \cite{Sh3}, Section 4, which states the following: There exists a {\it non-increasing sequence of non-empty closed sets} $\{E_i\}_{i\in\mZ}$, $E_{i+1}\subset E_i\subset E$, $i\in\mZ,$ such that for every $i\in\mZ$ the following conditions are satisfied:\medskip
%----------------------------------------------------------
\par (i). The points of the set $E_i$ are $2^i$-separated, i.e.,
%----------------------------------------------------------
\bel{SEP-E}
\|z-z'\|\ge 2^i~~~\text{for every}~~~z,z'\in E_i\,;
\ee
%----------------------------------------------------------
%@@@@@@@@@@@@@@@@@@@@@@@@@@@@@@@@@@@@@@@@@@@@@@@@@@@@@@@@@@
%----------------------------------------------------------
\par (ii). $E_i$ is a $2^{i+1}$-net in $E$, i.e., %----------------------------------------------------------
\bel{AP-E}
\text{for every}~~x\in E~~\text{there exists}~~z\in E_i~~\text{such that}~~\|x-z\|\le 2^{i+1}\,.
\ee
%----------------------------------------------------------
\smallskip
%----------------------------------------------------------
%@@@@@@@@@@@@@@@@@@@@@@@@@@@@@@@@@@@@@@@@@@@@@@@@@@@@@@@@@@
%@@@@@@@@@@@@@@@@@@@@@@@@@@@@@@@@@@@@@@@@@@@@@@@@@@@@@@@@@@
%@@@@@@@@@@@@@@@@@@@@@@@@@@@@@@@@@@@@@@@@@@@@@@@@@@@@@@@@@@
%@@@@@@@@@@@@@@@@@@@@@@@@@@@@@@@@@@@@@@@@@@@@@@@@@@@@@@@@@@
%@@@@@@@@@@@@@@@@@@@@@@@@@@@@@@@@@@@@@@@@@@@@@@@@@@@@@@@@@@
%----------------------------------------------------------
\par Let us apply this statement to the points $A_L$ and $B_L$. Since $E_{i_L-2}$ is a $2^{i_L-1}$-net in $E$, there exist points $\tA_L,\tB_L\in E_{i_L-2}$ such that
%----------------------------------------------------------
\bel{Y5}
\|A_L-\tA_L\|\le 2^{i_L-1}~~~\text{and}~~~\|B_L-\tB_L\|\le 2^{i_L-1}.
\ee
%----------------------------------------------------------
Since $\tA_L,\tB_L\in E_{i_L-2}$, by \rf{AP-E}, $\|\tA_L-\tB_L\|\ge 2^{i_L-2}$.\medskip
%----------------------------------------------------------
%@@@@@@@@@@@@@@@@@@@@@@@@@@@@@@@@@@@@@@@@@@@@@@@@@@@@@@@@@@
%@@@@@@@@@@@@@@@@@@@@@@@@@@@@@@@@@@@@@@@@@@@@@@@@@@@@@@@@@@
%@@@@@@@@@@@@@@@@@@@@@@@@@@@@@@@@@@@@@@@@@@@@@@@@@@@@@@@@@@
%@@@@@@@@@@@@@@@@@@@@@@@@@@@@@@@@@@@@@@@@@@@@@@@@@@@@@@@@@@
%@@@@@@@@@@@@@@@@@@@@@@@@@@@@@@@@@@@@@@@@@@@@@@@@@@@@@@@@@@
%----------------------------------------------------------
\par Prove that $\tA_L\ne\tB_L$ and
%----------------------------------------------------------
\bel{GW-1}
\{\tA_L,\tB_L\}\cap(E_{i_L-2}\setminus E_{i_L+2})\ne\emp.
\ee
%----------------------------------------------------------
\par In fact, by \rf{Y5},
%----------------------------------------------------------
$$
\|A_L-B_L\|\le \|A_L-\tA_L\|+\|\tA_L-\tB_L\|+\|\tB_L-B_L\|
\le 2^{i_L-1}+\|\tA_L-\tB_L\|+2^{i_L-1}
$$
%----------------------------------------------------------
so that
%----------------------------------------------------------
$$
\|\tA_L-\tB_L\|\ge \|A_L-B_L\|- 2^{i_L}.
$$
%----------------------------------------------------------
But, by \rf{IL-VL}, $\|A_L-B_L\|>2^{i_L}$ proving that $\tA_L\ne \tB_L$.
\medskip
%----------------------------------------------------------
%@@@@@@@@@@@@@@@@@@@@@@@@@@@@@@@@@@@@@@@@@@@@@@@@@@@@@@@@@@
%@@@@@@@@@@@@@@@@@@@@@@@@@@@@@@@@@@@@@@@@@@@@@@@@@@@@@@@@@@
%@@@@@@@@@@@@@@@@@@@@@@@@@@@@@@@@@@@@@@@@@@@@@@@@@@@@@@@@@@
%----------------------------------------------------------
\par Prove that the set $\{\tA_L,\tB_L\}\nsubseteq E_{i_L+2}$. In fact, if $\{\tA_L,\tB_L\}\subset E_{i_L+2}$ then $\|\tA_L-\tB_L\|\ge 2^{i_L+2}$, see \rf{SEP-E}. But, by \rf{IL-VL} and \rf{Y5},
%----------------------------------------------------------
$$
\|\tA_L-\tB_L\|\le \|\tA_L-A_L\|+\|A_L-B_L\|+\|B_L-\tB_L\|
\le 2^{i_L-1}+2^{i_L+1}+2^{i_L-1}=3\cdot 2^{i_L}<2^{i_L+2},
$$
%----------------------------------------------------------
a contradiction.
%----------------------------------------------------------
\par Since $\tA_L,\tB_L\in E_{i_L-2}$, the statement \rf{GW-1} follows. Thus there exists a point $C_L\in E$  such that
%----------------------------------------------------------
$$
C_L\in\{\tA_L,\tB_L\}\cap(E_{i_L-2}\setminus E_{i_L+2}).
$$
%----------------------------------------------------------
\par Note that, by \rf{Y5}, $C_L\in [V_L]_{\delta}$
with $\delta=2^{i_L-1}$.  Here given $\delta>0$ the sign $[\cdot]_{\delta}$ denotes the open $\delta$-neighborhood of a set. Hence, by \rf{IL-VL},
%----------------------------------------------------------
\bel{GW-3}
C_L\in [V_L]_{\ve}\cap(E_{i_L-2}\setminus E_{i_L+2})
\ee
%----------------------------------------------------------
with $\ve=\tfrac12\diam V_L$.
%----------------------------------------------------------
\medskip
%----------------------------------------------------------
%@@@@@@@@@@@@@@@@@@@@@@@@@@@@@@@@@@@@@@@@@@@@@@@@@@@@@@@@@@
%@@@@@@@@@@@@@@@@@@@@@@@@@@@@@@@@@@@@@@@@@@@@@@@@@@@@@@@@@@
%@@@@@@@@@@@@@@@@@@@@@@@@@@@@@@@@@@@@@@@@@@@@@@@@@@@@@@@@@@
%----------------------------------------------------------
\par We turn to definition of $\PRL(L)$ whenever $L$ is a bounded lacuna satisfying the following condition:
%----------------------------------------------------------
\bel{DF-SG1}
\diam\QL\le\sigma\,\diam V_L
\ee
%----------------------------------------------------------
where
%----------------------------------------------------------
\bel{SIG}
\sigma:=33\tau.
\ee
%----------------------------------------------------------
Recall that $\tau$ is the constant from part (ii) of Proposition \reff{M-LAC}.
%@@@@@@@@@@@@@@@@@@@@@@@@@@@@@@@@@@@@@@@@@@@@@@@@@@@@@@@@@@
%----------------------------------------------------------
\par In this case we define $\PRL(L)$ by
%----------------------------------------------------------
\bel{DF-CL}
\PRL(L):=C_L.
\ee
%----------------------------------------------------------
\par In particular, by \rf{A-DL}, each {\it elementary lacuna} $L\in \hL_E$, see \rf{EL-LC}, satisfies inequality \rf{DF-SG1}, so that
%----------------------------------------------------------
\bel{DF-ELCL}
\PRL(L):=C_L~~~\text{for every elementary lacuna}~~ L\in \hL_E.
\ee
%----------------------------------------------------------
\par It remains to define $\PRL(L)$ whenever $L$ is a true lacuna satisfying inequality
%----------------------------------------------------------
\bel{Y6}
\diam\QL>\sigma\,\diam V_L.
\ee
%----------------------------------------------------------
\par First suppose that $L$ is a true {\it bounded lacuna}. In this case the cube $\QL$ is well defined, see \rf{DF-QLU}. Let $j_L\in\mZ$ be an integer such that
%----------------------------------------------------------
\bel{Y7}
2^{j_L-1}<\tfrac{1}{\sigma}\,\diam \QL\le 2^{j_L}.
\ee
%----------------------------------------------------------
Recall that, by \rf{YY},
%----------------------------------------------------------
$$
40\,\diam \QL\le \dist(V_L,E\setminus V_L)
$$
%----------------------------------------------------------
so that, by \rf{Y7},
%----------------------------------------------------------
\bel{Y8}
\dist(V_L,E\setminus V_L)\ge 40\sigma\,(\diam \QL/\sigma)\ge
40\sigma\,2^{j_L-1}> 2^{j_L+2}\,.
\ee
%----------------------------------------------------------
%@@@@@@@@@@@@@@@@@@@@@@@@@@@@@@@@@@@@@@@@@@@@@@@@@@@@@@@@@@
%----------------------------------------------------------
\par On the other hand, by \rf{Y6} and \rf{Y7},
%----------------------------------------------------------
\bel{Y9}
\diam V_L<\tfrac{1}{\sigma}\,\diam \QL\le 2^{j_L}.
\ee
%----------------------------------------------------------
In particular, by \rf{Y8} and \rf{Y9}, $\diam V_L<\dist(V_L,E\setminus V_L)$ which implies that
%----------------------------------------------------------
\bel{CL-VL}
C_L\in V_L.
\ee
%---------------------------------------------------------
In fact, by \rf{GW-3}, $C_L\in [V_L]_{\ve}$
with $\ve=\tfrac12\diam V_L$ so that $C_L\notin E\setminus V_L$.\medskip
%----------------------------------------------------------
\par Let us consider now the set $E_{j_L}\subset E$.
We know that $\|z-z'\|\ge 2^{j_L}$ for all $z,z'\in E_{j_L}$, and that for each $x\in E$ there exists $z\in E_{j_L}$ such that $\|x-z\|\le 2^{j_L+1}$.
%--------------------------------------------------------
\par Prove that the set $V_L\cap E_{j_L}$ is a singleton. Let us fix a point $x_0\in V_L$. Then there exists a point $z_0\in E_{j_L}$ such that $\|x_0-z_0\|\le 2^{j_L+1}$. Prove that $z_0\in V_L$.
%--------------------------------------------------------
\par In fact, if $z_0\in E\setminus V_L$, then, by \rf{Y8}, %----------------------------------------------------------
$$
\|x_0-z_0\|\ge \dist(V_L,E\setminus V_L)>2^{j_L+2},
$$
%----------------------------------------------------------
a contradiction.
%--------------------------------------------------------
\par Prove that $\{z_0\}=V_L\cap E_{j_L}$. In fact, if there exists a point $z_1\in V_L\cap E_{j_L}$, $z_1\ne z_0$, then $\|z_0-z_1\|\ge 2^{j_L}$ so that
%----------------------------------------------------------
$$
\diam V_L\ge \|z_0-z_1\|\ge 2^{j_L}
$$
%----------------------------------------------------------
which contradicts \rf{Y9}.
%--------------------------------------------------------
\par We denote this unique point of the intersection $V_L\cap E_{j_L}$ by $D_L$; thus
%----------------------------------------------------------
\bel{DL-VL}
\{D_L\}=V_L\cap E_{j_L}.
\ee
%----------------------------------------------------------
%\medskip
%----------------------------------------------------------
%@@@@@@@@@@@@@@@@@@@@@@@@@@@@@@@@@@@@@@@@@@@@@@@@@@@@@@@@@@
%@@@@@@@@@@@@@@@@@@@@@@@@@@@@@@@@@@@@@@@@@@@@@@@@@@@@@@@@@@
%@@@@@@@@@@@@@@@@@@@@@@@@@@@@@@@@@@@@@@@@@@@@@@@@@@@@@@@@@@
%----------------------------------------------------------
\par We are now in a position to define $\PRL(L)$ for an arbitrary bounded lacuna satisfying inequality \rf{Y6}. In this case we put
%----------------------------------------------------------
\bel{DEF-PR2}
\PRL(L):=\left \{
%----------------------------------------------------------
\begin{array}{ll}
D_L,& \text{if}~~~D_L\in E_{j_L}\setminus E_{j_L+k},\smallskip\\
C_L,& \text{otherwise}.
\end{array}
%----------------------------------------------------------
\right.
\ee
%----------------------------------------------------------
Here
%----------------------------------------------------------
\bel{DEF-K}
k:=[\log_2(360\sigma)]+2.
\ee
%----------------------------------------------------------
%\smallskip
%----------------------------------------------------------
%@@@@@@@@@@@@@@@@@@@@@@@@@@@@@@@@@@@@@@@@@@@@@@@@@@@@@@@@@@
%@@@@@@@@@@@@@@@@@@@@@@@@@@@@@@@@@@@@@@@@@@@@@@@@@@@@@@@@@@
%@@@@@@@@@@@@@@@@@@@@@@@@@@@@@@@@@@@@@@@@@@@@@@@@@@@@@@@@@@
%----------------------------------------------------------
\par It remains to define $\PRL(L)$ for an {\it unbounded lacuna} $L\in\LE$. By Proposition \reff{P-2L}, such a lacuna exists if and only if $E$ is a bounded set. Furthermore, by part (ii) of this proposition, such a lacuna is unique, and $V_L=E$.
%----------------------------------------------------------
\par In this case we define $\PRL(L)$ by the formula \rf{DF-CL}. Thus
%----------------------------------------------------------
\bel{DF-UNB}
\PRL(L):=C_L~~~\text{provided the lacuna}~L~\text{is unbounded}.
\ee
%----------------------------------------------------------
\medskip
%----------------------------------------------------------
%@@@@@@@@@@@@@@@@@@@@@@@@@@@@@@@@@@@@@@@@@@@@@@@@@@@@@@@@@@
%@@@@@@@@@@@@@@@@@@@@@@@@@@@@@@@@@@@@@@@@@@@@@@@@@@@@@@@@@@
%@@@@@@@@@@@@@@@@@@@@@@@@@@@@@@@@@@@@@@@@@@@@@@@@@@@@@@@@@@
%----------------------------------------------------------
\par We have defined the projector $\PRL$ on all of the family $\LE$ of lacunae of the set $E$. Prove that this mapping satisfies properties (i)-(iii) of the theorem.\medskip
%----------------------------------------------------------
\par {\it Proof of part (i) of the theorem.} Let us prove (i) with any $\tgm\ge 180$. By formulae \rf{DF-CL}, \rf{DEF-PR2} and \rf{DF-UNB}, $\PRL(L)\in\{C_L,D_L\}$ for every lacuna $L\in\LE$.
Note that $D_L\in V_L$. Since $V_L=(90 Q)\cap E$ for each $Q\in L$, we conclude that $V_L\subset 90 Q$ proving that $D_L\in 90 Q$.
%----------------------------------------------------------
\par By \rf{GW-3}, $C_L$ belongs to the $\ve$-neighborhood of $V_L$ with $\ve=\tfrac12\diam V_L$. But $V_L\subset 90 Q$ so that $C_L$ belongs to the $\tilde{\ve}$-neighborhood of $90Q$ with $\tilde{\ve}=45\diam Q$. Hence, $C_L\in \tgm Q$ provided $\tgm\ge 180$.
%----------------------------------------------------------
\smallskip
%----------------------------------------------------------
%@@@@@@@@@@@@@@@@@@@@@@@@@@@@@@@@@@@@@@@@@@@@@@@@@@@@@@@@@@
%@@@@@@@@@@@@@@@@@@@@@@@@@@@@@@@@@@@@@@@@@@@@@@@@@@@@@@@@@@
%@@@@@@@@@@@@@@@@@@@@@@@@@@@@@@@@@@@@@@@@@@@@@@@@@@@@@@@@@@
%----------------------------------------------------------
\par {\it Proof of part (ii) of the theorem.} Let $L,L'\in\LE$ be two distinct lacunae such that their ``projections'' $\PRL(L)$ and $\PRL(L')$ are distinct as well. Thus $L\ne L'$ and $\PRL(L)\ne \PRL(L')$. Prove that under these conditions inequality \rf{DM-QQP} holds with some absolute constant $\tgm>0$.
%----------------------------------------------------------
\par First we note that, by part (1) of Lemma \reff{Wadd}, %----------------------------------------------------------
\bel{Q4P}
\frac{1}{4}\le\diam Q/\diam Q'\le 4\,.
\ee
%----------------------------------------------------------
\par By part (i) of Proposition \reff{M-LAC}, either $L$ or $L'$ is an elementary lacuna. Thus, without loss of generality, we may assume that $L'$ {\it is an elementary lacuna}.
%----------------------------------------------------------
\par Then, by \rf{A-DL},
%----------------------------------------------------------
\bel{QP-VL}
\diam Q'\le 2\,\diam V_{L'}.
\ee
%----------------------------------------------------------
We also recall that, by \rf{D-VL},
%----------------------------------------------------------
\bel{V-90L}
\diam V_L\le 90\diam Q~~~\text{and}~~~\diam V_{L'}\le 90\diam Q'\,.
\ee
%----------------------------------------------------------
\par We know that in this case $\PRL(L')=C_{L'}$.
See \rf{DF-ELCL}.
%----------------------------------------------------------
\par Recall that
%----------------------------------------------------------
\bel{CLP}
C_{L'}\in E_{i_{L'}-2}\setminus E_{i_{L'}+2}
\ee
%----------------------------------------------------------
where $i_{L'}$ is an  integer such that
%----------------------------------------------------------
\bel{IL-VP}
2^{i_{L'}}<\diam V_{L'}\le 2^{i_{L'}+1}.
\ee
%----------------------------------------------------------
See \rf{GW-3} and \rf{IL-VL}.
%----------------------------------------------------------
\par The next lemma shows that, under certain restriction on the cube $Q$ the inequality \rf{DM-QQP} holds.
%---------------------------------------------------------- %@@@@@@@@@@@@@@@@@@@@@@@@@@@@@@@@@@@@@@@@@@@@@@@@@@@@@@@@@@
%@@@@@@@@@@@@@@@@@@@@@@@@@@@@@@@@@@@@@@@@@@@@@@@@@@@@@@@@@@
%@@@@@@@@@@@@@@@@@@@@@@@@@@@@@@@@@@@@@@@@@@@@@@@@@@@@@@@@@@
%@@@@@@@@@@@@@@@@@@@@@@@@@@@@@@@@@@@@@@@@@@@@@@@@@@@@@@@@@@
%----------------------------------------------------------
\begin{lemma}\lbl{TWO-CL} Suppose that $\PRL(L)=C_L$ and $\diam Q\le\theta \diam V_L$ where $\theta$ is a positive constant. Then
%----------------------------------------------------------
$$
\diam Q+\diam Q'\le C\,\theta\,\|\PRL(L)-\PRL(L')\|
$$
%----------------------------------------------------------
where $C$ is an absolute constant.
%----------------------------------------------------------
\end{lemma}
%----------------------------------------------------------
%@@@@@@@@@@@@@@@@@@@@@@@@@@@@@@@@@@@@@@@@@@@@@@@@@@@@@@@@@@
%@@@@@@@@@@@@@@@@@@@@@@@@@@@@@@@@@@@@@@@@@@@@@@@@@@@@@@@@@@
%@@@@@@@@@@@@@@@@@@@@@@@@@@@@@@@@@@@@@@@@@@@@@@@@@@@@@@@@@@
%@@@@@@@@@@@@@@@@@@@@@@@@@@@@@@@@@@@@@@@@@@@@@@@@@@@@@@@@@@
%----------------------------------------------------------
\par {\it Proof.} Let $m\in\mZ$ be an integer such that $2^{m-1}<\theta\le 2^m$. Then, by \rf{IL-VP}, \rf{V-90L} and \rf{Q4P},
%----------------------------------------------------------
$$
2^{i_{L'}}\le\diam V_{L'}\le 90\diam Q'\le 360\diam Q\le 360\,\theta\diam V_L.
$$
%----------------------------------------------------------
Hence, by \rf{IL-VL},
%----------------------------------------------------------
$$
2^{i_{L'}}\le 360\,\theta\, 2^{i_{L}+1}
\le 2^{i_{L}+m+9}
$$
%----------------------------------------------------------
proving that $j:=i_{L'}-m-9\le i_{L}$.
%----------------------------------------------------------
\par Since $C_{L}\in E_{i_L}$, $C_{L'}\in E_{i_{L'}}$, we conclude that $C_{L},C_{L'}\in E_{j}$.  But $C_L=\PRL(L)$,
$C_{L'}=\PRL(L')$, and $\PRL(L)\ne\PRL(L')$ so that $C_L$ and $C_{L'}$ are two distinct points of $E_j$. Therefore, by \rf{SEP-E}, these two points are $2^j$-separated, i.e.,
%----------------------------------------------------------
$$
\|\PRL(L)-\PRL(L')\|\ge 2^j=2^{-(m+10)}\, 2^{i_{L'}+1}\ge
2^{-11}\,2^{i_{L'}+1}/\theta.
$$
%----------------------------------------------------------
Hence, by \rf{IL-VP},
%----------------------------------------------------------
$$
\diam V_{L'}\le 2^{11}\,\theta\,\|\PRL(L)-\PRL(L')\|
$$
%----------------------------------------------------------
so that, by \rf{QP-VL},
%----------------------------------------------------------
$$
\diam Q'\le 2^{12}\theta\,\|\PRL(L)-\PRL(L')\|\,.
$$
%----------------------------------------------------------
Finally, by \rf{Q4P},
%----------------------------------------------------------
$$
\diam Q+\diam Q'\le 5\diam Q'\le 2^{15}\theta\,\|\PRL(L)-\PRL(L')\|
$$
%----------------------------------------------------------
proving the lemma.\bx
%----------------------------------------------------------
\bigskip
%@@@@@@@@@@@@@@@@@@@@@@@@@@@@@@@@@@@@@@@@@@@@@@@@@@@@@@@@@@
%@@@@@@@@@@@@@@@@@@@@@@@@@@@@@@@@@@@@@@@@@@@@@@@@@@@@@@@@@@
%@@@@@@@@@@@@@@@@@@@@@@@@@@@@@@@@@@@@@@@@@@@@@@@@@@@@@@@@@@
%@@@@@@@@@@@@@@@@@@@@@@@@@@@@@@@@@@@@@@@@@@@@@@@@@@@@@@@@@@
%----------------------------------------------------------
\par We turn to the proof of the property (ii) of the theorem in general case. We do this in several steps. First we prove \rf{DM-QQP} for the lacuna $L$ satisfying inequality \rf{DF-SG1}. In this case, by \rf{DF-CL}, $\PRL(L)=C_L$ so that the conditions of Lemma \reff{TWO-CL} are satisfied. This lemma implies inequality \rf{DM-QQP} with a constant $\tgm=C\,\sigma$ where $C>0$ is an absolute constant.
%----------------------------------------------------------
\par Let now $L$ be an unbounded lacuna . By part (i) of Proposition \reff{P-2L}, in this case the set $E$ is bounded. Furthermore, $L$ coincides with the {\it unique} unbounded lacuna $L^{max}$ which is a true lacuna such that $V_L=V_{L^{max}}=E$.
%----------------------------------------------------------
\par By \rf{Q4P} and \rf{QP-VL},
%----------------------------------------------------------
$$
\diam Q\le 4\diam Q'\le 8\diam V_{L'}\le 8\diam E
$$
%----------------------------------------------------------
so that $\diam Q\le 8\diam V_L$. Furthermore, by \rf{DF-UNB}, $\PRL(L)=C_L$. Thus the conditions of Lemma \reff{TWO-CL} are satisfied so that inequality \rf{DM-QQP} for the unbounded lacuna $L$ holds with some absolute constant $\tgm>0$.\medskip
%----------------------------------------------------------
%@@@@@@@@@@@@@@@@@@@@@@@@@@@@@@@@@@@@@@@@@@@@@@@@@@@@@@@@@@
%@@@@@@@@@@@@@@@@@@@@@@@@@@@@@@@@@@@@@@@@@@@@@@@@@@@@@@@@@@
%@@@@@@@@@@@@@@@@@@@@@@@@@@@@@@@@@@@@@@@@@@@@@@@@@@@@@@@@@@
%@@@@@@@@@@@@@@@@@@@@@@@@@@@@@@@@@@@@@@@@@@@@@@@@@@@@@@@@@@
%----------------------------------------------------------
\par It remains to prove \rf{DM-QQP} for a {\it bounded} lacuna $L$ satisfying inequality \rf{Y6}. By part (ii) of Proposition \reff{M-LAC}, it suffices to consider two case.
%----------------------------------------------------------
\smallskip
%@@@@@@@@@@@@@@@@@@@@@@@@@@@@@@@@@@@@@@@@@@@@@@@@@@@@@@@@@@
%@@@@@@@@@@@@@@@@@@@@@@@@@@@@@@@@@@@@@@@@@@@@@@@@@@@@@@@@@@
%@@@@@@@@@@@@@@@@@@@@@@@@@@@@@@@@@@@@@@@@@@@@@@@@@@@@@@@@@@
%@@@@@@@@@@@@@@@@@@@@@@@@@@@@@@@@@@@@@@@@@@@@@@@@@@@@@@@@@@
%----------------------------------------------------------
\par {\it The first case: for the cubes $Q$ and $\QL$ the following inequality
%----------------------------------------------------------
\bel{F1Q}
\diam \QL\le \tau\,\diam Q
\ee
%----------------------------------------------------------
holds.}
%----------------------------------------------------------
\par Prove that
%----------------------------------------------------------
\bel{JLK}
j_L+2<i_{L'}<j_L+k-2.
\ee
%----------------------------------------------------------
(Recall that $j_L$, $i_{L'}$ and $k$ are defined by \rf{Y7}, \rf{IL-VP} and \rf{DEF-PR2} respectively.)
%----------------------------------------------------------
\par We begin with the proof of the inequality $j_L+2<i_{L'}$. By \rf{Y7}, \rf{F1Q} and \rf{Q4P},
%----------------------------------------------------------
$$
2^{j_L-1}\le \diam \QL/\sigma\le \tau\diam Q/\sigma
\le 4\tau\,\diam Q'/\sigma
$$
%----------------------------------------------------------
so that, by \rf{QP-VL} and \rf{IL-VP},
%----------------------------------------------------------
$$
2^{j_L-1}\le 2\tau\,\diam V_{L'}/\sigma\le 2\tau 2^{i_{L'}+1}/\sigma\,.
$$
%----------------------------------------------------------
Since $\sigma>32\tau$, see \rf{SIG}, we obtain the required inequality $j_L+2<i_{L'}$.
%----------------------------------------------------------
\par Prove that $i_{L'}< j_L+k-2$. By \rf{IL-VP}, \rf{V-90L}, \rf{Q4P} and \rf{Y7},
%----------------------------------------------------------
$$
2^{i_{L'}}< \diam V_{L'}\le 90\diam Q'
\le 360\diam Q\le 360\diam \QL\le 360\sigma\,2^{j_L}\,.
$$
%----------------------------------------------------------
Since $360\sigma\le 2^{k-2}$, see \rf{SIG} and \rf{DEF-K}, the required inequality $i_{L'}< j_L+k-2$ follows.\medskip
%----------------------------------------------------------
%@@@@@@@@@@@@@@@@@@@@@@@@@@@@@@@@@@@@@@@@@@@@@@@@@@@@@@@@@@
%@@@@@@@@@@@@@@@@@@@@@@@@@@@@@@@@@@@@@@@@@@@@@@@@@@@@@@@@@@
%@@@@@@@@@@@@@@@@@@@@@@@@@@@@@@@@@@@@@@@@@@@@@@@@@@@@@@@@@@
%@@@@@@@@@@@@@@@@@@@@@@@@@@@@@@@@@@@@@@@@@@@@@@@@@@@@@@@@@@
%----------------------------------------------------------
\par We recall that, by \rf{DF-ELCL}, $\PRL(L')=C_{L'}$ where $C_{L'}$ is a point satisfying \rf{CLP}. Prove that $C_{L'}\in E\setminus V_L$.
%----------------------------------------------------------
\par Suppose that this is not true, i.e., $C_{L'}\in V_L$. By \rf{CLP}, $C_{L'}\in E_{i_{L'}-2}$. But, by \rf{JLK}, $i_{L'}-2> j_L$ so that $C_{L'}\in E_{j_L}\cap V_L$ (because $E_{i_{L'}-2}\subset E_{j_L}$). On the other hand, by \rf{DL-VL}, the set $E_{j_L}\cap V_L=\{D_L\}$ is a singleton so that $C_{L'}=D_L$.
%----------------------------------------------------------
\par Note that, by \rf{CLP}, $C_{L'}\notin E_{i_{L'}+2}$. Since $i_{L'}+2< j_L+k$, see \rf{JLK}, $ E_{j_L+k}\subset E_{i_{L'}+2}$ proving that $\{D_L\}=C_{L'}\in E_{j_L}\setminus E_{j_L+k}$.
%----------------------------------------------------------
\par But, by definition \rf{DEF-PR2}, in this case $\PRL(L)=D_L$ so that
%----------------------------------------------------------
$$
\PRL(L)=D_L=C_{L'}=\PRL(L').
$$
%----------------------------------------------------------
This contradicts our assumption that
$\PRL(L)\ne\PRL(L')$ proving the required imbedding  $$C_{L'}=\PRL(L')\in E\setminus V_L.$$
%----------------------------------------------------------
%@@@@@@@@@@@@@@@@@@@@@@@@@@@@@@@@@@@@@@@@@@@@@@@@@@@@@@@@@@
%@@@@@@@@@@@@@@@@@@@@@@@@@@@@@@@@@@@@@@@@@@@@@@@@@@@@@@@@@@
%@@@@@@@@@@@@@@@@@@@@@@@@@@@@@@@@@@@@@@@@@@@@@@@@@@@@@@@@@@
%@@@@@@@@@@@@@@@@@@@@@@@@@@@@@@@@@@@@@@@@@@@@@@@@@@@@@@@@@@
%----------------------------------------------------------
\par Let us now estimate from below the distance between the points $\PRL(L)$ and $\PRL(L')$. We note that, by \rf{CL-VL} and \rf{DL-VL}, $C_L,D_L\in V_L$. But the point $\PRL(L)\in\{C_L,D_L\}$ proving that $\PRL(L)\in V_L$ as well. Since $\PRL(L')\in E\setminus V_L$, we have
%----------------------------------------------------------
$$
\|\PRL(L)-\PRL(L')\|\ge \dist(V_L,E\setminus V_L)
$$
%----------------------------------------------------------
so that, by part (i) of Proposition \reff{P-SM},
%----------------------------------------------------------
$$
40\diam \QL\le\|\PRL(L)-\PRL(L')\|.
$$
%----------------------------------------------------------
Since $\diam Q\le\diam\QL$, by \rf{Q4P},
%----------------------------------------------------------
$$
\diam Q+\diam Q'\le 5\diam Q\le 5\diam\QL\le \|\PRL(L)-\PRL(L')\|
$$
%----------------------------------------------------------
proving \rf{DM-QQP} with $\tgm=1$.
%----------------------------------------------------------
\bigskip
%----------------------------------------------------------
%@@@@@@@@@@@@@@@@@@@@@@@@@@@@@@@@@@@@@@@@@@@@@@@@@@@@@@@@@@
%@@@@@@@@@@@@@@@@@@@@@@@@@@@@@@@@@@@@@@@@@@@@@@@@@@@@@@@@@@
%@@@@@@@@@@@@@@@@@@@@@@@@@@@@@@@@@@@@@@@@@@@@@@@@@@@@@@@@@@
%@@@@@@@@@@@@@@@@@@@@@@@@@@@@@@@@@@@@@@@@@@@@@@@@@@@@@@@@@@
%----------------------------------------------------------
\par {\it The second case:}
%----------------------------------------------------------
\bel{F2Q}
\diam Q\le \tau\,\diam Q_L\,.
\ee
%----------------------------------------------------------
\par Clearly, in this case $\diam V_L>0$ so that, by part (ii) of Proposition \reff{P-SM},
%----------------------------------------------------------
$$
\diam Q_L\le \gamma_1\diam V_L~~~\text{for some absolute constant}~~~ \gamma_1>0.
$$
%----------------------------------------------------------
Hence
%----------------------------------------------------------
\bel{TG-1}
\diam Q\le \tau\,\gamma_1\diam V_L.
\ee
%----------------------------------------------------------
\par By definition \rf{DEF-PR2}, $\PRL(L)=C_L$ provided the point $D_L\notin E_{j_L}\setminus E_{j_L+k}$. In this case, by \rf{TG-1}, the conditions of Lemma \reff{TWO-CL} are satisfied. By this lemma inequality \rf{DM-QQP} holds with a constant $\tgm=C\,\tau\,\gamma_1$ where $C>0$ is an absolute constant.
%----------------------------------------------------------
\par Let now $D_L\in E_{j_L}\setminus E_{j_L+k}$ so that, by \rf{DEF-PR2}, $\PRL(L)=D_L$. Prove that in this case $i_{L'}\le j_L+m$ where $m:=[\log_2\gamma_2]$ and $\gamma_2:=360\tau\,\gamma_1$.
%----------------------------------------------------------
\par By \rf{IL-VP}, \rf{V-90L} and \rf{Q4P},
%----------------------------------------------------------
$$
2^{i_{L'}}\le \diam V_{L'}\le 90\diam Q'\le 360\diam Q
$$
%----------------------------------------------------------
so that, by \rf{F2Q} and \rf{YG},
%----------------------------------------------------------
$$
2^{i_{L'}}\le 360\,\tau\diam Q_L\le 360\,\tau\gamma_1\diam V_L=\gamma_2\diam V_L\le 2^{j_L+m}
$$
%----------------------------------------------------------
proving the required inequality $i_{L'}\le j_L+m$.
%----------------------------------------------------------
\par Recall that $\PRL(L)=D_L\in E_{j_L}$ and $\PRL(L')=C_{L'}\in E_{i_{L'}-2}$, see \rf{CLP}. Hence
%----------------------------------------------------------
$$
\PRL(L),\,\PRL(L')\in E_{i_{L'}-m}
$$
%---------------------------------------------------------
so that $\PRL(L)$ and $\PRL(L')$ are $2^{i_{L'}-m}$-separated points, i.e.,
%----------------------------------------------------------
$$
\|\PRL(L)-\PRL(L')\|\ge 2^{i_{L'}-m}.
$$
%----------------------------------------------------------
We also recall that, by \rf{IL-VP}, $\diam V_{L'}\le 2^{i_{L'}+1}$ so that
%----------------------------------------------------------
$$
\diam V_{L'}\le 2^{m+1}\|\PRL(L)-\PRL(L')\|\le 4\gamma_2 \|\PRL(L)-\PRL(L')\|.
$$
%----------------------------------------------------------
Combining this inequality with \rf{QP-VL}, we obtain
%----------------------------------------------------------
$$
\diam Q'\le 8\gamma_2 \|\PRL(L)-\PRL(L')\|
$$
%----------------------------------------------------------
so that, by \rf{Q4P},
%----------------------------------------------------------
$$
\diam Q+\diam Q'\le 5\diam Q'\le  40\gamma_2 \|\PRL(L)-\PRL(L')\|.
$$
%----------------------------------------------------------
\par This completes the proof of part (ii) of the theorem.
%----------------------------------------------------------
\bigskip
%----------------------------------------------------------
%@@@@@@@@@@@@@@@@@@@@@@@@@@@@@@@@@@@@@@@@@@@@@@@@@@@@@@@@@@
%@@@@@@@@@@@@@@@@@@@@@@@@@@@@@@@@@@@@@@@@@@@@@@@@@@@@@@@@@@
%@@@@@@@@@@@@@@@@@@@@@@@@@@@@@@@@@@@@@@@@@@@@@@@@@@@@@@@@@@
%@@@@@@@@@@@@@@@@@@@@@@@@@@@@@@@@@@@@@@@@@@@@@@@@@@@@@@@@@@
%----------------------------------------------------------
\par {\it Proof of part (iii) of the theorem.} Let $A\in E$ and let %----------------------------------------------------------
$$
\PRC(A):=\{L\in\LE: \PRL(L)=A\}.
$$
%----------------------------------------------------------
\par By Proposition \reff{P-2L}, the set $\PRC(A)$ contains at most one {\it unbounded} lacuna so that, without loss of generality, we may assume that all lacunae from $\PRC(A)$ are bounded.
%----------------------------------------------------------
\par Let us also note that if $A$ is an isolated point of $E$, then there exists a {\it unique} lacuna $L_A\in\LE$
such that $\{A\}=V_{L_A}$. (Of course, $L_A$ is a true lacuna.) Conversely, if $\{A\}=V_{L}$ for some $L\in\LE$, then $A$ is an isolated point of $E$ so that $L=L_A$.
This elementary remark enables us to assume that $\diam V_L>0$ for each lacuna $L\in\PRC(A)$.
%----------------------------------------------------------
\par Now, by definition \rf{DF-CL} and \rf{DEF-PR2}, the point $\PRL(L)$ coincides either with the point $C_L$, see \rf{GW-3}, or with the point $D_L$, see \rf{DL-VL}.
%----------------------------------------------------------
\par Note that if $\PRL(L)=C_L$, by \rf{GW-3},
%----------------------------------------------------------
$$
\PRL(L)\in E_{i_L-2}\setminus E_{i_L+2}.
$$
%----------------------------------------------------------
(Recall that $i_L\in\mZ$ is an integer determined by inequalities \rf{IL-VL}.) Hence, by \rf{Y6} and \rf{D-VL},  $\diam Q_L\sim 2^{i_L}$ with absolute constants in this equivalence.
%----------------------------------------------------------
\par Let now $\PRL(L)=D_L$ for some lacuna $L\in\PRC(A)$. Then, by \rf{DEF-PR2},
%----------------------------------------------------------
$$
\PRL(L)\in E_{j_L}\setminus E_{j_L+k}.
$$
%----------------------------------------------------------
See \rf{DEF-K}. Furthermore, by \rf{Y7}, $\diam \QL\sim 2^{j_L}$ with absolute constants.
%----------------------------------------------------------
\par Summarizing these properties of $\PRL(L)$ we conclude that for each lacuna $L\in\PRC(A)$ there exists an integer $m_L\in\mZ$ and a cube $K_L\in L$ such that
%----------------------------------------------------------
\bel{PR-EM}
\PRL(L)\in E_{m_L}\setminus E_{m_L+k}
\ee
%----------------------------------------------------------
and
%----------------------------------------------------------
\bel{ML-C1}
\tfrac{1}{C_1}\,2^{m_L}\le \diam K_L\le C_1\,2^{m_L}.
\ee
%----------------------------------------------------------
Here $C_1>0$ is an absolute constant, and $k$ is defined by \rf{DEF-K}.
%----------------------------------------------------------
\par Now let $L,L'\in\PRC(A)$, i.e., $\PRL(L)=\PRL(L')=A$. Then, by \rf{PR-EM},
%----------------------------------------------------------
$$
A\notin E_{m_L+k}~~~\text{and}~~~A\in E_{m_{L'}}
$$
%----------------------------------------------------------
so that $E_{m_L+k}\nsupseteq E_{m_{L'}}$. Hence $E_{m_L+k}\subset E_{m_{L'}}$ proving that $m_{L'}\le m_L+k$. In the same fashion we prove that $m_{L}\le m_{L'}+k$, so that $|m_{L}- m_{L'}|\le k$.
%----------------------------------------------------------
\par We also note that, by inequality \rf{ML-C1} (which we apply to $L'$), we have
%----------------------------------------------------------
$$
\tfrac{1}{C_1}\,2^{m_{L'}}\le \diam K_L\le C_1\,2^{m_{L'}}.
$$
%----------------------------------------------------------
Hence,
%----------------------------------------------------------
$$
\diam K_L\le C_1\,2^{m_{L}}\le C_1 2^k 2^{m_{L'}}
\le C_1^2 2^k \diam K_{L'}.
$$
%----------------------------------------------------------
In  the same way we obtain that $\diam K_{L'}\le C_1^2 2^k \diam K_{L}$. Thus
%----------------------------------------------------------
\bel{DKL}
C_1^{-2} 2^{-k}\le \diam K_L/\diam K_{L'}\le C_1^2 2^k.
\ee
%----------------------------------------------------------
\par We recall that $K_L\in L$ and $K_{L'}\in L'$ so that $K_L\ne K_{L'}$ provided $L\ne L'$. Thus the family
%----------------------------------------------------------
$$
\Kc_A:=\{K_L: L\in\PRC(A)\}
$$
%----------------------------------------------------------
consists of pairwise disjoint cubes. We also note that, by part (i) of the theorem (proven earlier),
%----------------------------------------------------------
$$
A\in \tgm K_L~~~\text{for every}~~~L\in\PRC(A)
$$
%----------------------------------------------------------
with $\tgm=180$. Hence,
%----------------------------------------------------------
$$
\bigcup_{K_L\in\Kc_A} K_L\subset (2\tgm) K^{max}
$$
%----------------------------------------------------------
where $K^{max}$ is a cube from $\Kc_A$ of the maximal diameter.
%----------------------------------------------------------
\par Note that the cubes of the family $\Kc_A$ are pairwise disjoint, and the diameters of these cubes are equivalent to the diameter of the cube $K^{max}$, see \rf{DKL}. Consequently, the number of these cubes, the quantity $\# \Kc_A$, is bounded by a constant $C=C(n)$ depending only on $n$. But $\#\,\Kc_A=\#\,\PRC(A)$ which proves part (iii) of the theorem.
%----------------------------------------------------------
\par The proof of Theorem \reff{L-PE} is complete.\bx
%----------------------------------------------------------
%@@@@@@@@@@@@@@@@@@@@@@@@@@@@@@@@@@@@@@@@@@@@@@@@@@@@@@@@@@
%@@@@@@@@@@@@@@@@@@@@@@@@@@@@@@@@@@@@@@@@@@@@@@@@@@@@@@@@@@
%@@@@@@@@@@@@@@@@@@@@@@@@@@@@@@@@@@@@@@@@@@@@@@@@@@@@@@@@@@
%@@@@@@@@@@@@@@@@@@@@@@@@@@@@@@@@@@@@@@@@@@@@@@@@@@@@@@@@@@
%----------------------------------------------------------
\bigskip
\par Let us note several simple and useful properties of the projector $\PRL$ for a finite set $E$.
%----------------------------------------------------------
\begin{proposition} \lbl{PRL-EF} Let $E$ be a finite subset of $\RN$ and let $x\in E$. There exists a unique lacuna $L^{(x)}$ such that $V_{L^{(x)}}=\{x\}$. Furthermore,  $L^{(x)}$ is a true lacuna, and
%----------------------------------------------------------
\bel{PRL-LX}
\PRL\left(L^{(x)}\right)=x\,.
\ee
%----------------------------------------------------------
%@@@@@@@@@@@@@@@@@@@@@@@@@@@@@@@@@@@@@@@@@@@@@@@@@@@@@@@@@@
%----------------------------------------------------------
\end{proposition}
%----------------------------------------------------------
%@@@@@@@@@@@@@@@@@@@@@@@@@@@@@@@@@@@@@@@@@@@@@@@@@@@@@@@@@@
%@@@@@@@@@@@@@@@@@@@@@@@@@@@@@@@@@@@@@@@@@@@@@@@@@@@@@@@@@@
%@@@@@@@@@@@@@@@@@@@@@@@@@@@@@@@@@@@@@@@@@@@@@@@@@@@@@@@@@@
%@@@@@@@@@@@@@@@@@@@@@@@@@@@@@@@@@@@@@@@@@@@@@@@@@@@@@@@@@@
%----------------------------------------------------------
\par {\it Proof.} We define $\LX$ by
%----------------------------------------------------------
\bel{DF-LX}
\LX:=\{Q\in W_E: (90Q)\cap E=\{x\}\}\,.
\ee
%----------------------------------------------------------
\par Let
%----------------------------------------------------------
$
\teps:=\dist(x,E\setminus\{x\})/180
$
%----------------------------------------------------------
and let $0<\ve<\teps$. Then each cube $Q\in W_E$ such that $Q\subset Q(x,\ve)$ belongs to $\LX$. In fact,
%----------------------------------------------------------
$$
(90Q)\subset Q(x,90\ve)\subset Q(x,90\teps)=Q(x,\tfrac12\dist(x,E\setminus\{x\}))\,.
$$
%----------------------------------------------------------
Hence $(90Q)\cap E=\{x\}$ proving that $Q\in \LX$.
%----------------------------------------------------------
\par In particular, $\LX\ne\emp$. Furthermore, every cube $Q\in\LX$ is {\it a lacunary cube}, i.e., \rf{L-PR} is satisfied. In fact, by \rf{9Q-INT}, $(9Q)\cap E\ne\emp$. Since $(90Q)\cap E=\{x\}$, we have
%----------------------------------------------------------
$$
(10Q)\cap E=(90Q)\cap E=\{x\},
$$
%----------------------------------------------------------
i.e., \rf{L-PR} holds. This equality also shows that $\LX$ is {\it a true lacuna}, i.e., an equivalence class with respect to the binary relation $\sim$, see \rf{BIN-R}, on the family $LW_E$ of all lacunary cubes. In addition, $V_{L^{(x)}}=(90Q)\cap E=\{x\}$, see \rf{D-VL}. Also note that the uniqueness of $\LX$ directly follows from definitions \rf{DF-LX} and \rf{D-VL}.
%----------------------------------------------------------
\par Prove \rf{PRL-LX}. By part (i) of Theorem \reff{L-PE}, $\PRL(\LX)\in (\tgm\,Q)\cap E$ for every $Q\in \LX$ where $\tgm\ge1$ is an absolute constant.
%----------------------------------------------------------
\par Let $\ve_0:=\teps/(2\tgm)$, and let $Q\subset Q(x,\ve_0)$. We know that $Q\in\LX$. But
%----------------------------------------------------------
$$
\tgm Q\subset Q(x,\tgm\ve_0)=Q(x,\teps/2)
$$
%----------------------------------------------------------
so that
%----------------------------------------------------------
$$
(\tgm Q)\cap E\subset Q(x,\teps/2)\cap E=\{x\}\,.
$$
%----------------------------------------------------------
Hence $\PRL(\LX)\in (\tgm\,Q)\cap E=\{x\}$ proving the proposition.\bx
%----------------------------------------------------------
%@@@@@@@@@@@@@@@@@@@@@@@@@@@@@@@@@@@@@@@@@@@@@@@@@@@@@@@@@@
%@@@@@@@@@@@@@@@@@@@@@@@@@@@@@@@@@@@@@@@@@@@@@@@@@@@@@@@@@@
%@@@@@@@@@@@@@@@@@@@@@@@@@@@@@@@@@@@@@@@@@@@@@@@@@@@@@@@@@@
%@@@@@@@@@@@@@@@@@@@@@@@@@@@@@@@@@@@@@@@@@@@@@@@@@@@@@@@@@@
%----------------------------------------------------------
\medskip
\par Theorem \reff{L-PE} enables us to compare the number
of lacunae with the number of points of a finite set $E\subset\RN$.
%----------------------------------------------------------
\begin{proposition} Let $E$ be a finite subset of $\RN$. Then
%----------------------------------------------------------
$$
\# E\le \# \LE\le C(n)\,\# E.
$$
%----------------------------------------------------------
%@@@@@@@@@@@@@@@@@@@@@@@@@@@@@@@@@@@@@@@@@@@@@@@@@@@@@@@@@@
%----------------------------------------------------------
\end{proposition}
%----------------------------------------------------------
%@@@@@@@@@@@@@@@@@@@@@@@@@@@@@@@@@@@@@@@@@@@@@@@@@@@@@@@@@@
%@@@@@@@@@@@@@@@@@@@@@@@@@@@@@@@@@@@@@@@@@@@@@@@@@@@@@@@@@@
%@@@@@@@@@@@@@@@@@@@@@@@@@@@@@@@@@@@@@@@@@@@@@@@@@@@@@@@@@@
%@@@@@@@@@@@@@@@@@@@@@@@@@@@@@@@@@@@@@@@@@@@@@@@@@@@@@@@@@@
%----------------------------------------------------------
\par {\it Proof.} By Proposition \reff{PRL-EF}, $E\ni x\to \LX\in \LE$ is a one-to-one mapping, so that $\# E\le \# \LE$. On the other hand, by  part (iii) of Theorem \reff{L-PE}, the projector $\PRL:\LE\to E$ is ``almost'' one-to-one, i.e., for each $x\in E$ the number of its sources $\{L\in\LE: \PRL(L)=x\}$ is bounded by a constant $C=C(n)$. Hence $\#\LE\le C(n)\,\# E$ proving the proposition.\bx
%----------------------------------------------------------
%@@@@@@@@@@@@@@@@@@@@@@@@@@@@@@@@@@@@@@@@@@@@@@@@@@@@@@@@@@
%@@@@@@@@@@@@@@@@@@@@@@@@@@@@@@@@@@@@@@@@@@@@@@@@@@@@@@@@@@
%@@@@@@@@@@@@@@@@@@@@@@@@@@@@@@@@@@@@@@@@@@@@@@@@@@@@@@@@@@
%@@@@@@@@@@@@@@@@@@@@@@@@@@@@@@@@@@@@@@@@@@@@@@@@@@@@@@@@@@
%@@@@@@@@@@@@@@@@@@@@@@@@@@@@@@@@@@@@@@@@@@@@@@@@@@@@@@@@@@
%@@@@@@@@@@@@@@@@@@@@@@@@@@@@@@@@@@@@@@@@@@@@@@@@@@@@@@@@@@
%@@@@@@@@@@@@@@@@@@@@@@@@@@@@@@@@@@@@@@@@@@@@@@@@@@@@@@@@@@
%----------------------------------------------------------
\bigskip
%----------------------------------------------------------
\par {\bf 6.3. The graph $\GE$ and its properties.}\medskip
%----------------------------------------------------------
\addtocontents{toc}{~~~~6.3. The graph $\GE$ and its properties.\hfill \thepage\par}
%----------------------------------------------------------
\par The lacunary projector constructed in Theorem \reff{L-PE} generates a certain (undirected) graph with vertices in $E$ which we denote by $\GE$.
%----------------------------------------------------------
%@@@@@@@@@@@@@@@@@@@@@@@@@@@@@@@@@@@@@@@@@@@@@@@@@@@@@@@@@@
%@@@@@@@@@@@@@@@@@@@@@@@@@@@@@@@@@@@@@@@@@@@@@@@@@@@@@@@@@@
%@@@@@@@@@@@@@@@@@@@@@@@@@@@@@@@@@@@@@@@@@@@@@@@@@@@@@@@@@@
%----------------------------------------------------------
%@@@@@@@@@@@@@@@@@@@@@@@@@@@@@@@@@@@@@@@@@@@@@@@@@@@@@@@@@@
\begin{definition}\lbl{GRAPH-GE} {\em We define the graph $\GE$ as follows:\smallskip
%----------------------------------------------------------
\par $\bullet$~ The set of vertices of the graph $\GE$ coincides with $E$;
%----------------------------------------------------------
%@@@@@@@@@@@@@@@@@@@@@@@@@@@@@@@@@@@@@@@@@@@@@@@@@@@@@@@@@@
%----------------------------------------------------------
\par $\bullet$~ Two distinct vertices $A,A'\in E$, $A\ne A'$, are joined by an edge in $\GE$ (we write $A\lr A'$) if there exist {\it contacting lacunae} $L,L'\in\LE$,
$L\lcr L'$, such that
%---------------------------------------------------------
$$
A=\PRL(L)~~~\text{and}~~~A'=\PRL(L').
$$
%----------------------------------------------------------
}
%----------------------------------------------------------
\end{definition}
%----------------------------------------------------------
\smallskip
%----------------------------------------------------------
\par More specifically, notation $A\lr A'$ means that  $A\ne A'$ and there exist lacunae $L,L'\in\LE$ and cubes $Q\in L$ and $Q'\in L'$ such that $Q\cap Q'\ne\emp$, $A=\PRL(L)$ and $A'=\PRL(L')$. See Definition \reff{CONT-L}.\medskip
%----------------------------------------------------------
\par Let us note two important properties of the graph $\GE$.
%----------------------------------------------------------
%@@@@@@@@@@@@@@@@@@@@@@@@@@@@@@@@@@@@@@@@@@@@@@@@@@@@@@@@@@
%@@@@@@@@@@@@@@@@@@@@@@@@@@@@@@@@@@@@@@@@@@@@@@@@@@@@@@@@@@
%----------------------------------------------------------
\begin{proposition}\lbl{PROP-GE} For every closed set $E\subset\RN$ the graph $\GE$ has the following properties:
%----------------------------------------------------------
\par (i). $\GE$ is a $\gamma$-sparse graph where $\gamma=\gamma(n)\ge1$ is a constant depending only on $n$;
%----------------------------------------------------------
\smallskip
\par (ii). $\deg_{\GE}(x)\le C(n)$ for each vertex $x$ of the graph $\GE$.
%----------------------------------------------------------
\end{proposition}
%----------------------------------------------------------
%@@@@@@@@@@@@@@@@@@@@@@@@@@@@@@@@@@@@@@@@@@@@@@@@@@@@@@@@@@
%@@@@@@@@@@@@@@@@@@@@@@@@@@@@@@@@@@@@@@@@@@@@@@@@@@@@@@@@@@
%----------------------------------------------------------
\par {\it Proof.} (i). We let $\ED$ denote the family of edges of the graph $\GE$. Let $u\in\ED$ and let $A_u$ and $A'_u$ be the ends of $u$. (We write $u=(A_u,A'_u)$.) Thus $A_u\ne A'_u$ and $A_u\lr A'_u$ in $\GE$.
%----------------------------------------------------------
\par By Definition \reff{GRAPH-GE}, there exit distinct contacting lacunae $L_u,L'_u\in\LE$, $L_u\lcr L'_u$, and cubes $Q_u\in L_u, Q'_u\in L'_u$ such that
%---------------------------------------------------------
\bel{AU}
A_u=\PRL(L_u)~~~\text{and}~~~A'_u=\PRL(L'_u)
\ee
%----------------------------------------------------------
and $Q_u\cap Q'_u\ne\emp$.
%----------------------------------------------------------
\par Since $A_u=\PRL(L_u)\ne\PRL(L'_u)=A'_u,$
by part (ii) of Theorem \reff{L-PE},
%----------------------------------------------------------
\bel{DQ-A}
\diam Q_u+\diam Q'_u\le \tgm\,\|A_u-A'_u\|
\ee
%----------------------------------------------------------
where $\tgm\ge 1$ is an absolute constant. Note that, by part (i) of this theorem, $A_u\in \tgm Q$ and
$A'_u\in \tgm Q'$. Since $Q_u\cap Q'_u\ne\emp$, by part (i) of Lemma \reff{Wadd}, $Q'_u\subset 5\,\tgm Q$. Hence,
%----------------------------------------------------------
\bel{5-QA}
A_u,A'_u\subset 5\,\tgm Q.
\ee
%----------------------------------------------------------
Furthermore, by \rf{DQ-A},
%----------------------------------------------------------
\bel{C-DS}
\diam Q_u\le \tgm\,\|A_u-A'_u\|.
\ee
%----------------------------------------------------------
\par Now let us consider a mapping $T:\ED\to W_E$ defined by
%---------------------------------------------------------
$$
T(u):=Q_u,~~~u\in\ED.
$$
%----------------------------------------------------------
\par In general, this mapping is not one-to-one. Prove that $T$ is ``almost'' one-to-one, i.e., for each $Q\in W_E$ the set of its sources
%---------------------------------------------------------
$$
\TC(Q):=\{u\in\ED: T(u)=Q\}
$$
%----------------------------------------------------------
has the cardinality
%----------------------------------------------------------
\bel{CSQ}
\#\,\TC(Q)\le M=M(n).
\ee
%----------------------------------------------------------
\par Let $u=(A_u,A'_u)\in \TC(Q)$. In other words, $u\in\ED$ is an edge of the graph $\GE$ with the ends at points $A_u$ and $A'_u$.
%----------------------------------------------------------
\par Let $L^{(Q)}\in\LE$ be a lacuna
containing $Q$. Then, by \rf{AU},
%----------------------------------------------------------
$$
A_u=\PRL(L^{(Q)})~~~\text{for every}~~~u\in \TC(Q)
$$
%----------------------------------------------------------
proving that {\it $A_u$ depend only on $Q$} and does not depend on $u\in \TC(Q)$. We denote this common value of $A_u$ by $A^{(Q)}$.
%----------------------------------------------------------
\par Thus $u=(A^{(Q)}, A'_u)$ for every $u\in \TC(Q)$. We know that $A'_u=\PRL(L'_u)$ where $L'_u\in\LE$ is a lacuna contacting to $L^{(Q)}$ ($L'_u\lcr L^{(Q)}$), see Definition \reff{CONT-L}. But, by Proposition \reff{D-LQN}, the number of such lacunae is bounded by a constant $M=M(n)$ proving the required estimate \rf{CSQ}.\medskip
%----------------------------------------------------------
\par Let $M_Q:=\#\,\TC(Q)$. We know that $M_Q\le M$ for every $Q\in W_E$. Let us enumerate the elements of $\TC(Q)$:
%----------------------------------------------------------
$$
\TC(Q)=\{u_Q^{(1)},...,u_Q^{(M_Q)}\}.
$$
%----------------------------------------------------------
\par Consider a partition of $Q$ into $M^n$ equal cubes $\{H_Q^{(1)},...,H_Q^{(M^n)}\}$ of diameter $\diam Q/M$.
To each edge $u=u_Q^{(i)}\in\TC(Q)$ we assign a cube
%----------------------------------------------------------
$$
K_u:=\tfrac12\, H_Q^{(i)}.
$$
%----------------------------------------------------------
\par Let $u=(A_u,A_u')$, i.e., the points $A_u,A'_u\in E$ are the ends of the edge $u$. Clearly, $\diam K_u\le\diam Q$ so that, by \rf{C-DS},
%----------------------------------------------------------
$$
\diam K_u\le\tgm\|A_u-A'_u\|.
$$
%---------------------------------------------------------
It is also clear that $Q\subset \lambda K_u$ for some constant $\lambda=\lambda(n)\ge 1$. Hence, by \rf{5-QA}, $A_u,A'_u\in 5\lambda\,\tgm K_u$.
%----------------------------------------------------------
\par Finally, it remains to note that the cubes of the family
%----------------------------------------------------------
$$
\{\tfrac12 H_Q^{(i)}: Q\in W_E, 1\le i\le M_Q\}
$$
%---------------------------------------------------------
are pairwise disjoint proving that the cubes of the family $\{\Kc_u: u\in\ED\}$ are pairwise disjoint as well.
%----------------------------------------------------------
\par These properties of the family $\{\Kc_u: u\in\ED\}$ show that the graph $\GE$ satisfies conditions of Definition \reff{DF-GRPR} with a constant $\gamma=5\lambda\,\tgm$ proving that $\GE$ is a $\gamma$-sparse graph.\medskip
%----------------------------------------------------------
%@@@@@@@@@@@@@@@@@@@@@@@@@@@@@@@@@@@@@@@@@@@@@@@@@@@@@@@@@@
%@@@@@@@@@@@@@@@@@@@@@@@@@@@@@@@@@@@@@@@@@@@@@@@@@@@@@@@@@@
%@@@@@@@@@@@@@@@@@@@@@@@@@@@@@@@@@@@@@@@@@@@@@@@@@@@@@@@@@@
%----------------------------------------------------------
\par (ii). Let $A\in E$. We have to prove that the set $\AD(A):=\{A'\in E: A'\lr A\}$ of vertices adjacent to $A$ consists of at most $C=C(n)$ elements.
%----------------------------------------------------------
\par By Definition \reff{GRAPH-GE}, for every $A'\in \AD(A)$ there exist contacting lacunae $L,L'\in\LE$ ($L\lcr L'$) such that $A=\PRL(L)$ and $A'=\PRL(L')$ so that
$\#\AD(A)$ is bounded by the cardinality of the set of lacunae
%---------------------------------------------------------
$$
I=\{L'\in\LE: A'=\PRL(L')~\text{and}~~\exists~L\in\LE~\text{such that}~ L\lcr L'~\text{and}~\PRL(L)=A\}.
$$
%----------------------------------------------------------
\par But, by part (iii) of Theorem \reff{L-PE},
$\#\{L\in\Lc_E:\PRL(L)=A\}\le C_1(n)$, and, by Proposition \reff{D-LQN}, $\#\,\{L'\in\LE:L'\lcr L\}\le C_2(n)$. Hence
%---------------------------------------------------------
$$
\deg_{\GE}(A)=\#\AD(A)\le\#I\le C_1(n)C_2(n)
$$
%----------------------------------------------------------
proving part (ii) of the proposition.\smallskip
%----------------------------------------------------------
\par The proof of the proposition is complete.\bx\bigskip
%----------------------------------------------------------
%@@@@@@@@@@@@@@@@@@@@@@@@@@@@@@@@@@@@@@@@@@@@@@@@@@@@@@@@@@
%@@@@@@@@@@@@@@@@@@@@@@@@@@@@@@@@@@@@@@@@@@@@@@@@@@@@@@@@@@
%@@@@@@@@@@@@@@@@@@@@@@@@@@@@@@@@@@@@@@@@@@@@@@@@@@@@@@@@@@
%@@@@@@@@@@@@@@@@@@@@@@@@@      @@@@@@@@@@@@@@@@@@@@@@@@@@@
%@@@@@@@@@@@@@@@@@@@@@@@          @@@@@@@@@@@@@@@@@@@@@@@@@
%@@@@@@@@@@@@@@@@@@@@@              @@@@@@@@@@@@@@@@@@@@@@@
%@@@@@@@@@@@@@@@@@@@     SECTION 7    @@@@@@@@@@@@@@@@@@@@@
%@@@@@@@@@@@@@@@@@@@@@              @@@@@@@@@@@@@@@@@@@@@@@
%@@@@@@@@@@@@@@@@@@@@@@@          @@@@@@@@@@@@@@@@@@@@@@@@@
%@@@@@@@@@@@@@@@@@@@@@@@@@      @@@@@@@@@@@@@@@@@@@@@@@@@@@
%@@@@@@@@@@@@@@@@@@@@@@@@@@@@@@@@@@@@@@@@@@@@@@@@@@@@@@@@@@
%@@@@@@@@@@@@@@@@@@@@@@@@@@@@@@@@@@@@@@@@@@@@@@@@@@@@@@@@@@
%@@@@@@@@@@@@@@@@@@@@@@@@@@@@@@@@@@@@@@@@@@@@@@@@@@@@@@@@@@
%----------------------------------------------------------
\SECT{7. The lacunary extension operator: a proof of the  variational criterion.}{7}
%----------------------------------------------------------
\addtocontents{toc}{7. The variational criterion: a proof of Theorem \reff{EX-TK}.\hfill\thepage\par\VST}
%----------------------------------------------------------
%@@@@@@@@@@@@@@@@@@@@@@@@@@@@@@@@@@@@@@@@@@@@@@@@@@@@@@@@@@
%----------------------------------------------------------
\par {\bf 7.1. The variational criterion: necessity.}\medskip
%----------------------------------------------------------
\addtocontents{toc}{~~~~7.1. The variational criterion: necessity. \hfill \thepage\par}
%----------------------------------------------------------
\par As we have mentioned above, in \cite{Sh2} we have proved a criterion which provides a cha\-rac\-terization of the trace space $\LOP|_E$ it terms of certain local oscillations of functions on subsets of the set $E$. Theorem \reff{EX-TK} which we prove in this section  refines and generalizes this criterion to the case of jet-spaces generated by $L^m_p(\RN)-$functions, $m\ge 1$, $p>n$.
%----------------------------------------------------------
%@@@@@@@@@@@@@@@@@@@@@@@@@@@@@@@@@@@@@@@@@@@@@@@@@@@@@@@@@@
%@@@@@@@@@@@@@@@@@@@@@@@@@@@@@@@@@@@@@@@@@@@@@@@@@@@@@@@@@@
%----------------------------------------------------------
\smallskip
\par {\it Proof.} {\it (Necessity.)} We prove a slightly more general result which immediately implies the necessity part of Theorem \reff{EX-TK}.
%----------------------------------------------------------
%@@@@@@@@@@@@@@@@@@@@@@@@@@@@@@@@@@@@@@@@@@@@@@@@@@@@@@@@@@
%@@@@@@@@@@@@@@@@@@@@@@@@@@@@@@@@@@@@@@@@@@@@@@@@@@@@@@@@@@
%@@@@@@@@@@@@@@@@@@@@@@@@@@@@@@@@@@@@@@@@@@@@@@@@@@@@@@@@@@
%----------------------------------------------------------
%@@@@@@@@@@@@@@@@@@@@@@@@@@@@@@@@@@@@@@@@@@@@@@@@@@@@@@@@@@
\begin{proposition}\lbl{N-VC} Let $n<p<\infty$. Let $F\in C^{m-1}(\RN)\cap\LMP$ and let $P_x=T^{m-1}_x[F]$, $x\in E$. %----------------------------------------------------------
\par Then for every constant $\gamma\ge 1$, every family $\{Q_i:i\in I\}$ of pairwise disjoint cubes in $\RN$, every collection of points $x_i,y_i\in (\gamma Q_i)\cap E$, $i\in I$, and every multiindex $\beta, |\beta|\le m-1,$ the following inequality
%----------------------------------------------------------
$$
\sbig_{i\in I}\,\,\frac{|D^\beta P_{x_i}(x_i)
-D^\beta P_{y_i}(x_i)|^p}
{(\diam Q_i)^{(m-|\beta|)p-n}} \le C\,\|F\|^p_{\LMP}
$$
%----------------------------------------------------------
holds. Here $C$ is a constant depending only on $m,n,p$ and $\gamma$.
%----------------------------------------------------------
\end{proposition}
%----------------------------------------------------------
%@@@@@@@@@@@@@@@@@@@@@@@@@@@@@@@@@@@@@@@@@@@@@@@@@@@@@@@@@@
%@@@@@@@@@@@@@@@@@@@@@@@@@@@@@@@@@@@@@@@@@@@@@@@@@@@@@@@@@@
%@@@@@@@@@@@@@@@@@@@@@@@@@@@@@@@@@@@@@@@@@@@@@@@@@@@@@@@@@@
%----------------------------------------------------------
\par {\it Proof.} Let $q=(n+p)/2$. Let $Q$ be a cube in $\RN$ and let $K:=\gamma Q$. Fix two points $x,y\in K$. Recall that $Q_{xy}=Q(x,\|x-y\|)$ so that $Q_{xy}\subset 2K$. By inequality \rf{SP-BTQ}, for every multiindex $\beta$, $|\beta|\le m-1,$
%----------------------------------------------------------
\be
|D^{\beta}P_x(x)-D^{\beta}P_y(x)|&\le& C\,\|x-y\|^{m-|\beta|}\left(\frac{1}{|Q_{xy}|}
\intl_{Q_{xy}}(\nabla ^mF(u))^qdu\right)^{\frac{1}{q}}\nn\\
&\le&
C\,\|x-y\|^{m-|\beta|-\tfrac{n}{q}}
\left(\,\intl_{Q_{xy}}(\nabla ^mF(u))^qdu\right)^{\frac{1}{q}}\nn.
\ee
%----------------------------------------------------------
\par Since $|\beta|\le m-1$ and $n<q$, we have   $m-|\beta|-\tfrac{n}{q}>0$. Hence
%----------------------------------------------------------
$$
|D^{\beta}P_x(x)-D^{\beta}P_y(x)|\le C\,
(\diam Q)^{m-|\beta|} \left(\frac{1}{|2K|}
\intl_{2K}(\nabla ^mF(u))^qdu\right)^{\frac{1}{q}}
$$
%----------------------------------------------------------
where $C$ is a constant depending only on $m,n,p$ and $\gamma$. By this inequality,
%----------------------------------------------------------
$$
|D^{\beta}P_x(x)-D^{\beta}P_y(x)|^p\le C\,
(\diam Q)^{p(m-|\beta|)} \left(\Mc[\nabla ^m F^q](z)\right)^{\frac{p}{q}}
$$
%----------------------------------------------------------
for arbitrary $z\in Q$. Integrating this inequality over $Q$ (with respect to $z$) we obtain the following:
%----------------------------------------------------------
$$
\frac{|D^{\beta}P_x(x)-D^{\beta}P_y(x)|^p}
{(\diam Q)^{(m-|\beta|)p-n}}\le C\,\intl_{Q}
\left(\Mc[\nabla ^m F^q](z)\right)^{\frac{p}{q}}\,dz.
$$
%----------------------------------------------------------
\par Hence,
%----------------------------------------------------------
\be
I_\beta:=\sbig_{i\in I}\,\,\frac{|D^\beta P_{x_i}(x_i)
-D^\beta P_{y_i}(x_i)|^p}
{(\diam Q_i)^{(m-|\beta|)p-n}}& \le&
C\,\sbig_{i\in I}\,
\intl_{Q_i}
\left(\Mc[\nabla ^m F^q](z)\right)^{\frac{p}{q}}\,dz\nn\\
&\le& C\,\intl_{\RN}
\left(\Mc[\nabla ^m F^q](z)\right)^{\frac{p}{q}}\,dz\nn
\ee
%----------------------------------------------------------
so that, by the Hardy-Littlewood maximal theorem,
%----------------------------------------------------------
$$
I_\beta\le C\,\intl_{\RN}
(\nabla ^m F)^p(z)\,dz\,\le C\,\|F\|^p_{\LMP}
$$
%----------------------------------------------------------
proving the proposition.\bx
%----------------------------------------------------------
%@@@@@@@@@@@@@@@@@@@@@@@@@@@@@@@@@@@@@@@@@@@@@@@@@@@@@@@@@@
%@@@@@@@@@@@@@@@@@@@@@@@@@@@@@@@@@@@@@@@@@@@@@@@@@@@@@@@@@@
%@@@@@@@@@@@@@@@@@@@@@@@@@@@@@@@@@@@@@@@@@@@@@@@@@@@@@@@@@@
%@@@@@@@@@@@@@@@@@@@@@@@@@@@@@@@@@@@@@@@@@@@@@@@@@@@@@@@@@@
%@@@@@@@@@@@@@@@@@@@@@@@@@@@@@@@@@@@@@@@@@@@@@@@@@@@@@@@@@@
%@@@@@@@@@@@@@@@@@@@@@@@@@@@@@@@@@@@@@@@@@@@@@@@@@@@@@@@@@@
%----------------------------------------------------------
\par Let us prove the necessity part of the theorem. Let  $\VP=\{P_x: x\in E\}$ be a family of polynomials of degree at most $m-1$ indexed by points of $E$. Suppose there exists a $C^{m-1}$-function $F\in\LMP$ such that $T_{x}^{m-1}[F]=P_{x}$ for every $x\in E$.
%----------------------------------------------------------
\par Let $\gamma\ge 1$ and let $\{\{x_i,y_i\}: i=1,...,k\}$ be an arbitrary finite $\gamma$-sparse collection of two point subsets of $E$. Then, by Definition \reff{DF-PR}, there exists a family $\{Q_i, i=1,...,k\}$ of pairwise disjoint cubes in $\RN$ satisfying condition \rf{XY-Q}, so that, by Proposition \reff{N-VC},
%----------------------------------------------------------
$$
\smed_{i=1}^k\,\,\smed_{|\alpha|\le m-1}\,\,\frac{|D^\alpha P_{x_i}(x_i)-D^\alpha P_{y_i}(x_i)|^p}
{(\diam Q_i)^{(m-|\alpha|)p-n}}\le C\,\|F\|^p_{\LMP}.
$$
%----------------------------------------------------------
Here $C$ is a constant depending only on $m,n,p,$ and $\gamma$.
%----------------------------------------------------------
\par But, by \rf{XY-Q}, $\diam Q_i\le\gamma\|x_i-y_i\|$, $i=1,...,k$, so that
%----------------------------------------------------------
$$
\smed_{i=1}^k\,\,\smed_{|\alpha|\le m-1}\frac{|D^\alpha P_{x_i}(x_i)-D^\alpha P_{y_i}(x_i)|^p}
{\|x_i-y_i\|^{(m-|\alpha|)p-n}}\le \gamma^{mp-n}\,C\,\|F\|^p_{\LMP}\,.
$$
%----------------------------------------------------------
Hence $\Nc_{m,p,E}(\VP)\le \gamma^{m-n/p}\,C^{1/p}\,\|F\|_{\LMP}$, see \rf{N-P-NEW}, and the proof of the necessity is complete.
%----------------------------------------------------------
%@@@@@@@@@@@@@@@@@@@@@@@@@@@@@@@@@@@@@@@@@@@@@@@@@@@@@@@@@@
%@@@@@@@@@@@@@@@@@@@@@@@@@@@@@@@@@@@@@@@@@@@@@@@@@@@@@@@@@@
%@@@@@@@@@@@@@@@@@@@@@@@@@@@@@@@@@@@@@@@@@@@@@@@@@@@@@@@@@@
%@@@@@@@@@@@@@@@@@@@@@@@@@@@@@@@@@@@@@@@@@@@@@@@@@@@@@@@@@@
%@@@@@@@@@@@@@@@@@@@@@@@@@@@@@@@@@@@@@@@@@@@@@@@@@@@@@@@@@@
%@@@@@@@@@@@@@@@@@@@@@@@@@@@@@@@@@@@@@@@@@@@@@@@@@@@@@@@@@@
%@@@@@@@@@@@@@@@@@@@@@@@@@@@@@@@@@@@@@@@@@@@@@@@@@@@@@@@@@@
%@@@@@@@@@@@@@@@@@@@@@@@@@@@@@@@@@@@@@@@@@@@@@@@@@@@@@@@@@@
%@@@@@@@@@@@@@@@@@@@@@@@@@@@@@@@@@@@@@@@@@@@@@@@@@@@@@@@@@@
%@@@@@@@@@@@@@@@@@@@@@@@@@@@@@@@@@@@@@@@@@@@@@@@@@@@@@@@@@@
%----------------------------------------------------------
\bigskip\medskip
%----------------------------------------------------------
\par {\bf 7.2. The variational criterion: sufficiency.} \medskip
%----------------------------------------------------------
\addtocontents{toc}{~~~~7.2. The variational criterion: a proof of the sufficiency. \hfill \thepage\par}
%----------------------------------------------------------
\par Let $\gamma:=10^4\tgm$ where $\tgm$ is the constant from Theorem \reff{L-PE}. Let  $\VP=\{P_x\in\PMRN:x\in E\}$ be a Whitney $(m-1)$-field on $E$.
%----------------------------------------------------------
\par Let $\lambda:=\Nc_{m,p,E}(\VP)^p$, see \rf{N-P-NEW}.  Suppose that $\lambda<\infty$. Then, by \rf{N-P-NEW}, for every finite $\gamma$-sparse collection $\{\{x_i,y_i\}: i=1,...,k\}$ of two point subsets of $E$ the following inequality
%----------------------------------------------------------
\bel{N-P-NEW-1}
\smed_{i=1}^k\,\,
\smed_{|\alpha|\le m-1}
\frac{|D^\alpha P_{x_i}(x_i)-D^\alpha P_{y_i}(x_i)|^p}
{\|x_i-y_i\|^{(m-|\alpha|)p-n}}\le \lambda
\ee
%----------------------------------------------------------
holds.
%----------------------------------------------------------
%@@@@@@@@@@@@@@@@@@@@@@@@@@@@@@@@@@@@@@@@@@@@@@@@@@@@@@@@@@
%@@@@@@@@@@@@@@@@@@@@@@@@@@@@@@@@@@@@@@@@@@@@@@@@@@@@@@@@@@
%----------------------------------------------------------
\smallskip
\par Let us introduce {\it the lacunary modification of the Whitney extension method.} Let $L\in\LE$ be a lacuna. For every cube $Q\in L$ we put
%----------------------------------------------------------
\bel{LM-WM-DEF}
a_Q:=\PRL(L).
\ee
%----------------------------------------------------------
Here $\PRL:\LE\to E$ is the ``projector'' from Theorem \reff{L-PE}.
%----------------------------------------------------------
\par Note that, by property (i) of this theorem,
%----------------------------------------------------------
\bel{LM-WM}
a_Q=\PRL(L)\in \tgm Q~~~\text{for every}~~Q\in L.
\ee
%----------------------------------------------------------
Hence, by \rf{IN-GQ}, for each $Q\in L$
%----------------------------------------------------------
\bel{5-5N}
a_Q~~~\text{is a}~~\vs-\text{nearest point to}~~Q~~
\text{with}~~\vs=(\tgm+1)/2.
\ee
%----------------------------------------------------------
%@@@@@@@@@@@@@@@@@@@@@@@@@@@@@@@@@@@@@@@@@@@@@@@@@@@@@@@@@@
%@@@@@@@@@@@@@@@@@@@@@@@@@@@@@@@@@@@@@@@@@@@@@@@@@@@@@@@@@@
%----------------------------------------------------------
\par Then we construct a function $F$ on $\RN$ using the Whitney extension formula \rf{DEF-F}.
%----------------------------------------------------------
%@@@@@@@@@@@@@@@@@@@@@@@@@@@@@@@@@@@@@@@@@@@@@@@@@@@@@@@@@@
%@@@@@@@@@@@@@@@@@@@@@@@@@@@@@@@@@@@@@@@@@@@@@@@@@@@@@@@@@@
%----------------------------------------------------------
\smallskip
\par Let us prove that the function $F$ satisfies the following conditions:\smallskip
%----------------------------------------------------------
\par (i). $F$ is a $C^{m-1}$-function such that  $T_x^{m-1}[F]=P_x$ for every $x\in E$\,;\smallskip
%----------------------------------------------------------
\par (ii). $F\in\LMP$ and
%----------------------------------------------------------
$$
\|F\|_{\LMP}\le C\,\lambda^{\frac1p}
$$
%----------------------------------------------------------
where $C$ is a constant depending only on $m,n$ and $p$.
\medskip
%-------------------------------------------------------
\par Prove (i). For every multiindex $\beta, |\beta|\le m-1,$ and every $x,y\in E$, by \rf{N-P-NEW-1},
%----------------------------------------------------------
$$
\frac{|D^\beta P_{x}(x)-D^\beta P_{y}(x)|^p}
{\|x-y\|^{(m-|\beta|)p-n}} \le\lambda.
$$
%----------------------------------------------------------
Hence
%----------------------------------------------------------
$$
|D^\beta P_{x}(x)-D^\beta P_{y}(x)|
\le C\,\lambda^{\frac1p}\|x-y\|^{m-|\beta|-1}\cdot
\|x-y\|^{1-\frac{n}{p}}.
$$
%----------------------------------------------------------
Since $n<p$,
%----------------------------------------------------------
$$
D^\beta P_{x}(x)-D^\beta P_{y}(x)
=o(\|x-y\|^{m-|\beta|-1})~~\text{as}~~y\to x,~y\in E.
$$
%----------------------------------------------------------
\par Thus the Whitney $(m-1)$-field $\VP=\{P_x\in\PMRN:x\in E\}$ satisfies the hypothesis of the Whitney extension theorem \cite{W1}. Recall that in Section 4 we have constructed the function $F$ using a modification of the Whitney extension method where for each $Q\in W_E$ the point $a_Q$ is $\vs$-nearest to $Q$ with $\vs=5.5$. See \rf{5-5N}. As we have noted in Section 5, see a remark after \rf{W-EF}, such a modification provides the statement (i) as well.
%----------------------------------------------------------
%@@@@@@@@@@@@@@@@@@@@@@@@@@@@@@@@@@@@@@@@@@@@@@@@@@@@@@@@@@
%@@@@@@@@@@@@@@@@@@@@@@@@@@@@@@@@@@@@@@@@@@@@@@@@@@@@@@@@@@
%----------------------------------------------------------
\smallskip
\par We turn to the proof of statement (ii). Let us fix a multiindex  $\beta$ of order $|\beta|=m-1$ and prove that %----------------------------------------------------------
\bel{DB-ST}
D^\beta F\in \LOP~~~\text{and}~~~\|D^\beta F\|_{\LOP}\le C\,\lambda^{\frac1p}.
\ee
%----------------------------------------------------------
The proof of this inequality is based on two auxiliary statements. The first of them is the following combinatorial
%---------------------------------------------------------
%@@@@@@@@@@@@@@@@@@@@@@@@@@@@@@@@@@@@@@@@@@@@@@@@@@@@@@@@@@
%@@@@@@@@@@@@@@@@@@@@@@@@@@@@@@@@@@@@@@@@@@@@@@@@@@@@@@@@@@
\begin{theorem}(\cite{BrK,Dol})\lbl{TFM} Let $\Ac=\{Q\}$ be a collection of cubes in $\RN$ with covering
multiplicity $ \MP(\Ac)<\infty$. Then $\Ac$ can be
partitioned into at most $N=2^{n-1}(\MP(\Ac)-1)+1$ families of disjoint cubes.
%----------------------------------------------------------
\end{theorem}
%----------------------------------------------------------
%@@@@@@@@@@@@@@@@@@@@@@@@@@@@@@@@@@@@@@@@@@@@@@@@@@@@@@@@@@
%@@@@@@@@@@@@@@@@@@@@@@@@@@@@@@@@@@@@@@@@@@@@@@@@@@@@@@@@@@
%@@@@@@@@@@@@@@@@@@@@@@@@@@@@@@@@@@@@@@@@@@@@@@@@@@@@@@@@@@
%----------------------------------------------------------
\par Recall that {\it covering multiplicity} $\MP(\Ac)$ of a family of cubes $\Ac$ is the minimal positive integer $M$ such that every point $x\in\RN$ is covered by at most $M$ cubes from $\Ac$.
%----------------------------------------------------------
%@@@@@@@@@@@@@@@@@@@@@@@@@@@@@@@@@@@@@@@@@@@@@@@@@@@@@@@@@@
%@@@@@@@@@@@@@@@@@@@@@@@@@@@@@@@@@@@@@@@@@@@@@@@@@@@@@@@@@@
\par The second auxiliary result which we are needed for the proof of the sufficiency is a certain variational description of the space $\LOP$ in terms of local oscillations. This result follows from a description of Sobolev spaces obtained in \cite{Br}.; see there \S\,4, subsection $3^{\circ}$.
%----------------------------------------------------------
%@@@@@@@@@@@@@@@@@@@@@@@@@@@@@@@@@@@@@@@@@@@@@@@@@@@@@@@@@@
%@@@@@@@@@@@@@@@@@@@@@@@@@@@@@@@@@@@@@@@@@@@@@@@@@@@@@@@@@@
%----------------------------------------------------------
\begin{theorem}\lbl{CR-SOB} Let $p>n$ and let $\tau>0$. Let $G$ be a continuous function on $\RN$ satisfying the following condition: There exists a constant $A>0$ such that for every finite family $\{Q_i: i=1,...,k\}$ of pairwise disjoint equal cubes in $\RN$ of diameter $\diam Q_i\le\tau$ and every $x_i\in Q_i$ the following inequality
%----------------------------------------------------------
$$
\sbig_{i=1}^k\,\frac{|\,G(x_i)-G(c_{Q_i})|^p}{(\diam Q_i)^{p-n}} \le A
$$
%----------------------------------------------------------
holds. Then $G\in\LOP$ and $\|G\|_{\LOP}\le C(n,p)\,A^{\frac1p}$.
%----------------------------------------------------------
%@@@@@@@@@@@@@@@@@@@@@@@@@@@@@@@@@@@@@@@@@@@@@@@@@@@@@@@@@@
%@@@@@@@@@@@@@@@@@@@@@@@@@@@@@@@@@@@@@@@@@@@@@@@@@@@@@@@@@@
%----------------------------------------------------------
\end{theorem}
%----------------------------------------------------------
%@@@@@@@@@@@@@@@@@@@@@@@@@@@@@@@@@@@@@@@@@@@@@@@@@@@@@@@@@@
%@@@@@@@@@@@@@@@@@@@@@@@@@@@@@@@@@@@@@@@@@@@@@@@@@@@@@@@@@@
%----------------------------------------------------------
\par {\it Proof.} For the case $\tau=\infty$ this criterion follows from a description of Sobolev spaces obtained in \cite{Br}.; see there \S\,4, subsection $3^{\circ}$.
%----------------------------------------------------------
\par Prove the result for $0<\tau<\infty$. Using a result from \cite{Br} related to an atomic decomposition of the modulus of smoothness in $L_p$, see there \S\,2, subsection $2^{\circ}$, Theorem 4, and the theorem's hypo\-the\-sis, we conclude that the first modulus of continuity of $G$ in  $\LPRN$, the function $\Omega_1(G,t)_{\LPRN}$, satisfies the following inequality:
%----------------------------------------------------------
$$
\Omega_1(G,t)_{\LPRN}\le C(n,p)\,A^{\frac1p}\,t~~~\text{for every}~~
t\in(0,\tau].
$$
%----------------------------------------------------------
But the function $\Omega_1(G,t)_{\LPRN}/t$ is a quasi-monotone function on $\R_+$, i.e.,
%----------------------------------------------------------
$$
\Omega_1(G,t_2)_{\LPRN}/t_2\le C\,
\Omega_1(G,t_1)_{\LPRN}/t_1
~~~\text{for every}~~0<t_1<t_2,
$$
%----------------------------------------------------------
where $C>1$ is an absolute constant.
%----------------------------------------------------------
\par Hence,
%----------------------------------------------------------
$$
\Omega_1(G,t)_{\LPRN}/t\le C\,
\Omega_1(G,\tau)_{\LPRN}/\tau\le C(n,p)\,A^{\frac1p}
~~~\text{for every}~~t>\tau
$$
%----------------------------------------------------------
so that
%----------------------------------------------------------
$$
\Omega_1(G,t)_{\LPRN}\le C(n,p)\,A^{\frac1p}\,t~~~\text{for every}~~t>0.
$$
%----------------------------------------------------------
\par It is shown in \cite{Br}, \S\,4, Subsection $3^{\circ}$, that this property implies that $G\in\LOP$ and that $\|G\|_{\LOP}\le C(n,p)\,A^{\frac1p}$. The proof of  the lemma is complete.\bx
\bigskip
%----------------------------------------------------------
%@@@@@@@@@@@@@@@@@@@@@@@@@@@@@@@@@@@@@@@@@@@@@@@@@@@@@@@@@@
%@@@@@@@@@@@@@@@@@@@@@@@@@@@@@@@@@@@@@@@@@@@@@@@@@@@@@@@@@@
%@@@@@@@@@@@@@@@@@@@@@@@@@@@@@@@@@@@@@@@@@@@@@@@@@@@@@@@@@@
%@@@@@@@@@@@@@@@@@@@@@@@@@@@@@@@@@@@@@@@@@@@@@@@@@@@@@@@@@@
%@@@@@@@@@@@@@@@@@@@@@@@@@@@@@@@@@@@@@@@@@@@@@@@@@@@@@@@@@@
%@@@@@@@@@@@@@@@@@@@@@@@@@@@@@@@@@@@@@@@@@@@@@@@@@@@@@@@@@@
%----------------------------------------------------------
\par We turn to the proof of the statement \rf{DB-ST}. We will do this by applying Theorem \reff{CR-SOB} to a function $G:=D^\beta F$. Basing on the following three lemmas we show that this function satisfies the hypothesis of Theorem \reff{CR-SOB}.\smallskip
%----------------------------------------------------------
\par Let $\Qc=\{Q_1,...,Q_k\}$ be a family of pairwise disjoint equal cubes in $\RN$. Let $c_i:=c_{Q_i}$, $i=1,...,k,$ be the center of the cube $Q_i$.
%----------------------------------------------------------
%@@@@@@@@@@@@@@@@@@@@@@@@@@@@@@@@@@@@@@@@@@@@@@@@@@@@@@@@@@
%@@@@@@@@@@@@@@@@@@@@@@@@@@@@@@@@@@@@@@@@@@@@@@@@@@@@@@@@@@
%----------------------------------------------------------
\begin{lemma}\lbl{L-BG} Suppose that %----------------------------------------------------------
\bel{D-L40}
\dist(c_i,E)\le 40\,\diam Q_i~~~\text{for all}~~~i=1,...,k.
\ee
%----------------------------------------------------------
\par Then for every $\beta, |\beta|=m-1$, and every $x_i\in Q_i$ the following inequality
%----------------------------------------------------------
$$
\sbig_{i=1}^k\,
\frac{|D^\beta F(x_i)-D^\beta F(c_i)|^p}
{(\diam Q_i)^{p-n}}
\le C\,\lambda\,.
$$
%----------------------------------------------------------
Here $C>0$ is a constant depending only on $m,n,$ and $p$. %----------------------------------------------------------
\end{lemma}
%----------------------------------------------------------
%@@@@@@@@@@@@@@@@@@@@@@@@@@@@@@@@@@@@@@@@@@@@@@@@@@@@@@@@@@
%@@@@@@@@@@@@@@@@@@@@@@@@@@@@@@@@@@@@@@@@@@@@@@@@@@@@@@@@@@
%@@@@@@@@@@@@@@@@@@@@@@@@@@@@@@@@@@@@@@@@@@@@@@@@@@@@@@@@@@
%----------------------------------------------------------
\par {\it Proof.} First we will prove the lemma for the case $c_i\in E$, $i=1,...,k$. Let
%----------------------------------------------------------
\bel{I1-DF}
I_1(F):=\sbig\,\left\{
\frac{|D^\beta F(x_i)-D^\beta F(c_i)|^p}
{(\diam Q_i)^{p-n}}: i\in\{1,...,k\},\,x_i\in E,~x_i\ne c_i\right\}.
\ee
%----------------------------------------------------------
Prove that $I_1(F)\le \lambda$.
%----------------------------------------------------------
\par Let us consider a collection of two point sets
%----------------------------------------------------------
$$
\Ac:=\{\{x_i,c_i\}: i=1,...,k,~x_i\in E, x_i\ne c_i\}.
$$
%----------------------------------------------------------
Prove that $\Ac$ is a $1$-sparse collection, see Definition \reff{DF-PR}. In fact, since $x_i,c_i\in Q_i$, there exists a cube $K_i$ such that $x_i,c_i\in K_i\subset Q_i$ and $\diam K_i=\|x_i-c_i\|$. Since the cubes $\{Q_i\}$ are pairwise disjoint, the same is true for the cubes $\{K_i\}$ proving that $\Ac$ is $1$-sparse. By \rf{J-PE}, $D^\beta F(x_i)=D^\beta P_{x_i}(x_i)$ and
$D^\beta F(c_i)=D^\beta P_{c_i}(c_i)$. But $D^\beta P_{x_i}$ is a constant function whenever $|\beta|=m-1$ so that $D^\beta F(x_i)=D^\beta P_{x_i}(c_i)$. Hence, by the theorems hypothesis (with $|\beta|=m-1$), see \rf{N-P-NEW-1},
%----------------------------------------------------------
\bel{I1-FEG}
I_1(F)\le\sbig\,\left\{
\frac{|D^\beta P_{x_i}(c_i)-D^\beta P_{c_i}(c_i)|^p}
{\|x_i-c_i\|^{p-n}}: i\in\{1,...k\},\,x_i\in E,~x_i\ne c_i\right\}\le \lambda.
\ee
%----------------------------------------------------------
\par Let
%----------------------------------------------------------
\bel{I2-DF}
I_2(F):=\sbig\,\left\{
\frac{|D^\beta F(x_i)-D^\beta F(c_i)|^p}
{(\diam Q_i)^{p-n}}: i\in\{1,...k\},\,x_i\in \RN\setminus E\right\}.
\ee
%----------------------------------------------------------
\par Let
%----------------------------------------------------------
$$
\Qc:=\{Q_i: x_i\in \RN\setminus E\}.
$$
%----------------------------------------------------------
This family consist of at most $k$ pairwise disjoint cubes. %----------------------------------------------------------
\par We recall that by $\tgm$ we denote the constant from Theorem \reff{L-PE}. Let us introduce a family of cubes
%----------------------------------------------------------
$$
\tQc:=\{5\tgm\,Q: Q\in\Qc\}.
$$
%----------------------------------------------------------
Clearly, the covering multiplicity $\MC(\tQc)\le C(n)$ so that, by Theorem \reff{TFM}, this family can be partitioned into at most $N(n)$ subfamilies of pairwise disjoint cubes. Therefore in this proof without loss of generality we can assume that the family $\tQc$ itself consists of pairwise disjoint cubes.
%----------------------------------------------------------
\par Fix a cube $Q\in \Qc$. Thus $Q=Q_i$ for some $i\in\{1,...,k\}$. Let $x_Q=x_i$. Thus $x_Q\notin E$ and $c_Q\in E$ for every $Q\in\Qc$. Note also that
$D^\beta F(c_Q)=D^\beta P_{c_Q}(c_Q)$.
%----------------------------------------------------------
\par Since $x_Q\notin E$, there exists a Whitney cube
%----------------------------------------------------------
\bel{KQ-DF1}
K_Q\in W_E~~~\text{which contains}~~~x_Q.
\ee
%----------------------------------------------------------
Recall that given a lacuna $L\in\LE$ and a cube $H\in L$ we have $a_H=\PRL(L)$, see \rf{LM-WM-DEF}. We also recall that by $T(K_Q)$ we denote the family of Whitney's cubes intersecting $K_Q$. See \rf{TK}.
%----------------------------------------------------------
\par For the sake of brevity let us introduce the following notation:
%----------------------------------------------------------
$$
S(Q):=|D^{\beta}F(x_Q)-D^{\beta}F(c_Q)|
$$
%----------------------------------------------------------
and
%----------------------------------------------------------
$$
V(Q):=|D^{\beta}P_{c_Q}(c_Q)-D^{\beta}P_{a_{K_Q}}(c_Q)|.
$$
%----------------------------------------------------------
Also given $H\in T(K_Q)$ and a multiindex $\xi$ with $|\xi|\le m-1$ let
%----------------------------------------------------------
$$
L(\xi:H,Q):=|D^{\xi}P_{a_{H}}(a_{K_Q})
-D^{\xi}P_{a_{K_Q}}(a_{K_Q})|.
$$
%----------------------------------------------------------
\par Then, by Lemma \reff{C-2},
%----------------------------------------------------------
$$
S(Q)\le
C\,\left\{V(Q)+
\sbig\limits_{H\in T(K_Q)}\,\,\sbig\limits_{|\xi|\le m-1}\frac{L(\xi:H,Q)}{(\diam K_Q)^{m-1-|\xi|}}\right\}.
$$
%----------------------------------------------------------
Since $\#T(K_Q)\le N(n),$ see Lemma \reff{Wadd}, we have
%----------------------------------------------------------
$$
\frac{S(Q)^p}
{(\diam Q)^{p-n}}\le
\frac{C\,V(Q)^p}
{(\diam Q)^{p-n}}+
C\sbig\limits_{H\in T(K_Q)}\,\,
\sbig\limits_{|\xi|\le m-1}
\,\left(\frac{\diam K_Q}{\diam Q}\right)^{p-n}\frac{L(\xi:H,Q)^p}
{(\diam K_Q)^{(m-|\xi|)p-n}}.
$$
%----------------------------------------------------------
%@@@@@@@@@@@@@@@@@@@@@@@@@@@@@@@@@@@@@@@@@@@@@@@@@@@@@@@@@@
%@@@@@@@@@@@@@@@@@@@@@@@@@@@@@@@@@@@@@@@@@@@@@@@@@@@@@@@@@@
%@@@@@@@@@@@@@@@@@@@@@@@@@@@@@@@@@@@@@@@@@@@@@@@@@@@@@@@@@@
%@@@@@@@@@@@@@@@@@@@@@@@@@@@@@@@@@@@@@@@@@@@@@@@@@@@@@@@@@@
%@@@@@@@@@@@@@@@@@@@@@@@@@@@@@@@@@@@@@@@@@@@@@@@@@@@@@@@@@@
%@@@@@@@@@@@@@@@@@@@@@@@@@@@@@@@@@@@@@@@@@@@@@@@@@@@@@@@@@@
%----------------------------------------------------------
\par Prove that $a_{K_Q}\in 5\tgm\,Q$. In fact, since $x_Q\in K_Q\cap Q$, we have
%----------------------------------------------------------
\bel{KQ-Q}
\diam K_Q\le 4\dist(K_Q,E)\le 4\dist(x_Q,E)
\le 4\|x_Q-c_Q\|\le 2\diam Q
\ee
%----------------------------------------------------------
Hence, for every $y\in K_Q$,
%----------------------------------------------------------
$$
\|y-c_Q\|\le\|y-x_Q\|+\|x_Q-c_Q\|\le \diam K_Q+r_Q\le 2\diam Q+r_Q=5r_Q
$$
%----------------------------------------------------------
proving that
%----------------------------------------------------------
\bel{K-INQ}
K_Q\subset 5Q.
\ee
%----------------------------------------------------------
%@@@@@@@@@@@@@@@@@@@@@@@@@@@@@@@@@@@@@@@@@@@@@@@@@@@@@@@@@@
%----------------------------------------------------------
On the other hand, by \rf{LM-WM}, $a_{K_Q}\in \tgm K_Q$ so that
%----------------------------------------------------------
\bel{50-Q}
a_{K_Q}\in 5\tgm\,Q.
\ee
%----------------------------------------------------------
%@@@@@@@@@@@@@@@@@@@@@@@@@@@@@@@@@@@@@@@@@@@@@@@@@@@@@@@@@@
%----------------------------------------------------------
\par In particular, by inequality \rf{KQ-Q},
%----------------------------------------------------------
$$
\frac{S(Q)^p}
{(\diam Q)^{p-n}}\le
C\,\left\{\frac{V(Q)^p}
{(\diam Q)^{p-n}}+
\sbig\limits_{H\in T(K_Q)}\,\,
\sbig\limits_{|\xi|\le m-1}
\,\frac{L(\xi:H,Q)^p}
{(\diam K_Q)^{(m-|\xi|)p-n}}\right\}.
$$
%----------------------------------------------------------
%@@@@@@@@@@@@@@@@@@@@@@@@@@@@@@@@@@@@@@@@@@@@@@@@@@@@@@@@@@
%@@@@@@@@@@@@@@@@@@@@@@@@@@@@@@@@@@@@@@@@@@@@@@@@@@@@@@@@@@
%@@@@@@@@@@@@@@@@@@@@@@@@@@@@@@@@@@@@@@@@@@@@@@@@@@@@@@@@@@
%@@@@@@@@@@@@@@@@@@@@@@@@@@@@@@@@@@@@@@@@@@@@@@@@@@@@@@@@@@
%----------------------------------------------------------
\par Prove that $a_H\in\gamma K_Q$ whenever $H\in T(K_Q)$. Since $H\cap K_Q\ne\emp$, by Lemma \reff{Wadd}, we have $\diam H\le 4\diam K_Q$ so that $H\subset 9 K_Q$. On the other hand, by Theorem \reff{L-PE}, $a_H\in\tgm H$. Hence
%----------------------------------------------------------
\bel{AH-KQ}
a_H\in 9\tgm\, K_Q~~~\text{for every}~~~H\in T(K_Q)
\ee
%---------------------------------------------------------
proving that
%----------------------------------------------------------
\bel{A-TK}
a_H\in\gamma K_Q.
\ee
%----------------------------------------------------------
(Recall that $\gamma:=10^4\tgm$, see the beginning of this subsection.)
%----------------------------------------------------------
%@@@@@@@@@@@@@@@@@@@@@@@@@@@@@@@@@@@@@@@@@@@@@@@@@@@@@@@@@@
%@@@@@@@@@@@@@@@@@@@@@@@@@@@@@@@@@@@@@@@@@@@@@@@@@@@@@@@@@@
%@@@@@@@@@@@@@@@@@@@@@@@@@@@@@@@@@@@@@@@@@@@@@@@@@@@@@@@@@@
%@@@@@@@@@@@@@@@@@@@@@@@@@@@@@@@@@@@@@@@@@@@@@@@@@@@@@@@@@@
%----------------------------------------------------------
\par By $H_Q$ we denote a cube $H\in T(K_Q)$ for which the quantity
%----------------------------------------------------------
\bel{HQ-DF}
\sbig\limits_{|\xi|\le m-1}
\,\frac{L(\xi:H,Q)^p}
{(\diam K_Q)^{(m-|\xi|)p-n}}~~~\text{is maximal on}~~T(K_Q).
\ee
%----------------------------------------------------------
Since $\#T(K_Q)\le N(n),$ we obtain the following inequality
%----------------------------------------------------------
$$
\frac{S(Q)^p}
{(\diam Q)^{p-n}}\le
C\,\left\{\frac{V(Q)^p}
{(\diam Q)^{p-n}}+
\sbig\limits_{|\xi|\le m-1}
\,\frac{L(\xi:H_Q,Q)^p}
{(\diam K_Q)^{(m-|\xi|)p-n}}\right\}.
$$
%----------------------------------------------------------
%@@@@@@@@@@@@@@@@@@@@@@@@@@@@@@@@@@@@@@@@@@@@@@@@@@@@@@@@@@
%@@@@@@@@@@@@@@@@@@@@@@@@@@@@@@@@@@@@@@@@@@@@@@@@@@@@@@@@@@
%----------------------------------------------------------
Hence
%----------------------------------------------------------
\be
I_2(F)&:=&\sbig_{Q\in\Qc}\frac{S(Q)^p}
{(\diam Q)^{p-n}}\nn\\
&\le&
C\,\left\{\sbig_{Q\in\Qc}\frac{V(Q)^p}
{(\diam Q)^{p-n}}+
\sbig\limits_{|\xi|\le m-1}\sbig_{Q\in\Qc}
\frac{L(\xi:H_Q,Q)^p}
{(\diam K_Q)^{(m-|\xi|)p-n}}\right\}\nn\\
&=&
C\,\{I_2^{(1)}(F)+I_2^{(2)}(F)\}.
\nn
\ee
%----------------------------------------------------------
\par Prove that
%----------------------------------------------------------
\bel{I21-DF1}
I_2^{(1)}(F):=\sbig_{Q\in\Qc}\frac{V(Q)^p}
{(\diam Q)^{p-n}}\le \lambda.
\ee
%----------------------------------------------------------
In fact, we know that $a_{K_Q}\in \tQ:=5\tgm\, Q$ and $c_Q=c_{\tQ}$, and the squares of the family  $\tQc=\{\tQ\}$ are pairwise disjoint. This enables us to use the same approach as in the proof of the inequality $I_1(F)\le\lambda$. This implies the required estimate $I_2^{(1)}(F)\le \lambda$.
%----------------------------------------------------------
\par Prove that $I_2^{(2)}(F)\le C\,\lambda$.
%----------------------------------------------------------
\par Let $\Kc:=\{K_Q: Q\in\Qc\}$. We know that $K_Q\subset 5Q\subset \tQ=5\tgm\,Q$ so that the cubes of the family $\Kc$ are pairwise disjoint. We also recall that, by \rf{A-TK}, $a_{H_Q}, a_{K_Q}\in\gamma K_Q$ for every $Q\in\Qc$. Furthermore, we know that $H_Q\cap K_Q\ne\emp$ and that $a_{K_Q}=\PRL(L)$, $a_{H_Q}=\PRL(L')$ where $L$ and $L'$ are lacunae containing $K_Q$ and $H_Q$ respectively. Therefore, by part (ii) of Theorem \reff{L-PE},
%----------------------------------------------------------
$$
\diam K_Q\le\tgm\,\|a_{K_Q}-a_{H_Q}\|\le
\gamma\,\|a_{K_Q}-a_{H_Q}\|
$$
%----------------------------------------------------------
provided $a_{K_Q}\ne a_{H_Q}$. Thus the family$\{\{a_{K_Q},a_{H_Q}\}: Q\in\Qc\}$ of two point subsets of $E$ is $\gamma$-sparse. See Definition \reff{DF-PR}. (We also note that $\diam K_Q\sim\|a_{K_Q}-a_{H_Q}\|$.)
%----------------------------------------------------------
\par Recall that
%----------------------------------------------------------
\bel{I22-DF2}
I_2^{(2)}(F):=\sbig\limits_{|\xi|\le m-1}\,\,
\sbig_{Q\in\Qc}\,
\frac{|D^{\xi}P_{a_{H_Q}}(a_{K_Q})
-D^{\xi}P_{a_{K_Q}}(a_{K_Q})|^p}
{(\diam K_Q)^{(m-|\xi|)p-n}}.
\ee
%----------------------------------------------------------
Then, by assumption \rf{N-P-NEW-1},
%----------------------------------------------------------
$$
I_2^{(2)}(F)\sim
\sbig\limits_{|\xi|\le m-1}\,\,
\sbig_{Q\in\Qc}\, \frac{|D^{\xi}P_{a_{H_Q}}(a_{K_Q})
-D^{\xi}P_{a_{K_Q}}(a_{K_Q})|^p}
{\|a_{K_Q}-a_{H_Q}\|^{(m-|\xi|)p-n}}\le \lambda
$$
%----------------------------------------------------------
%@@@@@@@@@@@@@@@@@@@@@@@@@@@@@@@@@@@@@@@@@@@@@@@@@@@@@@@@@@
%@@@@@@@@@@@@@@@@@@@@@@@@@@@@@@@@@@@@@@@@@@@@@@@@@@@@@@@@@@
%----------------------------------------------------------
proving the required inequality $I_2^{(2)}(F)\le C\,\lambda$. \medskip
%----------------------------------------------------------
%@@@@@@@@@@@@@@@@@@@@@@@@@@@@@@@@@@@@@@@@@@@@@@@@@@@@@@@@@@
%@@@@@@@@@@@@@@@@@@@@@@@@@@@@@@@@@@@@@@@@@@@@@@@@@@@@@@@@@@
%----------------------------------------------------------
\par We turn to the proof of the lemma in the general case, i.e., for an arbitrary family of equal cubes $\{Q_i:i=1,...,k\}$ satisfying inequality \rf{D-L40}. Let $u_i$ be a nearest point on $E$ to the square $Q_i$, $i=1,...,k$. By \rf{D-L40},
%----------------------------------------------------------
$$
\|u_i-c_i\|\le \dist(c_i,E)\le 40\,\diam Q_i.
$$
%----------------------------------------------------------
so that
%----------------------------------------------------------
\bel{Q-IT}
\|u_i-y\|\le 41\,\diam Q_i~~~\text{for every}~~~y\in Q_i.
\ee
%----------------------------------------------------------
\par Let
%----------------------------------------------------------
\bel{Q-H1}
\OV_i:=Q(u_i,R_i)~~~\text{where}~~~R_i=41\,\diam Q_i.
\ee
%----------------------------------------------------------
Then, by \rf{Q-IT}, $\OV_i\supset Q_i$. Also it can be readily seen that $\OV_i\subset \gamma' Q_i$ with $\gamma'=122$. Since the cubes $\{Q_i\}$ are pairwise disjoint, the family $\hat{\Qc}=\{\OV_i:i=1,...,k\}$ has covering multiplicity $ \MP(\hat{\Qc})<C(n)$. Therefore, by Theorem \reff{TFM}, $\hat{\Qc}$ can be partitioned into at most $N(n)$ subfamilies each consisting of pairwise disjoint cubes. This enables us to assume that the family $\hat{\Qc}$ itself consists of pairwise disjoint cubes.
%----------------------------------------------------------
\par We obtain
%----------------------------------------------------------
\be
J&:=&
\sbig_{i=1}^k\,
\frac{|D^\beta F(x_i)-D^\beta F(c_i)|^p}
{(\diam Q_i)^{p-n}}\nn\\
&\le& C\,\left\{\,\sbig_{i=1}^k\,
\frac{|D^\beta F(x_i)-D^\beta F(u_i)|^p}
{(\diam \OV_i)^{p-n}}+
\sbig_{i=1}^k\,
\frac{|D^\beta F(c_i)-D^\beta F(u_i)|^p}
{(\diam \OV_i)^{p-n}}\right\}\nn\\
&=& C\{J_1+J_2\}.\nn
%----------------------------------------------------------
\ee
%----------------------------------------------------------
\par But now $u_i=c_{\OV_i}\in E$, so that the problem is reduced to the case proven above. Thus $J_1,J_2\le C\lambda$ where $C=C(m,n,p)$. Hence $J\le  C\,\lambda$ proving the lemma.\bx
%----------------------------------------------------------
%@@@@@@@@@@@@@@@@@@@@@@@@@@@@@@@@@@@@@@@@@@@@@@@@@@@@@@@@@@
%@@@@@@@@@@@@@@@@@@@@@@@@@@@@@@@@@@@@@@@@@@@@@@@@@@@@@@@@@@
%@@@@@@@@@@@@@@@@@@@@@@@@@@@@@@@@@@@@@@@@@@@@@@@@@@@@@@@@@@
%----------------------------------------------------------
\smallskip
%---------------------------------------------------------
%@@@@@@@@@@@@@@@@@@@@@@@@@@@@@@@@@@@@@@@@@@@@@@@@@@@@@@@@@@
%@@@@@@@@@@@@@@@@@@@@@@@@@@@@@@@@@@@@@@@@@@@@@@@@@@@@@@@@@@
\begin{lemma}\lbl{L-SMC} Suppose that
%----------------------------------------------------------
\bel{D-S80}
\dist(c_i,E)>40\,\diam Q_i,~~~i=1,...,k.
\ee
%----------------------------------------------------------
Then for every $\beta, |\beta|=m-1$, and every $x_i\in Q_i$ the following inequality
%----------------------------------------------------------
$$
\sbig_{i=1}^k\,
\frac{|D^\beta F(x_i)-D^\beta F(c_i)|^p}
{(\diam Q_i)^{p-n}}
\le C\,\lambda
$$
%----------------------------------------------------------
holds. Here $C=C(m,n,p)$.
%----------------------------------------------------------
\end{lemma}
%----------------------------------------------------------
%@@@@@@@@@@@@@@@@@@@@@@@@@@@@@@@@@@@@@@@@@@@@@@@@@@@@@@@@@@
%@@@@@@@@@@@@@@@@@@@@@@@@@@@@@@@@@@@@@@@@@@@@@@@@@@@@@@@@@@
%@@@@@@@@@@@@@@@@@@@@@@@@@@@@@@@@@@@@@@@@@@@@@@@@@@@@@@@@@@
%----------------------------------------------------------
\par {\it Proof.} For every cube $Q\in\Qc$, by \rf{D-S80}, $\dist(c_{Q},E)> 40\diam Q>0$ so that $Q\subset \RN\setminus E$. Let $K_Q$ be a Whitney cube which contains $c_Q$. Prove that $Q\subset K^*_Q=\tfrac98 K_Q.$
%----------------------------------------------------------
\par In fact,
%----------------------------------------------------------
$$
\diam Q<\tfrac1{40}\dist(c_{Q},E)\le \tfrac1{40}\{\diam K_Q+\dist(K_{Q},E)\}
$$
%----------------------------------------------------------
so that
%----------------------------------------------------------
\bel{18-KQ}
\diam Q\le \tfrac1{40}
\{\diam K_Q+4\diam K_Q\}=\tfrac1{8}\diam K_Q\,.
\ee
%----------------------------------------------------------
Hence for every $z\in Q$ we have
%----------------------------------------------------------
\be
\|z-c_{K_Q}\|&\le& \|z-c_Q\|+\|c_Q-c_{K_Q}\|\le \tfrac12\diam Q+\tfrac12\diam K_Q\nn\\
&\le& \tfrac12\cdot\tfrac1{8}\diam K_Q+\tfrac12\diam K_Q=\left(\tfrac1{8}+1\right) r_{K_Q}\nn
\ee
%----------------------------------------------------------
so that $Q\subset\tfrac98 K_Q=K^*_Q$. In particular, $x_Q=x_i\in K^*_Q$ provided $Q=Q_i$.
%----------------------------------------------------------
%@@@@@@@@@@@@@@@@@@@@@@@@@@@@@@@@@@@@@@@@@@@@@@@@@@@@@@@@@@
%@@@@@@@@@@@@@@@@@@@@@@@@@@@@@@@@@@@@@@@@@@@@@@@@@@@@@@@@@@
%@@@@@@@@@@@@@@@@@@@@@@@@@@@@@@@@@@@@@@@@@@@@@@@@@@@@@@@@@@
%----------------------------------------------------------
\par Let $K\in W_E$ and let
%----------------------------------------------------------
$$
\Qc(K):=\{Q\in\Qc: c_Q\in K\}.
$$
%----------------------------------------------------------
Let
%----------------------------------------------------------
$$
\Kc:=\{K\in W_E: \Qc(K)\ne\emp\}.
$$
%----------------------------------------------------------
%@@@@@@@@@@@@@@@@@@@@@@@@@@@@@@@@@@@@@@@@@@@@@@@@@@@@@@@@@@
%----------------------------------------------------------
\par Then for each $K\in \Kc$, by Lemma \reff{C-4},
%----------------------------------------------------------
%@@@@@@@@@@@@@@@@@@@@@@@@@@@@@@@@@@@@@@@@@@@@@@@@@@@@@@@@@@
%----------------------------------------------------------
$$
|D^{\beta}F(x_Q)-D^{\beta}F(c_Q)|\le C
\|x_Q-c_Q\|\,
\sbig\limits_{H\in\, T(K)}\,\,\sbig\limits_{|\xi|\le m-1}
 \,\frac{|D^{\xi}P_{a_H}(a_K)-D^{\xi}P_{a_K}(a_K)|}{(\diam K)^{m-|\xi|}}.
$$
%----------------------------------------------------------
\par By \rf{A-TK},
%----------------------------------------------------------
\bel{IN-A}
a_K,a_H\in \gamma K~~~\text{for every}~~~H\in T(K).
\ee
%----------------------------------------------------------
\par Now we have
%----------------------------------------------------------
\be
I_K&:=&
\sbig_{Q\in\Qc(K)}\,
\frac{|D^\beta F(x_Q)-D^\beta F(c_Q)|^p}
{(\diam Q)^{p-n}}\nn\\
&\le&
C\left\{\sbig_{Q\in\Qc(K)}
\left(\frac{\|x_Q-c_Q\|}{\diam Q}\right)^p|Q|\right\}
\left\{\sbig\limits_{H\in T(K)}\,\,
\sbig\limits_{|\xi|\le m-1}
\frac{|D^{\xi}P_{a_H}(a_K)-D^{\xi}P_{a_K}(a_K)|}{(\diam K)^{m-|\xi|}}\right\}^p\nn\\
&\le&
C\,\left\{\sbig_{Q\in\Qc(K)}\,
|Q|\right\}
\left\{\sbig\limits_{H\in\, T(K)}\,\,
\sbig\limits_{|\xi|\le m-1}
\,\frac{|D^{\xi}P_{a_H}(a_K)-D^{\xi}P_{a_K}(a_K)|}{(\diam K)^{m-|\xi|}}\right\}^p\nn\\
&\le&
C\,|K|\left\{\sbig\limits_{H\in\, T(K)}\,\,
\sbig\limits_{|\xi|\le m-1}
\,\frac{|D^{\xi}P_{a_H}(a_K)-D^{\xi}P_{a_K}(a_K)|}{(\diam K)^{m-|\xi|}}\right\}^p.\nn
\ee
%----------------------------------------------------------
%@@@@@@@@@@@@@@@@@@@@@@@@@@@@@@@@@@@@@@@@@@@@@@@@@@@@@@@@@@
%@@@@@@@@@@@@@@@@@@@@@@@@@@@@@@@@@@@@@@@@@@@@@@@@@@@@@@@@@@
%----------------------------------------------------------
\par Since $\#T(K)\le N(n)$, see Lemma \reff{Wadd}, we obtain
%----------------------------------------------------------
$$
I_K\le
C\,\sbig\limits_{H\in\, T(K)}\,\,
\sbig\limits_{|\xi|\le m-1}
\,\frac{|D^{\xi}P_{a_H}(a_K)-D^{\xi}P_{a_K}(a_K)|^p}
{(\diam K)^{(m-|\xi|)p-n}}.
$$
%----------------------------------------------------------
%@@@@@@@@@@@@@@@@@@@@@@@@@@@@@@@@@@@@@@@@@@@@@@@@@@@@@@@@@@
%@@@@@@@@@@@@@@@@@@@@@@@@@@@@@@@@@@@@@@@@@@@@@@@@@@@@@@@@@@
%----------------------------------------------------------
Let $\KM\in T(K)$ be a cube such that the quantity
%----------------------------------------------------------
\bel{KM-DF}
\sbig\limits_{|\xi|\le m-1}
\,\frac{|D^{\xi}P_{a_H}(a_K)-D^{\xi}P_{a_K}(a_K)|^p}
{(\diam K)^{(m-|\xi|)p-n}}~~\text{takes the maximal value on}~~T(K).
\ee
%----------------------------------------------------------
Then
%----------------------------------------------------------
$$
I_K\le
C\,\sbig\limits_{|\xi|\le m-1}
\,\frac{|D^{\xi}P_{a_{\KM}}(a_K)-D^{\xi}P_{a_K}(a_K)|^p}
{(\diam K)^{(m-|\xi|)p-n}}.
$$
%----------------------------------------------------------
\par This enables us to estimate $I$ as follows:
%----------------------------------------------------------
$$
I:=\sbig_{i=1}^k\,
\frac{|D^\beta F(x_i)-D^\beta F(c_i)|^p}
{(\diam Q_i)^{p-n}}\le\smed_{K\in\, \Kc}\, I_K
$$
%----------------------------------------------------------
so that
%----------------------------------------------------------
\bel{L-FML}
I\le
C\,\sbig_{K\in \Kc}\,\,\,\sbig\limits_{|\xi|\le m-1}\,
\,\frac{|D^{\xi}P_{a_{\KM}}(a_K)-D^{\xi}P_{a_K}(a_K)|^p}
{(\diam K)^{(m-|\xi|)p-n}}.
\ee
%----------------------------------------------------------
\par As in the proof of the previous lemma, using Theorem \reff{TFM} we can assume that the cubes of the family $\Kc$ are pairwise disjoint. Furthermore, by \rf{LM-WM-DEF} and by part(ii) of Theorem \reff{L-PE},
%----------------------------------------------------------
$$
\diam K\le \tgm\,\|a_{\KM}-a_K\|\le \gamma\,\|a_{\KM}-a_K\|
$$
%----------------------------------------------------------
provided $a_{\KM}\ne a_K$. (Recall that $\gamma=10^4\tgm$.) In particular, this inequality and \rf{IN-A} imply the following:
%----------------------------------------------------------
\bel{DK-A2}
\frac{1}{\gamma}\,\diam K\le \|a_{\KM}-a_K\|\le \gamma \diam K.
\ee
%----------------------------------------------------------
Hence
%----------------------------------------------------------
\bel{IN-K1}
I\le
C\,\sbig_{K\in \Kc}\,\,\sbig\limits_{|\xi|\le m-1}\,
\,\frac{|D^{\xi}P_{a_{\KM}}(a_K)-D^{\xi}P_{a_K}(a_K)|^p}
{\|a_{\KM}-a_K\|^{(m-|\xi|)p-n}}.
\ee
%----------------------------------------------------------
\par Note that, by \rf{IN-A}, $a_{\KM},a_K\in\gamma\,K$. Combining this with \rf{DK-A2} and Definition \reff{DF-PR} we conclude that the family $\{\{a_{\KM},a_K\}: K\in\Kc\}$
of two point subsets of $E$ is $\gamma$-sparse. (Recall that $\Kc$ consists of pairwise disjoint cubes.) This enables us to apply assumption \rf{N-P-NEW-1} to the right hand side of the above inequality. By this assumption, $I\le C\,\lambda$, and the proof of the lemma is complete.\bx\medskip
%----------------------------------------------------------
%@@@@@@@@@@@@@@@@@@@@@@@@@@@@@@@@@@@@@@@@@@@@@@@@@@@@@@@@@@
%@@@@@@@@@@@@@@@@@@@@@@@@@@@@@@@@@@@@@@@@@@@@@@@@@@@@@@@@@@
%@@@@@@@@@@@@@@@@@@@@@@@@@@@@@@@@@@@@@@@@@@@@@@@@@@@@@@@@@@
%----------------------------------------------------------
\par Combining Lemma \reff{L-BG} and Lemma \reff{L-SMC} with the criterion \rf{CR-SOB}, we conclude that $F\in C^{m-1}(\RN)$, and for every multiindex $\beta$ of order $m-1$ the function $D^\beta F\in\LOP$ and $\|D^\beta F\|_{\LOP}\le C\,\lambda^{\frac1p}$. Since weak derivatives commute, the function $F\in\LMP$ and
$\| F\|_{\LMP}\le C\,\lambda^{\frac1p}$.\smallskip
%----------------------------------------------------------
%@@@@@@@@@@@@@@@@@@@@@@@@@@@@@@@@@@@@@@@@@@@@@@@@@@@@@@@@@@
%@@@@@@@@@@@@@@@@@@@@@@@@@@@@@@@@@@@@@@@@@@@@@@@@@@@@@@@@@@
%@@@@@@@@@@@@@@@@@@@@@@@@@@@@@@@@@@@@@@@@@@@@@@@@@@@@@@@@@@
%----------------------------------------------------------
\par The proof of Theorem \reff{EX-TK} is complete.\bx
%----------------------------------------------------------
\bigskip\medskip
%@@@@@@@@@@@@@@@@@@@@@@@@@@@@@@@@@@@@@@@@@@@@@@@@@@@@@@@@@@
%@@@@@@@@@@@@@@@@@@@@@@@@@@@@@@@@@@@@@@@@@@@@@@@@@@@@@@@@@@
%@@@@@@@@@@@@@@@@@@@@@@@@@@@@@@@@@@@@@@@@@@@@@@@@@@@@@@@@@@
%@@@@@@@@@@@@@@@@@@@@@@@@@@@@@@@@@@@@@@@@@@@@@@@@@@@@@@@@@@
%@@@@@@@@@@@@@@@@@@@@@@@@@@@@@@@@@@@@@@@@@@@@@@@@@@@@@@@@@@
%@@@@@@@@@@@@@@@@@@@@@@@@@@@@@@@@@@@@@@@@@@@@@@@@@@@@@@@@@@
%@@@@@@@@@@@@@@@@@@@@@@@@@@@@@@@@@@@@@@@@@@@@@@@@@@@@@@@@@@
%----------------------------------------------------------
\par {\bf 7.3. A refinement of the variational criterion for finite sets.}\medskip
%----------------------------------------------------------
\addtocontents{toc}{~~~~7.3. A refinement of the variational criterion for finite sets. \hfill \thepage\par}
%----------------------------------------------------------
\par In this subsection we prove Theorem \reff{EX-REFTK}.
Our proof is based on the following useful property of the lacunary extension operator.
%---------------------------------------------------------
%@@@@@@@@@@@@@@@@@@@@@@@@@@@@@@@@@@@@@@@@@@@@@@@@@@@@@@@@@@
%@@@@@@@@@@@@@@@@@@@@@@@@@@@@@@@@@@@@@@@@@@@@@@@@@@@@@@@@@@
\begin{proposition}\lbl{EXT-OUT} Let $E$ be a closed set and let $\VP=\{P_x:x\in E\}$ be a Whitney $(m-1)$-field on $E$. Let $F$ be the function obtained by the lacunary modification \rf{LM-WM-DEF} of the Whitney extension formula \rf{DEF-F}. Then the following inequality
%----------------------------------------------------------
$$
\smed_{|\alpha|=m}\,\, \|D^\alpha F\|_{L_p(\RN\setminus E)} \le C\,
\left\{\sbig_{x,y\in E,\, x\lr y}\,\,\,\,\sbig_{|\beta|\le m-1}\frac{|D^\beta P_{x}(x)-D^\beta P_{y}(x)|^p}
{\|x-y\|^{(m-|\beta|)p-n}}\right\}^{\frac1p}
%---------------------------------------------------------
$$
%----------------------------------------------------------
holds. Here the first sum is taken over all distinct points $x$ and $y$ in $E$ joined by an edge in the graph $\GE$. (In this section we use the notation $x\lr y$.)
%----------------------------------------------------------
\par The constant $C$ in this inequality depends only on $m,n,$ and $p$.
%----------------------------------------------------------
\end{proposition}
%----------------------------------------------------------
%@@@@@@@@@@@@@@@@@@@@@@@@@@@@@@@@@@@@@@@@@@@@@@@@@@@@@@@@@@
%@@@@@@@@@@@@@@@@@@@@@@@@@@@@@@@@@@@@@@@@@@@@@@@@@@@@@@@@@@
%@@@@@@@@@@@@@@@@@@@@@@@@@@@@@@@@@@@@@@@@@@@@@@@@@@@@@@@@@@
%----------------------------------------------------------
\par {\it Proof.} Let a cube $K\in W_E$ and let $u\in K$. Then, by Lemma \reff{MPO-C}, for every multiindex $\alpha$, $|\alpha|=m$,
%----------------------------------------------------------
$$
|D^{\alpha}F(u)|\le C\,
\smed\limits_{Q\in T(K),\, a_Q\ne a_K}
\,\,\,\smed\limits_{|\xi|\le m-1}\,\,
(\diam K)^{|\xi|-m} \,|D^{\xi}P_{a_Q}(a_K)-D^{\xi}P_{a_K}(a_K)|
$$
%----------------------------------------------------------
where $C=C(n,m)$. See \rf{5-5N}.
%----------------------------------------------------------
\par We recall that $F|_{\RN\setminus E}\in C^{\infty}(\RN\setminus E)$, $T(K)$ is the family of Whitney's cubes touching $K$, see \rf{TK}, and $a_Q$ is defined by \rf{LM-WM-DEF}. We also recall that $\#\,T(K)$ is bounded by $C(n)$, see Lemma \reff{Wadd}.
%----------------------------------------------------------
\par Integrating the latter inequality on $K$ (with respect to $u$) we obtain
%----------------------------------------------------------
$$
\intl_K |D^{\alpha}F(u)|^p\,du\le C
\smed\limits_{Q\in T(K),\, a_Q\ne a_K}
\,\,\,\smed\limits_{|\xi|\le m-1}\,\,
(\diam K)^{(|\xi|-m)p+n} \,|D^{\xi}P_{a_Q}(a_K)-D^{\xi}P_{a_K}(a_K)|^p
$$
%----------------------------------------------------------
where $C=C(n,m,p)$.
%----------------------------------------------------------
\par Let us fix a cube $Q\in T(K)$; thus $Q\cap K\ne\emp$. Let $L^{(Q)}$ and $L^{(K)}$ be lacunae containing $Q$ and $K$ respectively. We know that
%----------------------------------------------------------
$$
a_Q=\PRL(L^{(Q)})~~~\text{and}~~~a_K=\PRL(L^{(K)}),
$$
%----------------------------------------------------------
see \rf{LM-WM-DEF}, so that, by part (ii) of Theorem \reff{L-PE},
%----------------------------------------------------------
$$
\diam Q+\diam K\le \tgm\,\|a_Q-a_K\|
$$
%----------------------------------------------------------
provided $a_Q\ne a_K$. Recall that $\tgm$ is an absolute constant. Hence
%----------------------------------------------------------
$$
\intl_K |D^{\alpha}F(u)|^p\,du\le C\,
\smed\limits_{Q\in T(K),\, a_Q\ne a_K}
\,\,\smed\limits_{|\xi|\le m-1}\,\,
\frac{|D^\xi P_{a_Q}(a_K)-D^\xi P_{a_K}(a_K)|^p}
{\|a_Q-a_K\|^{(m-|\beta|)p-n}}
$$
%----------------------------------------------------------
so that, by Definition \reff{GRAPH-GE},
%----------------------------------------------------------
$$
\intl_K |D^{\alpha}F(u)|^p\,du\le C\,
\smed\limits_{Q\in W_E,\, a_Q\lr a_K}
\,\,\smed\limits_{|\xi|\le m-1}\,\,
\frac{|D^\xi P_{a_Q}(a_K)-D^\xi P_{a_K}(a_K)|^p}
{\|a_Q-a_K\|^{(m-|\beta|)p-n}}\,.
$$
%----------------------------------------------------------
\par This estimate implies the following inequality
%----------------------------------------------------------
$$
I:=\intl_{\RN\setminus E} |D^{\alpha}F(u)|^p\,du\le C\,\smed\limits_{K\in W_E}\,\,
\smed\limits_{Q\in W_E,\, a_Q\lr a_K}
\,\,\smed\limits_{|\xi|\le m-1}\,\,
\frac{|D^\xi P_{a_Q}(a_K)-D^\xi P_{a_K}(a_K)|^p}
{\|a_Q-a_K\|^{(m-|\beta|)p-n}}
$$
%----------------------------------------------------------
so that
%----------------------------------------------------------
$$
I\le C\,\smed\limits_{Q,K\in W_E,\, a_Q\lr a_K}
\,\,\smed\limits_{|\xi|\le m-1}\,\,
\frac{|D^\xi P_{a_Q}(a_K)-D^\xi P_{a_K}(a_K)|^p}
{\|a_Q-a_K\|^{(m-|\beta|)p-n}}\,.
$$
%----------------------------------------------------------
In turn, this inequality implies the following one:
%----------------------------------------------------------
\bel{F-DU}
I\le  C\,
\sbig_{x,y\in E,\, x\lr y}\,\#\,R(x,y)\sbig_{|\xi|\le m-1}\frac{|D^\xi P_{x}(x)-D^\xi P_{y}(x)|^p}
{\|x-y\|^{(m-|\xi|)p-n}}\,.
\ee
%----------------------------------------------------------
Here given $x,y\in E$, $x\lr y$,  by $R(x,y)$ we denote a subset of $W_E\times W_E$ consisting of all pairs $(Q,K)$ of contacting Whitney's cubes $Q$ and $K$, $Q\cap K\ne\emp$, such that there exist lacunae $L^{(Q)},L^{(K)}\in\LE$, $L^{(Q)}\ni Q$ and $L^{(K)}\ni K$, for which
%----------------------------------------------------------
$$
x=\PRL(L^{(Q)})~~~\text{and}~~~y=\PRL(L^{(K)})\,.
$$
%----------------------------------------------------------
\par But, by Theorem \reff{L-PE},
%----------------------------------------------------------
$$
\#\,\{L:\PRL(L)=x\}\le C_1(n)~~\text{and}~~\#\,\{L:\PRL(L)=y\}\le C_1(n).
$$
%----------------------------------------------------------
Furthermore, by Proposition \reff{D-LQN}, each lacuna contains at most $M(n)$ contacting cubes. Finally, each Whitney's cube has common points with at most $N(n)$ Whitney's cubes, see Lemma \reff{Wadd}. Hence,
%----------------------------------------------------------
$$
\#\,R(x,y)\le C_1(n)\,M(n)^2
$$
%----------------------------------------------------------
so that, by \rf{F-DU},
%----------------------------------------------------------
$$
I\le  C_2\,
\sbig_{x,y\in E,\, x\lr y}\,\,\,\sbig_{|\xi|\le m-1}\frac{|D^\xi P_{x}(x)-D^\xi P_{y}(x)|^p}
{\|x-y\|^{(m-|\xi|)p-n}}
$$
%----------------------------------------------------------
where $C_2:=C\,C_1(n)\,M(n)^2$.
%----------------------------------------------------------
\par The proof of the proposition is complete.\bx
%----------------------------------------------------------
\medskip
%----------------------------------------------------------
\par {\bf Proof of Theorem \reff{EX-REFTK}.} Let $E$ be a finite set and let $\VP=\{P_x:x\in E\}$ be a Whitney $(m-1)$-field on $E$. Let
%----------------------------------------------------------
$$
I:=\sbig_{x,y\in E,\, x\lr y}\,\,\,\,\sbig_{|\beta|\le m-1}\frac{|D^\beta P_{x}(x)-D^\beta P_{y}(x)|^p}
{\|x-y\|^{(m-|\beta|)p-n}}\,.
$$
%----------------------------------------------------------
\par Since $\GE$ is a $\gamma$-sparse graph with $\gamma=\gamma(n)$, see part (i) of Proposition \reff{PROP-GE}, by Definition \reff{DF-GRPR} and Proposition \reff{N-VC}, for every $C^{m-1}$-function $F\in \LMP$ such that $T^{m-1}_x[F]=P_x$ for all $x\in E$ the following inequality
%----------------------------------------------------------
$$
I\le  C(m,n,p)\,\|F\|^p_{\LMP}
$$
%----------------------------------------------------------
holds. This proves that $I\le C\,\PME$, see \rf{N-VP}.\smallskip
%----------------------------------------------------------
%@@@@@@@@@@@@@@@@@@@@@@@@@@@@@@@@@@@@@@@@@@@@@@@@@@@@@@@@@@
%@@@@@@@@@@@@@@@@@@@@@@@@@@@@@@@@@@@@@@@@@@@@@@@@@@@@@@@@@@
%@@@@@@@@@@@@@@@@@@@@@@@@@@@@@@@@@@@@@@@@@@@@@@@@@@@@@@@@@@
%----------------------------------------------------------
\par Prove that $\PME\le C\,I$. In fact, since $E$ is a finite set, the function $F$ obtained by the Whitney extension formula \rf{DEF-F} belongs to $C^{\infty}(\RN)$. Hence
%----------------------------------------------------------
$$
\|F\|_{\LMP}=\|F\|_{L^m_p(\RN\setminus E)}
=\smed_{|\alpha|=m}\,\,
\|D^\alpha F\|_{L_p(\RN\setminus E)}\,.
$$
%----------------------------------------------------------
\par Let $F$ be the function obtained by the lacunary modification of the Whitney extension formula. See \rf{LM-WM-DEF}. Then, by Proposition \reff{EXT-OUT}, $\|F\|^p_{L^m_p(\RN\setminus E)}\le C\,I$ so that $\|F\|^p_{\LMP}\le C\,I$. Hence, by \rf{N-VP}, $\PME\le C\,I$, and the proof of the theorem is complete.\bx
%----------------------------------------------------------
%@@@@@@@@@@@@@@@@@@@@@@@@@@@@@@@@@@@@@@@@@@@@@@@@@@@@@@@@@@
%@@@@@@@@@@@@@@@@@@@@@@@@@@@@@@@@@@@@@@@@@@@@@@@@@@@@@@@@@@
%@@@@@@@@@@@@@@@@@@@@@@@@@@@@@@@@@@@@@@@@@@@@@@@@@@@@@@@@@@
%@@@@@@@@@@@@@@@@@@@@@@@@@@@@@@@@@@@@@@@@@@@@@@@@@@@@@@@@@@
%@@@@@@@@@@@@@@@@@@@@@@@@@@@@@@@@@@@@@@@@@@@@@@@@@@@@@@@@@@
%----------------------------------------------------------
\medskip
\par  Let us note the following useful property of the graph $\GE$.
%----------------------------------------------------------
%@@@@@@@@@@@@@@@@@@@@@@@@@@@@@@@@@@@@@@@@@@@@@@@@@@@@@@@@@@
%@@@@@@@@@@@@@@@@@@@@@@@@@@@@@@@@@@@@@@@@@@@@@@@@@@@@@@@@@@
%@@@@@@@@@@@@@@@@@@@@@@@@@@@@@@@@@@@@@@@@@@@@@@@@@@@@@@@@@@
%----------------------------------------------------------
%@@@@@@@@@@@@@@@@@@@@@@@@@@@@@@@@@@@@@@@@@@@@@@@@@@@@@@@@@@
\begin{proposition} Let $E$ be a finite subset of $\RN$. Then $\GE$ is a connected graph.
%----------------------------------------------------------
\par Furthermore, for every $\hx,\hy\in E$ there is a path $\{z_0,...,z_m\}$ joining $\hx$ to $\hy$ in $\GE$ (i.e., $z_0=\hx$, $z_m=\hy$, and $z_i\lr z_{i+1}$, $i=0,...,m-1$) such that
%---------------------------------------------------------
$$
\smed_{i=0}^{m-1}\,\|z_i-z_{i+1}\|\le C(n)\,\|\hx-\hy\|.
$$
%----------------------------------------------------------
\end{proposition}
%----------------------------------------------------------
\par {\it Proof.} Let $f$ be a function on $E$. We know that $L^1_\infty(\RN)$ is isometrically isomorphic to $\Lip(\RN)$. On the other hand, by McShane's theorem \cite{McS}, $\|f\|_{\Lip(E)}=\|f\|_{\Lip(\RN)|_E}$ so that $\|f\|_{\Lip(E)}=\|f\|_{L^1_\infty(\RN)|_E}$. Hence, by \rf{P-INF1},
%----------------------------------------------------------
\bel{LIPE}
\|f\|_{\Lip(E)}\sim \sup_{x,y\in E,\, x\lr y}\, \frac{|f(x)-f(y)|}{\|x-y\|}
\ee
%----------------------------------------------------------
with the constants in this equivalence depending only on $n$.
%----------------------------------------------------------
\par Using this equivalence prove that $\GE$ is a connected graph. Otherwise there exists a non-empty subset $E'\subset E$, $E'\ne E$, such that $x\nleftrightarrow y$ for every $x\in E'$ and $y\in E\setminus E'$. Let $f:=\chi_{E'}$. Then $\|f\|_{\Lip(E)}>0$, while the right hand side of \rf{LIPE} equals $0$, a contradiction.
%----------------------------------------------------------
\par Prove the second statement of the proposition, i.e., the equivalence of the Euclidean and the  geodesic metrics in the graph $\GE$. As usual the geodesic metric $\DGE$ is defined by the formula
%----------------------------------------------------------
$$
\DGE(x,y):=\inf\,\smed_{i=0}^{m-1}\,\|z_i-z_{i+1}\|
$$
%---------------------------------------------------------
where the infimum is taken over all paths $\{z_0,...,z_m\}$ joining $x$ to $y$ in $\GE$. (Thus $z_0=x$, $z_m=y$, and $z_i\lr z_{i+1}$ for all $i=0,...,m-1$.)
%----------------------------------------------------------
\par Let $f(z):=\DGE(z,\hx)$, $z\in E$. Then for every $x,y\in E$, $x\lr y$,
%----------------------------------------------------------
$$
|f(x)-f(y)|=|\DGE(x,\hx)-\DGE(y,\hx)|\le \DGE(x,y)\le \|x-y\|
$$
%---------------------------------------------------------
so that, by \rf{LIPE}, $\|f\|_{\Lip(E)}\le C(n)$. Hence,  %----------------------------------------------------------
$$
\DGE(\hx,\hy)=|f(\hx)-f(\hy)|\le C(n)\, \|\hx-\hy\|
$$
%---------------------------------------------------------
proving the proposition.\bx
%----------------------------------------------------------
%@@@@@@@@@@@@@@@@@@@@@@@@@@@@@@@@@@@@@@@@@@@@@@@@@@@@@@@@@@
%@@@@@@@@@@@@@@@@@@@@@@@@@@@@@@@@@@@@@@@@@@@@@@@@@@@@@@@@@@
%@@@@@@@@@@@@@@@@@@@@@@@@@@@@@@@@@@@@@@@@@@@@@@@@@@@@@@@@@@
%@@@@@@@@@@@@@@@@@@@@@@@@@@@@@@@@@@@@@@@@@@@@@@@@@@@@@@@@@@
%----------------------------------------------------------
\begin{remark}\lbl{CR-TREE-PR}{\em Note that Theorem \reff{CR-TREE} easily follows from Theorems \reff{EX-TK} and \reff{EX-REFTK}. In fact, the necessity part of Theorem \reff{CR-TREE} and the inequality $\NW_{m,p,E}(\VP)\le C\,\PME$ directly follow from Definition \reff{DF-GRPR} and the necessity part of Theorem \reff{EX-TK}. Let us prove the sufficiency.
%---------------------------------------------------------
\par Let $\VP=\{P_x: x\in E\}$ be a Whitney $(m-1)$-field on $E$ such that $\NW_{m,p,E}(\VP)<\infty$. See \rf{NP-WIGL}. Let $\gamma\ge 1$ be the same as in Theorem \reff{EX-TK}, and let $A=\{\{x_i, y_i\} : i = 1, ..., k\}$ be an arbitrary finite $\gamma$-sparse collections of two point subsets of $E$. By $S$ we denote a subset of $E$ defined by
%----------------------------------------------------------
$$
S_A:=\bigcup_{i = 1, ..., k} \{x_i, y_i\}\,.
$$
%---------------------------------------------------------
\par Since $S_A$ is finite, by Theorem \reff{EX-REFTK}
there exists a connected tree $\Gamma_{S_A}$ whose set of vertices coincide with $S_A$ such that the conditions (i) and (ii) of this theorem are satisfied. Hence, by the assumption,  \rf{NP-WIGL} and the equivalence \rf{GE-NORM}, there exists a function $\tF\in \CMON$ which agrees with $\VP$ on $S_A$ such that $\|\tF\|_{\WMP}\le C\,\NW_{m,p,E}(\VP)$.
%---------------------------------------------------------
\par We again apply the necessity part of Theorem \reff{EX-TK} to the set $S_A$ and the function $\tF$, and obtain the following inequality:
%----------------------------------------------------------
$$
\left\{\,\smed_{i=1}^k\,\,
\smed_{|\alpha|\le m-1}\frac{|D^\alpha P_{x_i}(x_i)-D^\alpha P_{y_i}(x_i)|^p}
{\|x_i-y_i\|^{(m-|\alpha|)p-n}}\right\}^{1/p}
\le C\,\|\tF\|_{\WMP}\le C\,\NW_{m,p,E}(\VP).
$$
%----------------------------------------------------------
Since $A$ is arbitrary, by \rf{N-P-NEW}, $\Nc_{m,p,E}(\VP)\le C\,\NW_{m,p,E}(\VP)$ so that, by
\rf{EQV-JET}, $\PME\le C\,\NW_{m,p,E}(\VP)$. This completes the proof of Theorem \reff{CR-TREE}.\rbx
%----------------------------------------------------------
}
%----------------------------------------------------------
\end{remark}
%----------------------------------------------------------
%@@@@@@@@@@@@@@@@@@@@@@@@@@@@@@@@@@@@@@@@@@@@@@@@@@@@@@@@@@
%@@@@@@@@@@@@@@@@@@@@@@@@@@@@@@@@@@@@@@@@@@@@@@@@@@@@@@@@@@
%@@@@@@@@@@@@@@@@@@@@@@@@@@@@@@@@@@@@@@@@@@@@@@@@@@@@@@@@@@
%@@@@@@@@@@@@@@@@@@@@@@@@@@@@@@@@@@@@@@@@@@@@@@@@@@@@@@@@@@
%----------------------------------------------------------
%@@@@@@@@@@@@@@@@@@@@@@@@@@@@@@@@@@@@@@@@@@@@@@@@@@@@@@@@@@
%@@@@@@@@@@@@@@@@@@@@@@@@@@@@@@@@@@@@@@@@@@@@@@@@@@@@@@@@@@
%@@@@@@@@@@@@@@@@@@@@@@@@@@@@@@@@@@@@@@@@@@@@@@@@@@@@@@@@@@
%@@@@@@@@@@@@@@@@@@@@@@@@@      @@@@@@@@@@@@@@@@@@@@@@@@@@@
%@@@@@@@@@@@@@@@@@@@@@@@          @@@@@@@@@@@@@@@@@@@@@@@@@
%@@@@@@@@@@@@@@@@@@@@@              @@@@@@@@@@@@@@@@@@@@@@@
%@@@@@@@@@@@@@@@@@@@     SECTION 8    @@@@@@@@@@@@@@@@@@@@@
%@@@@@@@@@@@@@@@@@@@@@              @@@@@@@@@@@@@@@@@@@@@@@
%@@@@@@@@@@@@@@@@@@@@@@@          @@@@@@@@@@@@@@@@@@@@@@@@@
%@@@@@@@@@@@@@@@@@@@@@@@@@      @@@@@@@@@@@@@@@@@@@@@@@@@@@
%@@@@@@@@@@@@@@@@@@@@@@@@@@@@@@@@@@@@@@@@@@@@@@@@@@@@@@@@@@
%@@@@@@@@@@@@@@@@@@@@@@@@@@@@@@@@@@@@@@@@@@@@@@@@@@@@@@@@@@
%@@@@@@@@@@@@@@@@@@@@@@@@@@@@@@@@@@@@@@@@@@@@@@@@@@@@@@@@@@
%----------------------------------------------------------
\SECT{8. $\WMP$\,-jets on closed subsets of $\RN$.}{8}
%----------------------------------------------------------
\addtocontents{toc}{8. $\WMP$-jets on closed subsets of $\RN$. \hfill\thepage\par\VST}
%----------------------------------------------------------
%@@@@@@@@@@@@@@@@@@@@@@@@@@@@@@@@@@@@@@@@@@@@@@@@@@@@@@@@@@
%----------------------------------------------------------
\indent\par In this section we prove an analogue of Theorem \reff{EX-TK} for the normed Sobolev space, i.e., a variational criterion for jets generated by $\WMP$-functions.
%----------------------------------------------------------
\par Let $E$ be a closed subset of $\RN$, and let $\VP=\{P_x: x\in E\}$ be a family of polynomials of degree at most $m-1$ indexed by points of $E$. We define the $W^m_p$-``norm'' of $\VP$ by
%---------------------------------------------------------
\bel{N-VP-W}
\PMEW:=\inf\left\{\|F\|_{\WMP}:F\in \WMP,\, T^{m-1}_x[F]=P_x~\text{for every}~x\in E\right\}.
\ee
%----------------------------------------------------------
\par In this section given $\ve>0$ we let $E_\ve$ denote  the $\ve$-neighborhood of $E$.
%----------------------------------------------------------
%@@@@@@@@@@@@@@@@@@@@@@@@@@@@@@@@@@@@@@@@@@@@@@@@@@@@@@@@@@
%@@@@@@@@@@@@@@@@@@@@@@@@@@@@@@@@@@@@@@@@@@@@@@@@@@@@@@@@@@
%@@@@@@@@@@@@@@@@@@@@@@@@@@@@@@@@@@@@@@@@@@@@@@@@@@@@@@@@@@
%@@@@@@@@@@@@@@@@@@@@@@@@@@@@@@@@@@@@@@@@@@@@@@@@@@@@@@@@@@
%@@@@@@@@@@@@@@@@@@@@@@@@@@@@@@@@@@@@@@@@@@@@@@@@@@@@@@@@@@
%@@@@@@@@@@@@@@@@@@@@@@@@@@@@@@@@@@@@@@@@@@@@@@@@@@@@@@@@@@
%@@@@@@@@@@@@@@@@@@@@@@@@@@@@@@@@@@@@@@@@@@@@@@@@@@@@@@@@@@
%@@@@@@@@@@@@@@@@@@@@@@@@@@@@@@@@@@@@@@@@@@@@@@@@@@@@@@@@@@
%----------------------------------------------------------
\begin{theorem}\lbl{EX-WTK} Let $m\in\N$, $p\in(n,\infty)$, and let $\veh>0$. Fix a number $\theta\ge 1$ and let $\VE:E_{\veh}\to E$ be a measurable mapping such that
%----------------------------------------------------------
\bel{T-EM}
\|\VE(x)-x\|\le \theta \dist(x,E)~~~\text{for every}~~~x\in E_{\veh}.
\ee
%----------------------------------------------------------
\par There exists a constant $\gmh\ge 1$ depending only on $\theta$ for which the following result holds:
%----------------------------------------------------------
\par Suppose we are given a family of polynomials $\VP=\{P_x\in\PMRN: x\in E\}$. There exists a $C^{m-1}$-function $F\in\WMP$ such that %----------------------------------------------------------
\bel{J-PE-W}
T_{x}^{m-1}[F]=P_{x}~~~~\text{for every}~~~~x\in E
\ee
%----------------------------------------------------------
if and only if the function
%----------------------------------------------------------
\bel{PT-X}
P^{(\VE)}:=P_{\VE(x)}(x),~~~~x\in E_{\veh},
\ee
%----------------------------------------------------------
belongs to $L_p(E_{\veh})$, and the quantity $\Nc^*(\VP)=\Nc^*(\VP:m,n,p,E,\veh)$ defined by
%----------------------------------------------------------
\bel{NRM-W-NEW}
\Nc^*(\VP):=\sup\,\,
\left\{\smed_{i=1}^k\,\,\smed_{|\alpha|\le m-1}\frac{|D^\alpha P_{x_i}(x_i)-D^\alpha P_{y_i}(x_i)|^p}
{(\diam Q_i)^{(m-|\alpha|)p-n}}\right\}^{1/p}
\ee
%----------------------------------------------------------
is finite.  Here the supremum is taken over all finite families $\{Q_i: i = 1,...,k\}$ of pairwise disjoint cubes contained in $E_{\veh}$, and all choices of points $x_i,y_i\in (\gmh Q_i)\cap E$.
%----------------------------------------------------------
\par Furthermore,
%---------------------------------------------------------
\bel{EQV-JET-WP}
\PMEW\sim \|P^{(\VE)}\|_{L_p(E_{\veh})}+\Nc^*(\VP).
\ee
%---------------------------------------------------------
The constants of this equivalence depend only on $m,n,p,{\veh}$ and $\theta$.
%----------------------------------------------------------
\end{theorem}
%----------------------------------------------------------
%@@@@@@@@@@@@@@@@@@@@@@@@@@@@@@@@@@@@@@@@@@@@@@@@@@@@@@@@@@
%@@@@@@@@@@@@@@@@@@@@@@@@@@@@@@@@@@@@@@@@@@@@@@@@@@@@@@@@@@
%@@@@@@@@@@@@@@@@@@@@@@@@@@@@@@@@@@@@@@@@@@@@@@@@@@@@@@@@@@
%@@@@@@@@@@@@@@@@@@@@@@@@@@@@@@@@@@@@@@@@@@@@@@@@@@@@@@@@@@
%----------------------------------------------------------
\par See also Remark \reff{D-CR}. Note that for the case $m=1$, i.e., for the space $W^1_p(\RN)$, this result has been proven in \cite{Sh2}.
\medskip
%----------------------------------------------------------
%@@@@@@@@@@@@@@@@@@@@@@@@@@@@@@@@@@@@@@@@@@@@@@@@@@@@@@@@@@
%@@@@@@@@@@@@@@@@@@@@@@@@@@@@@@@@@@@@@@@@@@@@@@@@@@@@@@@@@@
%@@@@@@@@@@@@@@@@@@@@@@@@@@@@@@@@@@@@@@@@@@@@@@@@@@@@@@@@@@
%@@@@@@@@@@@@@@@@@@@@@@@@@@@@@@@@@@@@@@@@@@@@@@@@@@@@@@@@@@
%----------------------------------------------------------
\par Now given a closed $E\subset\RN$ we define a special mapping $V_E:\RN\to E$ which will enable us to refine the result of Theorem \reff{EX-WTK}. This mapping is generated by the lacunary projector $\PRL$ constructed in Subsection 6.2, see Theorem \reff{L-PE}.
%----------------------------------------------------------
\par We define $V_E$ as follows. We put
%----------------------------------------------------------
\bel{PV-ONE}
V_E(x):=x~~~\text{for every}~~~x\in E.
\ee
%----------------------------------------------------------
\par Let now $x\in\RN\setminus E$ and let $Q\in W_E$ be a Whitney cube containing $x$. By $L^{(Q)}$ we denote the (unique) lacuna which contains $Q$. See Subsection 6.1. Then we define $V_E(x)$ by letting
%----------------------------------------------------------
\bel{P-VE}
V_E(x):=\PRL(L^{(Q)}).
\ee
%----------------------------------------------------------
\par Clearly, $V_E$ is a measurable mapping which is well defined on the set
%----------------------------------------------------------
$$
S:=E\,\bigcup\, \left(\,\bigcup_{Q\in W_E} \intr(Q)\right).
$$
%----------------------------------------------------------
It is also clear that the Lebesgue measure of the set $\RN\setminus S$ is zero, so that  $V_E$ is well defined a.e. on $\RN$. Note that, by property \rf{DQ-E} of Whitney's cubes and by property (i) of Theorem \reff{L-PE}, the mapping $V=V_E$ satisfies inequality \rf{T-EM} with some absolute constant $\theta$.
%----------------------------------------------------------
%@@@@@@@@@@@@@@@@@@@@@@@@@@@@@@@@@@@@@@@@@@@@@@@@@@@@@@@@@@
%@@@@@@@@@@@@@@@@@@@@@@@@@@@@@@@@@@@@@@@@@@@@@@@@@@@@@@@@@@
%@@@@@@@@@@@@@@@@@@@@@@@@@@@@@@@@@@@@@@@@@@@@@@@@@@@@@@@@@@
%@@@@@@@@@@@@@@@@@@@@@@@@@@@@@@@@@@@@@@@@@@@@@@@@@@@@@@@@@@
%----------------------------------------------------------
\begin{theorem}\lbl{WP-PR} Let $m\in\N$, $p\in(n,\infty)$, $\ve>0$. There exists an absolute constant $\gamma\ge 1$ for which the following result holds:
%----------------------------------------------------------
\par Suppose we are given a family $\VP=\{P_x: x\in E\}$
of polynomials of degree at most $m-1$ indexed by points of $E$. There exists a $C^{m-1}$-function $F\in\WMP$ such that \rf{J-PE-W} is satisfied if and only if the function
%----------------------------------------------------------
\bel{PVE-1}
P^{(V_E)}:=P_{V_E(x)}(x),~~~~x\in E_\ve\,,
\ee
%----------------------------------------------------------
belongs to $L_p(E_\ve)$ and the quantity $\Nc^{\flat}(\VP)=\Nc^{\flat}(\VP:m,n,p,\ve,E)$ defined by %----------------------------------------------------------
\bel{N-P-FL}
\Nc^{\flat}(\VP):=\sup\left\{\,\smed_{i=1}^k\,\,
\smed_{|\alpha|\le m-1}
\frac{|D^\alpha P_{x_i}(x_i)-D^\alpha P_{y_i}(x_i)|^p}
{\|x_i-y_i\|^{(m-|\alpha|)p-n}}\right\}^{1/p}
\ee
%----------------------------------------------------------
is finite.  Here the supremum is taken over all finite $\gamma$-sparse collections $\{\{x_i, y_i\} : i = 1, ..., k\}$ of two point subsets of $E$ with $\|x_i-y_i\|\le \ve$, $i=1,...,k$.
%----------------------------------------------------------
\par Furthermore,
%---------------------------------------------------------
\bel{H-JE}
\PMEW\sim \|P^{(V_E)}\|_{L_p(E_\ve)}+\Nc^{\flat}(\VP).
\ee
%---------------------------------------------------------
The constants of this equivalence depend only on $m,n,p,$ and $\ve$.
%----------------------------------------------------------
\end{theorem}
%----------------------------------------------------------
\par See also Remark \reff{D-CR1}.
%----------------------------------------------------------
%@@@@@@@@@@@@@@@@@@@@@@@@@@@@@@@@@@@@@@@@@@@@@@@@@@@@@@@@@@
%@@@@@@@@@@@@@@@@@@@@@@@@@@@@@@@@@@@@@@@@@@@@@@@@@@@@@@@@@@
%@@@@@@@@@@@@@@@@@@@@@@@@@@@@@@@@@@@@@@@@@@@@@@@@@@@@@@@@@@
%@@@@@@@@@@@@@@@@@@@@@@@@@@@@@@@@@@@@@@@@@@@@@@@@@@@@@@@@@@
%----------------------------------------------------------
\medskip
\par We turn to the proofs of Theorems \reff{EX-WTK} and \reff{WP-PR}.
\medskip
%----------------------------------------------------------
\par {\it (Necessity.)} The proof of the necessity relies on the following auxiliary lemma.
%----------------------------------------------------------
%@@@@@@@@@@@@@@@@@@@@@@@@@@@@@@@@@@@@@@@@@@@@@@@@@@@@@@@@@@
%@@@@@@@@@@@@@@@@@@@@@@@@@@@@@@@@@@@@@@@@@@@@@@@@@@@@@@@@@@
%@@@@@@@@@@@@@@@@@@@@@@@@@@@@@@@@@@@@@@@@@@@@@@@@@@@@@@@@@@
%@@@@@@@@@@@@@@@@@@@@@@@@@@@@@@@@@@@@@@@@@@@@@@@@@@@@@@@@@@
%----------------------------------------------------------
\begin{lemma}\lbl{LP-VE} Let $F\in\WMP\cap \CMON$. Then the function $\tP:=T^{m-1}_{V(x)}[F](x)$, $x\in E_\ve,$ belongs to $L_p(E_\ve)$ and
%----------------------------------------------------------
\bel{L-WE}
\|\tP\|_{L_p(E_\ve)}\le C\|F\|_{\WMP}.
\ee
%----------------------------------------------------------
Here $C=C(m,n,p,\theta,\ve)$.
%----------------------------------------------------------
\end{lemma}
%----------------------------------------------------------
%@@@@@@@@@@@@@@@@@@@@@@@@@@@@@@@@@@@@@@@@@@@@@@@@@@@@@@@@@@
%----------------------------------------------------------
\par {\it Proof.} By \rf{T-EM}, $V(x)=x$ on $E$ so that $\tP(x)=T^{m-1}_x[F](x)=F(x)$ for every $x\in E$. Hence
%---------------------------------------------------------
\bel{PE-W2}
\|\tP\|_{L_p(E_\ve)}\le \|F\|_{L_p(E_\ve)}+\|F-\tP\|_{L_p(E_\ve)}=\|F\|_{L_p(E_\ve)}+
\|F-\tP\|_{L_p(E_\ve\setminus E)}.
\ee
%---------------------------------------------------------
\par Let $q:=(p+n)/2$. Prove that for every $x\in E_\ve$ the following inequality
%---------------------------------------------------------
\bel{PW-F}
|F(x)-\tP(x)|\le C_1\,(\Mc[(\nabla^m F)^q](x))^{\frac1q}
\ee
%---------------------------------------------------------
holds. In fact, let $Q=Q(x,\|x-V(x)\|)$. Since $x\in E_\ve$, we have $\dist(x,E)\le \ve$ so that, by \rf{T-EM},
$\diam Q\le 2\theta\ve$.
%---------------------------------------------------------
\par Let us apply to $x$ and $y=V(x)$ the Sobolev-Poincar\'e inequality \rf{SP-BTQ} with $\beta=0$. By this inequality,
%---------------------------------------------------------
\be
|F(x)-\tP(x)|&=&|F(x)-T^{m-1}_{V(x)}[F](x)|\le C
\,(\diam Q)^{m}\left(\frac{1}{|Q|}
\intl_{Q}(\nabla^m F(u))^qdu\right)^{\frac{1}{q}}
\nn\\
&\le&
C\,(2\theta\ve)^{m}\left(\frac{1}{|Q|}
\intl_{Q}(\nabla^m F(u))^qdu\right)^{\frac{1}{q}}
\le C_1\,(\Mc[(\nabla^m F)^q](x))^{\frac1q}
\nn
\ee
%---------------------------------------------------------
proving \rf{PW-F}. Here $C_1$ is a constant depending only on $m,n,p,\theta,$ and $\ve$.
%---------------------------------------------------------
\par Hence,
%---------------------------------------------------------
$$
\|F-\tP\|_{L_p(E_\ve\setminus E)}\le C_1\,
\|(\Mc[(\nabla^m F)^q])^{\frac1q}\|
_{L_p(E_\ve\setminus E)}\le C_1\,
\|(\Mc[(\nabla^m F)^q])^{\frac1q}\|_{L_p(\RN)}.
$$
%---------------------------------------------------------
Since $q<p$, by the Hardy-Littlewood maximal theorem,
%---------------------------------------------------------
$$
\|F-\tP\|_{L_p(E_\ve\setminus E)}\le C_2\,
\|\nabla^m F\|_{L_p(\RN)}\le C_2\, \|F\|_{\WMP}.
$$
%---------------------------------------------------------
\par Combining this inequality with \rf{PE-W2} we  obtain the required inequality \rf{L-WE}.\bx
%----------------------------------------------------------
%@@@@@@@@@@@@@@@@@@@@@@@@@@@@@@@@@@@@@@@@@@@@@@@@@@@@@@@@@@
%@@@@@@@@@@@@@@@@@@@@@@@@@@@@@@@@@@@@@@@@@@@@@@@@@@@@@@@@@@
%@@@@@@@@@@@@@@@@@@@@@@@@@@@@@@@@@@@@@@@@@@@@@@@@@@@@@@@@@@
%@@@@@@@@@@@@@@@@@@@@@@@@@@@@@@@@@@@@@@@@@@@@@@@@@@@@@@@@@@
%----------------------------------------------------------
\smallskip
\par Now the necessity part of Theorem \reff{EX-WTK} and Theorem \reff{WP-PR} and the inequalities
%---------------------------------------------------------
$$
\|P^{(\VE)}\|_{L_p(E_{\veh})}+\Nc^*(\VP)
\le C\,\PMEW~~~\text{and}~~~
\|P^{(V_E)}\|_{L_p(E_\ve)}+\Nc^{\flat}(\VP)\le C\,\PMEW
$$
%--------------------------------------------------------
directly follow from Lemma \reff{LP-VE}, Proposition \reff{N-VC} and definition \rf{N-VP-W}.
%----------------------------------------------------------
%@@@@@@@@@@@@@@@@@@@@@@@@@@@@@@@@@@@@@@@@@@@@@@@@@@@@@@@@@@
%@@@@@@@@@@@@@@@@@@@@@@@@@@@@@@@@@@@@@@@@@@@@@@@@@@@@@@@@@@
%@@@@@@@@@@@@@@@@@@@@@@@@@@@@@@@@@@@@@@@@@@@@@@@@@@@@@@@@@@
%@@@@@@@@@@@@@@@@@@@@@@@@@@@@@@@@@@@@@@@@@@@@@@@@@@@@@@@@@@
%----------------------------------------------------------
\bigskip
\par {\it Proof of the sufficiency part of Theorem \reff{WP-PR}.} The proof of the sufficiency is based on a slight modification of the lacunary extension method suggested in Subsection 7.2.
%----------------------------------------------------------
\par Let $L\in\LE$ be a lacuna and let a cube $Q\in L$.
We recall that
%----------------------------------------------------------
$$
a_Q:=\PRL(L)~~~\text{and}~~~P^{(Q)}:=P_{a_Q}.
$$
%----------------------------------------------------------
See \rf{LM-WM-DEF}.
%----------------------------------------------------------
\par Let $\ve>0$ and let $\delta:=\td\ve$. Let $\VP=\{P_x: x\in E\}$ be a Whitney $(m-1)$-field on $E$. Given a cube $Q\in W_E$ we let  $P^{(Q,\ve)}$ denote a polynomial of degree at most $m-1$ defined by the following formula:
%----------------------------------------------------------
\bel{P-EP1}
P^{(Q,\ve)}:=\left \{
%----------------------------------------------------------
\begin{array}{ll}
P^{(Q)},& \text{if}~~\diam Q<\delta,\\
0,& \text{if}~~\diam Q\ge\delta.
\end{array}
%----------------------------------------------------------
\right.
\ee
%----------------------------------------------------------
%@@@@@@@@@@@@@@@@@@@@@@@@@@@@@@@@@@@@@@@@@@@@@@@@@@@@@@@@@@
%----------------------------------------------------------
\par Finally, we define the extension $F_\ve=F(x;\ve,\VP)$, $x\in\RN$, by
%----------------------------------------------------------
\bel{W-DEF-F}
F_\ve(x):=\left \{
%----------------------------------------------------------
\begin{array}{ll}
P_x(x),& x\in E,\smallskip\\
\sbig\limits_{Q\in W_E}\,\,
\varphi_Q(x)P^{(Q,\ve)}(x),& x\in\RN\setminus E.
\end{array}
%----------------------------------------------------------
\right.
\ee
%----------------------------------------------------------
%@@@@@@@@@@@@@@@@@@@@@@@@@@@@@@@@@@@@@@@@@@@@@@@@@@@@@@@@@@
%@@@@@@@@@@@@@@@@@@@@@@@@@@@@@@@@@@@@@@@@@@@@@@@@@@@@@@@@@@
%@@@@@@@@@@@@@@@@@@@@@@@@@@@@@@@@@@@@@@@@@@@@@@@@@@@@@@@@@@
%@@@@@@@@@@@@@@@@@@@@@@@@@@@@@@@@@@@@@@@@@@@@@@@@@@@@@@@@@@
%----------------------------------------------------------
\medskip
\par Let $F$ be the extension operator constructed in Subsection 7.2, i.e., the operator given by the Whitney extension formula
%----------------------------------------------------------
\bel{F-LAC}
F(x):=\left \{
%----------------------------------------------------------
\begin{array}{ll}
P_x(x),& x\in E,\smallskip\\
\sbig\limits_{Q\in W_E}\,\,
\varphi_Q(x)P^{(Q)}(x),& x\in\RN\setminus E.
\end{array}
%----------------------------------------------------------
\right.
\ee
%----------------------------------------------------------
%@@@@@@@@@@@@@@@@@@@@@@@@@@@@@@@@@@@@@@@@@@@@@@@@@@@@@@@@@@
%@@@@@@@@@@@@@@@@@@@@@@@@@@@@@@@@@@@@@@@@@@@@@@@@@@@@@@@@@@
%@@@@@@@@@@@@@@@@@@@@@@@@@@@@@@@@@@@@@@@@@@@@@@@@@@@@@@@@@@
%@@@@@@@@@@@@@@@@@@@@@@@@@@@@@@@@@@@@@@@@@@@@@@@@@@@@@@@@@@
%----------------------------------------------------------
\par Let us note two properties of the function $F_\ve$.
%----------------------------------------------------------
%@@@@@@@@@@@@@@@@@@@@@@@@@@@@@@@@@@@@@@@@@@@@@@@@@@@@@@@@@@
%@@@@@@@@@@@@@@@@@@@@@@@@@@@@@@@@@@@@@@@@@@@@@@@@@@@@@@@@@@
%@@@@@@@@@@@@@@@@@@@@@@@@@@@@@@@@@@@@@@@@@@@@@@@@@@@@@@@@@@
%@@@@@@@@@@@@@@@@@@@@@@@@@@@@@@@@@@@@@@@@@@@@@@@@@@@@@@@@@@
%----------------------------------------------------------
\begin{lemma}\lbl{EP-F} (i). If $\dist(x,E)< \delta/4$, then $F_\ve(x)=F(x)$;\smallskip
%----------------------------------------------------------
\par (ii). $\supp F_\ve\subset E_{\tau}$ where $\tau=20\delta$.
%----------------------------------------------------------
\end{lemma}
%----------------------------------------------------------
%@@@@@@@@@@@@@@@@@@@@@@@@@@@@@@@@@@@@@@@@@@@@@@@@@@@@@@@@@@
%----------------------------------------------------------
\par {\it Proof.} (i). Clearly, by \rf{W-DEF-F} and \rf{F-LAC}, $F_\ve(x)=F(x)$ on $E$.
%----------------------------------------------------------
\par Suppose that $0<\dist(x,E)< \delta/4$. Let
%---------------------------------------------------------
$$
A_x:=\{Q\in W_E: Q^*\ni x\}.
$$
%---------------------------------------------------------
Recall that $Q^*:=\frac98 Q$. Then, by property (c) of Lemma \reff{P-U} and by \rf{F-LAC} and \rf{W-DEF-F},
%----------------------------------------------------------
\bel{FEF}
F(x)=
\sbig\limits_{Q\in A_x}\,\,
\varphi_Q(x)P^{(Q)}(x)~~~\text{and}~~~
F_\ve(x)=\sbig\limits_{Q\in A_x}\,\,
\varphi_Q(x)P^{(Q,\ve)}(x).
%----------------------------------------------------------
\ee
%----------------------------------------------------------
\par Let $K\in W_E$ and let $x\in K$. Then, by \rf{DQ-E},
%---------------------------------------------------------
$$
\diam K\le\dist(K,E)\le\dist(x,E)<\delta/4.
$$
%---------------------------------------------------------
Since $Q^*\cap K\ne\emp$ for each $Q\in A_x$, by property (1) of Lemma \reff{Wadd}, $\diam Q\le 4\diam K$, so that
%---------------------------------------------------------
$$
\diam Q<\delta~~~\text{for every}~~~Q\in A_x.
$$
%---------------------------------------------------------
Hence, by \rf{P-EP1}, $P^{(Q,\ve)}=P^{(Q)}$ for all $Q\in A_x$. Combining this property with \rf{FEF} we obtain the statement of part (i) of the lemma.\medskip
%----------------------------------------------------------
%@@@@@@@@@@@@@@@@@@@@@@@@@@@@@@@@@@@@@@@@@@@@@@@@@@@@@@@@@@
%@@@@@@@@@@@@@@@@@@@@@@@@@@@@@@@@@@@@@@@@@@@@@@@@@@@@@@@@@@
%@@@@@@@@@@@@@@@@@@@@@@@@@@@@@@@@@@@@@@@@@@@@@@@@@@@@@@@@@@
%----------------------------------------------------------
\par (ii). Let $\dist(x,E)\ge 20\delta$. Prove that $F_\ve(x)=0$.
%----------------------------------------------------------
\par By \rf{P-EP1} and \rf{FEF}, it suffices to show that
$\diam Q\ge \delta$ for every $Q\in A_x$.
%----------------------------------------------------------
\par Let $K\in W_E$ and let $x\in K$. (Clearly, $K\in A_x$.) Then, $Q^*\cap K\ne\emp$ for each $Q\in A_x$ so that, by part (3) of Lemma \reff{Wadd}, $Q\cap K\ne\emp$.
Hence, by part (1) of this lemma, $\diam Q\ge \frac14\diam K$. But $K\in W_E$ so that
%---------------------------------------------------------
$$
20\delta\le\dist(x,E)\le \dist(K,E)+\diam K\le 5\diam K.
$$
%---------------------------------------------------------
Hence $\diam K\ge 4\delta$ proving the required inequality $\diam Q\ge \delta$.\bx
%----------------------------------------------------------
%@@@@@@@@@@@@@@@@@@@@@@@@@@@@@@@@@@@@@@@@@@@@@@@@@@@@@@@@@@
%@@@@@@@@@@@@@@@@@@@@@@@@@@@@@@@@@@@@@@@@@@@@@@@@@@@@@@@@@@
%@@@@@@@@@@@@@@@@@@@@@@@@@@@@@@@@@@@@@@@@@@@@@@@@@@@@@@@@@@
%----------------------------------------------------------
\medskip
\par Let $\gamma\ge 1$ be the same absolute constant as in Theorem \reff{EX-TK}. Suppose that the Whitney $(m-1)$-field $\VP=\{P_x:x\in E\}$ satisfies the conditions of the sufficiency part of Theorem \reff{WP-PR}. Thus:\smallskip
%----------------------------------------------------------
\par (i) The function $P^{(V_E)}$ belongs to $L_p(E_\ve)$. See \rf{PVE-1} and \rf{P-VE};
\medskip
%----------------------------------------------------------
\par (ii) Let $\lw:=\Nc^\flat(\VP)^p$, see\rf{N-P-FL}. Then  $\lw<\infty$ so that for every finite $\gamma$-sparse collection $\{\{x_i, y_i\} : i = 1, ..., k\}$ of two point subsets of $E$ with
%---------------------------------------------------------
\bel{XY-E}
\|x_i-y_i\|\le \ve,~~i=1,...,k,
\ee
%---------------------------------------------------------
the following inequality
%----------------------------------------------------------
\bel{N-P-GW}
\smed_{i=1}^k\,\,
\smed_{|\alpha|\le m-1}
\frac{|D^\alpha P_{x_i}(x_i)-D^\alpha P_{y_i}(x_i)|^p}
{\|x_i-y_i\|^{(m-|\alpha|)p-n}}\le \lw
\ee
%----------------------------------------------------------
holds.
%----------------------------------------------------------
\medskip
\par In Section 7 we have proved that the function $F$ defined by the formula \rf{F-LAC} belongs to $C^{m-1}(\RN)$ and agrees with $\VP$ on $E$, i.e., \rf{J-PE-W} holds. Since, by Lemma \reff{EP-F}, $F_\ve$ coincides with $F$ on an open neighborhood of $E$, and $F_\ve\in C^{\infty}(\RN\setminus E)$, the function $F_\ve$ has similar properties, i.e., $F_\ve\in C^{m-1}(\RN)$ and $F_\ve$ agrees with $\VP$ on $E$.
%----------------------------------------------------------
\par Prove that
%---------------------------------------------------------
\bel{IN-JE}
\PMEW\le C\,\{\lw^{\frac1p}+\|P^{(V_E)}\|_{L_p(E_\ve)}\}
\ee
%---------------------------------------------------------
where $C$ is a constant depending only on $m,n,p,$ and $\ve$.
%----------------------------------------------------------
\par We will estimate the $L^m_p$-seminorm of $F_\ve$ using the scheme of the proof suggested in Subsection 7. More specifically, in Lemmas \reff{L-1E} and \reff{L-2E} we will prove analogues of Lemmas \reff{L-BG} and \reff{L-SMC} for the function $F_\ve$. Finally, in Lemma \reff{L-3E} we will estimate the $L_p$-norm of $F_\ve$.
%----------------------------------------------------------
\par  Let $\dw:=\delta/\tgm$ where $\tgm$ is the constant from Theorem \reff{L-PE}. Let $\Qc=\{Q_1,...,Q_k\}$ be a family of pairwise disjoint equal cubes in $\RN$ such that
%---------------------------------------------------------
\bel{DM-ET}
\diam Q_i\le \dw~~~\text{for every}~~~i=1,...,k.
\ee
%---------------------------------------------------------
Let $c_i:=c_{Q_i}$, $i=1,...,k,$ be the centers of these cubes. %----------------------------------------------------------
%@@@@@@@@@@@@@@@@@@@@@@@@@@@@@@@@@@@@@@@@@@@@@@@@@@@@@@@@@@
%@@@@@@@@@@@@@@@@@@@@@@@@@@@@@@@@@@@@@@@@@@@@@@@@@@@@@@@@@@
%----------------------------------------------------------
\begin{lemma}\lbl{L-1E} Suppose that the family $\Qc$ satisfies the condition \rf{D-L40}. %----------------------------------------------------------
%@@@@@@@@@@@@@@@@@@@@@@@@@@@@@@@@@@@@@@@@@@@@@@@@@@@@@@@@@@
%----------------------------------------------------------
\par Then for every $\beta, |\beta|=m-1$, and every $x_i\in Q_i$ the following inequality
%----------------------------------------------------------
$$
\sbig_{i=1}^k\,
\frac{|D^\beta F_\ve(x_i)-D^\beta F_\ve(c_i)|^p}
{(\diam Q_i)^{p-n}}
\le C\,\lw
$$
%----------------------------------------------------------
holds. Here $C>0$ is a constant depending only on $m,n,p$ and $\ve$.
%----------------------------------------------------------
\end{lemma}
%----------------------------------------------------------
%@@@@@@@@@@@@@@@@@@@@@@@@@@@@@@@@@@@@@@@@@@@@@@@@@@@@@@@@@@
%@@@@@@@@@@@@@@@@@@@@@@@@@@@@@@@@@@@@@@@@@@@@@@@@@@@@@@@@@@
%@@@@@@@@@@@@@@@@@@@@@@@@@@@@@@@@@@@@@@@@@@@@@@@@@@@@@@@@@@
%----------------------------------------------------------
\par {\it Proof.} By Lemma \reff{EP-F}, $F_\ve(x)=F(x)$ for $x\in E_{\delta'}$ where $\delta'=\delta/4$. Hence
%---------------------------------------------------------
$$
D^\beta F_\ve(x)=D^\beta F(x)~~\text{for every}~~\beta, |\beta|=m-1,~~\text{and every}~~x\in E_{\delta'}.
$$
%--------------------------------------------------------- %@@@@@@@@@@@@@@@@@@@@@@@@@@@@@@@@@@@@@@@@@@@@@@@@@@@@@@@@
%----------------------------------------------------------
\par Following the proof of Lemma \reff{L-BG} we first consider the case where $c_i\in E$ for each $i=1,...,k$. Let $I_1(F_\ve)$ be the quantity defined by \rf{I1-DF}. Prove that $I_1(F_\ve)\le \lw$.
%----------------------------------------------------------
\par In fact, since  $x_i\in Q_i$ and $\diam Q_i\le \dw$, the point $x_i\in E_{\dw}$. Hence, by part (i) of Lemma \reff{EP-F}, $D^\beta F(x_i)=D^\beta F_\ve(x_i)$ and
$D^\beta F(c_i)=D^\beta F_\ve(c_i)$
for every $\beta$, $|\beta|=m-1$, and every $i=1,...,k$. Hence $I_1(F_\ve)=I_1(F)$.
%----------------------------------------------------------
\par Since $x_i,c_i\in Q_i$ and
%----------------------------------------------------------
$$
\diam Q_i\le\dw<\delta/4<\ve,
$$
%----------------------------------------------------------
we conclude that $\|x_i-c_i\|<\ve$ as it is required in \rf{XY-E}. Now, using the same argument as in the proof of Lemma \reff{L-BG}, we obtain the desired inequality
$I_1(F_\ve)\le \lw$. C.f. \rf{I1-FEG}.
%----------------------------------------------------------
\par Let $I_2(F_\ve)$ be the quantity defined by \rf{I2-DF}. Prove that $I_2(F_\ve)\le C\,\lw$.
%--------------------------------------------------------- %@@@@@@@@@@@@@@@@@@@@@@@@@@@@@@@@@@@@@@@@@@@@@@@@@@@@@@@@
%----------------------------------------------------------
\par Again we will follow the proof given in Lemma \reff{L-BG} for a similar estimate of $I_2(F)$.
%----------------------------------------------------------
\par Using the same argument as for the case of $I_1(F)$, we conclude that $I_2(F_\ve)=I_2(F)$. Then we literally follow the proof of Lemma \reff{L-BG} after the definition \rf{I2-DF} of $I_2(F)$. This leads us to the corresponding estimate
%----------------------------------------------------------
$$
I_2(F_\ve)\le C\,\{I_2^{(1)}(F_\ve)+I_2^{(2)}(F_\ve)\}
$$
%----------------------------------------------------------
where $I_2^{(j)}(F_\ve)$, $j=1,2$, are defined in the same way as $I_2^{(j)}(F)$. See \rf{I21-DF1} and \rf{I22-DF2}.
%----------------------------------------------------------
\par By repeating the argument of Lemma \reff{L-BG} we show that $I_2^{(1)}(F_\ve)\le \lw$.  The unique additional requirement which we have to check is the inequality $\|a_{K_Q}-c_Q\|\le\ve$. See \rf{XY-E}.
%----------------------------------------------------------
\par But we know that $a_{K_Q}\in 5\tgm\,Q$, see \rf{50-Q}, and $\diam Q\le\dw$, see \rf{DM-ET}.
Hence
%----------------------------------------------------------
\bel{EP1-D}
\|a_{K_Q}-c_Q\|\le 5\tgm\,\diam Q\le 5\tgm\,\dw= 5\,\delta.
\ee
%----------------------------------------------------------
Recall that $\delta=\td\ve$ so that $\|a_{K_Q}-c_Q\|\le\ve$.
\smallskip
%----------------------------------------------------------
\par We turn to estimates of $I_2^{(2)}(F_\ve)$. Following the proof of the inequality $I_2^{(2)}(F)\le C\,\lambda$ in Lemma \reff{L-BG}, we obtain the required inequality $I_2^{(2)}(F_\ve)\le C\,\lw$ provided
%----------------------------------------------------------
\bel{A-HQ}
\|a_{H_Q}-a_{K_Q}\|\le\ve.
\ee
%----------------------------------------------------------
See \rf{XY-E}. We recall that the cubes $K_Q$ and $H_Q$ are determined by \rf{KQ-DF1} and \rf{HQ-DF} respectively.
%----------------------------------------------------------
\par But in the proof of this lemma we show that $K_Q\subset 5Q$, $a_{H_Q}\in 9\tgm K_Q$, and $a_{K_Q}\in 5\tgm\,Q$. See \rf{K-INQ}, \rf{50-Q} and \rf{AH-KQ}. Hence $a_{H_Q},a_{K_Q}\in 45\tgm Q$ so that
%----------------------------------------------------------
\bel{EP2-D}
\|a_{H_Q}-a_{K_Q}\|\le 45\tgm\,\diam Q\le 45\tgm\,\dw=45\,\delta\le \ve.
\ee
%---------------------------------------------------------
\par This proves the required inequality \rf{A-HQ}.
%----------------------------------------------------------
\par Finally we turn to the general case, i.e., to an arbitrary family of $\{Q_i: i=1,...,k\}$ of equal cubes of diameter at most $\dw$ satisfying inequality \rf{D-L40}. We again follow the proof of this part of  Lemma \reff{L-BG}. This enables us to reduce the proof to the case of the family $\hat{\Qc}=\{\OV_i:i=1,...,k\}$ of cubes defined by \rf{Q-H1}. These cubes are centered at $E$ so the proof is reduced to the previous cases proven above.
%---------------------------------------------------------
\par Also we know that
%----------------------------------------------------------
$$
\diam \OV_i\le \gamma'\diam Q_i\le  \gamma'\dw
$$
%---------------------------------------------------------
where $\gamma':=122$. Thus the proof is actually reduced to the known case but with $\dw'=\gamma'\dw$ instead of $\dw$.
%---------------------------------------------------------
\par We recall that $\dw=\delta/\tgm$. We use equality in inequalities \rf{EP1-D} and \rf{EP2-D} to prove that the left hand side of these inequalities are bounded by $\ve$. Clearly, the same is true whenever the constant $\dw$ in these inequalities is replaced by $\dw'$ satisfying the equality $\dw'=\gamma'\delta/\tgm$.
In fact, in this case the right hand sides of \rf{EP1-D} and \rf{EP2-D} are bounded by $5\gamma'\delta$ and $45\gamma'\delta$ respectively. Since $\delta=10^{-5}\ve$,
in the both cases the right hand side does not exceed $\ve$ which proves the lemma.\bx\medskip
%----------------------------------------------------------
%@@@@@@@@@@@@@@@@@@@@@@@@@@@@@@@@@@@@@@@@@@@@@@@@@@@@@@@@@@
%----------------------------------------------------------
\par The following properties of polynomials are well known.
%----------------------------------------------------------
\begin{lemma}\lbl{PLM-M} Let $P\in\PMP$ and let $1\le p<\infty$. Let $Q,\tQ$ be two cubes in $\RN$ such that $Q\subset\tQ$. Then %----------------------------------------------------------
%@@@@@@@@@@@@@@@@@@@@@@@@@@@@@@@@@@@@@@@@@@@@@@@@@@@@@@@@@@
%----------------------------------------------------------
$$
\sup_{\tQ}|P|\le C\,(\diam\tQ/\diam Q)^m\,\sup_{Q}|P|
$$
%----------------------------------------------------------
and
%----------------------------------------------------------
%@@@@@@@@@@@@@@@@@@@@@@@@@@@@@@@@@@@@@@@@@@@@@@@@@@@@@@@@@@
%----------------------------------------------------------
$$
\sup_{\tQ}|P|\le C\,\left\{\frac{1}{|Q|}\intl_Q|P(x)|^p\,dx\right\} ^{\frac1p}.
$$
%----------------------------------------------------------
\par Furthermore, for every $\xi,|\xi|\le m$, the following inequality
%----------------------------------------------------------
$$
\sup_{Q}|D^\xi P|\le C\,(\diam Q)^{-|\xi|}\,\sup_{Q}|P|
$$
%----------------------------------------------------------
holds. Here $C$ is a constant depending only on $m$ and $n$.
%----------------------------------------------------------
\end{lemma}
%----------------------------------------------------------
%@@@@@@@@@@@@@@@@@@@@@@@@@@@@@@@@@@@@@@@@@@@@@@@@@@@@@@@@@@
%@@@@@@@@@@@@@@@@@@@@@@@@@@@@@@@@@@@@@@@@@@@@@@@@@@@@@@@@@@
%@@@@@@@@@@@@@@@@@@@@@@@@@@@@@@@@@@@@@@@@@@@@@@@@@@@@@@@@@@
%@@@@@@@@@@@@@@@@@@@@@@@@@@@@@@@@@@@@@@@@@@@@@@@@@@@@@@@@@@
%----------------------------------------------------------
\begin{lemma}\lbl{L-2E} Suppose that the family $\Qc$ satisfies the condition \rf{D-S80}. %----------------------------------------------------------
%@@@@@@@@@@@@@@@@@@@@@@@@@@@@@@@@@@@@@@@@@@@@@@@@@@@@@@@@@@
%----------------------------------------------------------
\par Then for every $\beta, |\beta|=m-1$, and every $x_i\in Q_i$ the following inequality
%----------------------------------------------------------
$$
\sbig_{i=1}^k\,
\frac{|D^\beta F_\ve(x_i)-D^\beta F_\ve(c_i)|^p}
{(\diam Q_i)^{p-n}}
\le C\,\left\{\lw+\|P^{(V_E)}\|_{L_p(E_\ve)}^p\right\}
$$
%----------------------------------------------------------
holds. Here $C>0$ is a constant depending only on $m,n,p$ and $\ve$.
%----------------------------------------------------------
\end{lemma}
%----------------------------------------------------------
%@@@@@@@@@@@@@@@@@@@@@@@@@@@@@@@@@@@@@@@@@@@@@@@@@@@@@@@@@@
%@@@@@@@@@@@@@@@@@@@@@@@@@@@@@@@@@@@@@@@@@@@@@@@@@@@@@@@@@@
%@@@@@@@@@@@@@@@@@@@@@@@@@@@@@@@@@@@@@@@@@@@@@@@@@@@@@@@@@@
%----------------------------------------------------------
\par {\it Proof.} We literally follow the proof of Lemma  \reff{L-SMC} until the inequality \rf{IN-K1}. But before to continue the proof from this point let us make the following remark: Let $Q\in\Qc$ and, as in the proof of this lemma, let $K_Q$ be a cube such that $c_Q\in K_Q$. Note that if $Q\subset\RN\setminus E_\tau$, then, by part (ii) of Lemma \reff{EP-F}, $F_\ve|_Q\equiv 0$. This enables us to assume that $\dist(Q,E)\le\tau=20\delta$.
%----------------------------------------------------------
\par Now, by \rf{DQ-E} and \rf{18-KQ},
%----------------------------------------------------------
$$
\diam K_Q\le\dist(K_Q,E)\le \dist(c_Q,E)\le\dist(Q,E)+\diam Q\le\tau+\tfrac18\diam K_Q
$$
%---------------------------------------------------------
proving that $\diam K_Q\le 8\tau/7\le 24\,\delta$. Then, by part (1) of Lemma \reff{Wadd},
%----------------------------------------------------------
\bel{24-KD}
\diam H\le 4\cdot 24\,\delta=96\,\delta~~~\text{for every cube}~~H\in W_E~~\text{such that}~~H\cap K_Q\ne\emp.
\ee
%----------------------------------------------------------
We will use this estimate at the end of the proof.
%----------------------------------------------------------
\par We continue the proof of the lemma as follows. Let
%----------------------------------------------------------
$$
\Kc^{(1)}:=\{K\in\Kc:\diam K\le\dw\}~~~
\text{and let}~~~
\Kc^{(2)}:=\{K\in\Kc:\diam K>\dw\}.
$$
%---------------------------------------------------------
be a partition of the family $\Kc$. Then
%----------------------------------------------------------
\bel{FN-ST}
I:=\sbig_{i=1}^k\,
\frac{|D^\beta F_\ve(x_i)-D^\beta F_\ve(c_i)|^p}
{(\diam Q_i)^{p-n}}\le C\,\{I^{(1)}+I^{(2)}\}
\ee
%---------------------------------------------------------
where
%----------------------------------------------------------
$$
I^{(j)}:=\,\sbig_{K\in \Kc^{(j)}}
\,\,\sbig\limits_{|\xi|\le m-1}\,
\,\frac{|D^{\xi}P_{a_{\KM}}(a_K)-D^{\xi}P_{a_K}(a_K)|^p}
{\|a_{\KM}-a_K\|^{(m-|\xi|)p-n}},~~~j=1,2.
$$
%----------------------------------------------------------
Recall that $\KM$ is a cube touching $K$ and satisfying \rf{KM-DF}.
%----------------------------------------------------------
\par Let us prove that
%----------------------------------------------------------
\bel{I1-LW}
I^{(1)}\le C\lw~~~\text{with}~~~ C=C(m,n,p,\ve).
\ee
%---------------------------------------------------------
In fact, literally following the proof of Lemma \reff{L-SMC} after the inequality \rf{L-FML}, we obtain \rf{I1-LW} provided $\|a_{\KM}-a_K\|\le \ve$. But, by \rf{AH-KQ}, the points $a_{\KM},a_K\in 9\tgm\,K$ so that
%----------------------------------------------------------
$$
\|a_{\KM}-a_K\|\le 9\tgm\,\diam K\le 9\tgm\,\dw=9\delta\le\ve.
$$
%---------------------------------------------------------
%@@@@@@@@@@@@@@@@@@@@@@@@@@@@@@@@@@@@@@@@@@@@@@@@@@@@@@@@@
%----------------------------------------------------------
\par Let us estimate $I^{(2)}$. Let $K\in \Kc^{(2)}$. Since $\|a_{\KM}-a_K\|\sim \diam K$, see \rf{DK-A2},
%----------------------------------------------------------
$$
I^{(2)}\le C\,\,\sbig_{K\in \Kc^{(2)}}
\,\,\sbig\limits_{|\xi|\le m-1}\,
\,\frac{|D^{\xi}P_{a_{\KM}}(a_K)-D^{\xi}P_{a_K}(a_K)|^p}
{(\diam K)^{(m-|\xi|)p-n}}.
$$
%----------------------------------------------------------
\par We know that $a_K\in \tK:=\tgm\,K$, see part (i) of Theorem \reff{L-PE}. Since $\diam \KM\le 4\diam K$ and $\KM\cap K\ne\emp$, the cube $\KM\subset 9 K$ so that $\KM\subset \tK$.
%----------------------------------------------------------
\par Now, by Lemma \reff{PLM-M}, for every $\xi$, $|\xi|\le m-1$, the following inequality
%----------------------------------------------------------
\be
|D^{\xi}P_{a_{\KM}}(a_K)|&\le& \sup_{\tK}|D^{\xi}P_{a_{\KM}}|
\le C\,\sup_{\KM}|D^{\xi}P_{a_{\KM}}|\le
C\,(\diam \KM)^{-|\xi|}\,\sup_{\KM}|P_{a_{\KM}}|\nn\\
&\le& C\,(\diam K)^{-|\xi|}\,
\left\{\frac{1}{|\KM|}\intl_{\KM}|P_{a_{\KM}}(x)|^p\,dx\right\} ^{\frac1p}\nn
\ee
%----------------------------------------------------------
holds. Hence,
%----------------------------------------------------------
$$
\frac{|D^{\xi}P_{a_{\KM}}(a_K)|^p}
{(\diam K)^{(m-|\xi|)p-n}}\le C^p\,(\diam K)^{-mp}\,
\intl_{\KM}|P_{a_{\KM}}(x)|^p\,dx.
$$
%----------------------------------------------------------
Since $K\in\Kc_2$, its diameter is at least $\eta$ so that
%----------------------------------------------------------
$$
\frac{|D^{\xi}P_{a_{\KM}}(a_K)|^p}
{(\diam K)^{(m-|\xi|)p-n}}\le C^p\,\eta^{-mp}\,
\intl_{\KM}|P_{a_{\KM}}(x)|^p\,dx\le C_1\,
\intl_{\KM}|P_{a_{\KM}}(x)|^p\,dx
$$
%----------------------------------------------------------
where $C_1$ is a constant depending only on $m,n,p,$ and $\ve$. (Recall that $\eta:=\delta/\tgm$ depends on $\ve$.)
\smallskip
%----------------------------------------------------------
\par By \rf{P-VE} and \rf{PVE-1}, $P_{a_{\KM}}(x)=P_{V_E(x)}(x)=P^{(V_E)}(x)$ for every $x\in \KM$, so that
%----------------------------------------------------------
$$
\frac{|D^{\xi}P_{a_{\KM}}(a_K)|^p}
{(\diam K)^{(m-|\xi|)p-n}}\le C_1\,
\intl_{\KM}|P^{(V_E)}(x)|^p\,dx.
$$
%----------------------------------------------------------
\par In the same way we prove that
%----------------------------------------------------------
$$
\frac{|D^{\xi}P_{a_{K}}(a_K)|^p}
{(\diam K)^{(m-|\xi|)p-n}}\le C_1\,
\intl_{K}|P^{(V_E)}(x)|^p\,dx.
$$
%----------------------------------------------------------
\par Let us now introduce a family $\Ac$ of cubes by letting
%----------------------------------------------------------
$$
\Ac:=\{K,\KM: K\in\Kc^{(2)}\}.
$$
%----------------------------------------------------------
Clearly, by Lemma \reff{Wadd}, the covering multiplicity of $\Ac$ is bounded by a constant $C=C(n)$.
%----------------------------------------------------------
\par We have
%----------------------------------------------------------
\be
I^{(2)}&\le& C\,\,\sbig_{K\in \Kc^{(2)}}
\,\,\sbig\limits_{|\xi|\le m-1}\,
\,\frac{|D^{\xi}P_{a_{\KM}}(a_K)|^p+|D^{\xi}P_{a_K}(a_K)|^p}
{(\diam K)^{(m-|\xi|)p-n}}\nn\\
&\le&
C_2\,\,\sbig_{K\in \Ac}\,\,
\intl_K|P_{a_K}(x)|^p\,dx=
C_2\,\,\sbig_{K\in \Ac}\,\,\intl_K|P^{(V_E)}(x)|^p\,dx
\nn
\ee
%----------------------------------------------------------
with $C_2=C_2(m,n,p,\ve)$. Since the covering multiplicity of $\Ac$ is bounded by $C(n)$, we obtain
%----------------------------------------------------------
$$
I^{(2)}\le C_2\,\,\intl_{U_{\Ac}}|P^{(V_E)}(x)|^p\,dx
$$
%----------------------------------------------------------
where $U_{\Ac}:=\cup\{K: K\in\Ac\}$.
%----------------------------------------------------------
\par Prove that $U_{\Ac}\subset E_\ve$. In fact, by \rf{24-KD}, $\diam K\le 96\,\delta$ for every $K\in\Ac$. Since each cube $K\in\Ac$ is a Whitney cube, for every $y\in K$ we have
%----------------------------------------------------------
\bel{IM-DT}
\dist(y,E)\le \dist(K,E)+\diam K\le 5\diam K\le 480\,\delta =480\cdot 10^{-5}\ve<\ve
\ee
%----------------------------------------------------------
proving that $K\subset E_\ve$. Hence $U_{\Ac}\subset E_\ve$ so that
%----------------------------------------------------------
$$
I^{(2)}\le C\,\,\intl_{U_{\Ac}}|P^{(V_E)}(x)|^p\,dx
\le C\,\,\intl_{E_\ve}|P^{(V_E)}(x)|^p\,dx=
\|P^{(V_E)}\|_{L_p(E_\ve)}^p.
$$
%----------------------------------------------------------
Combining this inequality with \rf{I1-LW} and \rf{FN-ST} we obtain the statement of the lemma.\bx
%----------------------------------------------------------
%@@@@@@@@@@@@@@@@@@@@@@@@@@@@@@@@@@@@@@@@@@@@@@@@@@@@@@@@@@
%@@@@@@@@@@@@@@@@@@@@@@@@@@@@@@@@@@@@@@@@@@@@@@@@@@@@@@@@@@
%@@@@@@@@@@@@@@@@@@@@@@@@@@@@@@@@@@@@@@@@@@@@@@@@@@@@@@@@@@
%@@@@@@@@@@@@@@@@@@@@@@@@@@@@@@@@@@@@@@@@@@@@@@@@@@@@@@@@@@
%@@@@@@@@@@@@@@@@@@@@@@@@@@@@@@@@@@@@@@@@@@@@@@@@@@@@@@@@@@
%----------------------------------------------------------
\begin{corollary}\lbl{COR-FE} The function $F_\ve$ belongs to the space $\LMP$ and its seminorm in this space satisfies the following inequality:
%---------------------------------------------------------
\bel{NM-FEL}
\|F_\ve\|_{\LMP}\le C\,\{\lw^{\frac1p}+\|P^{(V_E)}\|_{L_p(E_\ve)}\}.
\ee
%---------------------------------------------------------
Here $C$ is a constant depending only on $m,n,p,$ and $\ve$.
%---------------------------------------------------------
\end{corollary}
%----------------------------------------------------------
%@@@@@@@@@@@@@@@@@@@@@@@@@@@@@@@@@@@@@@@@@@@@@@@@@@@@@@@@@@
%@@@@@@@@@@@@@@@@@@@@@@@@@@@@@@@@@@@@@@@@@@@@@@@@@@@@@@@@@@
%@@@@@@@@@@@@@@@@@@@@@@@@@@@@@@@@@@@@@@@@@@@@@@@@@@@@@@@@@@
%@@@@@@@@@@@@@@@@@@@@@@@@@@@@@@@@@@@@@@@@@@@@@@@@@@@@@@@@@@
%@@@@@@@@@@@@@@@@@@@@@@@@@@@@@@@@@@@@@@@@@@@@@@@@@@@@@@@@@@
%----------------------------------------------------------
\par {\it Proof.} Theorem \reff{CR-SOB}, Lemma \reff{L-1E} and Lemma \reff{L-2E} imply the following property of $F_\ve$: the function $F_\ve$ is a $C^{m-1}$-function such that for every multiindex $\beta$ of order $m-1$ the function $D^\beta F_\ve\in\LOP$ and
%----------------------------------------------------------
$$
\|D^\beta F_\ve\|_{\LOP}\le C\,\{\lw^{\frac1p}+\|P^{(V_E)}\|_{L_p(E_\ve)}\}.
$$
%----------------------------------------------------------
Since weak derivatives commute, the function $F_\ve\in\LMP$ and inequality \rf{NM-FEL} holds.\bx
%----------------------------------------------------------
%@@@@@@@@@@@@@@@@@@@@@@@@@@@@@@@@@@@@@@@@@@@@@@@@@@@@@@@@@@
%@@@@@@@@@@@@@@@@@@@@@@@@@@@@@@@@@@@@@@@@@@@@@@@@@@@@@@@@@@
%@@@@@@@@@@@@@@@@@@@@@@@@@@@@@@@@@@@@@@@@@@@@@@@@@@@@@@@@@@
%@@@@@@@@@@@@@@@@@@@@@@@@@@@@@@@@@@@@@@@@@@@@@@@@@@@@@@@@@@
%@@@@@@@@@@@@@@@@@@@@@@@@@@@@@@@@@@@@@@@@@@@@@@@@@@@@@@@@@@
%----------------------------------------------------------
\begin{lemma}\lbl{L-3E} The following inequality
%----------------------------------------------------------
\bel{LP-I7}
\|F_\ve\|_{\LPRN}\le C\, \|P^{(V_E)}\|_{L_p(E_\ve)}
\ee
%----------------------------------------------------------
holds. Here $C=C(m,n,p,\ve)$.
%----------------------------------------------------------
\end{lemma}
%----------------------------------------------------------
%@@@@@@@@@@@@@@@@@@@@@@@@@@@@@@@@@@@@@@@@@@@@@@@@@@@@@@@@@@
%@@@@@@@@@@@@@@@@@@@@@@@@@@@@@@@@@@@@@@@@@@@@@@@@@@@@@@@@@@
%@@@@@@@@@@@@@@@@@@@@@@@@@@@@@@@@@@@@@@@@@@@@@@@@@@@@@@@@@@
%----------------------------------------------------------
\par {\it Proof.} We recall that, by \rf{PV-ONE}, $V_E(x)=x$ on $E$ so that, by \rf{W-DEF-F}, $P^{(V_E)}(x)=F_\ve(x)$, $x\in E$. Also, by Lemma \reff{EP-F}, $\supp F_\ve\subset E_\tau$ with $\tau=20\,\delta$. Hence
%----------------------------------------------------------
\bel{F-MJ}
\|F_\ve\|_{\LPRN}^p=\|P^{(V_E)}\|_{L_p(E)}^p+
\|F_\ve\|_{L_p(\RN\setminus E)}^p=\|P^{(V_E)}\|_{L_p(E)}^p+
\|F_\ve\|_{L_p(E_\tau\setminus E)}^p.
\ee
%----------------------------------------------------------
\par Fix a point $y\in E_\tau\setminus E$. Let $y\in K$ where $K\in W_E$ is a Whitney cube. Let us apply Lemma \reff{C-A} to $y$, the cube $K$, the function $F_\ve$ and $\alpha=0$. By this lemma,
%----------------------------------------------------------
$$
|F(y)-P_{a_K}(y)|\le C\,
\smed\limits_{Q\in T(K)}\,\,
\smed\limits_{|\xi|\le m-1}\,\,
(\diam K)^{|\xi|} \,|D^{\xi}P_{a_Q}(a_K)-D^{\xi}P_{a_K}(a_K)|
$$
%----------------------------------------------------------
so that
%----------------------------------------------------------
$$
|F(y)-P_{a_K}(y)|\le C\,
\smed\limits_{Q\in T(K)}\,\,
\smed\limits_{|\xi|\le m-1}\,\,
(\diam K)^{|\xi|} \,|D^{\xi}P_{a_Q}(a_K)|\,.
$$
%----------------------------------------------------------
\par We know that $a_K\in \tK:=\tgm K$, see Theorem \reff{L-PE}, so that, by Lemma \reff{PLM-M}, for every $Q\in T(K)$
%----------------------------------------------------------
$$
(\diam K)^{|\xi|}\,|D^{\xi}P_{a_Q}(a_K)|
\le (\diam K)^{|\xi|}\sup_{\tK} |D^{\xi}P_{a_Q}|\le
C\,\sup_{\tK} |P_{a_Q}|.
$$
%----------------------------------------------------------
\par Since $Q\subset 9K$ (recall that $Q$ and $K$ are touching Whitney cubes) and $\tgm>9$, we have $Q\subset \tK$. Again, applying Lemma \reff{PLM-M} we get
%----------------------------------------------------------
$$
(\diam K)^{|\xi|}\,|D^{\xi}P_{a_Q}(a_K)|
\le C\,\sup_{\tK} |P_{a_Q}|\le C\,\sup_{Q} |P_{a_Q}|
\le C\,
\left\{\frac{1}{|Q|}\intl_{Q}|P_{a_{Q}}(x)|^p\,dx\right\} ^{\frac1p}
$$
%----------------------------------------------------------
provided $Q\in T(K)$ and $|\xi|\le m-1$. The same lemma implies the following inequality
%----------------------------------------------------------
$$
|P_{a_K}(y)|\le C\,
\left\{\frac{1}{|K|}\intl_{K}|P_{a_{K}}(x)|^p\,dx\right\} ^{\frac1p}.
$$
%----------------------------------------------------------
Hence, for every $y\in K$
%----------------------------------------------------------
$$
|F_\ve(y)|^p\le 2^p\{|F_\ve(y)-P_{a_K}(y)|^p+|P_{a_K}(y)|^p\}
\le C\,\smed\limits_{Q\in T(K)}\,\,
\left\{\frac{1}{|Q|}\intl_{Q}|P_{a_{Q}}(x)|^p\,dx\right\} ^{\frac1p}
$$
%----------------------------------------------------------
so that
%----------------------------------------------------------
\be
\intl_{K}|F_\ve(y)|^p\,dy&\le&
C\,\smed\limits_{Q\in T(K)}\,\,
\frac{|K|}{|Q|}\intl_{Q}|P_{a_{Q}}|^p\,dx\nn\\
&\le&
C\,\smed\limits_{Q\in T(K)}\,\,
\intl_{Q}|P_{a_{Q}}|^p\,dx=
C\,\smed\limits_{Q\in T(K)}\,\,
\intl_{Q}|P^{(V_E)}|^p\,dx.\nn
\ee
%----------------------------------------------------------
\par Let us introduce a family $\Ac$ of cubes by letting
%----------------------------------------------------------
$$
\Ac:=\{Q\in T(K): K\cap E_\tau\ne\emp\}.
$$
%----------------------------------------------------------
\par Prove that $Q\subset E_\ve$ for every $Q\in\Ac$. In fact, let $K\cap E_\tau\ne\emp$, and let $Q\in T(K)$. Then $\diam(K,E)\le\tau$, so that $\diam K\le\dist(K,E)\le\tau$. Therefore for each $Q\in T(K)$
%----------------------------------------------------------
$$
\diam Q\le 4\diam K\le 4\tau=80\delta.
$$
%----------------------------------------------------------
Using the same idea as in the proof of \rf{IM-DT}, we obtain the required inclusion $Q\subset E_\ve$.
%----------------------------------------------------------
\par We also note that the covering multiplicity of the family $\Ac$ is bounded by a constant $C=C(n)$. This easily follows from Lemma \reff{Wadd}.
%----------------------------------------------------------
\par Finally, we obtain
%----------------------------------------------------------
$$
\intl_{E_\tau\setminus E}|F_\ve(y)|^p\,dy
\le
C\,\smed\limits_{Q\in \Ac}\,\,
\intl_{Q}|P^{(V_E)}(x)|^p\,dx\le
C\,\,
\intl_{U_{\Ac}}|P^{(V_E)}(x)|^p\,dx
$$
%----------------------------------------------------------
where $U_{\Ac}:=\cup\{Q:Q\in \Ac\}$.
%----------------------------------------------------------
\par Since $Q\subset E_\ve$ for each $Q\in\Ac$, we conclude that $U_{\Ac}\subset E_\ve$ so that
%----------------------------------------------------------
$$
\intl_{E_\tau\setminus E}|F_\ve(y)|^p\,dy
\le
C\,\,
\intl_{E_\ve}|P^{(V_E)}(x)|^p\,dx.
$$
%----------------------------------------------------------
Combining this inequality with inequality \rf{F-MJ}, we obtain \rf{LP-I7} proving the lemma.\bx
%----------------------------------------------------------
%@@@@@@@@@@@@@@@@@@@@@@@@@@@@@@@@@@@@@@@@@@@@@@@@@@@@@@@@@@
%@@@@@@@@@@@@@@@@@@@@@@@@@@@@@@@@@@@@@@@@@@@@@@@@@@@@@@@@@@
%@@@@@@@@@@@@@@@@@@@@@@@@@@@@@@@@@@@@@@@@@@@@@@@@@@@@@@@@@@
%@@@@@@@@@@@@@@@@@@@@@@@@@@@@@@@@@@@@@@@@@@@@@@@@@@@@@@@@@@
%----------------------------------------------------------
\par We are in a position to finish the proof of the sufficiency part of Theorem \reff{WP-PR}. We know that the function $F_\ve$ agrees with the Whitney $(m-1)$-field $\VP=\{P_x: x\in E\}$ so that, by \rf{N-VP-W}, $$\PMEW\le\|F_\ve\|_{\WMP}.$$
%----------------------------------------------------------
\par It is well known, see, e.g. \cite{M}, p. 21, that for every $F\in\WMP$ the following equivalence
%----------------------------------------------------------
$$
\|F\|_{\WMP}\sim \|F\|_{\LMP}+\|F\|_{\LPRN}
$$
%----------------------------------------------------------
holds with constants depending only on $m,n,$ and $p$. Hence,
%----------------------------------------------------------
$$
\PMEW\le C\,\{\|F_\ve\|_{\LMP}+\|F_\ve\|_{\LPRN}\}.
$$
%----------------------------------------------------------
Combining this inequality with Corollary \reff{COR-FE} and Lemma \reff{L-3E}, we obtain \rf{IN-JE}.
%----------------------------------------------------------
\par The proof of Theorem \reff{WP-PR} is complete.\bx
\bigskip
%----------------------------------------------------------
%@@@@@@@@@@@@@@@@@@@@@@@@@@@@@@@@@@@@@@@@@@@@@@@@@@@@@@@@@@
%@@@@@@@@@@@@@@@@@@@@@@@@@@@@@@@@@@@@@@@@@@@@@@@@@@@@@@@@@@
%@@@@@@@@@@@@@@@@@@@@@@@@@@@@@@@@@@@@@@@@@@@@@@@@@@@@@@@@@@
%@@@@@@@@@@@@@@@@@@@@@@@@@@@@@@@@@@@@@@@@@@@@@@@@@@@@@@@@@@
%@@@@@@@@@@@@@@@@@@@@@@@@@@@@@@@@@@@@@@@@@@@@@@@@@@@@@@@@@@
%@@@@@@@@@@@@@@@@@@@@@@@@@@@@@@@@@@@@@@@@@@@@@@@@@@@@@@@@@@
%----------------------------------------------------------
\par {\it Proof of the sufficiency part of Theorem \reff{EX-WTK}.} Let $\gmh:=\tgm+12\theta$ where $\tgm$ is the constant from inequality \rf{PR-GM}. Let $\gamma$ be the same constant as in Theorem \reff{WP-PR}. (Recall that $\gamma$ coincides with the constant from Theorem \reff{EX-TK}.)\smallskip
%----------------------------------------------------------
\par Suppose that a Whitney $(m-1)$-field $\VP=\{P_x:x\in E\}$ satisfies the conditions of the sufficiency part of Theorem \reff{EX-WTK}. Thus:\smallskip
%----------------------------------------------------------
\par (a) The function $P^{(V)}$ belongs to $L_p(E_\ve)$. See \rf{PT-X};
\medskip
%----------------------------------------------------------
\par (ii) Let $\lambda:=\Nc^*(\VP)^p$, see\rf{NRM-W-NEW}. Then  $\lambda<\infty$ so that for every finite family $\{Q_i: i = 1,...,k\}$ of pairwise disjoint cubes contained in $E_{\veh}$, and every choice of points $x_i,y_i\in (\gmh Q_i)\cap E$ the following inequality

%----------------------------------------------------------
\bel{NRM-W}
\smed_{i=1}^k\,\,\smed_{|\alpha|\le m-1}\frac{|D^\alpha P_{x_i}(x_i)-D^\alpha P_{y_i}(x_i)|^p}
{(\diam Q_i)^{(m-|\alpha|)p-n}}\le \lambda
\ee
%----------------------------------------------------------
holds.\smallskip
%----------------------------------------------------------
\par Prove that $\VP=\{P_x: x\in E\}$ satisfies the hypothesis of the sufficiency part of Theorem \reff{WP-PR} with
%----------------------------------------------------------
\bel{EP-TG}
\ve:=\veh/(2\gamma^2).
\ee
%----------------------------------------------------------
More specifically, we will prove that the function $P^{(V_E)}$ defined by \rf{PVE-1} belongs to $L_p(E_\ve)$, and for every finite $\gamma$-sparse collection $\{\{x_i, y_i\} : i = 1, ..., k\}$ of two point subsets of $E$ with $\|x_i-y_i\|\le \ve$, $i=1,...,k$, the inequality \rf{N-P-GW} holds with
%---------------------------------------------------------
\bel{LM-WD}
\lw:=\gamma^{mp-n}\lcap.
\ee
%---------------------------------------------------------
%@@@@@@@@@@@@@@@@@@@@@@@@@@@@@@@@@@@@@@@@@@@@@@@@@@@@@@@@@
%----------------------------------------------------------
\par Furthermore, we will show that
%---------------------------------------------------------
\bel{PV-LPE}
\|P^{(V_E)}\|_{L_p(E_\ve)}\le C\, \{\|P^{(V)}\|_{L_p(E_{\veh})}+\lcap^{\frac1p}\}.
\ee
%---------------------------------------------------------
Here $C>0$ is a constant depending only on $m,n,p,\veh$ and $\theta$.
%----------------------------------------------------------
\medskip
\par Let $\{\{x_i, y_i\} : i = 1, ..., k\}$ be a finite $\gamma$-sparse collection of two point subsets of $E$ with $\|x_i-y_i\|\le \ve$, $i=1,...,k$. Then, by Definition \reff{DF-PR},  there exists a collection $\{Q_i: i=1,...k\}$ of pairwise disjoint cubes in $\RN$ such that $x_i,y_i\in \gamma Q_i$ and
%---------------------------------------------------------
\bel{DM-QI}
\diam Q_i\le \gamma \|x_i-y_i\|, ~~i=1,...,k.
\ee
%---------------------------------------------------------
Hence, $\diam Q_i\le \gamma\ve$. Since $x_i\in\gamma Q_i$,
we have $\|x_i-c_{Q_i}\|\le \gamma\diam Q_i/2$, so that for every $y\in Q_i$ the following inequality
%----------------------------------------------------------
$$
\dist(y,E)\le \|y-c_{Q_i}\|+\dist(c_{Q_i},E)\le\diam Q_i/2+\|x_i-c_{Q_i}\|\le \gamma\diam Q_i\le \gamma^2\ve.
$$
%----------------------------------------------------------
\par Since $\ve<\veh/\gamma^2$, see \rf{EP-TG}, we have $\dist(y,E)<\veh$ proving that $Q_i\subset E_{\veh}$ for every $i=1,...k$. Therefore, by the assumption, the inequality \rf{NRM-W} holds for the collection $\{\{x_i, y_i\} : i = 1, ..., k\}$.
%----------------------------------------------------------
\par Hence, by \rf{DM-QI},
%----------------------------------------------------------
\be
I&:=&\smed_{i=1}^k\,\,\smed_{|\alpha|\le m-1}\frac{|D^\alpha P_{x_i}(x_i)-D^\alpha P_{y_i}(x_i)|^p}
{\|x_i-y_i\|^{(m-|\alpha|)p-n}}
\nn\\
&\le&
\smed_{i=1}^k\,\,\gamma^{(m-|\alpha|)p-n}
\smed_{|\alpha|\le m-1}\frac{|D^\alpha P_{x_i}(x_i)-D^\alpha P_{y_i}(x_i)|^p}
{(\diam Q_i)^{(m-|\alpha|)p-n}},
\nn
\ee
%----------------------------------------------------------
so that, by \rf{NRM-W}, $I\le \gamma^{mp-n}\,\lcap=\lw$.
%----------------------------------------------------------
\par This proves inequality \rf{N-P-GW} with $\lw$ and $\ve$ defined by \rf{LM-WD} and \rf{EP-TG} respectively.
\medskip
%----------------------------------------------------------
\par Prove inequality \rf{PV-LPE}. We know that $P^{(V_E)}|_E=P^{(V)}|_E=f_{\VP}$ where
%---------------------------------------------------------
$$
f_{\VP}(x):=P_x(x),~~~x\in E.
$$
%---------------------------------------------------------
Hence,
%----------------------------------------------------------
$$
\|P^{(V_E)}\|^p_{L_p(E_\ve)}=\|f_{\VP}\|^p_{L_p(E)}+ \|P^{(V_E)}\|^p_{L_p(E_\ve\setminus E)}\le \|P^{(V)}\|^p_{L_p(E_{\veh})}+ \|P^{(V_E)}\|^p_{L_p(E_\ve\setminus E)}\,.
$$
%----------------------------------------------------------
\par Let
%----------------------------------------------------------
$$
\Ac_\ve:=\{Q\in W_E: Q\cap E_\ve\ne\emp\}.
$$
%----------------------------------------------------------
\par Note that for every $Q\in\Ac_\ve$ there exists a point $y\in Q$ such that $\dist(y,E)<\ve$. Hence we have  $\dist(Q,E)<\ve$ so that, by \rf{DQ-E},
%----------------------------------------------------------
\bel{RT-7}
\diam Q\le \dist(Q,E)<\ve~~~\text{for every}~~Q\in\Ac_\ve.
\ee
%----------------------------------------------------------
\par Now let $x\in Q$. Then
%----------------------------------------------------------
$$
\dist(x,E)\le \dist(y,E)+\|x-y\|\le \dist(y,E)+\diam Q.
$$
%----------------------------------------------------------
Since $Q\in W_E$,
%----------------------------------------------------------
$$
\dist(x,E)\le \dist(y,E)+\dist(Q,E)\le 2\dist(y,E)<2\ve
$$
%----------------------------------------------------------
proving that $Q\subset E_{2\ve}$. In particular,
%----------------------------------------------------------
\bel{Q-EE2}
Q\subset E_{\veh}~~~\text{for every}~~~Q\in\Ac_\ve.
\ee
%-------------------------------------------------------
See \rf{EP-TG}.
%-------------------------------------------------------
\par We obtain
%----------------------------------------------------------
$$
\|P^{(V_E)}\|^p_{L_p(E_\ve\setminus E)}\le
\smed_{Q\in\Ac_\ve}\,\,\|P^{(V_E)}\|^p_{L_p(Q)}\le
2^p\left\{\smed_{Q\in\Ac_\ve}\,\,\|P^{(V)}\|^p_{L_p(Q)}+I
\right\}
$$
%----------------------------------------------------------
where
%----------------------------------------------------------
$$
I:=\smed_{Q\in\Ac_\ve}\,\,\|P^{(V_E)}-P^{(V)}\|^p_{L_p(Q)}.
$$
%----------------------------------------------------------
\par Hence,
%----------------------------------------------------------
$$
\|P^{(V_E)}\|^p_{L_p(E_\ve\setminus E)}\le
2^p\left\{\|P^{(V)}\|^p_{L_p(\Uc_\ve)}+I\right\}
$$
%----------------------------------------------------------
where $\Uc_\ve:=\cup\{Q: Q\in\Ac_\ve\}$. By \rf{Q-EE2}, $\Uc_\ve\subset E_{\veh}$ so that
%----------------------------------------------------------
\bel{Y-4}
\|P^{(V_E)}\|^p_{L_p(E_\ve\setminus E)}\le
2^p\left\{\|P^{(V)}\|^p_{L_p(E_{\veh})}+I\right\}.
\ee
%----------------------------------------------------------
\smallskip
\par Prove that $I\le C(\ve)\,\lcap$.
%----------------------------------------------------------
\par Without loss of generality we may assume that the family $\Ac_\ve$ is finite, i.e.,
%----------------------------------------------------------
$$
\Ac_\ve=\{Q_i: i=1,...,k\}
$$
%----------------------------------------------------------
for some positive integer $k$. Furthermore, since the covering multiplicity of $\Ac_\ve$ is a bounded by a constant $C=C(n)$ (recall that $\Ac_\ve\subset W_E$), we may assume that the cubes of the family $\Ac_\ve$ are pairwise disjoint. See Theorem \reff{TFM}.
%----------------------------------------------------------
\par We recall that, by definition of $P^{(V_E)}$, see \rf{P-VE} and \rf{PVE-1}, for every cube $Q\in\Ac_\ve$
%----------------------------------------------------------
\bel{PVE-2}
P^{(V_E)}(x)=P_{a_Q}~~~\text{for every}~~~x\in Q.
\ee
%-------------------------------------------------------
Here $a_Q:=\PRL(L^{(Q)})$ and $L^{(Q)}$ is the (unique) lacuna containing $Q$. We know that $a_Q\in\tgm\,Q$, see \rf{PR-GM}.
%----------------------------------------------------------
\par Let $Q=Q_i\in\Ac_\ve$ and let
%----------------------------------------------------------
$$
\tau_Q=\tau_{Q_i}:= 2^{-\frac{i}{p}}\lcap^{\frac{1}{p}}\,|Q|^{-\frac{1}{p}},~~~~     i=1,...,k.
$$
%----------------------------------------------------------
\par By $x_Q=x_{Q_i}$ we denote a point in $Q$ such that
%----------------------------------------------------------
$$
\esssup_{Q}\,|P^{(V_E)}-P^{(V)}|\le
|P^{(V_E)}(x_Q)-P^{(V)}(x_Q)|+\tau_Q\,.
$$
%----------------------------------------------------------
Then
%----------------------------------------------------------
\be
\|P^{(V_E)}-P^{(V)}\|^p_{L_p(Q)}&\le& (|P^{(V_E)}(x_Q)-P^{(V)}(x_Q)|+\tau_Q)^p\,|Q|\nn\\
&\le&
2^p\,|P^{(V_E)}(x_Q)-P^{(V)}(x_Q)|^p\,|Q|+2^p\tau_Q^p\,|Q|.
\nn
\ee
%----------------------------------------------------------
Combining this inequality with \rf{PVE-2} and \rf{PT-X}, we obtain
%----------------------------------------------------------
\bel{P-14}
\|P^{(V_E)}-P^{(V)}\|^p_{L_p(Q)}\le
2^p\,|P_{a_Q}(x_Q)-P_{b_Q}(x_Q)|^p\,|Q|+2^{p-i}\lcap
\ee
%----------------------------------------------------------
provided $Q=Q_i$. Here $b_Q:=V(x_Q)$.
%----------------------------------------------------------
\par By \rf{T-EM}, $\|b_Q-x_Q\|\le\theta\dist(x_Q,E)$ so that
%----------------------------------------------------------
$$
\|b_Q-c_Q\|\le \|b_Q-x_Q\|+\|x_Q-c_Q\|\le
\theta\dist(x_Q,E)+\diam Q\le \theta\dist(Q,E)+2\diam Q.
$$
%----------------------------------------------------------
But $Q\in W_E$ so that $\dist(Q,E)\le 4\diam Q$. Hence
%----------------------------------------------------------
$$
\|b_Q-c_Q\|\le 4\theta\diam Q+2\diam Q\le 6\theta\diam Q
$$
%----------------------------------------------------------
proving that $b_Q\in 12\theta\, Q$. (Recall that $\theta\ge 1$.)
%----------------------------------------------------------
\par Let $H_Q:=P_{a_Q}-P_{b_Q}$. Since $H_Q\in\PMRN$,
%----------------------------------------------------------
$$
H_Q(x)=\smed_{|\alpha|\le m-1}\, \tfrac{1}{\alpha!} D^\alpha H_Q(a_Q)\,(x-a_Q)^\alpha,~~~~x\in\RN.
$$
%----------------------------------------------------------
Since $x_Q\in Q$ and $a_Q\in\tgm Q$, we have $\|x_Q-a_Q\|\le\tgm\diam Q$ so that
%----------------------------------------------------------
$$
|H_Q(x_Q)|\le \smed_{|\alpha|\le m-1}\,
|D^\alpha H_Q(a_Q)|\,\|x_Q-a_Q\|^{|\alpha|}
\le \tgm^m\,
\smed_{|\alpha|\le m-1}\,
|D^\alpha H_Q(a_Q)|\,(\diam Q)^{|\alpha|}
$$
%----------------------------------------------------------
\par Now, by \rf{P-14},
%----------------------------------------------------------
$$
\|P^{(V_E)}-P^{(V)}\|^p_{L_p(Q)}\le C\,
\smed_{|\alpha|\le m-1}\,
|D^\alpha H_Q(a_Q)|^p\,(\diam Q)^{|\alpha|p+n}+2^{p-i}\lcap
$$
%----------------------------------------------------------
where $Q=Q_i$. Hence,
%----------------------------------------------------------
\be
I&=&\smed_{Q\in\Ac_\ve}\,\,\|P^{(V_E)}-P^{(V)}\|^p_{L_p(Q)}
\nn\\
&\le&
C\,\smed_{Q\in\Ac_\ve}\,\smed_{|\alpha|\le m-1}\,
|D^\alpha H_Q(a_Q)|^p\,(\diam Q)^{|\alpha|p+n}
+2^p\,\lcap\,\smed_{i=1}^k\,2^{-i}\nn\,.
\ee
%----------------------------------------------------------
\par By \rf{RT-7}, $\diam Q<\ve$ so that
%----------------------------------------------------------
$$
I\le
C\ve^{mp}\,\smed_{Q\in\Ac_\ve}\,\smed_{|\alpha|\le m-1}\,
\frac{|D^\alpha H_Q(a_Q)|^p}{(\diam Q)^{(m-|\alpha|)p-n}}
+2^{p+1}\lcap.
$$
%----------------------------------------------------------
Recall that $a_Q\in\tgm Q$ and $b_Q\in (12\theta)Q$ so that $a_Q,b_Q\in (\tgm+12\theta)Q=\gmh\,Q$.
%----------------------------------------------------------
\par Hence, by \rf{NRM-W},
%----------------------------------------------------------
$$
I\le
C\ve^{mp}\,\lcap+2^{p+1}\lcap.
$$
%----------------------------------------------------------
proving the required inequality $I\le C(\ve)\lcap$.
%----------------------------------------------------------
\par Combining this inequality with inequality \rf{Y-4}, we conclude that inequality \rf{PV-LPE} holds.  Thus all conditions of the hypothesis of Theorem \reff{WP-PR} are satisfied (with $\lw$ defined by \rf{LM-WD}). By this theorem there exists a $C^{m-1}$-function $F\in\WMP$ which agrees with $\VP$ on $E$ (i.e., \rf{J-PE-W} holds). Furthermore, by \rf{H-JE},
%---------------------------------------------------------
$$
\PMEW\le C\,\{\|P^{(V_E)}\|_{L_p(E_\ve)}+\lw^{\frac1p}\}
$$
%---------------------------------------------------------
so that, by \rf{PV-LPE},
%---------------------------------------------------------
$$
\PMEW\le C\,\{\|P^{(V)}\|_{L_p(E_{\veh})}+\lcap^{\frac1p}+
(\gamma^{mp-n}\lcap)^{\frac1p}\}
$$
%---------------------------------------------------------
proving the required inequality
%---------------------------------------------------------
$$
\PMEW\le C\,\{\|P^{(V)}\|_{L_p(E_{\veh})}+\lcap^{\frac1p}\}.
$$
%---------------------------------------------------------
\par The proof of Theorem \reff{EX-WTK} is complete.\bx
%----------------------------------------------------------
%@@@@@@@@@@@@@@@@@@@@@@@@@@@@@@@@@@@@@@@@@@@@@@@@@@@@@@@@@@
%@@@@@@@@@@@@@@@@@@@@@@@@@@@@@@@@@@@@@@@@@@@@@@@@@@@@@@@@@@
%@@@@@@@@@@@@@@@@@@@@@@@@@@@@@@@@@@@@@@@@@@@@@@@@@@@@@@@@@@
%@@@@@@@@@@@@@@@@@@@@@@@@@@@@@@@@@@@@@@@@@@@@@@@@@@@@@@@@@@
%----------------------------------------------------------
\begin{remark} \lbl{D-CR}{\em The equivalence \rf{EQV-JET-WP} and Lemma \reff{PLM-M} imply the following ``discrete'' version of the result of Theorem \reff{EX-WTK}: Let $Q\in W_E$ be a Whitney cube with $\diam Q\le\veh$. Let $z_Q\in Q$ be an {\it arbitrary} point in $Q$ and let $t_Q:=\VE(z_Q)$. Finally, let
$f_{\VP}(x):=P_x(x)$, $x\in E$.
%----------------------------------------------------------
\par Then for every Whitney $(m-1)$-field  $\VP=\{P_x: x\in E\}$ on $E$ the following equivalence
%---------------------------------------------------------
$$
\PMEW\sim \,\|f_{\VP}\|_{L_p(E)}+\Nc^*(\VP)\,+\left
\{\,\sbig_{
\substack {Q\in\, W_E,\vspace*{1mm}\\ \diam Q\,\le\, {\veh}}}
\,\,\,\sbig_{|\alpha|\le m-1}\,\,\,
(\diam Q)^{|\alpha|p+n}\,
|D^\alpha P_{t_Q}(t_Q)|^p\,\right\}^{\frac1p}
$$
%---------------------------------------------------------
holds. The constants of this equivalence depend only on $m,n,p,\veh,$ and $\theta$.}\rbx
%----------------------------------------------------------
\end{remark}
%@@@@@@@@@@@@@@@@@@@@@@@@@@@@@@@@@@@@@@@@@@@@@@@@@@@@@@@@@@
%@@@@@@@@@@@@@@@@@@@@@@@@@@@@@@@@@@@@@@@@@@@@@@@@@@@@@@@@@@
%@@@@@@@@@@@@@@@@@@@@@@@@@@@@@@@@@@@@@@@@@@@@@@@@@@@@@@@@@@
%@@@@@@@@@@@@@@@@@@@@@@@@@@@@@@@@@@@@@@@@@@@@@@@@@@@@@@@@@@
%----------------------------------------------------------
\begin{remark} \lbl{D-CR1}{\em Using the result of Theorem \reff{WP-PR}, Lemma \reff{PLM-M} and properties of lacunae described in Subsections 6.1 and 6.2,  we also obtain a corresponding ``discrete'' version of the criterion \rf{H-JE}:
%----------------------------------------------------------
\par Let $L\in\LE$ be a lacuna, and let
$\diam L:=\sup\{\diam Q: Q\in L\}$. Let $s_L:=\PRL(L)$,
%----------------------------------------------------------
\par Then for every Whitney $(m-1)$-field
$\VP=\{P_x: x\in E\}$
%---------------------------------------------------------
$$
\PMEW\sim \,\|f_{\VP}\|_{L_p(E)}+\Nc^\flat(\VP)+
\,\left\{\sbig_{L\in\,\LE}
\,\,\,\sbig_{|\alpha|\le m-1}\,\,\,
\min\{\ve,\diam L\}^{|\alpha|p+n}\,
|D^\alpha P_{s_L}(s_L)|^p\,\right\}^{\frac1p}\,.
$$
%---------------------------------------------------------
The constants in this equivalence depend only on $m,n,p$ and $\ve$.}\rbx
%----------------------------------------------------------
\end{remark}
%----------------------------------------------------------
%@@@@@@@@@@@@@@@@@@@@@@@@@@@@@@@@@@@@@@@@@@@@@@@@@@@@@@@@@@
%@@@@@@@@@@@@@@@@@@@@@@@@@@@@@@@@@@@@@@@@@@@@@@@@@@@@@@@@@@
%@@@@@@@@@@@@@@@@@@@@@@@@@@@@@@@@@@@@@@@@@@@@@@@@@@@@@@@@@@
%@@@@@@@@@@@@@@@@@@@@@@@@@@@@@@@@@@@@@@@@@@@@@@@@@@@@@@@@@@
\par Our last result is an analogue of Theorem \reff{LIN-OP} for the normed Sobolev space $\WMP$, $p>n$.
\par Let $\ve=1$. We know that the operator $\JV(\WMP)|_E\ni\VP\to F_\ve$ defined by formula \rf{W-DEF-F} provides an almost optimal ``extensions'' of $(m-1)$-jets generated by Sobolev $\WMP$-functions. Since this operator is {\it linear}, we obtain the following
%----------------------------------------------------------
%@@@@@@@@@@@@@@@@@@@@@@@@@@@@@@@@@@@@@@@@@@@@@@@@@@@@@@@@@@
%@@@@@@@@@@@@@@@@@@@@@@@@@@@@@@@@@@@@@@@@@@@@@@@@@@@@@@@@@@
%@@@@@@@@@@@@@@@@@@@@@@@@@@@@@@@@@@@@@@@@@@@@@@@@@@@@@@@@@@
%@@@@@@@@@@@@@@@@@@@@@@@@@@@@@@@@@@@@@@@@@@@@@@@@@@@@@@@@@@
%@@@@@@@@@@@@@@@@@@@@@@@@@@@@@@@@@@@@@@@@@@@@@@@@@@@@@@@@@@
%@@@@@@@@@@@@@@@@@@@@@@@@@@@@@@@@@@@@@@@@@@@@@@@@@@@@@@@@@@
%----------------------------------------------------------
\begin{theorem} \lbl{LIN-WOP} For every closed subset $E\subset\RN$ and every $p>n$ there exists a continuous linear operator $\Fcw:\,\JV(\WMP)|_E\,\to\, \WMP$
%----------------------------------------------------------
such that for every Whitney $(m-1)$-field $$\VP=\{P_x:x\in E\}\in \JV(\WMP)|_E$$ the function $\Fcw(\VP)$ agrees with $\VP$ on $E$.
%----------------------------------------------------------
\par Furthermore, the operator norm of $\Fcw$ is bounded by a constant depending only on $m,n,$ and $p$.
%----------------------------------------------------------
\end{theorem}
%----------------------------------------------------------
%@@@@@@@@@@@@@@@@@@@@@@@@@@@@@@@@@@@@@@@@@@@@@@@@@@@@@@@@@@
%@@@@@@@@@@@@@@@@@@@@@@@@@@@@@@@@@@@@@@@@@@@@@@@@@@@@@@@@@@
%@@@@@@@@@@@@@@@@@@@@@@@@@@@@@@@@@@@@@@@@@@@@@@@@@@@@@@@@@@
%@@@@@@@@@@@@@@@@@@@@@@@@@@@@@@@@@@@@@@@@@@@@@@@@@@@@@@@@@@
%@@@@@@@@@@@@@@@@@@@@@@@@@@@@@@@@@@@@@@@@@@@@@@@@@@@@@@@@@@
%@@@@@@@@@@@@@@@@@@@@@@@@@@@@@@@@@@@@@@@@@@@@@@@@@@@@@@@@@@
%----------------------------------------------------------
%&&&&&&&&&&&&&&&&&&&&&&&&&&&&&&&&&&&&&&&&&&&&&&&&&&&&&&&&&&
%                                                         &
%                      REFERENCES                         &
%_________________________________________________________&
%&&&&&&&&&&&&&&&&&&&&&&&&&&&&&&&&&&&&&&&&&&&&&&&&&&&&&&&&&&

%@@@@@@@@@@@@@@@@@@@@@@@@@@@@@@@@@@@@@@@@@@@@@@@@@@@@@@@@@@

\begin{thebibliography}{ABCD}
%@@@@@@@@@@@@@@@@@@@@@@@@@@@@@@@@@@@@@@@@@@@@@@@@@@@@@@@@@@
%----------------------------------------------------------
\addtocontents{toc}{References \hfill \thepage\par}
%----------------------------------------------------------
%@@@@@@@@@@@@@@@@@@@@@@@@@@@@@@@@@@@@@@@@@@@@@@@@@@@@@@@@@@
%@@@@@@@@@@@@@@@@@@@@@@@@@@@@@@@@@@@@@@@@@@@@@@@@@@@@@@@@@@
%@@@@@@@@@@@@@@@@@@@@@@@@@@@@@@@@@@@@@@@@@@@@@@@@@@@@@@@@@@
%----------------------------------------------------------
%\bibitem {Ad} R. Adams, Sobolev Spaces, Academic Press, %New York, 1975.
%----------------------------------------------------------
\bibitem {BMP1} E. Bierstone, P. Milman, W. Pawlucki, Differentiable functions defined in closed sets. A problem of Whitney. Invent. Math. {151}, No. 2, (2003) 329--352.
%----------------------------------------------------------
\bibitem {Br} Yu. A. Brudnyi, Spaces that are definable by means of local approximations, Trudy Moscov. Math. Obshch. 24 (1971) 69--132; English transl.: Trans. Moscow Math. Soc. 24 (1974) 73--139.
%----------------------------------------------------------
\bibitem {BrK} Yu. A. Brudnyi, B. D. Kotljar, A certain problem of combinatorial geometry, Sibirsk. Mat. Zh. 11 (1970) 1171--1173; English transl.: Siberian Math. J. 11 (1970) 870--871.
%----------------------------------------------------------
\bibitem {BS1} Yu. Brudnyi, P. Shvartsman,
Generalizations  of Whitney  Extension Theorem,
Intern. Math.  Research Notices, No. 3, (1994) 129--139.
%----------------------------------------------------------
\bibitem {BS2} Yu. Brudnyi, P. Shvartsman,
The Whitney Problem  of  Existence of a Linear
Extension Operator, J. Geom. Anal., {7}, No. 4 (1997) 515--574.
%----------------------------------------------------------
\bibitem {BS3} Yu. Brudnyi, P. Shvartsman,
Whitney Extension Problem for Multivariate
$C^{1,\omega}$-functions, Trans. Amer. Math. Soc. {353} No. 6,  (2001) 2487--2512.
%----------------------------------------------------------
\bibitem {C1} A. P. Calder\'{o}n, Estimates for
    singular integral operators in terms of maximal
    functions, Studia Math. 44 (1972) 563--582.
%----------------------------------------------------------
\bibitem {CS}  A. P. Calder\'{o}n, R. Scott, Sobolev
    type inequalities for $p>0$,  Studia Math. 62 (1978)
    75--92.
%----------------------------------------------------------
\bibitem {Dol} V. L. Dolnikov, The partitioning of families of convex bodies, Sibirsk. Mat. Zh. 12 (1971) 664--667 (in Russian); English transl.: Siberian Math. J. 12 (1971) 473--475.
%----------------------------------------------------------
\bibitem {F2} C. Fefferman, A sharp form of Whitney extension theorem, Annals of Math. 161, No. 1 (2005) 509--577.
%----------------------------------------------------------
\bibitem {F-J} C. Fefferman, A Generalized Sharp Whitney Theorem for Jets, Rev. Mat. Iberoamericana 21, no.2, (2005) 577--688.
%----------------------------------------------------------
\bibitem {F4} C. Fefferman, Whitney extension problem for $C^m$, Annals of Math. 164, no. 1, (2006) 313--359.
%----------------------------------------------------------
\bibitem {F3} C. Fefferman,  $C^m$ Extension by Linear Operators, Annals of Math. 166, No. 3, (2007) 779--835.
%----------------------------------------------------------
\bibitem {F-IO} C. Fefferman, Extension of $C^{m,\omega}$-Smooth Functions by Linear Operators, Rev. Mat. Iberoamericana 25, No. 1, (2009) 1--48.
%----------------------------------------------------------
\bibitem {F-Bl} C. Fefferman, Whitney extension problems and interpolation of data, Bulletin A.M.S. 46, no. 2 (2009) 207--220.
%----------------------------------------------------------
\bibitem {F9} C. Fefferman, The $C^m$ norm of a function with prescribed jets I., Rev. Mat. Iberoamericana 26 (2010) 1075--1098.
%----------------------------------------------------------
\bibitem {FIL} C. Fefferman, A. Israel, G. K. Luli,
Sobolev extension by linear operators, J. Amer. Math. Soc. 27 (2014), 69--145.
%----------------------------------------------------------
\bibitem {G}  G. Glaeser, \'Etude de quelques algebres Tayloriennes, J. d'Analyse Math. {6} (1958) 1--125.
%----------------------------------------------------------
\bibitem {Guz} M. de Guzm\'{a}n, Differentiation of
    integrals in $\RN$, Lect. Notes in Math. 481,
    Springer-Verlag, 1975.
%----------------------------------------------------------
\bibitem {HHC} A. Herbert-Voss, M. J. Hirn, F. McCollum, Computing minimal interpolants in $C^{1,1}(\R^d)$, arXiv:1411.5668 (2016) 41 pp.
%----------------------------------------------------------
\bibitem {HL} M. J. Hirn, E. Le Gruyer, A general theorem of existence of quasi absolutely minimal Lipschitz extensions, Math. Ann. 359 (2014) 595--628.
%----------------------------------------------------------
\bibitem {Is} A. Israel, A Bounded Linear Extension Operator for $L^{2,p}(\R^2)$, Annals of Math.  178 (2013) 1--48.
%----------------------------------------------------------
\bibitem {LG} E. Le Gruyer, Minimal Lipschitz extensions to differentiable functions defined on a Hilbert space, GAFA 19 (2009) 1101--1118.
%----------------------------------------------------------
\bibitem {M}  V.G. Maz'ja,  Sobolev spaces,
     Springer-Verlag, Berlin, 1985, xix+486 pp.
%----------------------------------------------------------
\bibitem {MP} V. Maz'ya,  S. Poborchi, Differentiable
    Functions on Bad Domains, Word Scientific, River Edge,
    NJ, 1997.
%----------------------------------------------------------
\bibitem {McS} E. J. McShane, Extension of range of functions, Bull. Amer. Math. Soc. 40 no. 12, (1934) 837--842.
%----------------------------------------------------------
\bibitem {Sh1} P. Shvartsman, The Whitney extension problem and Lipschitz selections of set-valued mappings in jet-spaces, Trans. Amer. Math. Soc. 360, No. 10, (2008) 5529--5550.
%-----------------------------------------------------------
\bibitem {Sh2} P. Shvartsman, Sobolev $W^1_p$-spaces on closed subsets of $\RN$, Advances in Math. 220 (2009) 1842--1922.
%-----------------------------------------------------------
\bibitem {Sh3} P. Shvartsman, Sobolev $L^2_p$-functions on closed subsets of $\R^2$, Advances in Math. 252 (2014) 22--113.
%----------------------------------------------------------
\bibitem {Sh4} P. Shvartsman, On the sum of a Sobolev space and a weighted $L_p$-space, Advances in Math. 248 (2013) 155--228.
%----------------------------------------------------------
\bibitem {St} E. M. Stein, Singular integrals and
    differentiability properties of functions, Princeton
    Univ. Press, Princeton, New Jersey, 1970.
%----------------------------------------------------------
\bibitem {WEL} J. C. Wells, Differentiable functions on Banach spaces with Lipschitz derivatives, J. Differ. Geom. 8 (1973) 135--152.
%----------------------------------------------------------
\bibitem {W1}  H. Whitney, Analytic extension of differentiable functions defined in closed sets, Trans.  Amer. Math. Soc.  36 (1934) 63--89.
%----------------------------------------------------------
\bibitem {W2} H. Whitney, Differentiable functions defined in closed sets. I., Trans. Amer. Math. Soc. 36(1934) 369--387.
%----------------------------------------------------------
\bibitem {Z1} N. Zobin, Whitney's problem on    extendability of functions and an intrinsic metric,    Advances in Math. 133 (1998) 96--132.
%----------------------------------------------------------
\bibitem {Z2} N. Zobin, Extension of smooth functions from finitely connected planar domains, J. Geom. Anal. 9, no. 3, (1999) 489--509.
%----------------------------------------------------------
%@@@@@@@@@@@@@@@@@@@@@@@@@@@@@@@@@@@@@@@@@@@@@@@@@@@@@@@@@@
\end{thebibliography}
\end{document}